\documentclass[12pt,fleqn]{tsp}

\usepackage[hang,flushmargin]{footmisc}
\techscience{vol}{no}{IJNME RLT 90th birthday \red{v2.3.7 arXiv} {\scriptsize (submitted v1 on 2023.10.15)} - \today}{0}
\techscience{vol}{no}{\red{v2.3.10 arXiv}, \today\  \ding{228} IJNME, doi:\href{https://doi.org/10.1002/nme.7586}{10.1002/nme.7586} {\scriptsize (online as of 2024.10.17)}}{0}

\license{}
\def\thedoi{}
\def\thetype{ARTICLE}

\graphicspath{{./image/}}% Images path
\runningtitle{Manuscript Preparation for TSP}

\usepackage{amsmath}
\usepackage{mathtools}
\usepackage{amssymb}
\usepackage{amsthm}
\usepackage{upgreek}

\usepackage{pifont}

\usepackage[linesnumbered,lined,commentsnumbered]{algorithm2e}

\usepackage{accents}
\DeclareMathSymbol{\widehatsym}{\mathord}{largesymbols}{"62}
\DeclareMathSymbol{\widetildesym}{\mathord}{largesymbols}{"65}
\usepackage{stackengine}
\usepackage{scalerel}

\usepackage{bbm}

\usepackage{boondox-cal}
\usepackage[scr=boondoxo,scrscaled=1.05]{mathalfa}
\usepackage[table,xcdraw,dvipsnames]{xcolor}

\usepackage{listings}
\definecolor{codegreen}{rgb}{0,0.6,0}
\definecolor{codegray}{rgb}{0.5,0.5,0.5}
\definecolor{codepurple}{rgb}{0.58,0,0.82}
\definecolor{backcolour}{rgb}{0.95,0.95,0.92}
\lstdefinestyle{mystyle}{
	backgroundcolor=\color{backcolour},   
	commentstyle=\color{codegreen},
	keywordstyle=\color{magenta},
	numberstyle=\tiny\color{codegray},
	stringstyle=\color{codepurple},
	basicstyle=\ttfamily\scriptsize,
	breakatwhitespace=false,         
	breaklines=true,                 
	captionpos=b,                    
	keepspaces=true,                 
	numbers=left,                    
	numbersep=5pt,                  
	showspaces=false,                
	showstringspaces=false,
	showtabs=false,                  
	tabsize=2
}
\usepackage{graphicx,float,wrapfig}
\usepackage{subcaption}

\usepackage{tabularx,booktabs}
\newcolumntype{Y}{>{\centering\arraybackslash}X}

\usepackage{multirow}

\usepackage{times}

\usepackage[margin=1cm]{caption}

\usepackage{enumitem}
\usepackage{setspace}
\setlist[itemize,enumerate]{topsep=.1em,itemsep=.1em,parsep=0em,partopsep=0em}
\setlist*[itemize,enumerate]{first=\vskip\baselineskip\setstretch{1.0}\vspace{-\baselineskip}}

\usepackage{clesm_defs_1}
\newtheorem{rem}{Remark}[section]

\usepackage[title]{appendix}

\newcommand{\loss}{J}
\newcommand{\dtb}[1]{\hbox{$#1$}\kern-1.4ex\raise+1.8ex\hbox{\bigdot}\kern+1.4ex}

\newcommand{\dtbs}[1]{\hbox{\scriptsize $#1$}\kern-1.0ex\raise+1.4ex\hbox{$\centerdot$}\kern+1.0ex}

\newcommand{\x}{x}
\newcommand{\X}{X}

\newcommand{\learn}{\epsilon}

\newcommand{\M}{M}
\newcommand{\Ms}[1]{\M_{#1}}
\newcommand{\y}{y}
\newcommand{\Y}{Y}
\newcommand{\by}{\boldsymbol{\y}}

\newcommand{\slope}{s}

\newcommand{\ndt}[1]{\stackon[0.2pt]{$#1$}{\vstretch{0.55}{\hstretch{0.55}{\bullet}}}}

\newcommand{\nddt}[1]{\stackon[0.2pt]{$#1$}{\vstretch{0.55}{\hstretch{0.55}{\bullet\bullet}}}}

\newcommand{\Xb}{\karb \X}
\newcommand{\Yb}{\karb \Y}
\newcommand{\xb}{\karb \x}
\newcommand{\yb}{\karb \y}
\newcommand{\ti}{t}
\newcommand{\tb}{\karb \ti}
\newcommand{\ntip}{\tau} 

\newcommand{\sbv}{e}
\newcommand{\bsbv}[1]{\bkars{\sbv}{#1}}
\newcommand{\cbv}{t}
\newcommand{\bcbv}[1]{\bkars{\cbv}{#1}}
\let\bmbv\bcbv 

\newcommand{\pf}{\mathcal{b}}
\newcommand{\barf}[1]{\pf (#1)}
\let\psf\barf
\newcommand{\pfp}[1]{\karp{\pf}{#1}}
\newcommand{\bde}{\mathcal{d}}
\newcommand{\bwe}{\mathcal{w}}

\definecolor{myOrange}{RGB}{255,80,0}
\newcommand{\DDET}{\text{ \href{https://deepxde.readthedocs.io/en/latest/}{\color{myOrange}\underline{DDE}}-\href{https://deepxde.readthedocs.io/en/latest/user/installation.html?highlight=backend}{\color{myOrange}\underline{T}}}}
\newcommand{\JAX}{\href{https://jax.readthedocs.io/en/latest/index.html}{\color{Plum} \underline{JAX}}}
\newcommand{\DvJ}{\text{\DDET\ vs\ \JAX}}

\newcommand{\decayf}{\mathcal{f}_\text{decay}}

\newcommand{\dm}{\phi}
\newcommand{\bdm}{\bkar{\dm}}
\newcommand{\bdmz}{\kars{\bdm}{0}}

\newcommand{\du}{u}
\newcommand{\ub}{\karb \du}
\newcommand{\uc}{\check{\du}}
\newcommand{\up}[1]{\karp \du {(#1)}}
\newcommand{\ubd}{\ndt \ub}
\newcommand{\ubds}[1]{\kars{\ubd}{#1}}
\newcommand{\ubs}[1]{\kars{\ub}{#1}}
\newcommand{\ubp}[1]{\karp{\ub}{(#1)}}
\newcommand{\usta}{\kars{\ub}{\text{st}}}
\newcommand{\bub}{\bkar{\ub}}
\newcommand{\busta}{\bkars{\ub}{\text{st}}}
\newcommand{\dv}{v}
\newcommand{\vb}{\karb \dv}
\newcommand{\vp}[1]{\karp \dv {(#1)}}
\newcommand{\vbd}{\ndt \vb}
\newcommand{\vbp}[1]{\vb^{(#1)}}
\newcommand{\vsta}{\kars{\vb}{\text{st}}}

\newcommand{\ext}{e}
\newcommand{\extb}{\karb \ext}
\newcommand{\ope}{\mathscr{c}}
\newcommand{\opeb}{\karb \ope}
\newcommand{\iope}{\mathscr{k}}
\newcommand{\iopeb}{\karb \iope}
\newcommand{\F}{F}
\newcommand{\Fp}[1]{\F^{(#1)}}
\newcommand{\Fb}{\karb \F}
\newcommand{\Fbp}[1]{\Fb^{(#1)}}

\newcommand{\force}{f}
\newcommand{\df}{\force}
\newcommand{\dfs}[1]{\kars \force #1}
\newcommand{\dfsp}[2]{\karsp \force #1 {#2}}
\newcommand{\dfb}{\karb \force}
\newcommand{\dfbs}[1]{\kars \dfb #1}
\newcommand{\dfbhs}[1]{\kars{\widehat{\dfb}}{#1}}

\newcommand{\shear}{S}
\newcommand{\shrb}{\karb \shear}
\newcommand{\Sbs}[1]{\kars{\shrb}{#1}}
\newcommand{\shrbp}[1]{\karpb \shear {(#1)}}
\newcommand{\Sbsp}[2]{\karspb \shear #1 {(#2)}}

\newcommand{\N}{N}
\newcommand{\Ns}[1]{\kars{\N}{#1}}
\newcommand{\Nb}{\karb \N}
\newcommand{\Nbs}[1]{\kars{\Nb}{#1}}
\newcommand{\Nbsp}[2]{\karspb \N #1 {(#2)}}

\newcommand{\rtdn}{\mathcal{r}}
\newcommand{\sldn}{\mathcal{s}}

\renewcommand{\M}{M}
\newcommand{\Mp}[1]{\M^{(#1)}}
\newcommand{\Mb}{\karb \M}
\newcommand{\Mbp}[1]{\Mb^{(#1)}}

\newcommand{\p}{p}
\newcommand{\bpb}{\bkarb{\p}}
\newcommand{\pbs}[1]{\karsb \p #1}
\newcommand{\pbsd}[1]{\ndt {\karb \p}_{#1}}

\newcommand{\rc}{\theta}
\newcommand{\bRotM}{\bkar{\Lambda}}

\newcommand{\dpde}{\mathcal{D}}
\newcommand{\bdpde}{\boldsymbol{\dpde}}
\newcommand{\bdpdep}[1]{\karp{\bdpde}{#1}}
\newcommand{\dpdesp}[2]{\karsp{\dpde}{#1}{#2}}
\newcommand{\spde}{\mathcal{S}}
\newcommand{\bspde}{\boldsymbol{\spde}}
\newcommand{\bspdep}[1]{\karp{\bspde}{#1}}
\newcommand{\spdes}[1]{\kars{\spde}{#1}}
\newcommand{\spdesp}[2]{\karsp{\spde}{#1}{#2}}

\newcommand{\eslp}{\alpha}
\newcommand{\eslpp}[1]{\karp \eslp {(#1)}}
\newcommand{\eslpb}{\karb \eslp}
\newcommand{\eslpbp}[1]{\karp \eslpb {(#1)}}
\newcommand{\slps}[1]{\kars \slope #1}
\newcommand{\slpsp}[2]{\karsp \slope #1 {(#2)}}
\newcommand{\slpbs}[1]{\karsb \slope #1}
\newcommand{\slpbsp}[2]{\karspb \slope #1 {(#2)}}

\newcommand{\red}[1]{{\color{red} #1}}

\newcommand{\shiftPs}{\mathscr{s}}
\newcommand{\shiftPc}{\mathcal{c}}
\newcommand{\shiftPa}{\mathcal{a}}
\newcommand{\shiftPA}{\mathscr{A}}

\newcommand{\sstr}{\gamma}
\newcommand{\bsstr}{\boldsymbol{\sstr}}
\newcommand{\sstrs}[1]{\kars \sstr #1}
\newcommand{\cstr}{\Gamma}
\newcommand{\cstrs}[1]{\kars \cstr #1}

\newcommand{\wl}{\mathcal{w}}
\newcommand{\wls}[1]{\kars{\wl}{#1}}
\newcommand{\wlsp}[2]{\karsp{\wl}{#1}{#2}}

\usepackage[yyyymmdd,hhmmss]{datetime}

\newcommand{\abc}{}

\newcommand{\kr}[1]{#1}

\newcommand{\krb}[1]{\overline{#1}}

\newcommand{\karb}[1]{\abc \krb{#1} \abc{}}

\newcommand{\kars}[2]{\abc {\kr {#1}} _{#2}\abc{}}

\newcommand{\karsb}[2]{\abc {\krb {#1}} _{#2}\abc{}}

\newcommand{\karp}[2]{\abc {\kr {#1}} ^{#2}\abc{}}

\newcommand{\karpb}[2]{\abc {\krb {#1}} ^{#2}\abc{}}

\newcommand{\karsp}[3]{\abc  {\kr {#1}} _{#2} ^{#3}\abc{}}

\newcommand{\karspb}[3]{\abc  {\krb {#1}} _{#2} ^{#3}\abc{}}

\newcommand{\bkr}[1]{{\bfg{#1}}}

\newcommand{\bkrb}[1]{{\overline{\bfg{#1}}}}

\newcommand{\bkar}[1]{\abc \bkr{#1}\abc{}}

\newcommand{\bkarb}[1]{\abc \bkrb{#1}\abc{}}

\newcommand{\bkars}[2]{\abc {\bkr {#1}} _{#2}\abc{}}

\usepackage[misc]{ifsym}

\usepackage[pagebackref]{hyperref}
\hypersetup{
	colorlinks = true,
	linkbordercolor = {white},
	linkcolor=blue,
	citecolor=[RGB]{0,50,0},
	urlcolor=[RGB]{100,0,0},
}

\title{
	Partial-differential-algebraic equations of nonlinear dynamics by Physics-Informed Neural-Network: (I) Operator splitting and framework assessment
	\normalsize
		\vspace{12pt}
	\vskip -0.2cm
	\hskip -0.5ex
	\emph{Dedicated to Professor Robert L. Taylor for his 90th birthday.}
}
\author{
	Loc Vu-Quoc\thanks{
		$^{, \,\textrm{\Letter} \,}$Aerospace Engineering, University of Illinois at Urbana-Champaign, IL 61801, USA $\bullet$
		$^{\textrm{\Letter} \,}${vql@illinois.edu}
	}$^{\bm , \, \textrm{\Letter}}$ and \
	Alexander Humer\thanks{
		Institute of Technical Mechanics, Johannes Kepler University, A-4040 Linz, Austria $\bullet$
		{alexander.humer@jku.at}
	}\	 
}

\usepackage{xcolor}
\usepackage{mdframed}
\definecolor{mygray}{gray}{0.92}

\newsavebox\myVerb
\newenvironment{verbbox}{\lrbox\myVerb}{\endlrbox}
\newcommand*{\verbBox}{\usebox\myVerb}
\newsavebox\myVerbT
\newenvironment{verbboxT}{\lrbox\myVerbT}{\endlrbox}
\newcommand*{\verbBoxT}{\usebox\myVerbT}

\begin{document}
\maketitle

%\begin{center}
%	{\it Dedicated for the 90th birthday of Professor Robert L. Taylor.}
%\end{center}

\vskip -0.18cm
\begin{mdframed}[
	% style=mystyle,
	%	backgroundcolor=lightgray,
	backgroundcolor=mygray,
	topline=false,
	bottomline=false,
	leftline=false,
	rightline=false
	]
	
	{
		% \fontfamily{crimson}
		\fontsize{10pt}{12pt}\selectfont
		
		\noindent	
		%	{\bf ABSTRACT}
		{\bf ABSTRACT}
		%\newline
		
		\noindent
Several forms for constructing {novel} physics-informed neural-networks (PINNs) for the solution of partial-differential-algebraic equations (PDAEs) based on derivative operator splitting are proposed, using the nonlinear Kirchhoff rod as a prototype for demonstration.

%A formulation for the partial-differential-algebraic equations (PDAEs) of motion of geometrically-exact rods with no shear, known as the Kirchhoff-Love rod, is presented together with several forms for constructing physics-informed neural-networks for the solution of these PDAEs based on derivative operator splitting.

The present work is a natural extension of our review paper in \cite{vuquoc2023deep} aiming at both experts and first-time learners of both deep learning and PINN frameworks, among which the open-source DeepXDE (DDE) \cite{lu2021deepxde} is likely the most well documented framework with many examples.
Yet, we encountered some pathological problems (time shift, amplification, static solutions) and proposed {novel} methods to resolve them.   

Among these {novel} methods are the PDE forms, which evolve
from the lower-level form with fewer unknown dependent variables (e.g., displacements, slope, finite extension) to higher-level form with more dependent variables (e.g., forces, moments, momenta, etc.), in addition to those from lower-level forms.   {Traditionally}, the highest-level form, the balance-of-momenta form, is the starting point for (hand) deriving the lowest-level form through a tedious (and error prone) process of successive substitutions.  The next step in a finite element method is to discretize the lowest-level form upon forming a weak form and linearization with appropriate interpolation functions, followed by their implementation in a code and testing.  The time-consuming tedium in all of these steps could be bypassed by applying {the proposed novel} PINN directly to the highest-level form.  

We also developed a script based on \JAX, the High Performance Array Computing. For the axial motion of elastic bar, while our \JAX\ script did not show the pathological problems of \DDET\ (DDE with TensorFlow backend), 
% and is faster than \DDET\ for lower-level form, \JAX\ is much slower than \DDET\ for higher-level form.  
it is slower than \DDET.
Moreover, that \DDET\ itself being more efficient in higher-level form than in lower-level form makes working directly with higher-level form even more attractive in addition to the advantages mentioned further above.

Since coming up with an appropriate {learning-rate schedule for a good solution is more art than science}, we {systematically codified in detail our experience running optimization 
% (network training) 
through a normalization/standardization of the network-training process so readers can reproduce our results}.
%
% 24.7.16, code snippets
As the training speed is most likely related to the gradient computation, we provide snippets of our \DDET\ script and our \JAX\ script with detailed annotation.\footnote{
	\label{fn:future-update}
	The present update of Part I (full report)---which was shortened for publication in \cite{vuquoc2024PDAE-I}---includes our code snippets.
}

%\noindent
%\red{I AM HERE 23.10.10.}

		\vspace{0.25cm}
		\noindent
		{\bf KEYWORDS}
		\newline
		PINN, Physics-Informed Neural Network; Beam, small deformation, Euler-Bernoulli, large deformation, Kirchhoff rod, nonlinear boundary conditions, statics, dynamics.  
		
		\vspace{0.25cm}
		\noindent
		{\bf Subjects}, \href{https://arxiv.org/abs/2408.01914}{axXiv:2408.01914}: 
		Numerical Analysis (math.NA); Artificial Intelligence (cs.AI)
		\newline
		{\bf Mathematical Subject Classification (MSC classes):} \href{https://zbmath.org/static/msc2020.pdf}{MSC2020}
		\newline
		74-10 {\scriptsize (primary, Mathematical modeling or simulation for problems pertaining to mechanics of deformable solids)}; 
		\newline
		74H15 {\scriptsize (primary, Numerical approximation of solutions of dynamical problems in solid mechanics)}; 
		\newline
		74K10 {\scriptsize (primary, Rods (beams, columns, shafts, arches, rings, etc.))}; 
		\newline
		74B20 {\scriptsize (secondary, Nonlinear elasticity)}
		\newline
		{\bf ACM Computing Classification (ACM classes):} \href{https://www.acm.org/publications/computing-classification-system/1998/i}{ACM1998}
		\newline
		\href{https://www.acm.org/publications/computing-classification-system/1998/i.2}{I.2 {\scriptsize (Artificial Intelligence)}}; 
		\href{https://www.acm.org/publications/computing-classification-system/1998/i.2.6}{I.2.6 {\scriptsize (Learning)}}
		\newline
		\href{https://www.acm.org/publications/computing-classification-system/1998/i.6}{I.6 {\scriptsize (Simulation and Modeling)}}; 
		\href{https://www.acm.org/publications/computing-classification-system/1998/i.6.5}{I.6.5 {\scriptsize (Model Development)}}
			
	}
	
\end{mdframed}

% \newpage
% table of contents
%\vspace{0.1cm}
\hrule

%\vspace{1.5\baselineskip}
\vspace{0.5\baselineskip}
\noindent
%{\bf Table of Contents}
{\bf TABLE OF CONTENTS}

%\vspace{0.5\baselineskip}
\vspace{0.3\baselineskip}

% print the TOC using the internal command \@starttoc, which has
% the special character @, and therefore the command \newcommand must
% be enclosed between two commands \makeatletter and \makeatother
\makeatletter
\newcommand*{\toccontents}{\@starttoc{toc}}
\makeatother

% reduce the \baselineskip by half for the TOC
\setlength{\baselineskip}{0.5\baselineskip}

% print the TOC
\toccontents

% after printing the TOC, double \baselineskip back to its previous value
\setlength{\baselineskip}{2.0\baselineskip}

%\vspace{0.5cm}
\vspace{0.7\baselineskip}
\hrule

% nomenclature
% \input{01-nomenclature}

% introduction
% reset footnote counter so the first footnote number is 1. The reason is because
% footnote numbers 1 and 2 were used in the authors's affiliations.
\setcounter{footnote}{0}
\section{Introduction}
\noindent
% DAEs of nonlinear structural dynamics as exemplified by the Kirchhoff-Love rod.
%The partial-differential-algebraic equations (PDAEs) of motion for geometrically-exact rod with no shear considered here start from 
As a prototype of partial-differential-algebraic equations (PDAEs) for 
%\red{novel} 
novel  
physics-informed neural network (PINN) formulations, the geome\-trically-exact rod with no shear is considered here, starting with
the assumption often attributed to Kirchhoff (or Kirchhoff-Clebsch, or Kirchhoff-Love, or Kirchhoff-Clebsch-Love; see the review in \cite{dill1992kirchhoff} \cite{neukirch2021comparison} and the references therein): 
Plane cross section initially perpendicular to the undeformed centroidal line remains perpendicular to the deformed centroidal line.\footnote{
	Even though this assumption is part of the ``Kirchhoff-Love hypotheses'' in contemporary literature, neither Kirchhoff, nor Love, made these hypotheses \cite{dill1992kirchhoff}, in which the author wrote ``Major and minor authors alike have found the arguments of {\scriptsize\sc KIRCHHOFF} not  persuasive.'' 
	The method of derivation in the 1992 review paper \cite{dill1992kirchhoff} based on the first Piola-Kirchhoff stress tensor closely resembles, however, the derivation in 
	% Simo (1982) PhD thesis \cite{simo1982consistent}, 
	Simo (1982) \cite{simo1982consistent},
	% Simo \& Kelly JAM paper, 
	% Simo, J.C., K.H. Hjelmstad & R.L. Taylor [1984]. "Numerical formulations for the elasto-viscoplastic response of beams accounting for the effect of shear,"
	Simo, Heljmstad \& Taylor (1984) \cite{simo1984numerical}, and
	% Simo (1985) ADD REFS. 
	% Simo, J.C. [19B5], "A finite strain beam formulation. Part I: The three dimensional dynamic problem," Comp. Meth. Appl. Mech. Engrg., 49, 55-70.
	Simo (1985) \cite{simo1985finite}.
	% ``Contemporary treatments of the bending and twisting of rods are frequently  based on some set of hypotheses about the motion and the state of stress such as the following ones. (i) Cross-sections remain plane, undistorted, and normal to the axis of the rod. (ii) The transverse stress is zero. (iii) The bending moments and the twisting moment are proportional to the components of curvature and twist of the axis. Such assumptions are called the KIRCHHOFF-LOVE hypotheses. One even sees statements in the current literature about the mutual contradictions within the KIRCHHOFF-LOVE hypotheses. In fact, neither author made such assumptions.''
}  

Our formulation's starting point is similar to that in \cite{neukirch2021comparison}, then takes on a different direction as we are aiming at developing a finite deformable rod, without introducing approximations, culminating in a system of nonlinear partial differential equations, complemented by algebraic expressions with derivatives (AEDs), i.e., a PDAE system of the form
\begin{align}
	&
	\dpdesp{i}{(k)} (\bkar{\du} (\X , \ti), t) = 0 
	\ , \quad
	\bkar{\du} (\X , \ti) := \left[\{\partial_\X^p u_j , \partial_\X^{p-1} u_j , \ldots, \partial_\X^1 u_j\} , 
	\{\partial_\ti^q u_j , \partial_\X^{q-1} u_j , \ldots, \partial_\ti^1 u_j\} , 
	u_j\right]
	\ ,
	%	\nonumber
	%	\\
	%	&
	%	\text{with }
	%	k = 1,\ldots, 4 \ , \quad
	%	i = 1,\ldots, n_i^{(k)} , \quad j = 1,\ldots, n_j^{(k)} , \quad 
	%	p = 1,\ldots, n_p^{(k)} , \quad q = 1,\ldots, n_q^{(k)} 
	%	\ .
	\tag{\ref{eq:generic-PDE-aux-conditions}}
\end{align}
where $\dpdesp{i}{(k)}$ is a nonlinear differential operator, operating on $\bkar{\du} (\X , \ti)$, which collects the unknown dependent variables $\left\{u_j\right\}$ and their partial derivatives,
with $\left\{k, i, j, p, q\right\}$ being indices to be defined later in the paper; see Eq.~\eqref{eq:generic-PDE-aux-conditions} below.

An example of a PDAE system is the motion of the Kirchhoff-rod Form-1 (Section~\ref{sc:Kirchhoff-rod-Form-1}) Eq.~\eqref{eq:eom-u-bar} (PDE), Eq.~\eqref{eq:eom-v-bar} (PDE), Eq.~\eqref{eq:exact-slope-bar} (AED), Eq.~\eqref{eq:inverse-one-plus-e} (AED); these equations form an implicit system of nonlinear PDEs in terms of the four dependent variables $\{\ub , \vb , \eslpb\}$, $\iopeb$.  
While it is possible to eliminate two dependent variables $\{\eslpb, \iopeb\}$ by substituting Eq.~\eqref{eq:exact-slope-bar} and Eq.~\eqref{eq:inverse-one-plus-e} into Eq.~\eqref{eq:eom-u-bar} and Eq.~\eqref{eq:eom-v-bar}, the resulting equations, which can be put in generic form as,
\begin{align}
	\ndt \by (\X, t) = 
	% \mathcal{F} 
	\mathbb{F}
	(\by (\X, t), t)
	\ ,
\end{align}
where the partial space derivatives are included in the nonlinear operator $\mathbb{F}$,
become much more complex.

% from 08-recurrent.tex in our deep-learning review paper
%Some review papers on PINN are 
%\cite{cuomo2022scientific} 
%\cite{karniadakis2021physics}, with the latter being more general than \cite{cai2021fluid}, which was restricted to fluid mechanics, and touching on many different fields.
%Table~\ref{tb:PINN-frameworks} lists PINN frameworks that are currently actively developed, and a few selected \emph{solvers} among which are summarized below.

After deep learning achieved a sharp decrease in image classification error in 2012, and then by 2015 surpassed human performance \cite{vuquoc2023deep}, research in AI was revived after a long ``winter'' and took on a burgeoning turn, especially in deep learning.
Following this trend, as a method to solve PDEs without a need for discretization as in traditional methods (such as finite-difference, finite-element, boundary-element methods), physics-informed neural network (PINN) was reintroduced \cite{karniadakis2021physics} \cite{cuomo2022scientific} after a long dormant period since first proposed in the late 1990s \cite{lagaris1998artificial} \cite{lagaris2000neural}.
%
% from 08-recurrent.tex in our deep-learning review paper
% In laying out the roadmap for ``Simulation Intelligence'' (SI) the authors of \cite{lavin2021simulation} considered PINN as a key player in the first of the nine SI ``motifs,'' called ``Multi-physics \& multi-scale modeling.''
%
PINN is given a key role in ``Multi-physics \& multi-scale modeling,'' the first of the nine motifs in the ``Simulation Intelligence'' roadmap proposed in \cite{lavin2021simulation}.
For nonlinear problems, it is recently reported that PINN is 37 times faster than the traditional finite-difference method 
\cite{bazmara2023physics}.
Since PINN converts all problems, even linear ones, to a nonlinear optimization problem, the more a problem is 
%\red{nonlinear}, 
nonlinear,
the better for the application of PINN.

The present work is a natural extension of our review paper in \cite{vuquoc2023deep} aiming at both experts and first-time learners of both deep learning and PINN frameworks, among which the open-source DeepXDE (DDE) \cite{lu2021deepxde} is likely the most well documented framework with many examples and a forum for discussion among users.
So naturally, we chose to implement our 
%\red{novel} 
novel
PINN formulations in DDE with TensorFlow backend, abbreviated as \DDET.
Yet, we encountered some pathological problems (time shift, amplification, static solutions, as exemplified through the motion of an elastic bar) and propose 
%\red{novel} 
novel
methods to resolve them.   

Among these 
%\red{novel} 
novel 
PINN formulations are the 
%\red{different} 
different
PDE forms, obtained by splitting the derivative operators, evolving
from the lower-level form with fewer unknown dependent variables (e.g., displacements, slope, finite extension; see the Kirchhoff-rod Form 1 in Section~\ref{sc:Kirchhoff-rod-Form-1}) to higher-level form with more dependent variables (e.g., forces, moments, momenta, etc.), in addition to those from lower-level forms (see Kirchhoff-rod Form 4 in Section~\ref{sc:Kirchhoff-rod-Form-4}).   

%\red{Traditionally}, 
Traditionally,
the highest-level form, the balance-of-momenta form, is the starting point for (hand) deriving the lowest-level form through a tedious (and error prone) process of successive substitutions.  The next step in a finite element method is to discretize the lowest-level form, after forming the weak form and constructing a linearization, with appropriate interpolation functions, followed by their implementation and testing in a FE code \cite{zienkiewicz2013finite} \cite{zienkiewicz2014solid} \cite{zienkiewicz2013fluid}.  There are a number of finite-element frameworks to choose from, e.g., the Berkeley code \href{http://feap.berkeley.edu/wiki/index.php?title=FEAP_Wiki_Main_Page}{FEAP} by Taylor \cite{taylor2020feap} 
% (see \href{http://feap.berkeley.edu/wiki/index.php?title=FEAP_Wiki_Main_Page}{FEAP}) 
and \href{https://www.ngsolve.org}{Netgen/NGSolve} \cite{ngsolve.2014}. 
%(see \href{https://www.ngsolve.org}{Netgen/NGSolve}).

The time-consuming tedium of all of these 
%\red{traditional} 
traditional
steps could be bypassed by our 
%\red{novel application of} 
novel application of
PINN directly to the highest-level momentum form.\footnote{
	\label{fn:barrier functions}
	To filter the static solutions with near-zero acceleration encountered when using \DDET, we introduce various forms of barrier functions, 
	%\red{a novelty in PINN formulations}, 
	a novelty in PINN formulations,
	with room for improvements.  
%	A full-length report, with much more details than in the shortened version by half that is the present paper, can be found on arXiv.org \cite{vuquoc2024pinn.beam.2d.part1}. 
	(Appendix~\ref{app:barrier-functions}).
}

Perplexed with the number of pathological problems encountered with linear dynamic problems when using \DDET, as mentioned above, in parallel with developing 
%\red{novel} 
novel
methods to circumvent the problems in \DDET,
we also developed a script based on \JAX, the High Performance Array Computing. 

For the axial motion of elastic bar, while our \JAX\ script did not show the pathological problems of \DDET\ (DDE with TensorFlow backend), 
%and is faster than \DDET\ for lower-level Form 1, 
% by 10\%, 
% \JAX\ is much slower than \DDET\ for higher-level Form 3 (or 4).
% by 32\%.
it is slower than \DDET\ in all PDE forms tested.  
Moreover, that \DDET\ itself being more efficient in higher-level form than in lower-level form makes working directly with higher-level form even more attractive in addition to the advantages mentioned further above.

Arriving at an appropriate learning-rate schedule for a category of problems is more art than science, and has to be done by trial-and-error.  Because of the large number of parameters in the training process, we systematically codify our experience running optimization (training process) through a normalization/standardization of the network-training process and the result presentation so readers could understand and reproduce our results.  
%
% 24.7.16, code snippets
The training speed, or computational efficiency, is most likely related to how the gradient is computed in different software (e.g., \DDET\ vs \JAX), we provide snippets of our \DDET\ script and \JAX\ script with detailed explanatory annotation.\footnote{
	See Footnote~\ref{fn:future-update}.
}
%The training speed, or computational efficiency, is most likely related to how the gradient is computed in different software (e.g., \DDET\ vs \JAX).

%\red{Large computational efficiency is recently obtained with the proposed novel PINN formulations for the transverse motion---by more than 400\%, four hundred percent---and the nonlinear Kirchhoff rod---by more than 1000\%, one thousand percent---to be reported in a follow-up publication.}
%\red{Large computational efficiency is recently obtained with the proposed novel PINN formulations for the transverse motion of the Euler-Bernoulli beam (by more than 400\%, four hundred percent) and for the geometrically-exact Kirchhoff rod (by more than 1000\%, one thousand percent) to be reported in a follow-up publication.}
Large computational efficiency is recently obtained with the proposed novel PINN formulations for the transverse motion of the Euler-Bernoulli beam and for the geometrically-exact Kirchhoff rod to be reported in a follow-up publication.
%Forward and inverse problem formulation and results using the fully nonlinear Kirchhoff rod will be reported in a follow-up publication.

% beam models
\section{PDAEs: Beam formulations with no shear}
\noindent
%\red{PDEAs:}
% We begin by deriving the geometrically-exact beam with no shear, i.e., the Kirchhoff rod, followed by pointing out the inconsistency of such formulation based on geometrically-exact beam with shear.  Then we introduce approximations to recover the classical linear Euler-Bernoulli beam.
{As a particular case of the general partial-differential-algebraic equations (PDAEs) represented by Eq.~\eqref{eq:generic-PDE-aux-conditions}}, consider the geometrically-exact beam with no shear, i.e., the Kirchhoff rod, a derivation for which is first provided followed by pointing out the inconsistency of such formulation based on geometrically-exact beam with shear.
%A derivation of the geometrically-exact beam with no shear, i.e., the Kirchhoff rod, is first provided, followed by pointing out the inconsistency of such formulation based on geometrically-exact beam with shear.  
Then approximations are introduced to recover the classical linear Euler-Bernoulli beam with axial deformation in preparation for an implementation in a PINN framework, beginning with \DDET, then with \JAX.

% p.D237
\begin{figure}[tph]
	\centering
	\includegraphics[width=0.80\textwidth]{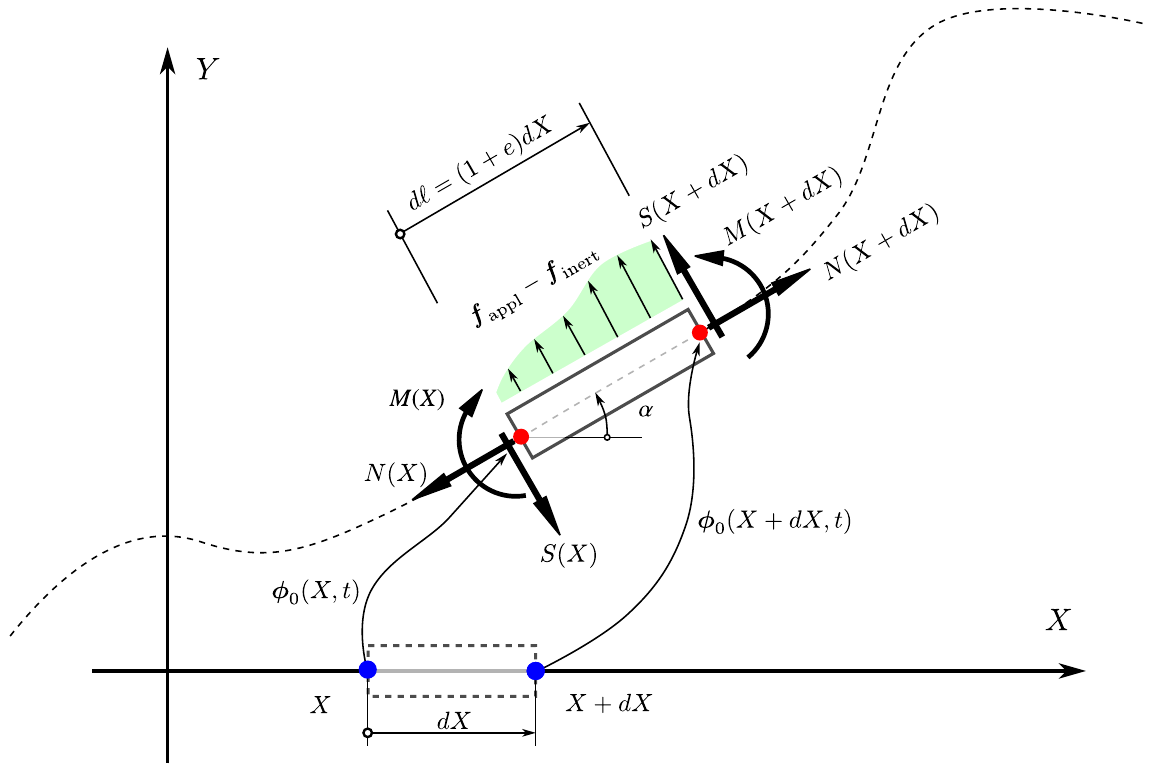}
	\caption{
		{\it Geometrically-exact beam without shear deformation.} 
		Section~\ref{sc:Kirchhoff-rod-derivation}.
		Shear forces introduced for equilibrium. Section~\ref{sc:inconsistency-Kirchhoff-rod-Euler-Bernoulli-beam}: Inconsistency in Kirchhoff rod (large deformation) and Euler-Bernoulli (small deformation) beam theories.  
		Figure~\ref{fig:geom-exact-with-shear}, geometrically-exact beam with shear deformation.
	}
	\label{fig:geom-exact-no-shear}
\end{figure}

%\subsection{Large deformation, Kirchhoff rod}
\subsection{Geometrically-exact beam with no shear, Kirchhoff rod}
\label{sc:Kirchhoff-rod-derivation}
\noindent
All equations are first derived in dimensional form, then subsequently converted to non-dimensional form.

\subsubsection{Kinematics, finite extension, constitutive relations}
%\noindent
%Kinematics
\noindent
Figure~\ref{fig:geom-exact-no-shear} provides a pictorial definition of the finite extension $\ext$ expressed as
\begin{align}
	% p.D224
	\ext = \frac{d \ell}{dL} - 1 
	= \frac{d \ell}{d \X} - 1 , \text{ with } d L = d \X 
	\ , \text{ and }
	d \ell = (1 + \ext) d\X
	\ ,
	\label{eq:finite-extension-1}
\end{align}
% \noindent
with the following kinematics of deformation $(x, y) = \phi_0 (\X, t)$ such that
\begin{align}
	% p.D224, D229, D230
	\x = \X + \du(\X, \ti) \Rightarrow
	d \x = d \X + d \du 
	= d \ell 
	% \cos \slps{\dv} 
	\cos \eslp
	\ , \quad
	\y = \dv(\X, \ti) \Rightarrow
	d \y = d \dv 
	= d \ell 
	% \sin \slps{\dv}
	\sin \eslp 
	\ ,
	\label{eq:kinematics}
\end{align}
where $\{\du (\X, t), \dv(\X, t) \}$ are the axial and transversal displacements, respectively.
% \noindent
The gradients of the spatial coordinates $(x, y)$ with respect to the material coordinate $\X$ are then:
\begin{align}
	\frac{d \x}{d \X} = 1 + \up{1}  
	= \frac{d \ell 
		% \cos \slps{\dv}
		\cos \eslp
	}{d L} 
	= (1 + \ext) 
	% \cos \slps{\dv}
	\cos \eslp 
	\ , \quad
	\frac{d \y}{d \X} = \vp{1} 
	= \frac{d \ell 
		% \cos \slps{\dv}
		\cos \eslp
	}{d L} 
	= (1 + \ext) 
	% \sin \slps{\dv}
	\sin \eslp 
	\ ,
	\label{eq:grad-spatial-coord}
\end{align}
% \noindent
which give rise to the expression for the
slope $\eslp$ of the centroidal line, together with the relations involving finite extension $\ext$ and the displacements $(\du , \dv)$:
\begin{align}
	&
	\tan \eslp = \frac{\vp{1}}{1 + \up{1}}
	\ , \quad 
	% p.D231
	\ext 
	\cos \eslp 
	= 1 + \up{1} - \cos \eslp
	\ , \quad
	\ext 
	\sin \eslp 
	= \vp{1} -	\sin \eslp
	\ .
	%\\
	\label{eq:slope-exact}
\end{align}
The last two expressions in Eq.~\eqref{eq:slope-exact}$_{2,3}$ provide then the extension $\ext$ in terms of the displacements $(\du, \dv)$ and the slope $\eslp$:
\begin{align}
	&
	e = \left[ \left( 1 + \up{1} - \cos \eslp \right)^2 + \left( \vp{1} - \sin \eslp \right)^2 \right]^{1/2}
	\ .
	\label{eq:finite-extension-2}
\end{align}
% \noindent
Using Eq.~\eqref{eq:grad-spatial-coord}, some convenient variables involving the extension $\ext$ and the displacement derivatives $(\up{1}$, $\vp{1})$ are introduced:
\begin{align}
	\ope^2 := (1+e)^2 = \left[ 1 + \up{1} \right]^2 + \left[ \vp{1} \right]^2
	\Rightarrow
	\ope =
	\left\{ 
		\left[ 1 + \up{1} \right]^2 + \left[ \vp{1} \right]^2 
	\right\}^{1/2}
	=: \iope^{-1}
	\ .
	\label{eq:define-ope}
\end{align}
%\noindent
The traditional elastic axial and bending constitutive relations are:
\begin{align}
	\N = EA \, e
	\ , 
	\quad
	\M = EI \eslpp{1}
	\ .
	\label{eq:elasticity-axial-bending}
\end{align}

\subsubsection{Balance of momenta, non-dimensional form}
\noindent
Referring to Figure~\ref{fig:geom-exact-no-shear}, the equilibrium of forces and moment of the depicted infinitesimal segment leads to the following expressions.

\noindent
Balance of linear momentum along $\X$ axis:
% p.236
\begin{align}
	\sum F_\X = 0 = 
	&
	\left[ \left. (\N \cos \eslp) \right|_{\X + d \X}
	- \left. (\N \cos \eslp) \right|_{\X} \right]
	- \left[  \left. (\shear \sin \eslp) \right|_{\X + d \X}
	- \left. (\shear \sin \eslp) \right|_{\X}  \right]
	\nonumber
	\\
	&
	+ \dfsp{\X}{\text{appl}} d\X - \dfsp{\X}{\text{inert}} d\X
	= d (\N \cos \eslp) - d (\shear \sin \eslp) + \dfs{\X} d\X - (\rho A d\X) \nddt \du 
	\ .
	\label{eq:bal-linear-mom-X}
\end{align}
% \noindent
Balance of linear momentum along $\Y$ axis:
% p.237
\begin{align}
	\sum F_\Y = 0 = 
	&
	\left[ \left. (\N \sin \eslp) \right|_{\X + d \X}
	- \left. (\N \sin \eslp) \right|_{\X} \right]
	- \left[  \left. (\shear \cos \eslp) \right|_{\X + d \X}
	- \left. (\shear \cos \eslp) \right|_{\X}  \right]
	\nonumber
	\\
	&
	+ \dfsp{\Y}{\text{appl}} d\X - \dfsp{\Y}{\text{inert}} d\X
	= d (\N \sin \eslp) - d (\shear \cos \eslp) + \dfs{\Y} d\X - (\rho A d\X) \nddt \dv 
	\ .
	\label{eq:bal-linear-mom-Y}
\end{align}
% \noindent
Balance of angular momentum:
\begin{align}
	\sum \left. \M \right|_{a} = 0 =
	- \M (\X) + \M (\X + d\X) 
	+ \shear (\X + d\X)  d \ell
	+ \Ms{\text{appl}} 
	- \Ms{\text{inert}}
	\ ,
\end{align}
\begin{align}
	\text{with }
	\Ms{\text{appl}} = O(d\ell^2)
	\ , \quad
	\Ms{\text{inert}} = \rho \, I_A \, d\ell \, \nddt{\slps{\dv}}
	\ .
\end{align}
Neglecting higher-order terms in $d\ell$ and the area moment of intertia $I_A$ of the cross section $A$ as customarily done for thin beams without shear deformation, and taking the limit $d\X \rightarrow 0$, we obtain
\begin{align}
	0 = \frac{d \M}{d\X} \frac{d \X}{d\ell} + \shear
	\Rightarrow
	0 = \frac{d \M}{d\X} + (1 + \ext) \shear = \Mp{1} + \ope \shear
	\ .
	\label{eq:moment-shear-LARGE-deformation}
\end{align}

\noindent
{\bf Non-dimensionalization}
\begin{align}
	% p.D225
	&
	\text{Coordinates and displacements, using Eq.~\eqref{eq:kinematics}:}
	\nonumber
	\\
	&
	\Xb = \frac{\X}{L}
	\ , \quad
	\xb = \frac{\x}{L} = \Xb + \ub
	\ , \quad
	\ub = \frac{\du}{L}
	\ , \quad
	\yb = \frac{\y}{L} = \frac{\dv}{L} = \vb \ ,
	\label{eq:non-dim-coord-disp}
	\\
	&
	\left[ \Xb \right] = \left[ \xb \right] = \left[ \ub \right] = \left[ \yb \right] = \left[ \vb \right] = 1
	\ .
	% \\
\end{align}
\begin{align}
	&
	\text{Normal force and shear force:}
	\nonumber
	\\
	&
	% p.D225
	\Nb = \frac{L^2}{EI} \N , \text{ with }
	\left[ \frac{L^2}{EI} \right] = \frac{\mathscr L^2}{(\mathscr{F/L^2}) \mathscr{L}^4} = \frac{1}{\mathscr{F}} , \text{ and } \left[ N \right] = \mathscr{F} , \text{ so } \left[ \Nb \right] = 1
	\label{eq:non-dim-normal-force}
	\ ,
	\\
	&
	\shrb = \frac{L^2}{EI} \shear, \text{ and } \left[ \shrb \right] = 1
	\ .
	\label{eq:non-dim-shear}
	% \\
\end{align}
In Eq.~\eqref{eq:non-dim-normal-force}, $\left[ N \right] = \mathscr{F}$ means the axial force $N$ has the dimension of force, denoted by $\mathscr{F}$, and similarly for length $\mathscr{L}$, mass $\mathscr{M}$, and time $\mathscr{T}$ as used in Eq.~\eqref{eq:non-dim-distributed-force-1} below.
\begin{align}
	&
	\text{Distributed force using Eq.~\eqref{eq:non-dim-normal-force}:}
	\nonumber
	\\
	&
	\dfb = \frac{L^3}{EI} \df , \text{ with } \left[ \dfb \right] = \left[ \frac{L^3}{EI} \right] \left[ \df \right] = \frac{\mathscr{L}}{\mathscr{F}} \frac{\mathscr{F}}{\mathscr{L}} = 1
	% \\
\end{align}
\begin{align}
	&
	\text{Moment:}
	\nonumber
	\\
	&
	% p.D229
	\Mb 
	= \frac{L}{EI} \M , \text{ since } \left[ \Mb \right] 
	= \left[ \frac{L}{EI} \right] \left[ \M \right] = \frac{\mathscr{L}}{\mathscr{F}} \frac{\mathscr{F}}{\mathscr{L}} = 1
	\ .
	\label{eq:non-dim-moment}
	% \\
\end{align}
\begin{align}
	&
	\text{Time:}
	\nonumber
	\\
	&
	% p.D226
	\ntip := L^2 \sqrt{(\rho A) / (EI)} ; \text{ since }
	% p.D226 
	\left[ \frac{\rho A}{EI} \right] = \frac{\mathscr{M} / \mathscr{L}}{\mathscr{F} \mathscr{L}^2}
	= \frac{\mathscr{M} / \mathscr{L}}{(\mathscr{M} \mathscr{L} / \mathscr{T}^2) \mathscr{L}^2}
	= \frac{\mathscr{T}^2}{\mathscr{L}^4} , \text{ thus }
	\left[ \ntip \right] = \mathscr{T}
	\label{eq:non-dim-distributed-force-1}
	\ ,
	\\
	&
	% p.D226
	\tb = \frac{t}{\ntip}  ,  \text{ and } \left[ \tb \right] = 1
	\ .
	\label{eq:non-dim-distributed-force}
\end{align}

\noindent
{\bf Balance of momenta}
\newline
\noindent
Space derivative of a generic force $\F (\X)$:
\begin{align}
	\Fp{1} = \frac{d \F(\X)}{d \X} 
	=
	\frac{d\Xb}{d\X}  \frac{d }{d\Xb} \left( \frac{EI}{L^2} \Fb(\Xb) \right)
	= \frac{1}{L} \frac{EI}{L^2} \frac{d \Fb(\Xb)}{d\Xb}
	= \frac{EI}{L^3} \Fbp{1}
	\ .
	\label{eq:space-deriv-generic-force}
\end{align}

\noindent
Time derivatives of a generic displacement function $x(\ti)$ with length dimension $\mathscr{L}$:
\begin{align}
	&
	\ndt x = \frac{d x}{d t} = \frac{d\tb}{dt} \frac{d [L \xb (\tb)]}{d\tb} = \frac{L}{\ntip} \frac{d\xb (\tb)}{d\tb} = \frac{L}{\ntip} \ndt \xb
	\ ,
	\\
	&
	\nddt \x 
	= \frac{L}{\ntip} \frac{d\tb}{dt} \frac{d^2 \xb (\tb)}{(d\tb)^2} 
	= \frac{L}{\ntip^2} \nddt \xb = \frac{EI}{\rho A L^3} \nddt \xb	
	\ .
	\label{eq:time-deriv-generic-displacement}
\end{align}

\noindent
\emph{Generic balance of linear momentum} in non-dimensional form using Eq.~\eqref{eq:non-dim-distributed-force}, Eq.~\eqref{eq:space-deriv-generic-force} and Eq.~\eqref{eq:time-deriv-generic-displacement}:
\begin{align}
	\Fp{1} + \df = \rho A \, \nddt \x
	\Rightarrow
	\Fbp{1} + \dfb = \nddt \xb 
	\ .
	\label{eq:non-dim-generic-bal-linear-mom}
\end{align}

\noindent
Space derivative of displacements and extension using Eq.~\eqref{eq:non-dim-coord-disp}:
\begin{align}
	\frac{\partial \du(\X, \ti)}{\partial \X} =:  	
	\up{1} = \frac{\partial \Xb}{\partial \X} \frac{\partial \left[L \ub(\Xb, \tb) \right]}{\partial \Xb}
	= \frac{\partial \ub(\Xb, \tb)}{\partial \Xb} =: \ubp{1}
	, \text{ and }
	\vp{1} = \vbp{1}
	\ .
	\label{eq:space-deriv-disp}
\end{align}
Similarly, using Eq.~\eqref{eq:space-deriv-disp} in Eqs.~\eqref{eq:slope-exact}-\eqref{eq:finite-extension-2} and Eq.~\eqref{eq:define-ope}, we have:
\begin{align}
	&
	\eslp = \eslpb , \quad
	\ext = \extb , \quad
	\ext 
	\cos \eslp 
	= \extb \cos \eslpb 
	, \quad
	\ext 
	\sin \eslp  
	= \extb \sin \eslpb
	, \quad
	%	\\
	%	&
	\ope = \opeb , \quad
	\iope = \iopeb
	\ .
	\label{eq:extension-bar}
\end{align}
% \noindent
Using Eq.~\eqref{eq:non-dim-generic-bal-linear-mom} in Eq.~\eqref{eq:bal-linear-mom-X} and Eq.~\eqref{eq:bal-linear-mom-Y}, together with Eq.~\eqref{eq:extension-bar}$_1$, we have the balance of linear momenta along $\Xb$ and $\Yb$ as:
\begin{align}
	\left( \Nb \cos \eslpb \right)^{(1)} 
	- \left( \shrb \sin \eslpb \right)^{(1)}
	+ \dfbs{\Xb} = \nddt \ub , \quad
	\left( \Nb \sin \eslpb \right)^{(1)} 
	+ \left( \shrb \cos \eslpb \right)^{(1)}
	+ \dfbs{\Yb} = \nddt \vb \ .
	\label{eq:non-dim-bal-linear-mom}
\end{align}
% \noindent
Space derivative of moment from Eq.~\eqref{eq:non-dim-moment} and similar to Eq.~\eqref{eq:space-deriv-generic-force}, we have:
\begin{align}
	\Mp{1} = \frac{EI}{L^2} \Mbp{1}
	\ .
	\label{eq:space-deriv-moment}
\end{align}

\noindent
\emph{Balance of angular momentum}. Eq.~\eqref{eq:moment-shear-LARGE-deformation} in non-dimensional form based on Eq.~\eqref{eq:space-deriv-moment}, Eq.~\eqref{eq:non-dim-shear}, and Eq.~\eqref{eq:extension-bar}$_5$:
\begin{align}
	 \Mp{1} + \ope \shear = 0 \Rightarrow \Mbp{1} + \opeb \shrb = 0 
%	 , \text{ with }
%	\ope = 1 + \ext = 1 + \extb = \opeb
	\ ,
	\label{eq:non-dim-bal-angular-mom}
	%\label{eq:non-dim-moment-shear-LARGE-deformation}
\end{align}
which $\opeb$ can be expressed in terms of the space derivative of the displacement components $(\ub, \vb)$ as defined in Eq.~\eqref{eq:define-ope}, using Eq.~\eqref{eq:space-deriv-disp}.

\subsubsection{Constitutive relations, non-dimensional form}
%\noindent
%{\bf Constitutive relations}
\noindent
Non-dimensional constitutive relation for axial deformation:
\begin{align}
	\N = EA \ext
	\Rightarrow
	\frac{EI}{L^2} \Nb = EA \extb
	\Rightarrow
	\Nb = \frac{A L^2}{I} \extb = \sldn \extb
	\ ,
	\label{eq:elasticity-axial}
\end{align}
with
non-dimensional slenderness parameter $\sldn$ and its inverse being the rotundness $\rtdn$
\begin{align}
	\sldn := \frac{A L^2}{I} 
	\Leftrightarrow
	\rtdn := \sldn^{-1} = \frac{I}{A L^2}
	\ .
	\label{eq:slenderness-rotundness}
\end{align}
% \noindent
The space derivative of the slope $\eslp$ of the centroidal line is:
\begin{align}
	\eslpp{1} = \frac{\partial \eslp}{\partial \X}
	= \frac{\partial \Xb}{\partial \X} \frac{\partial \eslpb}{\partial \Xb}
	= \frac{1}{L} \eslpbp{1}
	\ .
	\label{eq:space-deriv-slope-centroid-line}
\end{align}
% \noindent
Non-dimensional constitutive relations for bending deformation using Eq.~\eqref{eq:elasticity-axial-bending}, Eq.~\eqref{eq:non-dim-moment}, and Eq.~\eqref{eq:space-deriv-slope-centroid-line}:
\begin{align}
	\M = EI \eslpp{1}
	\Rightarrow
	\Mb = \eslpbp{1}
	\ .
	\label{eq:elasticity-bending-exact}
\end{align}

\subsubsection{Equations of motion, non-dimensional form}
%\noindent
%{\bf Equations of motion}
\noindent
The first two terms in each of the non-dimensional balance of linear momentum in $\Xb$ and in $\Yb$ direction, respectively, as expressed in Eq.~\eqref{eq:non-dim-bal-linear-mom} can be written in terms of displacements as follows.  Using the axial constitutive relation Eq.~\eqref{eq:elasticity-axial} and the definition of $\eslp$ in Eq.~\eqref{eq:slope-exact} together with the expressions related to the finite extension $\ext$ defined in Eqs.~\eqref{eq:slope-exact}-\eqref{eq:finite-extension-2}, and again using Eq.~\eqref{eq:extension-bar}, we have:
\begin{align}
	&
	\Nb \cos \eslpb = \sldn \extb \cos \eslpb = \sldn \left[ 1 + \ubp{1} - \cos \eslpb \right]
	, \quad
	\Nb \sin \eslpb = \sldn \extb \sin \eslpb = \sldn \left[ \vbp{1} - \sin \eslpb \right]
	\ ,
	\\
	&
	\tan \eslpb = \frac{\vbp{1}}{1 + \ubp{1}}
	\ ,
	\label{eq:exact-slope-bar}
\end{align}
from which 
%the first term as a function of $(\Nb , \eslpb)$ with a space derivative, in each equation in Eq.~\eqref{eq:non-dim-bal-linear-mom} is, respectively:
the first term in each of the two equations in Eq.~\eqref{eq:non-dim-bal-linear-mom} is obtained respectively as a function of $(\Nb , \eslpb)$ with a space derivative as follows:
\begin{align}
	\left( \Nb \cos \eslpb \right)^{(1)} = \sldn \left[ \ubp{2} + \eslpbp{1} \sin \eslpb \right]
	, \quad
	\left( \Nb \sin \eslpb \right)^{(1)} = \sldn \left[ \vbp{2} - \eslpbp{1} \cos \eslpb \right]
	\ .
	\label{eq:deriv-terms-N-bar}
\end{align}

% \noindent
Next, using the balance of angular momentum, yielding the moment-shear relation with finite-extension effect, in Eq.~\eqref{eq:non-dim-bal-angular-mom}, the bending constitutive relation in Eq.~\eqref{eq:elasticity-bending-exact}, and the definition of $\eslp$ in Eq.~\eqref{eq:slope-exact}, and making use of Eq.~\eqref{eq:extension-bar}$_1$, we obtain:
\begin{align}
	&
	\shrb \sin \eslpb = - \iopeb \Mbp{1} \sin \eslpb = - \iopeb \eslpbp{2} \sin \eslpb
	, \quad
	\shrb \cos \eslpb = - \iopeb \eslpbp{2} \cos \eslpb
	\ ,
	\\
	&
	- \left( \shrb \sin \eslpb \right)^{(1)} 
	= - \left( - \iopeb \eslpbp{2} \sin \eslpb \right)^{(1)}
	= \left[ \iopeb \eslpbp{2} \right]^{(1)} \sin \eslpb 
	+ \left[\iopeb \eslpbp{2} \right] \eslpbp{1} \cos \eslpb
	\ ,
	\label{eq:deriv-terms-S-bar-A}
	\\
	&
	+ \left(\shrb \cos \eslpb \right)^{(1)}
	= \left( - \iopeb \eslpbp{2} \cos \eslpb \right)^{(1)}
	= - \left[ \iopeb \eslpbp{2} \right]^{(1)} \cos \eslpb 
	+ \left[\iopeb \eslpbp{2} \right] \eslpbp{1} \sin \eslpbp{1}
	\ .
	\label{eq:deriv-terms-S-bar-B}
\end{align}

\noindent
Using Eq.~\eqref{eq:deriv-terms-N-bar} for the derivative of the terms with $(\Nb , \eslpb)$, and Eqs.~\eqref{eq:deriv-terms-S-bar-A}-\eqref{eq:deriv-terms-S-bar-B} for the derivative of the terms with $(\shrb , \eslpb)$ in Eq.~\eqref{eq:non-dim-bal-linear-mom}, the non-dimensionalized balance of linear momenta, we obtain the equations of motion in terms of the displacements $(u,v)$.

\noindent
Axial equation of motion:
% Eq.(21), p.D242
\begin{align}
	\sldn 
	\left[ \ubp 2 + \eslpbp{1} \sin \eslpb \right]
	+ \left[ \iopeb \eslpbp{2} \right]^{(1)} \sin \eslpb 
	+ \left[\iopeb \eslpbp{2} \right] \eslpbp{1} \cos \eslpb + \dfbs{\X} 
	= \nddt \ub 
	\ ,
	\label{eq:eom-u-bar}
\end{align}

\noindent
Transversal equation of motion:
% Eq.(26), p.D243
\begin{align}
	\sldn 
	\left[ \vbp 2 - \eslpbp{1} \cos \eslpb \right]
	- \left[ \iopeb \eslpbp{2} \right]^{(1)} \cos \eslpb 
	+ \left[\iopeb \eslpbp{2} \right] \eslpbp{1} \sin \eslpb + \dfbs{\Y}
	= \nddt \vb 
	\ ,
	\label{eq:eom-v-bar}
\end{align}
where the slenderness $\sldn$ was defined in Eq.~\eqref{eq:slenderness-rotundness}$_1$, and $\iopeb$ in Eq.~\eqref{eq:extension-bar}, Eq.~\eqref{eq:space-deriv-disp}, and Eq.~\eqref{eq:define-ope}, reproduced here so all relevant terms in the above equations of motion are grouped together for convenience:
\begin{align}
	\sldn := \frac{A L^2}{I}
	\tag{\ref{eq:slenderness-rotundness}}
	\ ,
\end{align}
\begin{align}
	\iopeb =
	\left\{ 
	\left[ 1 + \ubp{1} \right]^2 + \left[ \vbp{1} \right]^2 
	\right\}^{-1/2}
	\ .
	\label{eq:inverse-one-plus-e}
\end{align}

\subsubsection{Input parameters: Boundary and initial conditions, distributed load}
\noindent
In addition to the slenderness $\sldn$ defined in Eq.~\eqref{eq:slenderness-rotundness} and the input distributed load functions $\dfbs{\X} (\Xb , \tb)$ and $\dfbs{\Y} (\Xb , \tb)$ for Eqs.~\eqref{eq:eom-u-bar}-\eqref{eq:eom-v-bar}, respectively, i.e.,
\begin{align}
	\dfbs{\X} (\Xb , \tb) = \dfbhs{\X} (\Xb , \tb)
	\ , \quad
	\dfbs{\Y} (\Xb , \tb) = \dfbhs{\Y} (\Xb , \tb)
	\ ,
	\label{eq:distributed-forces}
\end{align}
there are six input parameters in each of the two sets of boundary conditions considered here: One set for cantilever beam, and another for simply-supported beam.

\noindent
\emph{Cantilever beam}

\noindent
Clamped end: Three input parameters are the prescribed displacements and rotation $\left\{ \widehat{\ub}, \widehat{\vb},  \widehat{\eslpb} \right\}$, such that
\begin{align}
	\text{At } \Xb = 0, \quad 
	\ub (\Xb = 0, \tb) = \widehat{\ub}, \quad 
	\vb (0, \tb) = \widehat{\vb}, \quad
	\eslpb (0, \tb) = \widehat{\eslpb} \ .
	\label{eq:cantilever-clamped-end-bound-cond}
\end{align}

\noindent
Free end: Three input parameters, two algebraic expressions with derivatives (AEDs).
% https://www.britannica.com/science/algebraic-equation
% https://en.wikipedia.org/wiki/Differential-algebraic_system_of_equations
\begin{align}
	\text{At } \Xb = 1, \quad
	\sldn  \extb (1, \tb) = \widehat{\Nb} , \quad
	- \iopeb (1, \tb) \Mbp{1} (1, \tb) = - \iopeb (1, \tb) \eslpbp{2} (1, \tb) = \widehat{\shrb} , \quad
	\eslpbp{1} (1, \tb) = \widehat{\Mb}
	\ .
	\label{eq:cantilever-free-end-bound-cond}
\end{align}
The three input parameters in Eq.~\eqref{eq:cantilever-free-end-bound-cond} are the prescribed concentrated forces and moment $\{ \widehat{\Nb}$, $\widehat{\shrb}$, $\widehat{\Mb} \}$, and the two AEDs are Eqs.~\eqref{eq:cantilever-free-end-bound-cond}$_{1,2}$ due to the expression of $\extb$ as a result of 
%Eq.~\eqref{eq:finite-extension-2}, Eq.~\eqref{eq:space-deriv-disp}, Eq.~\eqref{eq:extension-bar}$_2$, 
Eqs.~\{\eqref{eq:slope-exact}, \eqref{eq:finite-extension-2}, \eqref{eq:space-deriv-disp}, \eqref{eq:extension-bar}$_2$\},
and the expression of $\iopeb$ in Eq.~\eqref{eq:inverse-one-plus-e}, 
with both $\extb$ and $\iopeb$ involving the first-order derivatives $\ubp{1}$ and $\vbp{1}$. 
These AEDs together with the differential equations from the balance of momenta form a system of differential-algebraic equations (DAEs).

\noindent
\emph{Simply-supported beam}
\begin{align}
	\text{At } \Xb = 0, \quad 
	\ub (0, \tb) = \widehat{\ub}_0 , \quad 
	\vb (0, \tb) = \widehat{\vb}_0 , \quad
	\eslpbp{1} (0, \tb) = \widehat{\Mb}_0 \ .
	\label{eq:pinned-left-end-bound-cond}
\end{align}
\begin{align}
	\text{At } \Xb = 1, \quad 
	\sldn \extb (1, \tb) = \widehat{\Nb}_1 , \quad
	\vb (1, \tb) = \widehat{\vb}_1 , \quad 
	\eslpbp{1} (0, \tb) = \widehat{\Mb}_1 \ .
	\label{eq:roller-right-end-bound-cond}
\end{align}
Because of the expression for $\extb$ in Eq.~\eqref{eq:roller-right-end-bound-cond}$_1$, the equations of motion of the simply-supported beam also form a system of DAEs.

\begin{figure}[tph]
	\centering
	%
	% 24.7.12, for arXiv replace blank spaces by underscores in filename
	%	\includegraphics[width=0.5\textwidth]{Figures/GE_beam_strains_-_ALL.pdf}
	\includegraphics[width=0.5\textwidth]{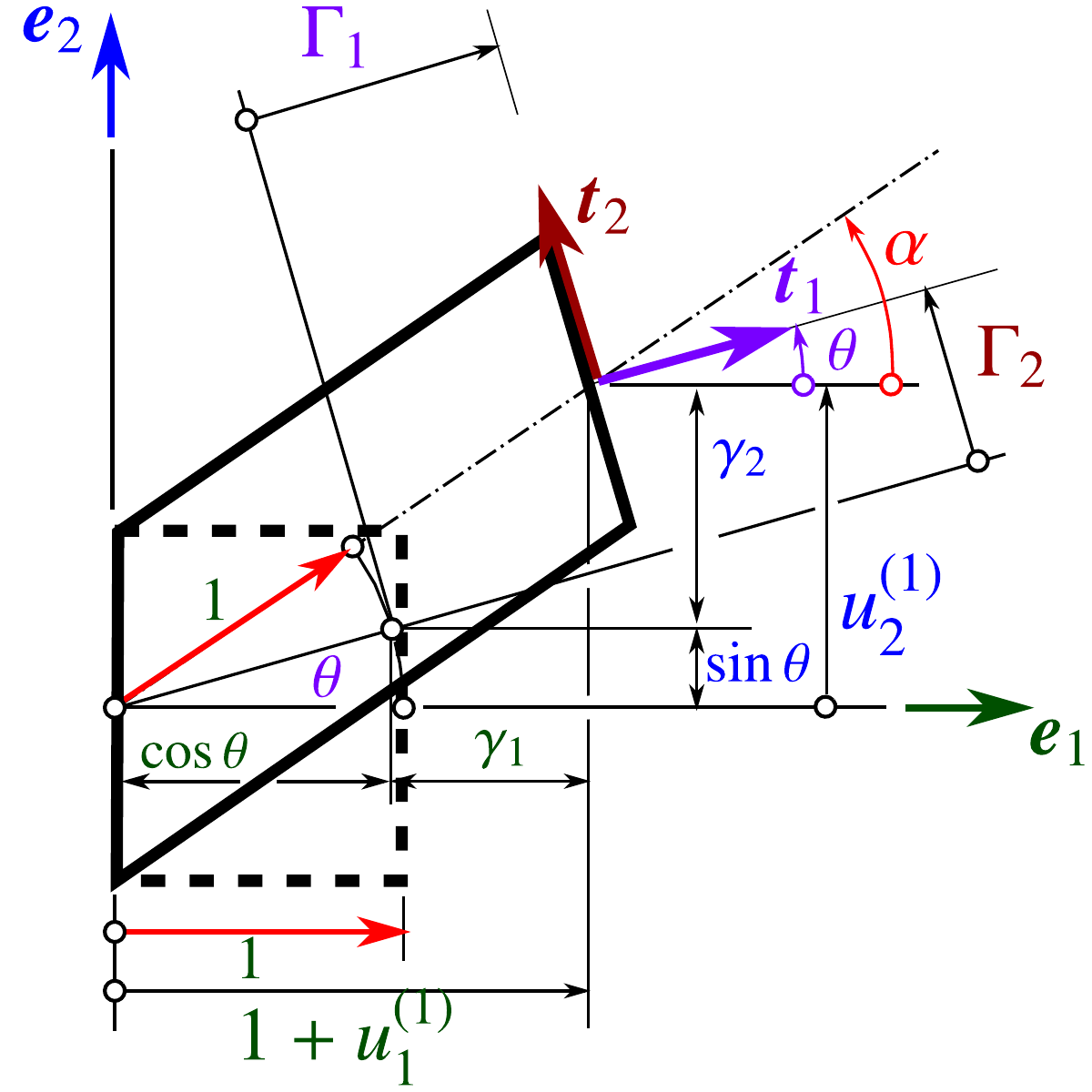}
	\caption{
		\emph{Geometrically-exact beam with shear deformation.}  
		The deformed configuration (solid line) is superposed on the initial configuration (dotted line) of unit length (which could be multiplied by $d\X$).  Shear deformation is the difference between the angle $\eslp$ of the deformed centroidal line and the rotation $\rc$ of the cross section. Spatial strains $\{\sstrs{1} , \sstrs{2}\}$ and material strains $\{\cstrs{1} , \cstrs{2}\}$, such that  
		$\bsstr = \sstrs{1} \bsbv{1} + \sstrs{2} \bsbv{2} =  \cstrs{1} \bmbv{1} + \cstrs{2} \bmbv{2}$.
		See \cite{simo1986dynamicsI} and Figure~\ref{fig:geom-exact-no-shear} for geometrically-exact beam with no shear.
		% Initial (reference) configuration: Dotted line.  Deformed configuration: Solid line.  Rotation of cross section: $\eslp$ 
	}
	\label{fig:geom-exact-with-shear}
\end{figure}

\begin{rem}
	{AEDs cannot be imposed exactly.}
	{\rm 
		Because ``complex'' boundary conditions are AEDs, such as Eq.~\eqref{eq:roller-right-end-bound-cond}$_1$, which involve derivatives of the dependent variables, they cannot be imposed exactly by hard constraints.  Only ``simple'' boundary conditions not involving derivatives of dependent variables, such as Eqs.\eqref{eq:roller-right-end-bound-cond}$_{2,3}$, can be imposed exactly.
	}
	$\hfill\blacksquare$	
\end{rem}

%\noindent
%{\color{red}  NOTE: 2023.05.26 - We should think about how to devise a network that can impose a hard constraint on AEDs.  To this end, we need to be able to evaluate the derivatives.  ENDNOTE}

\noindent
There are four prescribed functions for the four initial conditions:
\begin{align}
	\ub (\Xb , 0) = \widetilde{\ub}_0 (\Xb) , \quad
	\ubd (\Xb , 0) = \widetilde{\ubd}_0 (\Xb) \ ,
	\label{eq:initial-conditions-1}
	\\
	\vb (\Xb , 0) = \widetilde{\vb}_0 (\Xb) , \quad
	\vbd (\Xb , 0) = \widetilde{\vbd}_0 (\Xb) \ .
	\label{eq:initial-conditions-2}
\end{align}

\subsection{Inconsistency, reduction from geometrically-exact beam with shear}
\label{sc:inconsistency-Kirchhoff-rod-Euler-Bernoulli-beam}
\noindent
Deformation map of centroidal line \cite{simo1986dynamicsI}, which is valid for both geometrically-exact beam without shear (Figure~\ref{fig:geom-exact-no-shear}) and with shear (Figure~\ref{fig:geom-exact-with-shear}), is written as:
\begin{align}
	\bdmz (\X , \ti) = \left[\X + \du (\X, \ti)\right] \bsbv{1} 
	 + \dv (\X, \ti) \bsbv{2}
	 \ ,
	\label{eq:geom-exact-defm-map}
\end{align}
where $\{ \bsbv{1} , \bsbv{2} \}$ are the spatial basis vectors.
The material basis vectors $\{ \bmbv{1} , \bmbv{2} \}$ attached to the cross section are then in terms of the rotation $\rc$ of the cross section
\begin{align}
	\bcbv{1} (\X , \ti) = \cos \rc (\X , \ti) \cdot \bsbv{1} + \sin \rc (\X , \ti) \cdot \bsbv{2}
	\ , \quad
	\bcbv{2} (\X , \ti) = - \sin \rc (\X , \ti) \cdot \bsbv{1} + \cos \rc (\X , \ti) \cdot \bsbv{2}
	\ .
	\label{eq:material-basis-vectors}
\end{align}
\begin{align}
	\begin{Bmatrix}
		\bcbv{1}
		\\
		\bcbv{2}
	\end{Bmatrix}
	=
	\begin{bmatrix}
		\phantom{-} \cos \rc & \sin \rc
		\\
		- \sin \rc & \cos \rc
	\end{bmatrix}
	\begin{Bmatrix}
		\bsbv{1}
		\\
		\bsbv{2}
	\end{Bmatrix}
	=:
	\bRotM^T
	\begin{Bmatrix}
		\bsbv{1}
		\\
		\bsbv{2}
	\end{Bmatrix} 
	\ .
	\label{eq:matl-to-sptl-lambdaT}
\end{align}
Spatial strains $\{\sstrs{1}, \sstrs{2}\}$ and material strain $\{\cstrs{1} , \cstrs{2} \}$; see Figure~\ref{fig:geom-exact-with-shear}:
\begin{align}
	&
	\begin{Bmatrix}
		\sstrs{1}
		\\
		\sstrs{2}
	\end{Bmatrix}
	=
	\begin{Bmatrix}
		1 + \up{1} - \cos \rc
		\\
		\vp{1} - \sin \rc
	\end{Bmatrix}
%	=
%	\begin{Bmatrix}
%		(1+\ext) \cos \eslp - \cos \rc
%		\\
%		(1+\ext) \sin \eslp - \sin \rc
%	\end{Bmatrix}
	\ , 
	\quad
	\begin{Bmatrix}
		\cstrs{1}
		\\
		\cstrs{2}
	\end{Bmatrix}
	=
	\bRotM^T 
	\begin{Bmatrix}
		\sstrs{1}
		\\
		\sstrs{2}
	\end{Bmatrix}
	\ .
	\label{eq:spatial-material-strains}
\end{align}
Next, using Eq.~\eqref{eq:grad-spatial-coord} in Eq.~\eqref{eq:spatial-material-strains}$_1$, 
and then Eq.~\eqref{eq:matl-to-sptl-lambdaT} in Eq.~\eqref{eq:spatial-material-strains}$_2$,
we obtain
\begin{align}
	\begin{Bmatrix}
		\sstrs{1}
		\\
		\sstrs{2}
	\end{Bmatrix}
	=
	\begin{Bmatrix}
		(1+\ext) \cos \eslp - \cos \rc
		\\
		(1+\ext) \sin \eslp - \sin \rc
	\end{Bmatrix}
	\ , \quad
	\begin{Bmatrix}
		\cstrs{1}
		\\
		\cstrs{2}
	\end{Bmatrix}
	=
	\begin{Bmatrix}
		(1+\ext) \cos (\eslp - \rc) - 1
		\\
		(1+\ext) \sin (\eslp - \rc)
	\end{Bmatrix}
	\ ,
	\label{eq:spatial-material-strains-2}
\end{align}
and in the case where $\rc = \eslp$, we have
\begin{align}
	\begin{Bmatrix}
		\sstrs{1}
		\\
		\sstrs{2}
	\end{Bmatrix}
	=
	\ext
	\begin{Bmatrix}
		\cos \eslp
		\\
		\sin \eslp
	\end{Bmatrix}
	\ , \quad
	\begin{Bmatrix}
		\cstrs{1}
		\\
		\cstrs{2}
	\end{Bmatrix}
%	=
%	\bRotM^T 
%	\begin{Bmatrix}
%		\sstrs{1}
%		\\
%		\sstrs{2}
%	\end{Bmatrix}
	=
	\begin{Bmatrix}
		\ext
		\\
		0
	\end{Bmatrix}
	\ ,
	\label{eq:spatial-material-strains-3}
\end{align}
which indicates that there is zero shear deformation, and thus zero shear force $\shear$, since the material (axial and shear) forces $\{\Ns{1}, \Ns{2}\} = \{\N, \shear\}$ are computed from the linear constitutive relation \cite{simo1986dynamicsI} as follows:
\begin{align}
	\begin{Bmatrix}
		\Ns{1}
		\\
		\Ns{2}
	\end{Bmatrix}
	=
	\begin{bmatrix}
		EA & 0
		\\
		0 & GA_s
	\end{bmatrix}
	\begin{Bmatrix}
		\cstrs{1}
		\\
		\cstrs{2}
	\end{Bmatrix}
	=
	\begin{Bmatrix}
		EA \ext
		\\
		0
	\end{Bmatrix}
	=
	\begin{Bmatrix}
		\N
		\\
		\shear
	\end{Bmatrix}
	\ .
\end{align}
But a non-zero shear force $\shear$ was introduced in Figure~\ref{fig:geom-exact-no-shear} and in the balance of linear momenta to derive the equations of motion as done above. 

\begin{rem}
	{\rm
		The authors of \cite{simo1986dynamicsI} cited (among other references) \cite{reissner1972one}, where neither the dynamic example, nor the computational formulation, was considered. 
		The expressions for $\{ \cstrs{1}, \cstrs{2} \}^T$ in Eq.~\eqref{eq:spatial-material-strains-3}$_2$ are the same as those for $\{ \epsilon, \gamma \}^T$ in Eqs.~(13 a,b) in \cite{reissner1972one}, upon recognizing the following relations between the variables $(\eslp, \rc)$ in \cite{simo1986dynamicsI} and the variables $(\chi, \phi)$ in \cite{reissner1972one}: $\rc \equiv \phi$ and $(\eslp - \rc) \equiv \chi$.
		Specifically, even though the dynamic equilibrium equations (``Form 4'') were mentioned in \cite{reissner1972one}, the ``dynamic'' operators (accelerations) in \cite{reissner1972one} were implicitly hidden in the applied force and moment (d'Alembert principle), expressed either in the fixed inertial frame, or projected on the cross-section basis vectors, without the explicit equations of motion in terms of the displacements and rotation (``Form 1''), as given in \cite{simo1986dynamicsI}.  The simplicity of the acceleration components in a fixed inertial frame compared to those in a rotating (floating) frame was explicitly demonstrated in \cite{simo1986dynamicsI}. 
		Simulating flying flexible beams under large deformation and large overall motion was a breakthrough first reported in \cite{simo1986dynamicsI}. 
%		Hence while the ``dynamic'' ``Form 4'' in \cite{reissner1972one} was the result of a projection onto the rotating cross-section basis vectors $\{\bcbv{1} , \bcbv{2}\}$, and thus with hidden complex expressions for the acceleration components, the dynamic Form 4 (Section~\ref{sc:Kirchhoff-rod-Form-4}) considered here remains in the inertial frame, with simple expressions for the acceleration components.
	}
	\hfill$\blacksquare$
\end{rem}

%\noindent
%{\color{red} \ding{42} I AM HERE \today}

\subsection{Small deformation, Euler-Bernoulli beam}
\noindent
For small axial deformation $\ext \ll 1$, which includes inextensibility $\ext = 0$, the moment-shear relation in Eq.~\eqref{eq:moment-shear-LARGE-deformation} is approximated as
\begin{align}
	\ext \ll 1 
	\Rightarrow
	0 = \frac{d \M}{d\X} + \shear
	\ ,
	\label{eq:moment-shear-SMALL-deformation}
\end{align}
which is the classical relation between moment and shear.

\noindent
For small slope $\eslp$ and small axial strain $\up{1}$, the slope of the centroidal line can be approximated as the slope of the transverse displacement using 
%Eq.~\eqref{eq:finite-extension-2}, 
Eq.~\eqref{eq:slope-exact},
i.e., 
\begin{align}
	\eslp < 10^\circ \text{ (small angle) } \text{ and } \up{1} \ll 1
	\Rightarrow
	\eslp \approx \vp{1}
	= \slps{\dv}
	\Rightarrow
	\eslpb \approx \vbp{1} = \slpbs{\dv}
	\ .
	\label{eq:small-angle}
\end{align}

\noindent
From Eq.~\eqref{eq:space-deriv-disp}$_5$ on the space-derivative of the transverse displacement, the slope of the deformed centroidal line can be written as
\begin{align}
	\slps{v} (\X, \ti) := \frac{\partial \dv(\X, \ti)}{\partial \X} = \vp{1} (\X, \ti) 
	, \quad 
	\slpbs{v} (\Xb, \tb) := \frac{\partial \vb(\Xb, \tb)}{\partial \Xb} = \vbp{1} (\Xb, \tb) , \quad \slps{v} (\X, \ti) = \slpbs{v} (\Xb, \tb)
	\ ,
	\label{eq:non-dim-slope-v}
\end{align}
and thus the space derivative of the slope is:
\begin{align}
	\slpsp{v}{1} = \frac{\partial \vp{1}}{\partial \X}
	= \frac{\partial \Xb}{\partial \X} \frac{\partial \vbp{1}}{\partial \Xb}
	= \frac{1}{L} \slpbsp{\dv}{1}
	\ ,
	\label{eq:space-deriv-slope}
\end{align}
again using Eq.~\eqref{eq:space-deriv-disp}$_5$.

\noindent
The non-dimensional constitutive relations for bending deformation using Eq.~\eqref{eq:elasticity-axial-bending}, Eq.~\eqref{eq:non-dim-moment}, Eq.~\eqref{eq:elasticity-bending-exact}, Eq.~\eqref{eq:small-angle}, and Eq.~\eqref{eq:space-deriv-slope} are obtained as:
\begin{align}
	\M = EI \slpsp{\dv}{1} = \vp{2}
	\Rightarrow
	\Mb = \slpbsp{\dv}{1} = \vbp{2}
	\ ,
	\label{eq:elasticity-bending}
\end{align}
the last equation of which is the moment-curvature relation for the Euler-Bernoulli beam.  

\noindent
For the axial equation of motion in Eq.~\eqref{eq:eom-u-bar}, the small-angle approximation in Eq.~\eqref{eq:small-angle} together with the near-zero slope approximation:
\begin{align}
	&
	\eslp \approx 0 \Rightarrow \sin \eslp \approx 0 \text{ and } \cos \eslp \approx 1 \ ,
	\label{eq:zero-slope}
%	\\
%	&
%	\sldn \ubp{2} + \dfbs{\X} \approx \nddt \ub
%	\ .
%	\label{eq:Euler-Bernoulli-axial}
\end{align}
lead to the standard axial equation of motion with distributed load:
\begin{align}
%	&
%	\eslp \approx 0 \Rightarrow \sin \eslp \approx 0 \text{ and } \cos \eslp \approx 1 \ ,
%	\label{eq:zero-slope}
%	\\
%	&
	\sldn \ubp{2} + \dfbs{\X} \approx \nddt \ub
	\ .
	\label{eq:eom-euler-bernoulli-axial}
\end{align}
Without applied axial distributed load, the axial Eq.~\eqref{eq:eom-euler-bernoulli-axial} becomes the standard 1-D wave equation:
\begin{align}
	\sldn \ubp{2} = \nddt \ub 
	\Rightarrow
	\ubp{2} = \rtdn \nddt{\ub}
	\ ,
	\label{eq:wave-equation}
\end{align} 
with the {square-root of the slenderness, i.e., $\sqrt{\sldn}$, as the non-dimensionalized wave speed.  The more slender (i.e., the less rotund) the bar is, the faster the wave propagates along the axial direction.}

\noindent
For the transverse equation of motion in Eq.~\eqref{eq:eom-v-bar}, the approximations in Eq.~\eqref{eq:small-angle} (small angle) and Eq.~\eqref{eq:zero-slope} (near-zero slope) lead to: 
\begin{align}
	\sldn 
	\left[ \vbp 2 - \eslpbp{1} \cos \eslpb \right]
	\approx 0
	\ , \quad	
	\left[\iopeb \eslpbp{2} \right] \eslpbp{1} \sin \eslpb
	\approx 0
	\ , \quad
	- \left[ \iopeb \eslpbp{2} \right]^{(1)}  
	+ \dfbs{\Y}
	\approx \nddt{\vb}
	\ ,
\end{align}
which, upon using the approximation $\iopeb \approx 1$ together with the small-angle approximation in Eq.~\eqref{eq:small-angle}, leads to the standard Euler-Bernoulli beam equation of motion with fourth space derivative on the transverse displacement $\vb$:
\begin{align}
	- \vbp{4} + \dfbs{\Y} \approx \nddt{\vb}
	\ .
	\label{eq:eom-euler-bernoulli-transversal}
\end{align}
Eq.~\eqref{eq:eom-euler-bernoulli-axial} and Eq.~\eqref{eq:eom-euler-bernoulli-transversal} form the system of non-dimensionalized PDEs that describes the dynamics of small-deformation Euler-Bernoulli beam. 

%\noindent
%{\color{red} \ding{42} I AM HERE \today}
%SECTION NOT DONE

% PINN
\section{PINN formulations}
\noindent
Several forms of the equations of motion of the Kirchhoff rod---starting from the lower level (Form 1) to the higher level (Form 4)---are provided below with significant consequences on the corresponding PINN implementation and computational efficiency.  All these forms have nonlinear boundary conditions.
% \red{[NOTE: 2023.09.13, To write. ENDNOTE]}

\subsection{Inputs}
%\vspace{0.5cm}
\begin{mdframed}[
	% style=mystyle,
	% backgroundcolor=lightgray,
	backgroundcolor=mygray,
	% backgroundcolor=white,
	topline=true,
	bottomline=true,
	leftline=true,
	rightline=true
	]
	
	{
		% \fontfamily{crimson}
		% \fontsize{10pt}{12pt}\selectfont
		
		\noindent	
		{\bf Inputs:} 		
		\begin{tabularx}{\linewidth}{llX}
			$\displaystyle \sldn := \frac{A L^2}{I}$ \quad (\ref{eq:slenderness-rotundness})
			\ , 
			&
			$\dfbs{\X} (\Xb , \tb) = \dfbhs{\X} (\Xb , \tb)$  
			\ , 
			&  
			$\dfbs{\X} (\Xb , \tb) = \dfbhs{\X} (\Xb , \tb)$ \quad (\ref{eq:distributed-forces})
		\end{tabularx}

	}
	
\end{mdframed}

% 2023.06.28
% changed Form i to Form i+1 to be consistent with the code
\subsection{Form 1: No operator splitting}
\label{sc:Kirchhoff-rod-Form-1}
%\vspace{0.5cm}
\begin{mdframed}[
	% style=mystyle,
	% backgroundcolor=lightgray,
	backgroundcolor=mygray,
	% backgroundcolor=white,
	topline=true,
	bottomline=true,
	leftline=true,
	rightline=true
	]
	
	{
		% \fontfamily{crimson}
		% \fontsize{10pt}{12pt}\selectfont
		
		\noindent	
		{\bf Form 1:} No operator splitting		
		\begin{align}
			\sldn 
			\left[ \ubp 2 + \eslpbp{1} \sin \eslpb \right]
			+ \left[ \iopeb \eslpbp{2} \right]^{(1)} \sin \eslpb 
			+ \left[\iopeb \eslpbp{2} \right] \eslpbp{1} \cos \eslpb + \dfbs{\X} 
			= \nddt \ub 
			\ ,
			\tag{\ref{eq:eom-u-bar}}
		\end{align}
		\begin{align}
			\sldn 
			\left[ \vbp 2 - \eslpbp{1} \cos \eslpb \right]
			- \left[ \iopeb \eslpbp{2} \right]^{(1)} \cos \eslpb 
			+ \left[\iopeb \eslpbp{2} \right] \eslpbp{1} \sin \eslpb + \dfbs{\Y}
			= \nddt \vb 
			\ ,
			\tag{\ref{eq:eom-v-bar}}
		\end{align}
		
		\begin{tabularx}{\linewidth}{llX}
%			$\displaystyle \sldn := \frac{A L^2}{I}$ \quad (\ref{eq:slenderness-rotundness})
%			\ , 
%			&
			$\displaystyle \tan \eslpb = \frac{\vbp{1}}{1 + \ubp{1}}$ \quad (\ref{eq:exact-slope-bar})
			\ , 
			&  
			$\displaystyle \iopeb =
			\left\{ 
			\left[ 1 + \ubp{1} \right]^2 + \left[ \vbp{1} \right]^2 
			\right\}^{-1/2}$ \quad (\ref{eq:inverse-one-plus-e})
		\end{tabularx}		
		
		\noindent
		\emph{Dependent variables:}	4
		\newline
		$\{\ub , \vb , \eslpb\}$, $\iopeb$.

	}
	
\end{mdframed}
% Even though in principle,
% the number of dependent variables could be reduced from 4 to 3 by substituting the expression for $\iopeb$ in Eq.~\eqref{eq:inverse-one-plus-e} into the balance of linear momenta in Eqs.~\eqref{eq:eom-u-bar}-\eqref{eq:eom-v-bar}, but because of the space derivative on $\iopeb$, it would be simpler and computationally more efficient to use Form 1 above to avoid complex expressions and multiple computations of the derivatives $\ubp{1}$ and $\vbp{1}$.
In principle, the number of dependent variables could be reduced from 4 to 3 by substituting the expression for $\iopeb$ in Eq.~\eqref{eq:inverse-one-plus-e} into the balance of linear momenta in Eqs.~\eqref{eq:eom-u-bar}-\eqref{eq:eom-v-bar}. 
Due to the space derivative on $\iopeb$, it would be simpler and computationally more efficient to use Form 1 above to avoid complex expressions and multiple computations of the derivatives $\ubp{1}$ and $\vbp{1}$.

Similarly, while it is possible to further reduce the number of dependent variables from 3 to 2, because of the derivative $\eslpbp{1}$, it is simpler and more efficient to use Form 1 above, while paying attention to avoid multiple computation of $\sin \eslpb$ and $\cos \eslpb$.

%{\color{red} [NOTE: 
%	2023.08.01 - Substitution of the expressions for $\eslpb$ in Eq.~\eqref{eq:exact-slope-bar} and for $\iopeb$ in Eq.~\eqref{eq:inverse-one-plus-e} into the balance of linear momenta in Eqs.~\eqref{eq:eom-u-bar}-\eqref{eq:eom-v-bar} would lead to complex expression, with a loss of computational efficiency.
%	
%	2023.05.26 - It remains to demonstrate via numerical experiments whether this substitution would reduce efficiency, as predicted, without improving accuracy.  ENDNOTE]
%}

\begin{rem}
	\label{rm:time-shift-early-stopping}
	Form 1, pinned-pinned bar: Time shift and early stopping.
%	Time shift and early stopping: Form 1, pinned-pinned bar.
	{\rm
		A characteristic of Form 1 is the floating/shifting computed time history when there is convergence in the optimization process.  Time shift examples are shown in Figures~\ref{fig:23.7.22 R1 A-PP v1 midspan,shape100000} and \ref{fig:23.7.23 R1d A-PP v1 midspan,shape200000} for the axial motion of a pinned-pinned elastic bar; {see Eq.~\eqref{eq:axial-pinned-pinned-BCs} for the definition of the pinned boundary conditions}.
		
		Another characteristic of Form 1 is ``early stopping,'' i.e., the lowest loss value occurs well before the end of the optimization process (see, e.g., \cite{vuquoc2023deep}), which may diverge for large learning rates, and which could be \emph{stopped earlier}.
		See Remark~\ref{rm:early-stopping} and examples therein.
		% {\color{red} [NOTE: 2023.08.05, need examples. ENDNOTE]}
	}
	\hfill $\blacksquare$
\end{rem}

\begin{rem}
	\label{rm:static-solution-Form-1-pinned-free-1}
	Form 1, pinned-free bar: 
	Static solution.
	{\rm
		For sufficiently large learning rate, a static solution could manifest in the training process; see Remark~\ref{rm:static-solution-Form-1-pinned-free-2}.
		% For a pinned-free bar, a static solution could also be found in Form 2a and form 3. 
	}
	\hfill$\blacksquare$
\end{rem}

% 2023.06.28
% changed Form i to Form i+1 to be consistent with the code
\subsection{Form 2a: Split time derivatives}
\label{sc:Kirchhoff-rod-Form-2a}
%\vspace{0.5cm}
\begin{mdframed}[
	% style=mystyle,
	% backgroundcolor=lightgray,
	backgroundcolor=mygray,
	% backgroundcolor=white,
	topline=true,
	bottomline=true,
	leftline=true,
	rightline=true
	]
	
	{
		% \fontfamily{crimson}
		% \fontsize{10pt}{12pt}\selectfont
		
		\noindent	
		{\bf Form 2a:} Split only time derivatives to reduce to first order in time.		
		\begin{align}
			\sldn 
			\left[ \ubp 2 + \eslpbp{1} \sin \eslpb \right]
			+ \left[ \iopeb \eslpbp{2} \right]^{(1)} \sin \eslpb 
			+ \left[\iopeb \eslpbp{2} \right] \eslpbp{1} \cos \eslpb + \dfbs{\X} 
			= \pbsd{\X}
			\ , \quad
			\ndt \ub = \pbs{\X}
			\label{eq:eom-u-bar-form-2}
		\end{align}
		\begin{align}
			\sldn 
			\left[ \vbp 2 - \eslpbp{1} \cos \eslpb \right]
			- \left[ \iopeb \eslpbp{2} \right]^{(1)} \cos \eslpb 
			+ \left[\iopeb \eslpbp{2} \right] \eslpbp{1} \sin \eslpb + \dfbs{\Y}
			= \pbsd{\Y}
			\ , \quad
			\ndt \vb = \pbs{\Y}
			\label{eq:eom-v-bar-form-2}
		\end{align}
		
		\begin{tabularx}{\linewidth}{llX}
			$\displaystyle \tan \eslpb = \frac{\vbp{1}}{1 + \ubp{1}}$ \quad (\ref{eq:exact-slope-bar})
			\ , 
			&  
			$\displaystyle \iopeb =
			\left\{ 
			\left[ 1 + \ubp{1} \right]^2 + \left[ \vbp{1} \right]^2 
			\right\}^{-1/2}$ \quad (\ref{eq:inverse-one-plus-e})
		\end{tabularx}	
		
		\noindent
		\emph{Dependent variables:}	6
		\newline
		$\{\ub , \vb , \eslpb\}$, $\iopeb$, $\{\pbs{\X}, \pbs{\Y}\}$.
		
	}
	
\end{mdframed}
The number of dependent variables could be reduced from 6 to 4 by substituting the expressions for $\eslpb$ in Eq.~\eqref{eq:exact-slope-bar} and for $\iopeb$ in Eq.~\eqref{eq:inverse-one-plus-e} into the balance of linear momenta in Eqs.~\eqref{eq:eom-u-bar-form-2}-\eqref{eq:eom-v-bar-form-2} at the cost of repetitive evaluations and thus a loss of efficiency.
%\newline
%{\color{red} NOTE: 2023.05.26. It remains to demonstrate via numerical experiments whether this substitution would reduce efficiency, as predicted, without improving accuracy.  ENDNOTE}

\begin{rem}
	\label{rm:static-solution}
	Form 2a: Static solution.
	{\rm
		Under the right conditions, a static solution could occur in Form 2a for both the pinned-pinned bar and the pinned-free bar.
		Examples are given in Remark~\ref{rm:static-solution-avoid} with methods to avoid the static solution mentioned.
		%
%		\red{23.9.6 REWRITE:}
%		A characteristic of Form 2a is convergence toward a static (or quasi-static) solution for large learning rates, instead of ``early stopping'' as in Form 1 mentioned in Remark~\ref{rm:time-shift-early-stopping}.
		% {\color{red} [NOTE: 2023.08.05, need examples. ENDNOTE]} 
		%
		Static solution could also be found in Form 1 for the pinned-free bar; see
		Remark~\ref{rm:static-solution-Form-1-pinned-free-1}, Remark~\ref{rm:static-solution-Form-1-pinned-free-2}. 
	}
	\hfill $\blacksquare$
\end{rem}

\subsection{Form 2b: Split space derivatives}
\label{sc:Kirchhoff-rod-Form-2b}
%\vspace{0.5cm}
\begin{mdframed}[
	% style=mystyle,
	% backgroundcolor=lightgray,
	backgroundcolor=mygray,
	% backgroundcolor=white,
	topline=true,
	bottomline=true,
	leftline=true,
	rightline=true
	]
	
	{
		% \fontfamily{crimson}
		% \fontsize{10pt}{12pt}\selectfont
		
		\noindent	
		{\bf Form 2b:} Split only space derivatives to reduce to first order in space.		
		\begin{align}
			\sldn 
			\left[ \slpbsp{\du}{1} + \eslpbp{1} \sin \eslpb \right]
			- \shrbp{1}
			% + \left[ \iopeb \eslpbp{2} \right]^{(1)} 
			\sin \eslpb 
			- \shrb
			% + \left[\iopeb \eslpbp{2} \right] 
			\eslpbp{1} \cos \eslpb + \dfbs{\X} 
			= \nddt \ub
			\ , 
			\label{eq:eom-u-bar-form-2b}
		\end{align}
		\begin{align}
			\sldn 
			\left[ \slpbsp{\dv}{1} - \eslpbp{1} \cos \eslpb \right]
			+ \shrbp{1}
			% - \left[ \iopeb \eslpbp{2} \right]^{(1)} 
			\cos \eslpb
			- \shrb 
			% + \left[\iopeb \eslpbp{2} \right]  
			\eslpbp{1} \sin \eslpb + \dfbs{\Y}
			= \nddt \vb
			\ , \quad
			\label{eq:eom-v-bar-form-2b} 
		\end{align}		
		\begin{align}
			&
			\ubp{1} = \slpbs{\du} , \quad
			% &
			\eslpbp{1} = \Mb  , \quad 
			\label{eq:axial-slope-curvature-moment}
			\\
			&
			\vbp{1} = \slpbs{\dv}  , \quad
			% &
			\Mbp{1} = - \opeb \shrb  , 
			\label{eq:transversal-slope-moment-shear}
		\end{align}
		
		\hspace{-0.5cm}
		\begin{tabularx}{\linewidth}{llX}
			%			$\displaystyle \sldn := \frac{A L^2}{I}$ \quad (\ref{eq:slenderness-rotundness})
			%			\ , 
			&
			$\displaystyle \tan \eslpb = \frac{\vbp{1}}{1 + \ubp{1}}$ \quad (\ref{eq:exact-slope-bar})
			\ , 
			&  
			$\displaystyle 
			\opeb =
			\left\{ 
			\left[ 1 + \ubp{1} \right]^2 + \left[ \vbp{1} \right]^2 
			\right\}^{1/2}$ \quad \eqref{eq:define-ope}, \eqref{eq:extension-bar}
		\end{tabularx}	
		
		\noindent
		\emph{Dependent variables:}	7
		\newline
		$\{\ub , \vb , \eslpb\}$, $\opeb$, $\Mb$, $\{\slpbs{\du} , \slpbs{\dv}\}$
		
	}
	
\end{mdframed}

\subsection{Form 3: Split time and space derivatives}
\label{sc:Kirchhoff-rod-Form-3}
%\vspace{0.5cm}
\begin{mdframed}[
	% style=mystyle,
	% backgroundcolor=lightgray,
	backgroundcolor=mygray,
	% backgroundcolor=white,
	topline=true,
	bottomline=true,
	leftline=true,
	rightline=true
	]
	
	{
		% \fontfamily{crimson}
		% \fontsize{10pt}{12pt}\selectfont
		
		\noindent	
		{\bf Form 3:} Split both time and space derivatives to reduce to first order in both time space.		
		\begin{align}
			\sldn 
			\left[ \slpbsp{\du}{1} + \eslpbp{1} \sin \eslpb \right]
			- \shrbp{1}
			% + \left[ \iopeb \eslpbp{2} \right]^{(1)} 
			\sin \eslpb 
			- \shrb
			% + \left[\iopeb \eslpbp{2} \right] 
			\eslpbp{1} \cos \eslpb + \dfbs{\X} 
			= \pbsd{\X}
			\ , 
			\label{eq:eom-u-bar-form-3}
		\end{align}
		\begin{align}
			\sldn 
			\left[ \slpbsp{\dv}{1} - \eslpbp{1} \cos \eslpb \right]
			+ \shrbp{1}
			% - \left[ \iopeb \eslpbp{2} \right]^{(1)} 
			\cos \eslpb
			- \shrb 
			% + \left[\iopeb \eslpbp{2} \right]  
			\eslpbp{1} \sin \eslpb + \dfbs{\Y}
			= \pbsd{\Y}
			\ , \quad
			\label{eq:eom-v-bar-form-3} 
		\end{align}		
		\begin{align}
			&
			\ubp{1} = \slpbs{\du} , \quad
			% &
			\ndt \ub = \pbs{\X}  , \quad
			% &
			\eslpbp{1} = \Mb  , \quad 
			\label{eq:axial-slope-momentum-curvature-moment}
			\\
			&
			\vbp{1} = \slpbs{\dv}  , \quad
			% &
			\ndt \vb = \pbs{\Y}  , \quad
			% &
			\Mbp{1} = - \opeb \shrb  , 
		\end{align}
		
		\hspace{-0.5cm}
		\begin{tabularx}{\linewidth}{llX}
%			$\displaystyle \sldn := \frac{A L^2}{I}$ \quad (\ref{eq:slenderness-rotundness})
%			\ , 
			&
			$\displaystyle \tan \eslpb = \frac{\vbp{1}}{1 + \ubp{1}}$ \quad (\ref{eq:exact-slope-bar})
			\ , 
			&  
			$\displaystyle 
			\opeb =
			\left\{ 
			\left[ 1 + \ubp{1} \right]^2 + \left[ \vbp{1} \right]^2 
			\right\}^{1/2}$ \quad \eqref{eq:define-ope}, \eqref{eq:extension-bar}
		\end{tabularx}	
		
		\noindent
		\emph{Dependent variables:}	9
		\newline
		$\{\ub , \vb , \eslpb\}$, $\opeb$, $\Mb$, $\{\pbs{\X}, \pbs{\Y}\}$, $\{\slpbs{\du} , \slpbs{\dv}\}$
		
	}
	
\end{mdframed}
The number of dependent variables can be reduced from 9 to 8 by substituting the expression for $\opeb$ only once in the expression relating moment to shear in Eq.~\eqref{eq:axial-slope-momentum-curvature-moment}$_3$, without repetitive computation.
On the other hand, further reducing the number of dependent variable from 8 to 7 by substituting the expression for $\eslpb$ in Eq.~\eqref{eq:exact-slope-bar} into the balance of linear momenta Eqs.~\eqref{eq:eom-u-bar-form-3}-\eqref{eq:eom-v-bar-form-3} would be inefficient due to repetitive computation of $\eslpb$.
%\newline
%{\color{red} NOTE: 2023.05.26. It remains to demonstrate via numerical experiments whether this substitution would reduce efficiency, as predicted, without improving accuracy.  ENDNOTE}

\subsection{Form 4: Balance of momenta with no substitution}
\label{sc:Kirchhoff-rod-Form-4}

%\vspace{0.5cm}
\begin{mdframed}[
	% style=mystyle,
	% backgroundcolor=lightgray,
	backgroundcolor=mygray,
	% backgroundcolor=white,
	topline=true,
	bottomline=true,
	leftline=true,
	rightline=true
	]
	
	{
		% \fontfamily{crimson}
		% \fontsize{10pt}{12pt}\selectfont
		
		\noindent	
		{\bf Form 4:} Balance of momenta with first-order derivatives in both space and time.		
		\begin{align}
			\Nbsp{\X}{1} - \Sbsp{\X}{1}
			+ \dfbs{\X} 
			= \pbsd{\X} 
			\ , \quad
			\Nbsp{\Y}{1} + \Sbsp{\Y}{1}
			+\dfbs{\Y}
			= \pbsd{\Y}
			\label{eq:bal-lin-mom-form-4}
		\end{align}
		\begin{align}
			&
			\Nbs{\X} = \sldn \left[ 1 + \ubp{1} - \cos \eslpb \right] , \quad
			\Sbs{\X} = - \iopeb \Mbp{1} \sin \eslpb
			\ ,
			\label{eq:normal-force-shear-bar-X}
			\\
			&
			\Nbs{\Y} = \sldn \left[ \vbp{1} - \sin \eslpb \right] , \quad
			\Sbs{\Y} = - \iopeb \Mbp{1} \cos \eslpb
			\ ,
			\label{eq:normal-force-shear-bar-Y}
			\\
			&
			\eslpbp{1} = \Mb , \quad 
			\ubd = \pbs{\X} , \quad
			\vbd = \pbs{\Y}
			\ ,
		\end{align}
		
		\begin{tabularx}{\linewidth}{llX}
%			$\displaystyle \sldn := \frac{A L^2}{I}$ \quad (\ref{eq:slenderness-rotundness})
%			\ , 
			%&
			$\displaystyle \tan \eslpb = \frac{\vbp{1}}{1 + \ubp{1}}$ \quad (\ref{eq:exact-slope-bar})
			\ , 
			&  
			$\displaystyle \iopeb =
			\left\{ 
			\left[ 1 + \ubp{1} \right]^2 + \left[ \vbp{1} \right]^2 
			\right\}^{-1/2}$ \quad (\ref{eq:inverse-one-plus-e})
		\end{tabularx}
		
		\noindent
		\emph{Dependent variables:}	11
		\newline
		$\{\ub , \vb , \eslpb\}$, $\iopeb$, $\Mb$, $\{\pbs{\X}, \pbs{\Y}\}$, $\{\Nbs{\X}, \Nbs{\Y}, \Sbs{\X}, \Sbs{\Y}\}$

	}
	
\end{mdframed}
The number of dependent variables can be reduced from 11 to 9 by substituting the expressions for $\eslpb$ in Eq.~\eqref{eq:exact-slope-bar} and for $\iopeb$ in Eq.~\eqref{eq:inverse-one-plus-e} into the expressions for normal and shear force components $(\Nbs{\X} , \Nbs{\Y})$ and $(\Sbs{\X} , \Sbs{\Y})$ in Eq.~\eqref{eq:normal-force-shear-bar-X}$_2$ and Eq.~\eqref{eq:normal-force-shear-bar-Y}$_2$, respectively.  But the repetitive computations of the same quantities (e.g., $\eslpb$ in four different expressions, and $\iopeb$ in two different expressions) could make such substitution inefficient.
%\newline
%{\color{red} NOTE: 2023.05.23. It remains to demonstrate via numerical experiments whether this substitution would reduce efficiency, as predicted, without improving accuracy.  ENDNOTE}

%\noindent
%{\color{red} \ding{42} I AM HERE \today}
%\newline
%PINN PDE FORMS, boundary conditions, initial conditions DONE.

%\subsection{Loss functions}

% PINN formulation for generic PDEs
% \subsection{Generic PDEs, auxilliary conditions, loss function}
\section{Generic PDEs, auxilliary conditions, loss function}
\noindent
For short, boundary conditions and initial conditions are together referred to as \emph{auxilliary} conditions.  The system in all of the above forms $k = 1, \ldots, 4$---PDEs and the associated auxilliary conditions---can be generically written as a series of representative \emph{dynamic} nonlinear differential operators $\dpdesp{i}{(k)}$ as follows:
\begin{align}
	&
	\dpdesp{i}{(k)} (
	\{\partial_\X^p u_j , \partial_\X^{p-1} u_j , \ldots, \partial_\X^1 u_j\} , 
	\{\partial_\ti^q u_j , \partial_\X^{q-1} u_j , \ldots, \partial_\ti^1 u_j\} , 
	u_j ,
	\ti) 
	= \dpdesp{i}{(k)} (\bkar{\du} (\X , \ti), t)
	%	= \dpdesp{i}{(k)} (\X , \ti)
	= 0
	\ , 
	\nonumber
	\\
	&
	\text{with }
	k = 1,\ldots, 4 \ , \quad
	i = 1,\ldots, n_i^{(k)} , \quad j = 1,\ldots, n_j^{(k)} , \quad 
	p = 1,\ldots, n_p^{(k)} , \quad q = 1,\ldots, n_q^{(k)} 
	\ ,
	\label{eq:generic-PDE-aux-conditions}
\end{align}
where $u_j (\X, \ti)$ are the dependent variables, 
% functions of $(\X, \ti)$, 
and $\partial_\X^p u_j (\X, \ti)$ is the $p$th partial derivative of $u_j$ with respect to $\X$, and similarly for $\partial_\ti^q u_j (\X, \ti)$.
When $n_q^{(k)} = 0$, there are no terms with time derivatives, as in the static case.  For convenience, the following dynamic-operator array, which will be used later, is defined:
\begin{align}
	\bdpdep{(k)} = \left\{ \dpdesp{i}{(k)} , \ i=1, \ldots , n_i^{(k)} \right\}
	\ , \quad
	k = 0, \ldots, 3
	\ .
	\label{eq:array-operators-dynamics}
\end{align}

\noindent
The loss function for Form $k$ to minimize is
\begin{align}
	J^{(k)} 
	= \sum_{i=1}^{i=n_i^{(k)}} \wlsp{i}{(k)} \left[ \dpdesp{i}{(k)} \right]^2
	\ ,
	\label{eq:loss-function}
\end{align}
where $\wls{i}$ is the weight associated with the differential operator $\dpdesp{i}{(k)}$, for $i=1, \ldots, n_i^{(k)}$, playing the role of an ``error'' function.

% stochastic gradient descent
% \input{07-optimization}

% backpropagation method
\section{Numerical examples}
\label{sc:numerical-examples}
\noindent
At first, to set up PINN problems of structural dynamics, we opted to use the DeepXDE (abbreviated as ``DDE'') framework \cite{lu2021deepxde}, which was likely the most well-documented open-source software for PINN, with many examples to guide users and a forum for open discussions of problems users encountered.  {Specifically, since performance depends on the DDE version used, we provide this information.  In the present paper, we use DDE v1.9.2 or v1.9.3 with the TensorFlow backend, which together is abbreviated as, e.g., DDE-T v1.9.3, or simply \DDET, since DDE v1.9.3 is mostly used in this paper}.\footnote{
	{For example, Figure~\ref{fig:23.7.26 R5c A-PP F2a lr0.02x5 cycles=9 Nsteps=400000 regular N51 W32 H2 He init-1}, RunID 23726R5c-1, was obtained with ``deepxde-1.9.2 pyaml-23.7.0 scikit-optimize-0.9.0'', whereas Figure~\ref{fig:DEBUG v1.9.3 seed42, 23.9.4 R1 A-P bcrP F1J(2a,3) lr0.001 Cy1-9-NCA Nsteps400000 regular N51 W64 H4 Glorot init-1}, RunID 2394R1a-1, was obtained with ``deepxde-1.9.3 pyaml-23.9.1 scikit-optimize-0.9.0''.}  {In the follow-up paper, we will report a significant GPU time difference between DDE-T v1.10.0 and DDE-T v1.10.1.}
} 

In using \DDET, we encountered many problems such as shift and amplification in the dynamic solution, static-solution time history, in addition to the difficulty in designing appropriate learning-rate schedules, mostly by trial-and-error and accumulated experience that we want to convey to readers, to arrive at good solutions.  Along the way, we devised various strategies to solve these pathological problems, such as different forms of the governing PDE system to solve the shift and amplification problem, and the barrier function to avoid static solutions. Even then, a significant amount of time was devoted to get these methods to work, compared to our script based on the \JAX\ framework, the \emph{High-Performance Array Computing} library, which did not manifest any of these pathological problems; see Remarks~\ref{rm:learning-rate-schedule-4}, \ref{rm:static-solution-Form-1-DDE-T-v-JAX}.  
% \red{I AM HERE 2023.09.19.}

% ALEX: removed paragraph
% LOC: Keep this paragraph
Since \DDET\ is accessible to the public, the detailed documentation of the \DDET\ pathological problems, our proposed solutions and results, would not only be useful for the general readers, but also for the developers.

Because of the large number of parameters, resulting in many different cases, and to allow for convenient, frequent back-referencing in figure captions and elsewhere in the text, the sections below are organized as a series of remarks, each addressing a specific issue/topic {with immediate reference to the corresponding figures/results for illustration where applicable}.\footnote{
	{The remarks are important integral parts of the text.  It is recommended, particularly for readers not familiar with PINN and its training process, not to skip, but to read through all remarks so to be familiar with the corresponding illustrative examples on the remark topics. The ultimate objective is to solve nonlinear PDAEs, after developing an experience with the training process using \DDET\ and \JAX\ to solve standard linear PDEs using the proposed novel PINN formulations through a normalization/standardization presented in this paper.}
}

\begin{rem}
	\label{rm:colab-GPU}
	Colab GPU time.
	{\rm
		Even though we did use a dedicated Linux machine with Nvidia GPU (Graphics Processing Units) for script development, 
		all results presented here were obtained using \href{https://colab.research.google.com/?utm_source=scs-index}{Google Colab}'s Nvidia Tesla K80 GPU with 12 GB RAM (free).  See the GPU time in, e.g., Figure~\ref{fig:23.8.16 R1a A-PP F(1)2a lr0.01 cy5-CA Nsteps200000 regular N51 W64 H2 He init}.
	}
	\hfill$\blacksquare$
\end{rem}

\begin{rem}
	\label{rm:runid}
	RunID.
	{\rm
		At the end of the caption of a figure, there is a RunID, such as ``23723R1d-1,'' 
		% in Figure~\ref{fig:23.7.23 R1d A-PP v1 loss200000 shape25000} (Step 25,000) 
		to indicate which run produced that figure
		(two images): 
		``23723'' being the date 2023 Jul 23,  ``R1d'' being ``Run 1d,'' the order of the DeepXDE script file executed on that date, and ``-1'' being the first figure (two images) coming from that run, from which 
%		there are three more figures:   
%		Figure~\ref{fig:23.7.23 R1d A-PP v1 midspan,shape50000} with RunID 23723R1d-2 (Step 50,000),
%		Figure~\ref{fig:23.7.23 R1d A-PP v1 midspan,shape100000} with RunID 
%		23723R1d-3 (Step 100,000),
%		and
%		Figure~\ref{fig:23.7.23 R1d A-PP v1 midspan,shape200000} with RunID
%		23723R1d-4 (Step 200,000), the final result, which 
		%
		the final result in Figure~\ref{fig:23.7.23 R1d A-PP v1 midspan,shape200000} with RunID
		23723R1d-4 (Step 200,000)
		is presented in the paper body, Section~\ref{sc:axial-Form-1} to illustrate the time shift in Form 1 (Remark~\ref{rm:time-shift-Form-1}, Section~\ref{sc:axial-Form-1}).
		%, whereas the previous three figures (\ref{fig:23.7.23 R1d A-PP v1 loss200000 shape25000}, \ref{fig:23.7.23 R1d A-PP v1 midspan,shape50000}, \ref{fig:23.7.23 R1d A-PP v1 midspan,shape100000}) are relegated to  Appendix~\ref{app:time-shift-Form-1-Form-2a}. 
		
		A reason for introducing the RunID is because, even though ``stochastically reproducible,'' the results are ``deterministically irreproducible,'' 
		unless we run under deterministic mode, which is slower,
		as explained in Remark~\ref{rm:stochastic-reproducibility}. For this reason, we kept a large majority of the executed scripts (run under non-deterministic mode), the results from a small number of which are presented here.
	}
	\hfill$\blacksquare$
\end{rem}

\subsection{Network architecture, data-point grids}
\label{sc:networks-grids}

\begin{rem}
	\label{rm:parameter-names}
	% Network architecture and parameters.
	Computational-domain, network variables, dynamic results.
	{\rm
		To shortened the description, the following parameter symbols are used.
		For bars and beams, the domain for the space coordinate $x$ is the interval $[0, 1]$, the domain for time coordinate is $[0, T]$, and the rectangular space-time \emph{computational domain} is $[0, 1] \times [0, T]$, with N being the number of data points on one of its sides.  The data points, with a total number of NxN, are used to evaluate the loss functions.
		% For the loss function, N is the number of data points on one side of the rectangular computation space-time domain.  
		For the network architecture, 
		n\_inp is the number of inputs,
		W the width of a hidden layer,
		H the number of hidden layers, 
		and 
		n\_out the number of outputs.
		
		\emph{Network inputs:} For static problems, n\_inp=1 (for the data-point $x$-coordinates), whereas for dynamic problems, n\_inp=2 (for the data point $(x, t)$ coordinates). Since the value of n\_inp is clear from the context, it will not be listed in figure captions; see, e.g., Figure~\ref{fig:23.8.22 R3b A-PP F(1J)2a(3) lr0.001 cy5-CA cy6-NCA Nsteps250000 regular N51 W64 H4 He init}.
		
		\emph{Network outputs:} The number of dependent variables of the PDE(s) depends on the Form used.  For example, Form 4 of the Kirchhoff-Love rod (Section~\ref{sc:Kirchhoff-rod-Form-4}) has n\_out=11 outputs; see, e.g., Figure~\ref{fig:23.8.22 R3b A-PP F(1J)2a(3) lr0.001 cy5-CA cy6-NCA Nsteps250000 regular N51 W64 H4 He init} for Form 2a of a pinned-pinned bar, with n\_out=2.
		
		For results related to dynamic problems, the phrase ``\emph{time history}'' is omitted in figure captions for conciseness, since the meaning is clear from the context.  For example, in Figure~\ref{fig:23.8.22 R3b A-PP F(1J)2a(3) lr0.001 cy5-CA cy6-NCA Nsteps250000 regular N51 W64 H4 He init}, the shape (left) and the midspan displacement (right) are time histories, with the shape having ``x1'' as the space $x$ axis and ``x2'' the time $t$ axis, which is also the horizontal axis of the midspan displacement.
	}
	\hfill $\blacksquare$
\end{rem}

\begin{figure}[tph]
	\centering
	%
	% 24.7.12 replace blank spaces by underscores as required by arXiv
%	\includegraphics[width=0.5\textwidth]{Figures/23.7.21_R1_N51_W64_H4_model_summary_12737_params.PNG}
	\includegraphics[width=0.5\textwidth]{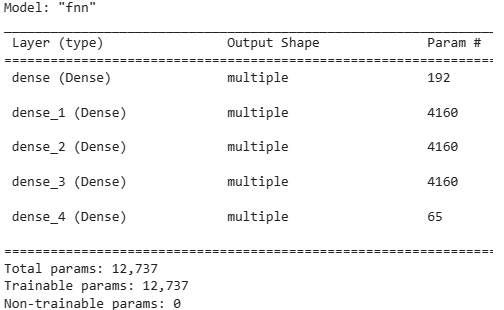}
	\caption{
		% DDE-T.
		% v1.9.3.
		\DDET.
		\emph{Pinned-pinned elastic bar, axial motion}.    
		Model summary.
		Feedforward neural network (fnn).
		Remark~\ref{rm:network-params}.  
		Section~\ref{sc:axial-Form-1},
		Form 1.
		Dense-connection layers, each between two consecutive layers of neurons.
		\emph{Network:}
		Remark~\ref{rm:parameter-names},
%		n\_inp=2 inputs $(x, t)$, 
%		W=64 neurons per hidden layer, 
%		H=4 hidden layers,
%		n\_out=1 output $u(x, t)$, 
		n\_inp=2, 
		W=64, 
		H=4,
		n\_out=1,
		12737 parameters.
		Six neuron layers (1 input, 4 hidden, 1 output), five connection layers (pairs of consecutive neuron layers).
		{\scriptsize (23721R1-1)}
	}
	\label{fig:model-summary-Form-1-N51W64H4-12737-params}
\end{figure}

\begin{rem}
	\label{rm:network-params}
	% Network parameters, dense connection. 
	Total parameters of dense networks.
	{\rm
		As an example, Form 1 of the axial motion of a pinned-pinned elastic bar (Section~\ref{sc:axial-Form-1}) and a network with n\_inp=2 inputs $(x, t)$, 
		W=64 neurons per hidden layer, 
		H=4 hidden layers,
		n\_out=1 output $u(x, t)$ lead to 12,737 parameters (weights and biases),  Figure~\ref{fig:model-summary-Form-1-N51W64H4-12737-params}, obtained as follows. 
		
		From the 2 inputs to the 64 neurons in hidden layer 1, there are 64*2 weights (connections) and 64 biases, thus in total (64*2 + 64) = 192 parameters.  From hidden layer 1 to hidden layer 2, there are 64*64 weights (connections) and 64 biases, thus in total (64*64 + 64) = 1460 parameters. 
		 
		Similarly for the two subsequent pairs of hidden layers, the number of weights and biases is 1460 per pair, thus in total (64*64 + 64)*3 = 12480 parameters for all three pairs of hidden layers. From the hidden layer 4 to the one output layer, there are 64*1 weights (connections) and 1 bias, thus in total (64 + 1) = 65 parameters.  Summing up, there are 12,737 parameters (weights and biases) in total, Figure~\ref{fig:model-summary-Form-1-N51W64H4-12737-params}. 
		
		In the model-summary output of the
		DeepXDE framework \cite{lu2021deepxde}, 
		Figure~\ref{fig:model-summary-Form-1-N51W64H4-12737-params},
		the term ``layer'' is used to designate ``layer of connections,'' which is formed by two ``layers of neurons'' as used here, and the type of connections mentioned above is termed ``dense,'' i.e., each neuron in the 1st layer is connected to all neurons in the 2nd layer. 
	}
	\hfill $\blacksquare$
\end{rem}

\begin{rem}
	\label{rm:data-point-grids}
	Data-point grids.
	{\rm
		For the evaluation of the loss functions. We essentially used three types of grids: Regular grid, fixed-random grid, and varying-random grid.  
		
		For a \emph{regular} grid with $N \in \{11, 21, 31, $ $ 41, {51} \}$ as the number of points\footnote{
			{ Meaning, the number $N$ of points could take any value in the set $\{11, 21, 31, $ $ 41, {51} \}$, with ${51}$ as the default number in the normalization (or standardization) of the presentation of our results.  To alleviate the notation, we also omit the overbar on $(\xb, \tb)$ to write simply as $(x, t)$.}
		} on, say one side of the space-time domain {$[0, 1] \times [0, T]$, containing the space-time point} $(x,t)$ with $x \in [0, 1]$ and $t \in [0, T]$, where for example $T \in \{2, 4, 8\}$, there {is} an even number of intervals on each side (space or time), i.e., $\{10, 20, 30, 40, 50\}$, respectively, with $N \times N$ as the total number of data points,  such that one data point is at the center of the space-time domain.  
		Figure~\ref{fig:23.8.22 R3b A-PP F(1J)2a(3) lr0.001 cy5-CA cy6-NCA Nsteps250000 regular N51 W64 H4 He init} shows the shape (left) and midspan-displacement (right) time histories of the axial motion of pinned-pinned elastic bar using Form 2a (Section~\ref{sc:axial-Form-2a}), a model with {N51} W64 H4, having 12,802 parameters, a regular grid, and init\_lr=0.001.
		{By default, N = 51} is used when the value of N is not mentioned in a figure caption. 
		Figure~\ref{fig:pinned-pinned-static-solution-1} shows a static solution obtained using \DDET\ with Form 2a, the same network and number of parameters and regular grid, but with an init\_lr=0.01 (ten times larger).
		Figure~\ref{fig:23.8.21 R1c A-PP F(1J,2a)3 lr0.04,0.03,0.02,0.01,0.005 cy5-VCA Nsteps200000 regular N51 W32 H2 He init-1} shows good shape and midspan displacement using Form 3 (Section~\ref{sc:axial-Form-3}), a model with N51 W32 H2, having 1,251 parameters, and a regular grid. 
		
%		\red{[NOTE: Incoporate into the text these  
%		Figures~\ref{fig:23.8.21 R1c A-PP F(1J,2a)3 lr0.04,0.03,0.02,0.01,0.005 cy5-VCA Nsteps200000 regular N51 W32 H2 He init-1}, \ref{fig:pinned-pinned-static-solution-1} on the use of regular grids.  ENDNOTE]}
		
		For a \emph{random} grid, {$N \times N$ data points are randomly distributed in the space-time $(x, t)$ domain, while ensuring that (1) there are $N$ randomly distributed points along the time $(t)$ axis for each of the two boundary locations $x = 0$ and $x = 1$, and (2) there are $N$ randomly distributed points on the space $(x)$ axis for the initial time $t = 0$.}   A \emph{fixed} random distribution is always generated with the seed 42 \cite{nafis2021story} so to obtain the same pseudo-random $N$ points for every execution.  For a \emph{varying} random grid, no seed is used.
		%
		% For a \emph{random} grid, the $N$ data points are first distributed randomly on each side of the space-time domain, then ``multiplied up'' to form a random grid.  
		% For a \emph{random} grid, $N \times N$ data points are randomly sampled from the space-time domain.  
		% A \emph{fixed} random distribution is always generated with the seed 42 \cite{nafis2021story} so to obtain the same pseudo-random $N$ points for every execution.  For a \emph{varying} random grid, no seed is used.

		As an example of a random grid, see Figure~\ref{fig:random-grid-NO-static-solution}, which corresponded to one of the lowest total loss 0.537e-06 for Form 2a.
		Two reruns of the same exact script that produced Figure~\ref{fig:random-grid-NO-static-solution} yielded the total loss of 0.626e-06 (RunID 2384R2b) and 0.518e-06 (RunID 2384R2c), both at Step 192,000, respectively.  
		Figure~\ref{fig:random-grid-NO-static-solution-2} shows that stopping earlier at Step 146,000 at which the total loss of 7.22e-06 was lowest in Cycle-4, with quasi-perfect midspan displacement and GPU time 442 sec, is a good trade-off.  To give more information, the loss function was shown in Figure~\ref{fig:random-grid-NO-static-solution-2}, since the random grid is similar to that in Figure~\ref{fig:random-grid-NO-static-solution}.
		An even better trade-off is given in Figure~\ref{fig:23.8.4 R2e.1 good trade-off}, with a total loss of 1.179e-06, a quasi-perfect midspan displacement with damping\%=0.1\% (See Remark~\ref{rm:shift-amplification-Form-1}) and GPU time 278 sec. 
		For an additional example of random grids, see also Figures~\ref{fig:23.8.24 R1 A-PP F(1J,2a)3 lr0.04,0.03,0.02,0.01,0.01,0.005 cy1-6-VCA cy7-9-NCA Nsteps400000 random N51 W32 H2 He init-1}-\ref{fig:23.8.24 R1 A-PP F(1J,2a)3 lr0.04,0.03,0.02,0.01,0.01,0.005 cy1-6-VCA cy7-9-NCA Nsteps400000 random N51 W32 H2 He init-3}.
		
		Switching from a regular grid to a \emph{random} grid, {while} keeping all other parameters the same, could lower the total loss; see Remark~\ref{rm:extension-cycle-5} on ``Extended learning-rate schedule'' (ELRS) for an example. 
	}
	\hfill $\blacksquare$
\end{rem}

% {\color{red} I AM HERE, \today}

% \emph{Learning-rate scheduling}.
\subsection{Training (optimization) learning-rate scheduling.}
\label{sc:optimization-learning-rate-scheduling}
%\subsection{Training, learning rate}
%\label{sc:training-learning-rate}

\begin{figure}[tph]
	\centering
	\includegraphics[width=0.9\textwidth]{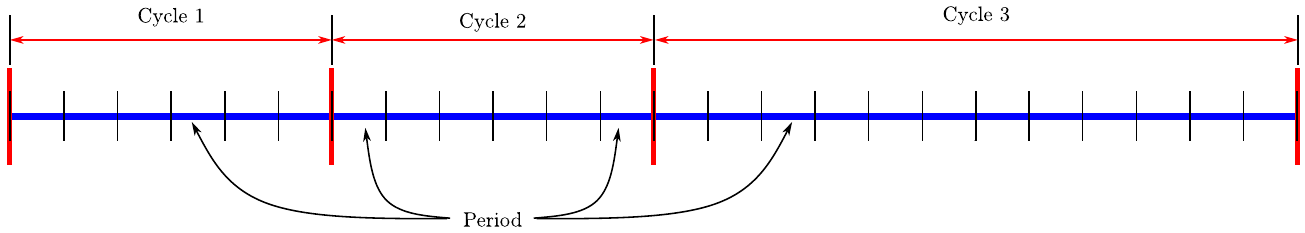}
	\caption{
		\emph{Optimization learning-rate scheduling}. % Section~\ref{sc:numerical-examples}. 
		Section~\ref{sc:optimization-learning-rate-scheduling}. 
		Cycles and periods.
	}
	\label{fig:cycles-periods}
\end{figure}

The Adam optimizer is used in all examples in the present paper,
with learning-rate decay over $N_\text{steps}$, the total number of steps, which is divided into $n_\text{cycles}$, the number of decay cycles, each of which contains $n\text{-steps}_\text{cycle}$, the number of steps per cycle.  This $n\text{-steps}_\text{cycle}$ may vary from cycle to cycle, and is subdivided into $n_\text{periods}$, the number of decay periods, each of which contains $n\text{-steps}_\text{period}$, the number of steps per period; see Figure~\ref{fig:cycles-periods}.  This $n\text{-steps}_\text{period}$ is generally fixed when compiling a model with decay, with the $n\text{-steps}_\text{cycle}$ chosen as a multiple of $n\text{-steps}_\text{period}$. 

After an initial-learning rate for a number of cycles was chosen, e.g., $\text{init\_lr} = 0.001$, 
in any given period (p) within a cycle, the learning rate $\learn_\text{(p-1)}$ at the end of the previous period (p-1) is decayed to the value $\learn_\text{(p)} = \decayf \cdot \learn_\text{(p-1)}$, with $\decayf$ being a decay factor less than one, following a decay function, such as the ``inverse time'' function.
% commented out equation number, 23.10.15
% in Eq.~\eqref{eq:learning-rate-schedule-3}.  

After each cycle and at the beginning of the subsequent cycle, the learning rate is reset to the initial-learning rate, similar to cyclic annealing, for which a detailed explanation can be found in \cite{vuquoc2023deep}.

As an example, we could set in a cycle $n\text{-steps}_\text{cycle} = 50,000$, $n\text{-steps}_\text{period} = 2,500$, so that $n_\text{periods} = 20$, and $\decayf = 0.9$.  

The number of cycles $n_\text{cycles}$ and the total number of steps $N_\text{steps}$ vary depending on the convergence properties of an execution.  With $n\text{-steps}_\text{cycle} = 50,000$ and $n\text{-steps}_\text{period} = 2,500$, we could set $n_\text{cycles} = 4$ and $N_\text{steps} = 200,000$ to ``run all'' these cycles one after another at once, or we could also set $n\text{-steps}_\text{cycle} = 25,000$  then run only the first cycle to examine the solution obtained to decide whether to stop (in the case of reaching a static solution in a dynamic problem, as shown in Figure~\ref{fig:barrier-1} right and Figure~\ref{fig:barrier-2} left) or to continue the optimization process, possibly with different values for $n\text{-steps}_\text{cycle}$, $n_\text{cycles}$, and $N_\text{steps}$.

\begin{figure}[tph]
	\includegraphics[width=0.49\textwidth]{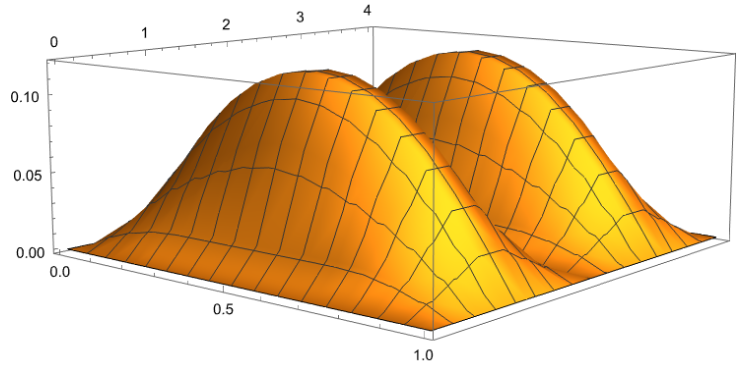}
	\includegraphics[width=0.49\textwidth]{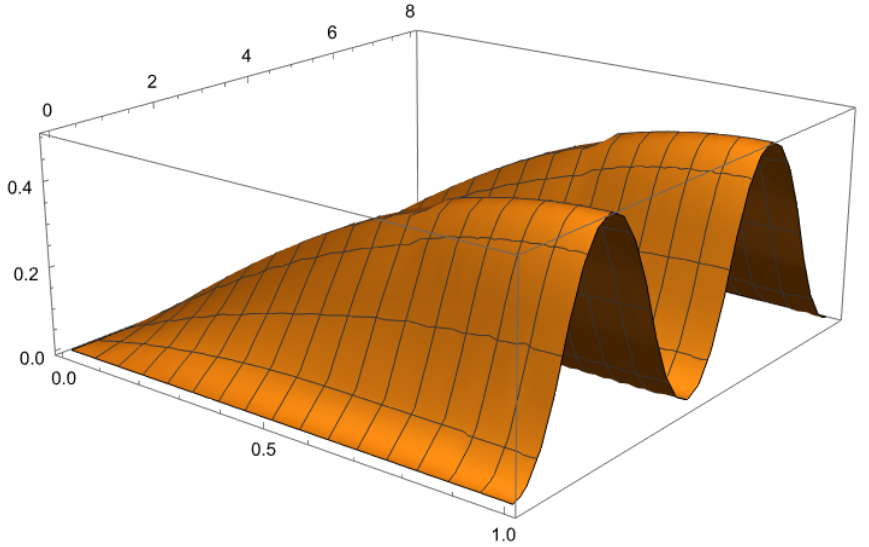}
	\caption{
		\emph{Axial motion of elastic bar.} Section~\ref{sc:axial-motion}. Mathematica solutions.  Eq.~\eqref{eq:eom-euler-bernoulli-axial} with slenderness $\sldn = 1$ and distributed load $\dfbs{\X} = 1/2$. 
		Two vibration periods for each set of boundary conditions.
		\emph{Left:} Pinned-pinned bar, with $\Xb \in \left[0, 1\right]$ and $\tb \in \left[0, 4\right]$.
		\emph{Right:}  Pinned-free bar, with $\Xb \in \left[0, 1\right]$ and $\tb \in \left[0, 8\right]$.
		% {\color{red} [NOTE: 2023.07.22, Alex, we need to values of the maximum displacement to compare with PINN results.  ENDNOTE]}
	}
	\label{fig:axial-Mathematica-solutions}
\end{figure}

% 23.9.20, OLD Fig.9 to move to section Form 2a
\begin{figure}[tph]
	\includegraphics[width=0.49\textwidth]{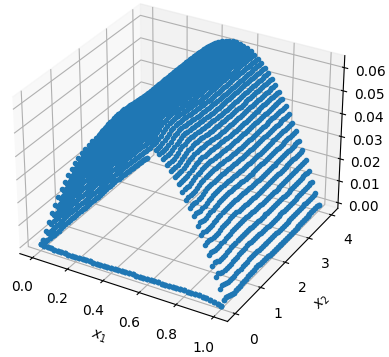}
	\includegraphics[width=0.49\textwidth]{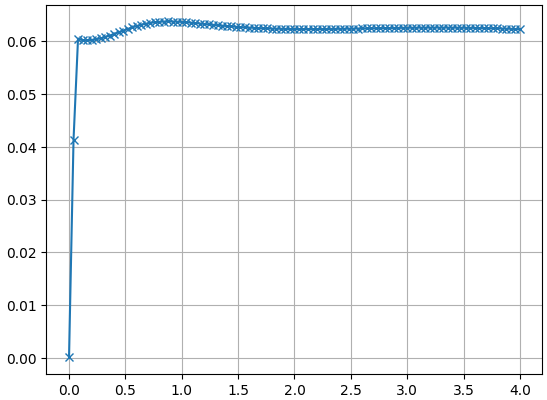}
	\caption{
		% A regular grid with CA may lead quickly to a static solution, as shown in Shape 50000 static 8.17 R10a and Center disp 50000 static 8.17 R10a.
		% A regular grid with NCA may lead to waves with damping (due to a low-capacity network) as shown in Shape 50000 waves, damping 8.17 R10b and Center disp 50000 waves, damping 8.17 R10b.
		\DDET.
		\emph{Pinned-pinned bar, cyclic annealing} \red{(CA)}. 
		\emph{Static}-shape (left), \emph{static} midspan displacement (right), Step 50,000.
		\red{$\star$}
		Section~\ref{sc:pinned-pinned-Form-2a}, \red{Form 2a}.
		Remarks~\ref{rm:static-solution} (Static solution), \ref{rm:static-solution-avoid} (How to avoid).
		\emph{Network:}
		\red{$\star$}
		Remarks~\ref{rm:parameter-names},~\ref{rm:data-point-grids}, 
		T=4,
		W=32,
		H=2, 
		n\_out=2, 
		He-uniform initializer, 
		\red{1,218} parameters,
		\emph{regular} grid.
		\emph{Training:} 
		\red{$\star$}
		Remark~\ref{rm:learning-rate-schedule-1},
		\red{LRS~1},
		init\_lr=\red{0.07},
		n-cycles=2,
		N\_steps=50,000.
		$\bullet$
		{\footnotesize
			Figure~\ref{fig:23.8.17 R10b A-PP F(1)2a lr0.07 cy2-NCA Nsteps50000 regular N51 W32 H2 He init}, NCA, same model and init\_lr=0.07, waves,  damping.
			$\triangleright$
			Figure~\ref{fig:axial-Mathematica-solutions}, reference solution to compare.
		}
		{\scriptsize (23817R10a-1)}
	}
	\label{fig:23.8.17 R10a A-PP F(1)2a lr0.07 cy2-CA Nsteps50000 regular N51 W32 H2 He init}
\end{figure}

% 23.9.20, OLD Fig.10 to move to section Form 2a
\begin{figure}[tph]
	\includegraphics[width=0.49\textwidth]{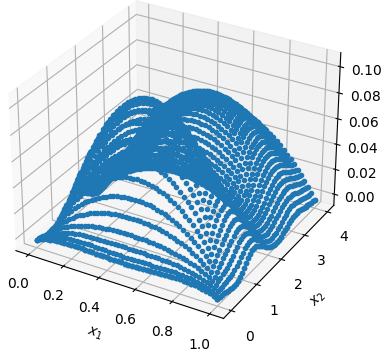}
	\includegraphics[width=0.49\textwidth]{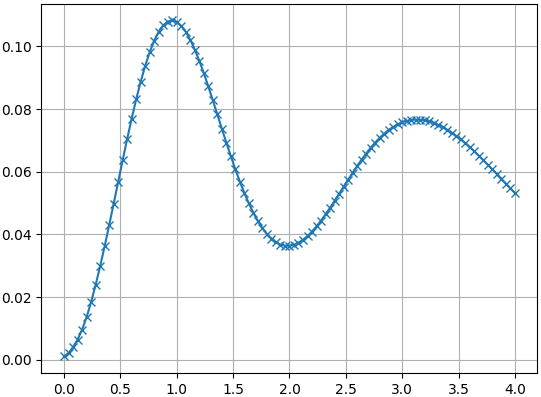}
	\caption{
		% A regular grid with CA may lead quickly to a static solution, as shown in Shape 50000 static 8.17 R10a and Center disp 50000 static 8.17 R10a.
		% A regular grid with NCA may lead to waves with damping (due to a low-capacity network) as shown in Shape 50000 waves, damping 8.17 R10b and Center disp 50000 waves, damping 8.17 R10b.
		\DDET.
		\emph{Pinned-pinned bar, No cyclic annealing} \red{(NCA)}. 
		Shape (left), midspan displacement (right), Step 50,000, \emph{waves with damping}.
		\red{$\star$}
		Section~\ref{sc:pinned-pinned-Form-2a}, \red{Form~2a}.
		Remark~\ref{rm:damping-low-capacity-models}, Damping.
		\emph{Network:}
		\red{$\star$}
		Remarks~\ref{rm:parameter-names},~\ref{rm:data-point-grids}, 
		T=4,
		W=32,
		H=2, 
		n\_out=2, 
		He-uniform initializer, 
		\red{1,218} parameters,
		\emph{regular} grid.
		\emph{Training:} 
		Remark~\ref{rm:learning-rate-schedule-3},
		\red{LRS~3},
		init\_lr=0.07,
		n-cycles=2,
		N\_steps=50,000.
		$\bullet$
		{\footnotesize
			Figure~\ref{fig:23.8.17 R10a A-PP F(1)2a lr0.07 cy2-CA Nsteps50000 regular N51 W32 H2 He init}, CA, same model and init\_lr=0.07, static solution.
			Figures~\ref{fig:23.8.21 R1 A-PP F(1)2a lr0.04,0.03,0.02,0.01,0.005 cy5-VCA Nsteps200000 regular N51 W32 H2 He init} (Form 2a, 1,218 parameters), \ref{fig:23.8.21 R1c A-PP F(1J,2a)3 lr0.04,0.03,0.02,0.01,0.005 cy5-VCA Nsteps200000 regular N51 W32 H2 He init-1}-\ref{fig:23.8.21 R1c A-PP F(1J,2a)3 lr0.04,0.03,0.02,0.01,0.005 cy5-VCA Nsteps200000 regular N51 W32 H2 He init-2} (Form 3, 1,251 parameters), VCA, low damping.
			Figures~\ref{fig:23.8.24 R1 A-PP F(1J,2a)3 lr0.04,0.03,0.02,0.01,0.01,0.005 cy1-6-VCA cy7-9-NCA Nsteps400000 random N51 W32 H2 He init-1}-\ref{fig:23.8.24 R1 A-PP F(1J,2a)3 lr0.04,0.03,0.02,0.01,0.01,0.005 cy1-6-VCA cy7-9-NCA Nsteps400000 random N51 W32 H2 He init-3}, \red{Form 3}, \red{1,251} parameters, random grid, \red{\emph{quasi-perfect}} solution.
		}
		{\scriptsize (23817R10b-1)}
	}
	\label{fig:23.8.17 R10b A-PP F(1)2a lr0.07 cy2-NCA Nsteps50000 regular N51 W32 H2 He init}
\end{figure}

% 23.9.20, OLD Fig.11 to move to section Form 2a
\begin{figure}[tph]
%	\includegraphics[width=0.49\textwidth]{Figures/23.8.15_R3b_A-PP_F_1_2a_lr0.03_cy5-CA_Nsteps200000_regular_N51_W64_H2_He_init_-_loss200000.png}
%	%_\includegraphics[width=0.49\textwidth]{Figures/23.8.15_R3b_A-PP_F_1_2a_lr0.03_cy5-CA_Nsteps200000_regular_N51_W64_H2_He_init_-_shape200000.png}
%	\includegraphics[width=0.49\textwidth]{Figures/23.8.15_R3b_A-PP_F_1_2a_lr0.03_cy5-CA_Nsteps200000_regular_N51_W64_H2_He_init_-_center_disp200000.png}
	%
	\includegraphics[width=0.49\textwidth]{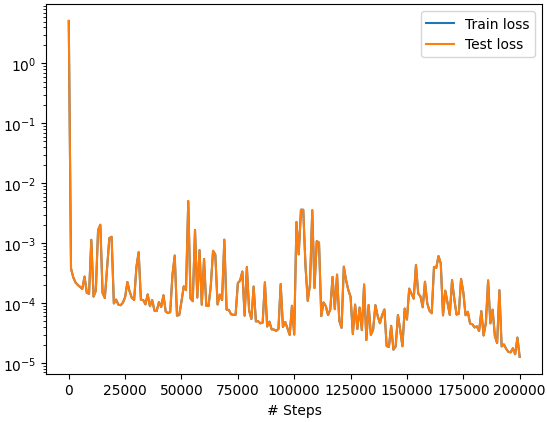}
	\includegraphics[width=0.49\textwidth]{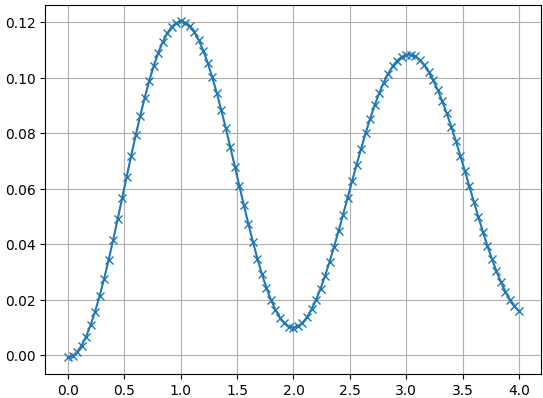}
	\caption{
		% CA led to waves with damping.
		% In 23.8.13 R4ab, 1218 params, NCA led to waves with damping, whereas CA led to the static solution.  In 23.8.15 R3ab, 4482 params, it's the opposite: NCA led to the static solution, whereas CA led to waves with damping ! 
		%
		\DDET.
		\emph{Pinned-pinned bar, cyclic annealing} (CA). 
		Loss function (left), midspan displacement, Step 200,000 (right), \emph{waves with damping}.
		\red{$\star$}
		Section~\ref{sc:pinned-pinned-Form-2a}, \red{Form 2a}.
		\emph{Network:}
		\red{$\star$}
		Remarks~\ref{rm:parameter-names},~\ref{rm:data-point-grids}, 
		T=4,
		\red{W=64},
		H=2, 
		n\_out=2, 
		He-uniform initializer, 
		\red{4,482} parameters,
		\red{$\star$}
		\red{\emph{regular}} grid.
		\emph{Training:} 
		\red{$\star$}
		Remark~\ref{rm:learning-rate-schedule-1},
		\red{LRS~1 (CA)},
		init\_lr=0.03.
		$\bullet$
		{\footnotesize
			Figure~\ref{fig:23.8.13 R4b A-PP F(1)2a lr0.04 cy5-CA Nsteps200000 regular N51 W32 H2 He init}, 1,218 parameters, CA, init\_lr=0.04, static solution.
			$\triangleright$
			Figure~\ref{fig:axial-Mathematica-solutions}, reference solution to compare.
		}
		% \red{I AM HERE 2023.08.29.}
		{\scriptsize (23815R3b-1)}
		$\bullet$
		{\footnotesize
			{The train-loss history in blue coincides with the test-loss history in orange, and is not visible in the plot.  The same situation occurs in subsequent loss-history figures.}
		}
	}
	\label{fig:23.8.15 R3b A-PP F(1)2a lr0.03 cy5-CA Nsteps200000 regular N51 W64 H2 He init}
\end{figure}

% 23.9.20, OLD Fig.12 to move to section Form 2a
\begin{figure}[tph]
	\includegraphics[width=0.49\textwidth]{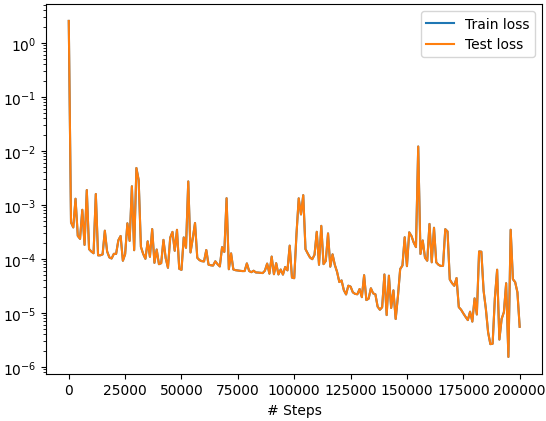}
	\includegraphics[width=0.49\textwidth]{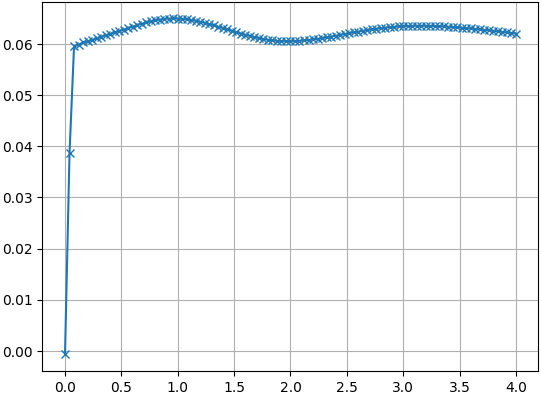}
	\caption{
		% CA led to waves with damping.
		% In 23.8.13 R4ab, 1218 params, NCA led to waves with damping, whereas CA led to the static solution.  In 23.8.15 R3ab, 4482 params, it's the opposite: NCA led to the static solution, whereas CA led to waves with damping ! 
		%
		\DDET.
		\emph{Pinned-pinned bar, cyclic annealing} (CA). 
		Loss function (left), midspan displacement, Step 200,000 (right), \emph{static solution}.
		\red{$\star$}
		Section~\ref{sc:pinned-pinned-Form-2a}, \red{Form 2a}.
		\emph{Network:}
		\red{$\star$}
		Remarks~\ref{rm:parameter-names},~\ref{rm:data-point-grids}, 
		T=4,
		\red{W=32},
		H=2, 
		n\_out=2, 
		He-uniform initializer, 
		\red{1,218} parameters,
		\red{$\star$}
		\mbox{\red{\emph{regular}}} grid.
		\emph{Training:} 
		\red{$\star$}
		Remark~\ref{rm:learning-rate-schedule-1},
		% Section~\ref{sc:optimization-learning-rate-scheduling}.
		\red{LRS~1 (CA)},
		init\_lr=0.04.
		$\bullet$
		{\footnotesize
			Figure~\ref{fig:23.8.15 R3b A-PP F(1)2a lr0.03 cy5-CA Nsteps200000 regular N51 W64 H2 He init}, 4,482 parameters, CA, init\_lr=0.03, waves, damping.
			Figure~\ref{fig:23.8.21 R1 A-PP F(1)2a lr0.04,0.03,0.02,0.01,0.005 cy5-VCA Nsteps200000 regular N51 W32 H2 He init}, VCA, Form 2a, same model, waves, small damping.
		}
		{\scriptsize (23813R4b-1)}
	}
	\label{fig:23.8.13 R4b A-PP F(1)2a lr0.04 cy5-CA Nsteps200000 regular N51 W32 H2 He init}
\end{figure}

\begin{rem}
	% \label{rm:learning-rate-decay-1}
	\label{rm:learning-rate-schedule-1}
	Learning-rate schedule 1
	{\rm 
		(LRS 1).
		\emph{Cyclic annealing} (CA) with \emph{same} learning rate.
		To allow for a systematic comparison of the PINN results, LRS 1 is set to consist of $n_\text{cycles} = 5$ cycles, with Cycles 1 and 2 having $n\text{-steps}_\text{cycle} = 25,000$ steps, totaling $50,000$ steps at the end of cycle 2, whereas Cycles 3, 4, 5, each has $n\text{-steps}_\text{cycle} = 50,000$ steps, making a total of $N_\text{steps} = 200,000$ steps over these 5 cycles.  Computational results (loss function, deformed shape, GPU time) would be gathered at the end of these 5 cycles.  
		Each period is kept at a constant $n\text{-steps}_\text{period} = 2,500$ steps, and has a learning-rate decay to $\decayf =$ 90\% 
		% ($\decayf = 0.9$) 
		by inverse time.\footnote{
			The learning rate at the end of a period is equal to 90\% of the learning rate at the beginning of that period.
		}
		
		Any deviation from the LRS 1 would be clearly indicated.  
		For example, we may run a case using only two cycles, n-cycles=2, init\_lr=0.07, and then stop the execution after N\_steps=50,000, such as in
		Figure~\ref{fig:23.8.17 R10a A-PP F(1)2a lr0.07 cy2-CA Nsteps50000 regular N51 W32 H2 He init}. 
		% {\color{red} FIGURE.  NOTE: 2023.08.20, NO, Figure~\ref{fig:23.7.22 R1 A-PP v1 midspan,shape100000} was obtained with NCA, not with CA.  Find a CA figure}. 
		
		Using \DDET, while CA led to waves with damping as in
		Figure~\ref{fig:23.8.15 R3b A-PP F(1)2a lr0.03 cy5-CA Nsteps200000 regular N51 W64 H2 He init} (4,482 parameters, init\_lr=0.03, regular grid), but CA could also lead to a static solution as in Figure~\ref{fig:23.8.17 R10a A-PP F(1)2a lr0.07 cy2-CA Nsteps50000 regular N51 W32 H2 He init} and Figure~\ref{fig:23.8.13 R4b A-PP F(1)2a lr0.04 cy5-CA Nsteps200000 regular N51 W32 H2 He init} (1,218 parameters, init\_lr=0.04, regular grid). 
		For the counterpart with ``no cyclic annealing'' (NCA), see Remark~\ref{rm:learning-rate-schedule-3} on LRS~3.
		
		It is important to note that the static solution is one of the pathological problems we encountered with using 
		% our \DDET\ script
		\DDET\ 
		(Figure~\ref{fig:DDE-T-Form-1-pinned-free-bar-static-solutions}), but did not appear when using 
		% our \JAX\ script
		\JAX\ 
		(Figure~\ref{fig:JAX-Form-1-pinned-free-bar-NO-static-solutions}); see Remarks~\ref{rm:learning-rate-schedule-4}, \ref{rm:static-solution-Form-1-DDE-T-v-JAX}.
		
		A characteristic of CA is the jump in the loss function after each resetting of the initial learning rate at the beginning of each cycle; see the left subfigures of Figures~\ref{fig:23.8.15 R3b A-PP F(1)2a lr0.03 cy5-CA Nsteps200000 regular N51 W64 H2 He init}-\ref{fig:23.8.13 R4b A-PP F(1)2a lr0.04 cy5-CA Nsteps200000 regular N51 W32 H2 He init}.
	}
	\hfill
	$\blacksquare$
\end{rem}

% 23.9.20, OLD Fig.13 to move to section Form 2a
\begin{figure}[tph]
	\includegraphics[width=0.49\textwidth]{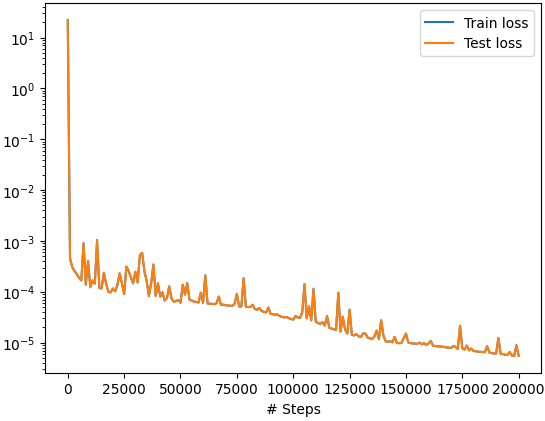}
	\includegraphics[width=0.49\textwidth]{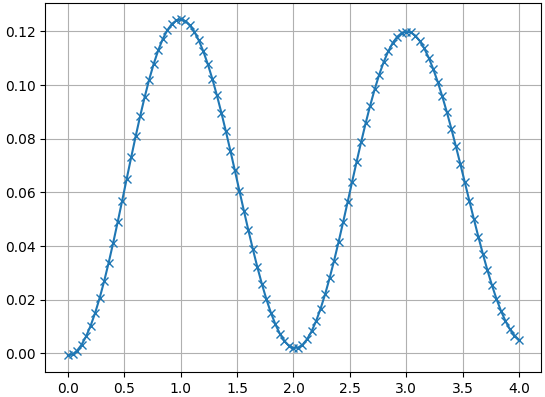}
	\caption{
		% CA led to waves with damping.
		% In 23.8.13 R4ab, 1218 params, NCA led to waves with damping, whereas CA led to the static solution.  In 23.8.15 R3ab, 4482 params, it's the opposite: NCA led to the static solution, whereas CA led to waves with damping ! 
		%
		% \red{I AM HERE, 23.9.7. Rewrite below.}
		\DDET.
		\emph{Pinned-pinned bar, Varying init\_lr cyclic annealing} (VCA). 
		Loss function (left), midspan displacement (right), Step 200,000, \emph{waves, small damping}.
		Section~\ref{sc:pinned-pinned-Form-2a}, {\color{red} Form~2a}.
		\red{$\star$}
		Remark~\ref{rm:damping-low-capacity-models}, \emph{Damping}.
		\emph{Network:}
		Remarks~\ref{rm:parameter-names},~\ref{rm:data-point-grids}, 
		T=4,
		W=32,
		H=2, 
		n\_out=2, 
		He-uniform initializer, 
		\red{1,218} parameters,
		\emph{regular} grid.
		\emph{Training:} 
		\red{$\star$}
		Remark~\ref{rm:learning-rate-schedule-2},
		\red{LRS~2 (VCA)},
		init\_lr=[0.04, 0.03, 0.02, 0.01, 0.005].
		$\bullet$
		{\footnotesize
			Figure~\ref{fig:23.8.13 R4b A-PP F(1)2a lr0.04 cy5-CA Nsteps200000 regular N51 W32 H2 He init}, CA, Form 2a, same model, init\_lr=0.04, static solution.
			$\triangleright$
			Figure~\ref{fig:23.8.21 R1c A-PP F(1J,2a)3 lr0.04,0.03,0.02,0.01,0.005 cy5-VCA Nsteps200000 regular N51 W32 H2 He init-1}, VCA, Form 3, same model and init\_lr sequence, waves, small damping.
			$\triangleright$
			Figure~\ref{fig:axial-Mathematica-solutions}, reference solution to compare.
		}
		{\scriptsize (23821R1a-1)}
	}
	\label{fig:23.8.21 R1 A-PP F(1)2a lr0.04,0.03,0.02,0.01,0.005 cy5-VCA Nsteps200000 regular N51 W32 H2 He init}
\end{figure}

\begin{rem}
	\label{rm:learning-rate-schedule-2}
	Learning-rate schedule 2
	{\rm
		(LRS 2).
		\emph{Varying initial-learning-rate cyclic annealing} (VCA).
		% with \emph{different} learning rates.
		Each cycle {has} its own learning rate, {which may (or may not) vary in the subsequent cycles}.
		An example would be to use the same parameters as in LRS~1 (Remark~\ref{rm:learning-rate-schedule-1}), except for the different initial learning rate at the beginning of each cycle as prescribed by the list init\_lr = [0.02, 0.02, 0.02, 0.01, 0.005], with decreasing initial learning rates, such that init\_lr = 0.02 is used for Cycle 1, 2, 3, init\_lr = 0.01 for Cycle 4, and init\_lr = 0.005 for Cycle 5.
		
		Figure~\ref{fig:23.8.21 R1 A-PP F(1)2a lr0.04,0.03,0.02,0.01,0.005 cy5-VCA Nsteps200000 regular N51 W32 H2 He init}, VCA with Form 2a, 
		Figures~\ref{fig:23.8.21 R1c A-PP F(1J,2a)3 lr0.04,0.03,0.02,0.01,0.005 cy5-VCA Nsteps200000 regular N51 W32 H2 He init-1}-\ref{fig:23.8.21 R1c A-PP F(1J,2a)3 lr0.04,0.03,0.02,0.01,0.005 cy5-VCA Nsteps200000 regular N51 W32 H2 He init-2}, VCA with Form 3,
		the same model as in Figure~\ref{fig:23.8.13 R4b A-PP F(1)2a lr0.04 cy5-CA Nsteps200000 regular N51 W32 H2 He init} (N51 W32 H2, 1,218 parameters, regular grid), and init\_lr sequence [0.04, 0.03, 0.02, 0.01, 0.005] led to waves with small damping.
	}
	\hfill
	$\blacksquare$
\end{rem}

% 23.9.20, OLD Fig.16 to move to section Form 2a
\begin{figure}[tph]
	\includegraphics[width=0.49\textwidth]{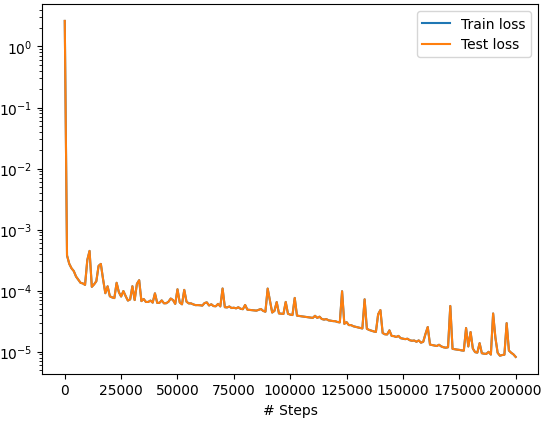}
	\includegraphics[width=0.49\textwidth]{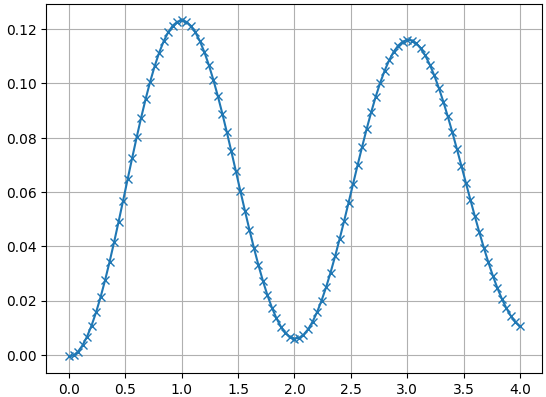}
	\caption{
		% CA led to waves with damping.
		% In 23.8.13 R4ab, 1218 params, NCA led to waves with damping, whereas CA led to the static solution.  In 23.8.15 R3ab, 4482 params, it's the opposite: NCA led to the static solution, whereas CA led to waves with damping ! 
		%
		% \red{I AM HERE, 23.9.8. Rewrite below.}
		\DDET.
		\emph{Pinned-pinned bar, NO cyclic annealing} (NCA).
		Loss function (left), midspan displacement, Step 200,000 (right), \emph{waves with damping}.
		\red{$\star$}
		Section~\ref{sc:pinned-pinned-Form-2a}, \red{Form 2a}.
		\emph{Network:}
		\red{$\star$}
		Remarks~\ref{rm:parameter-names},~\ref{rm:data-point-grids}, 
		T=4,
		W=32,
		H=2, 
		n\_out=2, 
		He-uniform initializer, 
		\red{1,218} parameters,
		\red{$\star$}
		\red{\emph{regular}} grid.
		\emph{Training:} 
		\red{$\star$}
		Remark~\ref{rm:learning-rate-schedule-3},
		\red{LRS~3 (NCA)},
		init\_lr=0.04.
		$\bullet$
		{\footnotesize
			Figure~\ref{fig:23.8.13 R4b A-PP F(1)2a lr0.04 cy5-CA Nsteps200000 regular N51 W32 H2 He init}, CA, same model, static solution.
			% Figure~\ref{fig:23.8.15 R3b A-PP F(1)2a lr0.03 cy5-CA Nsteps200000 regular N51 W64 H2 He init}, CA led to waves with damping.
			Figure~\ref{fig:23.8.15 R3a A-PP F(1)2a lr0.03 cy5-NCA Nsteps200000 regular N51 W64 H2 He init-1}, NCA, model with 4,482 parameters, static solution.
			$\triangleright$
			Figure~\ref{fig:axial-Mathematica-solutions}, reference solution to compare.
		}
		{\scriptsize (23813R4a-1)}
	}
	\label{fig:23.8.13 R4a A-PP F(1)2a lr0.04 cy5-NCA Nsteps200000 regular N51 W32 H2 He init}
\end{figure}

\begin{rem}
	\label{rm:learning-rate-schedule-3}
	Learning-rate schedule 3
	{\rm
		(LRS 3). \emph{No cyclic annealing} (NCA), i.e., no resetting of the learning rate at the beginning of each cycle, unlike LRS 1 in Remark~\ref{rm:learning-rate-schedule-1}.  
		The number of cycles are just breakpoints to plot intermediate results.
		For example, a training with n-cycles=1 and N\_steps=200,000 correspond to one breakpoint at Step 200,000, whereas a training with n-cycles=5 and N\_steps=200,000 has five breakpoints.
		The initial learning rate init\_lr set at the beginning of Cycle 1 decreases with inverse-time decay without being reset to init\_lr at the beginning of each subsequent cycle---unlike LRS~1 in Remark~\ref{rm:learning-rate-schedule-1} and LRS~2 in Remark~\ref{rm:learning-rate-schedule-2}---n-steps\_period=2500, $\decayf = 0.9$, and N\_steps set to either 25,000, 50,000, 100,000, or 200,000 steps for standardization.  
		%In other words, the learning rate is left to continue to decay over N\_steps without being reset to the initial learning rate at the beginning of each cycle, as done in LRS~1 or LRS~2.
		
		Figure~\ref{fig:23.8.17 R10b A-PP F(1)2a lr0.07 cy2-NCA Nsteps50000 regular N51 W32 H2 He init},
		model with N51 W32 H2, having 1,218 parameters,
		shows waves with damping in the midspan displacement at Step 50,000 when using NCA with init\_lr=0.07.
		
		Figure~\ref{fig:23.8.13 R4a A-PP F(1)2a lr0.04 cy5-NCA Nsteps200000 regular N51 W32 H2 He init}, NCA with init\_lr=0.04 and a model having 1,218 parameters led to waves with damping.
		Figure~\ref{fig:23.8.15 R3a A-PP F(1)2a lr0.03 cy5-NCA Nsteps200000 regular N51 W64 H2 He init-1}, NCA with init\_lr=0.03 and a model having 4,482 parameters led to static solution. 
		To avoid the static solution, reduce the initial learning rate init\_lr for the same model capacity, or reduce the model capacity (with a possible increase in initial learning rate as done in Figure~\ref{fig:23.8.13 R4a A-PP F(1)2a lr0.04 cy5-NCA Nsteps200000 regular N51 W32 H2 He init}); see Remark~\ref{rm:optimal-capacity}.
		% {\color{red} [NOTE: 2023.08.20, add the above figures to the remark on the relationship between model capacity and learning rate. ENDNOTE]}
	}
	\hfill
	$\blacksquare$
\end{rem}

% 23.9.20, OLD Fig.17 to move to section Form 2a
\begin{figure}[tph]
	\includegraphics[width=0.49\textwidth]{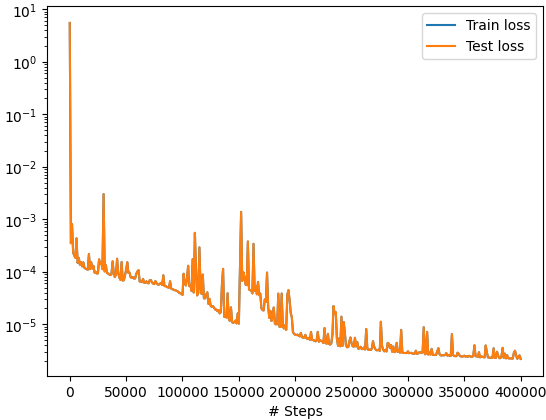}
	\includegraphics[width=0.49\textwidth]{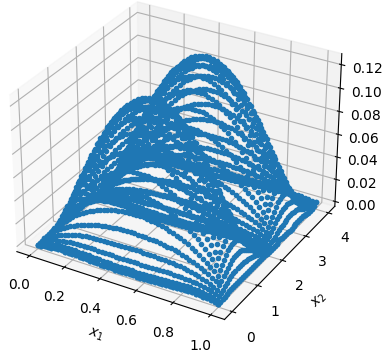}
	\caption{
		% \emph{Pinned-pinned bar, Cycle-5 extension}. 
		% \emph{Pinned-pinned bar, ELRS}. 
		% \red{I AM HERE, 23.9.8. Rewrite below.}
		\DDET.
		\emph{Pinned-pinned bar},
		\emph{Extended learning-rate schedule} (ELRS).
		Loss function (left), 
		shape, Step 400,000 (right).
		\red{$\star$}
		Section~\ref{sc:pinned-pinned-Form-2a}, \red{Form 2a}.
		\emph{Network:}
		\red{$\star$}
		Remarks~\ref{rm:parameter-names},~\ref{rm:data-point-grids}, 
		T=4,
		W=32,
		H=2, 
		n\_out=2, 
		He-uniform initializer, 
		\red{1,218} parameters,
		\red{$\star$}
		\red{\emph{regular}} grid.
		\emph{Training:} 
		\red{$\star$}
		Remarks~\ref{rm:learning-rate-schedule-1}
		(LRS~1),
		~\ref{rm:extension-cycle-5} \red{(ELRS)},
		init\_lr=0.02,
		% n-cycles=9,
		Cycles 1-5 (CA), 
		Cycles 6-9 (NCA),
		N\_steps=400,000.
		\red{$\star$}
		Lowest total loss \red{2.19e-06}, Step 400,000 (sum of 6 losses).
		% {\color{red} HERE}
		$\bullet$ 
		{\footnotesize
			% Figure~\ref{fig:23.7.26 R5c A-PP F2a lr0.02x5 cycles=9 Nsteps=400000 regular N51 W32 H2 He init-1}, loss function and shape time history at Step 400,000. 
			Figure~\ref{fig:23.7.26 R5c A-PP F2a lr0.02x5 cycles=9 Nsteps=400000 regular N51 W32 H2 He init-2}, midspan displacements, Steps 200,000 \& 300,000. 
			Figure~\ref{fig:23.7.26 R5c A-PP F2a lr0.02x5 cycles=9 Nsteps=400000 regular N51 W32 H2 He init-3}, velocity, \red{\emph{very-good}} midspan displacement, Step 400,000.
			$\triangleright$
			Figure~\ref{fig:axial-Mathematica-solutions}, reference solution to compare.
		}
		{\scriptsize (23726R5c-1)}
	}
	\label{fig:23.7.26 R5c A-PP F2a lr0.02x5 cycles=9 Nsteps=400000 regular N51 W32 H2 He init-1}
\end{figure}

% 23.9.20, OLD Fig.18 to move to section Form 2a
\begin{figure}[tph]
	\includegraphics[width=0.49\textwidth]{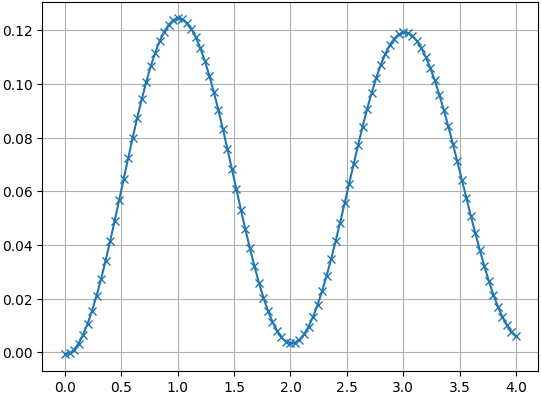}
	\includegraphics[width=0.49\textwidth]{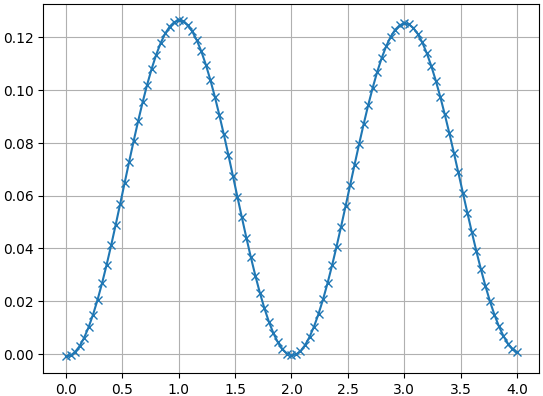}
	\caption{
		% \emph{Pinned-pinned bar, Cycle-5 extension}.
		% \red{I AM HERE, 23.9.8. Rewrite below.} 
		\DDET.
		\emph{Pinned-pinned bar, ELRS}. 
		Midspan displacements: Step 200,000 (left), Step 300,000 (right).
		Section~\ref{sc:pinned-pinned-Form-2a}, \red{Form 2a}.
		% \red{REWRITE below.}
		\emph{Network:}
		\red{$\star$}
		Remarks~\ref{rm:parameter-names},~\ref{rm:data-point-grids}, 
		T=4,
		W=32,
		H=2, 
		n\_out=2, 
		He-uniform initializer, 
		\red{1,218} parameters,
		\red{$\star$}
		\red{\emph{regular}} grid.
		\emph{Training:} 
		\red{$\star$}
		Remarks~\ref{rm:learning-rate-schedule-1}
		(LRS~1),
		~\ref{rm:extension-cycle-5} \red{(ELRS)},
		init\_lr=0.02,
		n-cycles=9,
		Cycles 1-5 (CA), 
		Cycles 6-9 (NCA),
		N\_steps=400,000.
		% {\color{red} HERE}
		$\bullet$ 
		{\footnotesize
			Figure~\ref{fig:23.7.26 R5c A-PP F2a lr0.02x5 cycles=9 Nsteps=400000 regular N51 W32 H2 He init-1}, loss function, shape, Step 400,000.
			%, lowest total loss. 
			% Figure~\ref{fig:23.7.26 R5c A-PP F2a lr0.02x5 cycles=9 Nsteps=400000 regular N51 W32 H2 He init-2}, midspan-displacement time histories at Step 200,000 and at Step 300,000. 
			Figure~\ref{fig:23.7.26 R5c A-PP F2a lr0.02x5 cycles=9 Nsteps=400000 regular N51 W32 H2 He init-3}, velocity, \red{\emph{very-good}} solution, Step 400,000.
		}
		{\scriptsize (23726R5c-2)}
	}
	\label{fig:23.7.26 R5c A-PP F2a lr0.02x5 cycles=9 Nsteps=400000 regular N51 W32 H2 He init-2}
\end{figure}

% 23.9.24, reran this case with results in 23.7.26 R5d
% DEBUG 23.7.26 R5d DDE A-PP F2a lr0.02 cy1-5-CA cy6-9-NCA Nsteps400000 regular N51 W32 H2 He init - center disp400000.png
% 23.9.20, OLD Fig.19 to move to section Form 2a
\begin{figure}[tph]
	\includegraphics[width=0.49\textwidth]{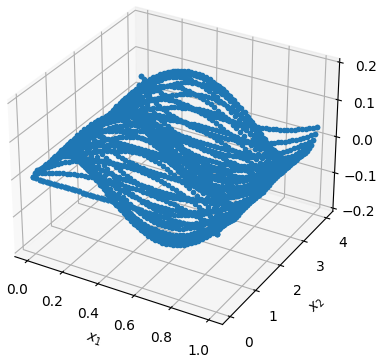}
	\includegraphics[width=0.49\textwidth]{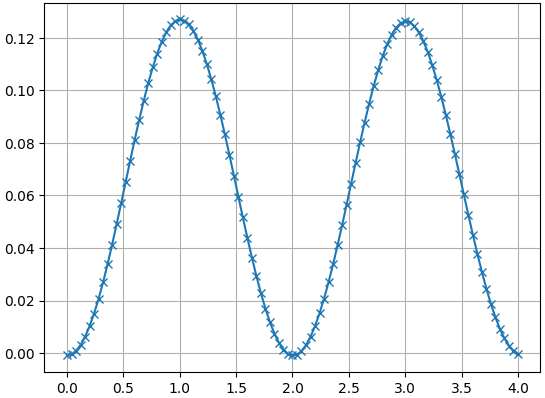}
	\caption{
		% \emph{Pinned-pinned bar, Cycle-5 extension}. 
		% \red{I AM HERE, 23.9.8. Rewrite below.}
		\DDET.
		\emph{Pinned-pinned bar, ELRS}.
		Velocity (left) and \red{\emph{very-good}} midspan displacement (right), Step 400,000.
		\red{$\star$}
		Section~\ref{sc:pinned-pinned-Form-2a}, \red{Form 2a}.
		% {\color{red} 2023.08.21, Write the section for Form 3.}
		\emph{Network:}
		\red{$\star$}
		Remarks~\ref{rm:parameter-names},~\ref{rm:data-point-grids},
		T=4, 
		W=32,
		H=2, 
		n\_out=2, 
		He-uniform initializer, 
		\red{1,218} parameters,
		\red{\emph{regular}} grid.
		\emph{Training:} 
		\red{$\star$}
		Remarks~\ref{rm:learning-rate-schedule-1},
		(LRS~1),
		~\ref{rm:extension-cycle-5} \red{(ELRS)},
		init\_lr=0.02,
		n-cycles=9,
		Cycles 1-5 (CA), 
		Cycles 6-9 (NCA),
		N\_steps=400,000.
		% Lowest total loss \red{2.19e-06}, Step 400,000 (sum of 6 losses).
		% {\color{red} HERE}
		$\bullet$ 
		{\footnotesize
			Figure~\ref{fig:23.7.26 R5c A-PP F2a lr0.02x5 cycles=9 Nsteps=400000 regular N51 W32 H2 He init-1}, loss function, shape, Step 400,000. 
			Figure~\ref{fig:23.7.26 R5c A-PP F2a lr0.02x5 cycles=9 Nsteps=400000 regular N51 W32 H2 He init-2}, midspan displacements, Steps 200,000 \&  300,000. 
			Figure~\ref{fig:random-grid-NO-static-solution}, Form 2a, 12,802 parameters, random grid, \emph{quasi-perfect} solution.
			% Figure~\ref{fig:23.7.26 R5c A-PP F2a lr0.02x5 cycles=9 Nsteps=400000 regular N51 W32 H2 He init-3}, velocity and midspan-displacement time histories at Step 400,000.
			$\triangleright$
			Figure~\ref{fig:axial-Mathematica-solutions}, reference solution to compare.
		}
		{\scriptsize (23726R5c-3)}
	}
	\label{fig:23.7.26 R5c A-PP F2a lr0.02x5 cycles=9 Nsteps=400000 regular N51 W32 H2 He init-3}
\end{figure}

\begin{rem}
	\label{rm:extension-cycle-5}
	% Extension of Cycle 5.
	Extended learning-rate schedule
	{\rm
		(ELRS).
		In some examples, N\_steps may be extended beyond 200,000, such as 400,000.  The extension could be in CA mode (Remarks~\ref{rm:learning-rate-schedule-1}-\ref{rm:learning-rate-schedule-2}) or in NCA mode (Remark~\ref{rm:learning-rate-schedule-3}). 
		
		In NCA mode, the learning rate simply continues to decay from the end of Cycle 5 without being reset. Doing so is equivalent to extend n-steps\_cycle for Cycle 5 from 50,000 steps to 250,000 steps.
		This Cycle 5 extension is applicable to both LRS 1 in Remark~\ref{rm:learning-rate-schedule-1} and LRS 2 in Remark~\ref{rm:learning-rate-schedule-2}.  
		% ALEX: removed paragraph
		An example of Cycle-5 extension for Form 2a and a low-capacity model N51 W32 H2 with 1,218 parameters is given in 
		Figure~\ref{fig:23.7.26 R5c A-PP F2a lr0.02x5 cycles=9 Nsteps=400000 regular N51 W32 H2 He init-1}, loss function and shape time history at Step 400,000; 
		Figure~\ref{fig:23.7.26 R5c A-PP F2a lr0.02x5 cycles=9 Nsteps=400000 regular N51 W32 H2 He init-2}, midspan-displacement time histories at Step 200,000 and at Step 300,000; 
		Figure~\ref{fig:23.7.26 R5c A-PP F2a lr0.02x5 cycles=9 Nsteps=400000 regular N51 W32 H2 He init-3}, velocity and \emph{quasi-perfect} midspan-displace\-ment time histories at Step 400,000.
		The loss function in Figure~\ref{fig:23.7.26 R5c A-PP F2a lr0.02x5 cycles=9 Nsteps=400000 regular N51 W32 H2 He init-1}, together with the midspan-displacement time histories in Figures~\ref{fig:23.7.26 R5c A-PP F2a lr0.02x5 cycles=9 Nsteps=400000 regular N51 W32 H2 He init-2}-\ref{fig:23.7.26 R5c A-PP F2a lr0.02x5 cycles=9 Nsteps=400000 regular N51 W32 H2 He init-3} gives an idea of the gradual improvement of the solution during the optimization process.
		The smallest total loss at Step 400,000 was 2.19e-06, which was the sum of six losses: PDE, momentum, two boundary conditions, and two initial conditions.
		
		% 23.8.24 R1 (random) vs 23.8.23 R6 (regular)
		% compare the Loss value F3 400000 = 1.25e-06, not SMALLEST on 23.8.24 from the script 23.8.24 R1 A-PP F(1J,2a)3 lr0.04,0.03,0.02,0.01,0.01,0.005 cy1-6-VCA cy7-9-NCA Nsteps400000 random N51 W32 H2 He init.ipynb using a random grid to the Loss value F3 400000 = 2.56e-06, 8.23 R6  from the script 23.8.23 R6 A-PP F(1J,2a)3 lr0.04,0.03,0.02,0.01,0.01,0.005 cy1-6-VCA cy7-9-NCA Nsteps400000 regular N51 W32 H2 He init.ipynb using a regular grid.
		Switching from regular grid to \emph{random} grid, maintaining all other parameters the same, could lower the total loss.
		For example, using Form 3 with the model N51 W32 H2 having 1,251 params, and init\_lr sequence [0.04,0.03,0.02,0.01,0.01,0.005] for Cycles 1-6 (VCA), and then Cycles 7-9 (NCA), using a regular grid yielded a total loss of 2.56e-06, the sum of seven losses for Form 3, whereas using a random grid yielded a smaller total loss of 1.25e-06, both at the same last Step 400,000.    
	}
	\hfill
	$\blacksquare$
\end{rem}

\begin{rem}
	\label{rm:learning-rate-schedule-piecewise-constant}
	\label{rm:learning-rate-schedule-4}
	Piecewise-constant learning-rate schedule.
	{\rm
		(LRS~4)
		While the inverse-time decay is used in 
		LRS~1 (Remark~\ref{rm:learning-rate-schedule-1}),
		LRS~2 (CA, Remark~\ref{rm:learning-rate-schedule-2}), 
		LRS~3 (NCA, Remark~\ref{rm:learning-rate-schedule-3}),
		in our PINN script written using the 
		% \href{https://jax.readthedocs.io/en/latest/index.html}{JAX} 
		\JAX\ framework, 
		% the \emph{High-Performance Array Computing} by GoogleMind,
		due to time constraint,\footnote{
			The main reason was the time constraint due to the approaching deadline for this special issue, as we only started developing our JAX script after spending a significant amount of time with DDE-T without satisfactory results, and had our JAX script in working order as the deadline was already looming close.   There is no technical reason for not implementing the inverse time decay or any other time decay methods.
		}
		we only implemented the
		\href{https://optax.readthedocs.io/en/latest/api.html#optax.piecewise_constant_schedule}{piecewise-constant schedule} of the 
		\href{https://optax.readthedocs.io/en/latest/index.html}{Optax} gradient processing and optimization library for JAX.
		
		An example of LRS~4 would be: init\_lr=0.003, with factor\_lr=[0.9, 0.8, 0.7, 1, 1], so that the learning rate decays as follows:
		0.003 in Cycle 1 (up to Step 25,000),
		(0.003 * 0.9) = 2.7e-03 in Cycle 2 (from Step 25,000+1 to 50,000),
		(2.7e-03 * 0.8) = 2.16e-03 in Cycle 3 (from 50,000+1 to 100,000),
		(2.16e-03 * 0.7) = 1.512e-03 in Cycle 4 (from 100,000+1 to 150,000),
		1.512e-03 in Cycle 5 (from 150,000+1 to 200,000),	
		1.512e-03 in Cycle 6 (from 200,000+1 to 250,000).	
		
		Hence LRS~4 is also a NCA schedule as LRS~3 (Remark~\ref{rm:learning-rate-schedule-3}), and can be extended beyond 200,000 steps as in the above example. 
		Figure~\ref{fig:JAX-Form-1-pinned-free-bar-NO-static-solutions} depicts the motion of a pinned-free bar under constant distributed axial load (Figure~\ref{fig:axial-Mathematica-solutions}) obtained with JAX and LRS~4.
		See also Remark~\ref{rm:static-solution-Form-1-DDE-T-v-JAX}, \DvJ, which exhibits none of the \DDET\ pathological problems.
	}
	\hfill$\blacksquare$
\end{rem}

\begin{rem}
	\label{rm:He-vs-Glorot}
	Initializer: ``He uniform'' vs. ``Glorot uniform.''
	{\rm
		% \red{NOTE: 23.9.6, Need to rewrite. ENDNOTE}
		The \emph{He-uniform} initializer is equivalent to a larger learning rate, compared to the \emph{Glorot-uniform} initializer.
		As a result, the He-uniform initializer may help the training converge faster, but could push the iterate to a static solution, depending on the PDE and its form.
		
		Consider 
		Form 2a of a pinned-free bar, using a network with N51 W64 H4, 12,802 parameters, and the learning-rate scheduling LRS~3 (no cyclic annealing), and init\_lr=0.005.
		Figure~\ref{fig:DEBUG v1.9.3 seed42, 23.9.5 R3a A-P bcrF F(1J)2a(3) lr0.005 cy1-3-NCA Nsteps100000 regular N51 W64 H4 He init-1} uses the \emph{He-uniform} initializer, with {regular} grid, yielded a \emph{static} solution.
		Switching from regular grid to {random} grid (Remark~\ref{rm:static-solution-avoid}), keeping the same network architecture and number of data points, i.e., 12,802 parameters, 
		Figure~\ref{fig:23.9.5 R3b A-P bcrF F(1J)2a(3) lr0.005 cy1-5-NCA Nsteps200000 random N51 W64 H4 He init-1}, also using the \emph{He-uniform} initializer, yielded a \emph{quasi-perfect} solution. 
		% See Remark~\ref{rm:static-solution-avoid} on how to avoid static solution.
		
		From the pre-static solution of a pinned-free bar shown in SubFigs~\ref{fig:23.9.5 R3d A-PF free-end disp197000}-\ref{fig:23.9.5 R3d A-PF free-end slope197000} obtained with Form 3, He-uniform initializer, and initial learning rate init\_lr=0.005, simply 
		switching from He-uniform to Glorot-uniform initializer, which is equivalent to reducing the learning rate (even though init\_lr remained at 0.005), resulted in a very-good free-end displacement shown in SubFigs~\ref{fig:23.9.5 R3d.3 A-PF shape184000}-\ref{fig:23.9.5 R3d.3 A-PF free-end disp184000}.
	}
	\hfill$\blacksquare$
\end{rem}

\begin{rem}
	\label{rm:rule-of-thumb-learning-rate}
	Rule of thumb for learning rate.
	{\rm
		A quick way to assess whether a learning rate would lead to the convergence to a solution (even though the converged solution may not be the one desired), is to run a script with \emph{zero applied force} to see how the optimization converges to the \emph{zero solution}.
	}
	\hfill
	$\blacksquare$
\end{rem}

% 23.9.20, OLD Fig.20 to move to section Form 2a
\begin{figure}[tph]
	\includegraphics[width=0.49\textwidth]{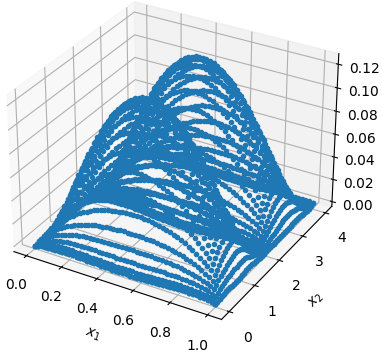}
	\includegraphics[width=0.49\textwidth]{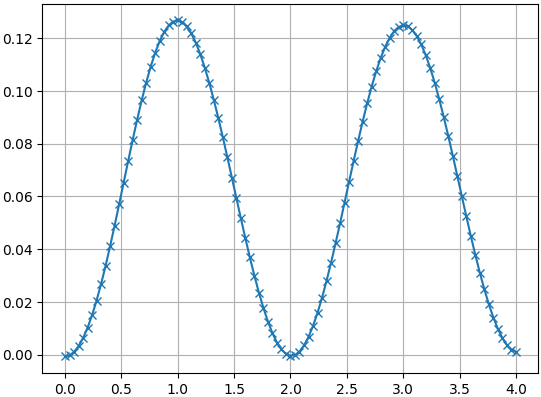}
	\caption{
		% Low-capacity network N51 W64 H2, 4,482 parameters.
		% \red{I AM HERE, 23.9.8. Rewrite below.}
		\DDET.
		\emph{Pinned-pinned bar, lower-capacity network 1}.
		Shape (left), midspan displacement (right), Step 200,000.
		\red{$\star$} 
		Remark~\ref{rm:optimal-capacity}, \red{\emph{Optimal capacity}}.
		Remark~\ref{rm:unstable-solution}, \red{\emph{Unstable solution}}.
		% Section~\ref{sc:networks-grids}. 
		Section~\ref{sc:pinned-pinned-Form-2a}, \red{Form 2a}.
		% Computational domain: $(x, t) \in [0, 1] \times [0, 4]$.
		\emph{Network:}
		\red{$\star$}
		Remarks~\ref{rm:parameter-names},~\ref{rm:data-point-grids},
		T=4,
		W=64,
		H=2, 
		n\_out=2, 
		He-uniform initializer, 
		\red{4,482} parameters,
		\red{$\star$}
		\emph{\emph{regular}} grid.
		\emph{Training:}
		Remark~\ref{rm:learning-rate-schedule-1} (LRS~1),
		% cyclic annealing,
		init\_lr=0.01.
		% n-cycles=5, N\_steps=200,000.
		% Figures~\ref{fig:pinned-pinned-static-solution-1}-\ref{fig:pinned-pinned-static-solution-2}, \emph{regular} grid can lead to static solution (same LRS 3 and init\_lr).
		\red{$\star$}
		Lowest loss value \red{2.12e-06}, Step 200,000.
		Total GPU time \red{392 sec}.
		$\bullet$
		{\footnotesize
			Figure~\ref{fig:23.8.22 R3b A-PP F(1J)2a(3) lr0.001 cy5-CA cy6-NCA Nsteps250000 regular N51 W64 H4 He init}, Form 2a, 12,802 parameters, regular grid, small damping, and
			Figure~\ref{fig:random-grid-NO-static-solution}, random grid, visually \emph{quasi-perfect} solution.
			Figure~\ref{fig:23.8.16 R2a A-PP F(1)2a lr0.01 cy5-CA Nsteps200000 regular N51 W32 H2 He init}, Form 2a, 1,218 parameters, regular grid, damping.
			Figures~\ref{fig:23.8.24 R1 A-PP F(1J,2a)3 lr0.04,0.03,0.02,0.01,0.01,0.005 cy1-6-VCA cy7-9-NCA Nsteps400000 random N51 W32 H2 He init-1}, \ref{fig:23.8.24 R1 A-PP F(1J,2a)3 lr0.04,0.03,0.02,0.01,0.01,0.005 cy1-6-VCA cy7-9-NCA Nsteps400000 random N51 W32 H2 He init-2}, \ref{fig:23.8.24 R1 A-PP F(1J,2a)3 lr0.04,0.03,0.02,0.01,0.01,0.005 cy1-6-VCA cy7-9-NCA Nsteps400000 random N51 W32 H2 He init-3}, \red{Form 3}, \red{1,251} parameters, visually \red{\emph{quasi-perfect}} solution.
			$\triangleright$
			Figure~\ref{fig:axial-Mathematica-solutions}, reference solution to compare.
		}
		{\scriptsize \mbox{(23816R1a-1)}}
	}
	\label{fig:23.8.16 R1a A-PP F(1)2a lr0.01 cy5-CA Nsteps200000 regular N51 W64 H2 He init}
\end{figure}

% 23.9.20, OLD Fig.21 to move to section Form 2a
\begin{figure}[tph]
	\includegraphics[width=0.49\textwidth]{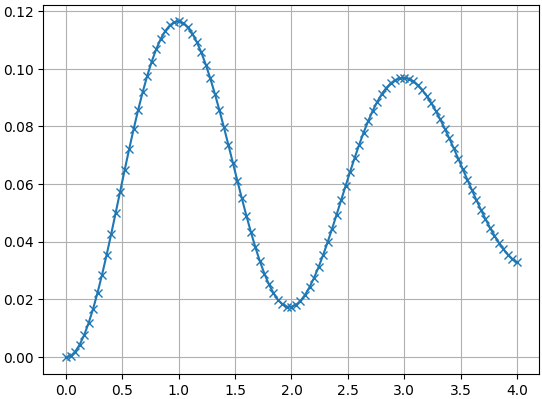}
	\includegraphics[width=0.49\textwidth]{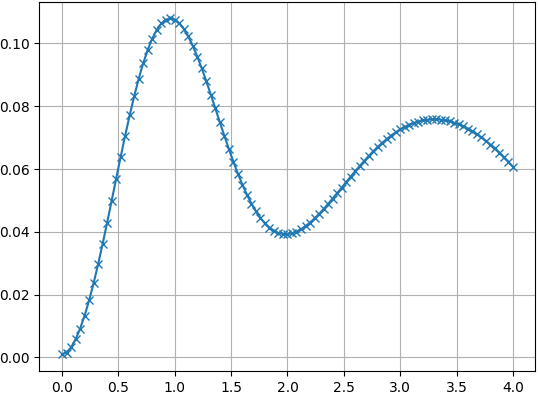}
	\caption{
		% 23.8.16 R2a and 23.8.16 R2b
		% \red{I AM HERE 2023.08.30.}
		% \red{I AM HERE, 23.9.8. Rewrite below.}
		\DDET.
		\emph{Pinned-pinned bar, lower-capacity network 2}. 
		Midspan displacement, \emph{damping}, Step 200,000.
		\red{$\star$}
		Remark~\ref{rm:optimal-capacity}, Optimal capacity. 
		% Section~\ref{sc:optimization-learning-rate-scheduling}
		Section~\ref{sc:pinned-pinned-Form-2a}, \red{Form 2a}.
		% Computational domain: $(x, t) \in [0, 1] \times [0, 4]$.
		\emph{Network:}
		Remark~\ref{rm:parameter-names},~\ref{rm:data-point-grids}, 
		T=4,
		W=32,
		H=2, 
		n\_out=2, 
		He-uniform initializer, 
		\red{1,218} parameters,
		\red{\emph{regular}} grid.
		\emph{Training:}
		init\_lr=0.01.
		\emph{Left:}
		Remark~\ref{rm:learning-rate-schedule-1}, LRS~1 (CA) {\scriptsize (23816R2a-1)};
		\emph{Right:}
		Remark~\ref{rm:learning-rate-schedule-3}, LRS~3 (NCA) {\scriptsize (23816R2b-1)}.
		% init\_lr=0.01,
		% n-cycles=5, N\_steps=200,000.
		$\bullet$
		{\footnotesize
			\red{\emph{Quasi-perfect solutions:}}
			Figure~\ref{fig:random-grid-NO-static-solution}, Form 2a, 12,802 parameters, random grid.
			Figure~\ref{fig:23.8.16 R1a A-PP F(1)2a lr0.01 cy5-CA Nsteps200000 regular N51 W64 H2 He init}, Form 2a, 4,482 parameters, regular grid; Remark~\ref{rm:unstable-solution}, \emph{Unstable} solution.
			Figures~\ref{fig:23.8.24 R1 A-PP F(1J,2a)3 lr0.04,0.03,0.02,0.01,0.01,0.005 cy1-6-VCA cy7-9-NCA Nsteps400000 random N51 W32 H2 He init-1}-\ref{fig:23.8.24 R1 A-PP F(1J,2a)3 lr0.04,0.03,0.02,0.01,0.01,0.005 cy1-6-VCA cy7-9-NCA Nsteps400000 random N51 W32 H2 He init-3}, \red{Form 3}, 1,251 parameters, random grid.
			Figure~\ref{fig:23.7.26 R5c A-PP F2a lr0.02x5 cycles=9 Nsteps=400000 regular N51 W32 H2 He init-3}, Form 2a, 1,218 parameters, regular grid.
		}	
	}
	\label{fig:23.8.16 R2a A-PP F(1)2a lr0.01 cy5-CA Nsteps200000 regular N51 W32 H2 He init}
\end{figure}

\begin{rem}
	\label{rm:optimal-capacity}
	Optimal network capacity.
	{\rm
		Figure~\ref{fig:random-grid-NO-static-solution} shows a \emph{quasi-perfect} solution using the network N51 W64 H4, with 12,802 parameters: Lowest total loss 0.537e-06, Step 192,000; Total GPU time 640 sec.  A question would be could a network with smaller capacity (fewer parameters) achieve the same (or similar) results?  Recall Remark~\ref{rm:learning-rate-schedule-3}, in which 
		Figure~\ref{fig:23.8.13 R4a A-PP F(1)2a lr0.04 cy5-NCA Nsteps200000 regular N51 W32 H2 He init}, NCA with init\_lr=0.04 and a model having 1,218 parameters led to waves with damping.
		% Figure~\ref{fig:23.8.15 R3a A-PP F(1)2a lr0.03 cy5-NCA Nsteps200000 regular N51 W64 H2 He init-1}, NCA with init\_lr=0.03 and a model having 4,482 parameters led to static solution.
		
		Figure~\ref{fig:23.8.16 R1a A-PP F(1)2a lr0.01 cy5-CA Nsteps200000 regular N51 W64 H2 He init} shows a much lower-capacity model with roughly one third the number of parameters, namely 4,482, and yet also produced similarly excellent results: Lowest loss value 2.12e-06, Step 200,000; Total GPU time 392 sec.  Could the model capacity be further reduced?
	
		Figure~\ref{fig:23.8.16 R2a A-PP F(1)2a lr0.01 cy5-CA Nsteps200000 regular N51 W32 H2 He init} shows a model with even lower capacity at 1,218 parameters, for which \emph{damping} in the time histories was clearly significant, less damping with cyclic annealing (Remark~\ref{rm:learning-rate-schedule-1}, LRS~1), and more damping without cyclic annealing (Remark~\ref{rm:learning-rate-schedule-3}, LRS~3).  
		But is it possible to reduce damping and lower the loss value with 1,218 parameters?  Yes, indeed, see Figure~\ref{fig:23.8.24 R1 A-PP F(1J,2a)3 lr0.04,0.03,0.02,0.01,0.01,0.005 cy1-6-VCA cy7-9-NCA Nsteps400000 random N51 W32 H2 He init-1}-\ref{fig:23.8.24 R1 A-PP F(1J,2a)3 lr0.04,0.03,0.02,0.01,0.01,0.005 cy1-6-VCA cy7-9-NCA Nsteps400000 random N51 W32 H2 He init-2}: 1,251 parameters; Lowest total loss 1.25e-06, Step 400,000;
		Total GPU time 598 sec; \emph{quasi-perfect} solution.
		
		Of the above three models, the second model with 4,482 parameters is the best for the axial motion of a pinned-pinned bar.
		See Remark~\ref{rm:damping-low-capacity-models} on how to avoid damping in low-capacity models using varying initial-learning-rate cyclic-annealing (VCA) (Remark~\ref{rm:learning-rate-schedule-2}) and extension of learning-rate schedule (ELRS) (Remark~\ref{rm:extension-cycle-5}).
		
		For the pinned-free bar, in Figure~\ref{fig:JAX-Form-1-pinned-free-bar-NO-static-solutions}, the network with N51 W32 H2, with 1,185 parameters, has the optimal capacity, providing low total loss (lowest among these four cases), quasi-perfect damping\% (See Remark~\ref{rm:shift-amplification-Form-1}) in the free-end displacement, and lowest GPU time.
	}
	\hfill $\blacksquare$
\end{rem}

% 23.9.20, OLD Fig.22 to move to section Form 2a
\begin{figure}[tph]
	\includegraphics[width=0.49\textwidth]{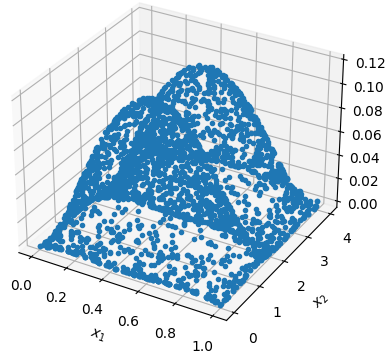}
	\includegraphics[width=0.49\textwidth]{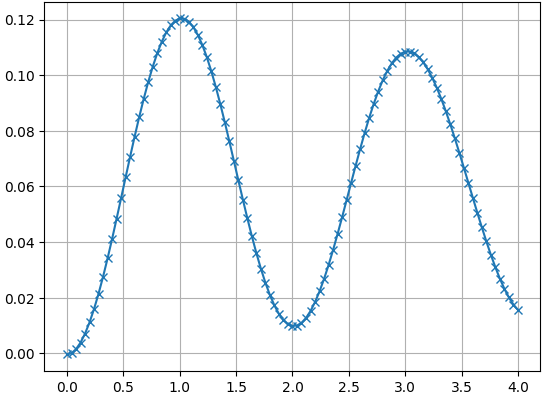}
	\caption{
		% {\color{red} Shape, Center disp, Step 200000.}
		% \red{I AM HERE 2023.09.03.}
		% \red{I AM HERE, 23.9.8. Rewrite below.}
		\DDET.
		\emph{Pinned-pinned bar, damping in low-capacity networks}.
		Shape (left), midspan displacement (right, \emph{damping}), Step 200,000.
		% {\color{red} TO REWRITE below for this case.} 
		\red{$\star$}
		Remark~\ref{rm:damping-low-capacity-models}, \red{\emph{Damping}}.
%		Section~\ref{sc:optimization-learning-rate-scheduling}.
		Section~\ref{sc:pinned-pinned-Form-2a}, \red{Form 2a}.
		% Computational domain: $(x, t) \in [0, 1] \times [0, 4]$.
		\emph{Network:}
		Remarks~\ref{rm:parameter-names},~\ref{rm:data-point-grids}, 
		T=4,
		W=32,
		H=2, 
		n\_out=2, 
		He-uniform initializer, 
		\red{1,218} parameters,
		\emph{random} grid.
		\emph{Training:}
		Remark~\ref{rm:learning-rate-schedule-1}, LRS~1, 
		% cyclic annealing,
		init\_lr=0.03.
		% n-cycles=5, N\_steps=200,000.
		$\bullet$
		{\footnotesize
			Figure~\ref{fig:23.8.17 R3a A-PP F(1)2a lr0.03 cy5-CA Nsteps200000 random N51 W32 H2 He init-2}, loss function, velocity, Step 200,000.
			Figure~\ref{fig:23.8.17 R2 - 23.8.17 R11}, midspan displacements, init\_lr=0.001 and init\_lr=0.1.
			$\triangleright$
			Figure~\ref{fig:axial-Mathematica-solutions}, reference solution to compare.
		}
		{\scriptsize (23817R3a-1)}
	}
	\label{fig:23.8.17 R3a A-PP F(1)2a lr0.03 cy5-CA Nsteps200000 random N51 W32 H2 He init-1}
\end{figure}

% 23.9.20, OLD Fig.23 to move to section Form 2a
\begin{figure}[tph]
	\includegraphics[width=0.49\textwidth]{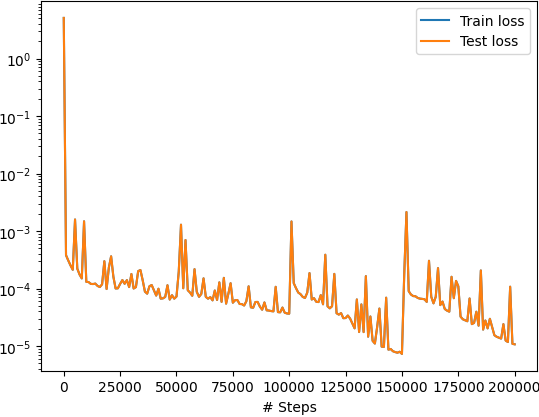}
	\includegraphics[width=0.49\textwidth]{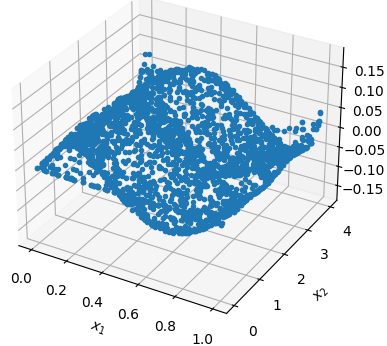}
	\caption{
		% Loss function, velocity, Step 200000.
		% \red{I AM HERE, 23.9.8. Rewrite below.}
		\DDET.
		\emph{Pinned-pinned bar, damping in low-capacity networks}.
		% {\color{red} I AM HERE 2023.09.03.} 
		Loss function (left), velocity (right), Step 200,000.
		Remark~\ref{rm:damping-low-capacity-models}, \emph{Damping}.
		% Section~\ref{sc:optimization-learning-rate-scheduling}. 
		\red{$\star$}
		Section~\ref{sc:pinned-pinned-Form-2a}, \red{Form 2a}.
		% Computational domain: $(x, t) \in [0, 1] \times [0, 4]$.
		\emph{Network:}
		Remarks~\ref{rm:parameter-names},~\ref{rm:data-point-grids}, 
		T=4,
		W=32,
		H=2, 
		n\_out=2, 
		He-uniform initializer, 
		\red{1,218} parameters,
		\emph{random} grid.
		\emph{Training:}
		\red{$\star$}
		Remark~\ref{rm:learning-rate-schedule-1}, \red{LRS~1}, 
		% cyclic annealing,
		init\_lr=0.03.
		% n-cycles=5, N\_steps=200,000.
		$\bullet$
		{\footnotesize
			Figure~\ref{fig:23.8.17 R3a A-PP F(1)2a lr0.03 cy5-CA Nsteps200000 random N51 W32 H2 He init-1}, shape, midspan displacement, Step 200,000.
			% Figure~\ref{fig:23.8.17 R3a A-PP F(1)2a lr0.03 cy5-CA Nsteps200000 random N51 W32 H2 He init-2}, loss function, velocity, Step 200,000.
			Figure~\ref{fig:23.8.17 R2 - 23.8.17 R11}, midspan displacements, init\_lr=0.001 \& init\_lr=0.1.
		}
		{\scriptsize (23817R3a-2)}
	}
	\label{fig:23.8.17 R3a A-PP F(1)2a lr0.03 cy5-CA Nsteps200000 random N51 W32 H2 He init-2}
\end{figure}

% 23.9.20, OLD Fig.24 to move to section Form 2a
\begin{figure}[tph]
	\includegraphics[width=0.49\textwidth]{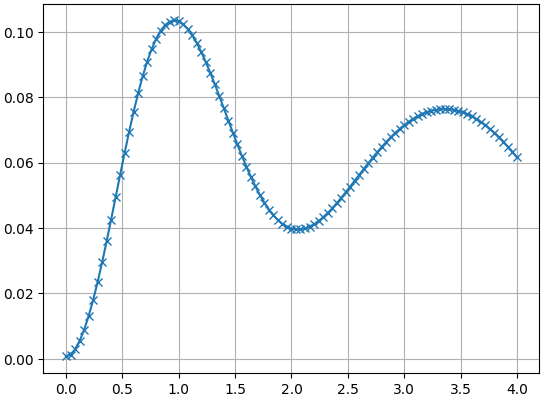}
	\includegraphics[width=0.49\textwidth]{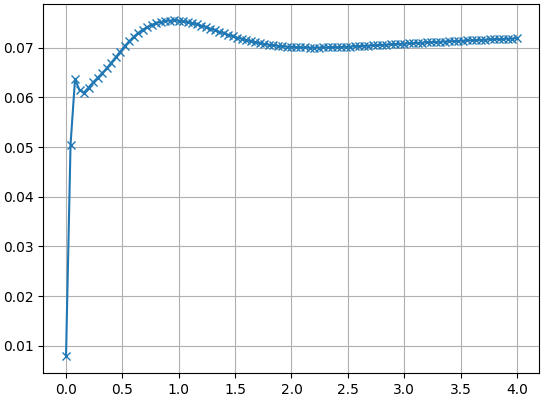}
	\caption{
		% Midspan displacement time histories.
		% Left: Waves, damping, 23.8.17 R2 A-PP F(1)2a lr0.001, Step 200,000.
		% Right: Static solution, 23.8.17 R11 A-PP F(1)2a lr0.1, Step 50,000.
		% {\color{red} I AM HERE 2023.09.03.}
		% \red{I AM HERE, 23.9.8. Rewrite below.}
		\DDET.
		\emph{Pinned-pinned bar, damping in low-capacity networks}.
		\emph{Left:} Waves, damping, Step 200,000. 
		\emph{Right:} Static solution, Step 50,000. 
		\red{$\star$}
		Remark~\ref{rm:damping-low-capacity-models}, \emph{Damping}.
		% Section~\ref{sc:optimization-learning-rate-scheduling}.
		Section~\ref{sc:pinned-pinned-Form-2a}, \red{Form~2a}.
		% Computational domain: $(x, t) \in [0, 1] \times [0, 4]$.
		\emph{Network:}
		Remarks~\ref{rm:parameter-names},~\ref{rm:data-point-grids}, 
		T=4,
		W=32,
		H=2, 
		n\_out=2, 
		He-uniform initializer, 1,218 parameters,
		\emph{random} grid (left), \emph{regular} grid (right).
		\emph{Training:}
		\red{$\star$}
		\emph{Left:}
		Remark~\ref{rm:learning-rate-schedule-1}, \red{LRS~1},
		init\_lr=\red{0.001}
		% n-cycles=5, N\_steps=200,000 
		{\scriptsize (23817R2-1)}.
		\emph{Right:}
		Remark~\ref{rm:learning-rate-schedule-3}, \red{LRS~3}, 
		init\_lr=\red{0.1},
		n-cycles=2, N\_steps=50,000 {\scriptsize (23817R11-1)}.
	}
	\label{fig:23.8.17 R2 - 23.8.17 R11}
\end{figure}

\begin{rem}
	\label{rm:damping-low-capacity-models}
	Damping in low-capacity models.
	{\rm
		For a pinned-pinned bar using Form 2a and a model with N51 W32 H2, having 1,218 parameters, stopping at N\_steps=200,000, the lowest damping was obtained around the initial learning rate init\_lr=0.03, which yielded a ratio of 1st peak to 2nd peak around 0.12 / 0.11 = 1.09; Figure~\ref{fig:23.8.17 R3a A-PP F(1)2a lr0.03 cy5-CA Nsteps200000 random N51 W32 H2 He init-1} (shape, midspan displacement, Step 200,000), Figure~\ref{fig:23.8.17 R3a A-PP F(1)2a lr0.03 cy5-CA Nsteps200000 random N51 W32 H2 He init-2} (loss function, velocity, Step 200,000).
		Several simulations using this same model and a wide range of initial learning rate, from 0.001 (Figure~\ref{fig:23.8.17 R2 - 23.8.17 R11}, left)
		% (23.8.17 R2, waves, damping) 
		to 0.1 (Figure~\ref{fig:23.8.17 R2 - 23.8.17 R11}, right), 
		% (23.8.17 R11, pre-static) 
		damping was inescapable in the results.
		Note that the damping with the smaller init\_lr=0.001 was more pronounced than the damping with the 30 times larger init\_lr=0.03.
		
		So at first, \emph{damping} appeared to be a characteristic of low-capacity models with N\_steps below 200,000 until we tried the varying initial-learning-rate cyclic annealing (VCA) mentioned in Remark~\ref{rm:learning-rate-schedule-2} and Figure~\ref{fig:23.8.21 R1 A-PP F(1)2a lr0.04,0.03,0.02,0.01,0.005 cy5-VCA Nsteps200000 regular N51 W32 H2 He init} (Form 2a, 1,218 parameters), and Figures~\ref{fig:23.8.21 R1c A-PP F(1J,2a)3 lr0.04,0.03,0.02,0.01,0.005 cy5-VCA Nsteps200000 regular N51 W32 H2 He init-1}-\ref{fig:23.8.21 R1c A-PP F(1J,2a)3 lr0.04,0.03,0.02,0.01,0.005 cy5-VCA Nsteps200000 regular N51 W32 H2 He init-2} (Form 3, 1,251 parameters), yielding much lower damping compared to that in Figure~\ref{fig:23.8.17 R3a A-PP F(1)2a lr0.03 cy5-CA Nsteps200000 random N51 W32 H2 He init-1} and Figure~\ref{fig:23.8.17 R2 - 23.8.17 R11}, left.
		% {\color{red} HERE}.
		
		Damping can further reduced after VCA scheduling by extension of learning-rate scheduling (ELRS) to, say, Step 400,000. 
		An example is given in
		Figures~\ref{fig:23.8.24 R1 A-PP F(1J,2a)3 lr0.04,0.03,0.02,0.01,0.01,0.005 cy1-6-VCA cy7-9-NCA Nsteps400000 random N51 W32 H2 He init-1}-\ref{fig:23.8.24 R1 A-PP F(1J,2a)3 lr0.04,0.03,0.02,0.01,0.01,0.005 cy1-6-VCA cy7-9-NCA Nsteps400000 random N51 W32 H2 He init-3}, Form 3, 1,251 parameters, random grid, \emph{quasi-perfect} solution. 
		See also Remark~\ref{rm:unstable-solution}, Unstable solution for another example of VCA-ELRS leading to \emph{quasi-perfect} solution.
		% {\color{red} [NOTE: 2023.08.24, add examples of From 2a and Form 3 ran on 23.8.23 here.  ENDNOTE]}
	}
	\hfill $\blacksquare$
\end{rem}

%\begin{rem}
%	\label{rm:learning-rate-grid-convergence-static-solution}
%	Learning rate, data-point grid, convergence, static solution.
%	{\rm
%		Hello
%	}
%\end{rem}

\begin{rem}
	\label{rm:stochastic-reproducibility}
	Stochastic reproducibility. 
	{\rm
		The results are only reproducible in the stochastic sense, i.e., the results could vary a little, even though qualitatively similar, when rerunning the same identical DeepXDE script, with exactly the same parameters. There are several sources of stochasticity: 
		(1) The stochastic optimization initializer, i.e, either 
		``Glorot uniform''
		% (Section~\ref{sc:xavier-glorot})
		or
		``He uniform'' \cite{vuquoc2023deep}
		% (Section~\ref{sc:kaiming-he}), 
		(2) the density and probabilistic distribution of the data points in the (random) grid for the evaluation of the loss function. 
		In other words, the results while stochastically reproducible are \emph{deterministically irreproducible}.
		% {\color{red} [NOTE: 2023.08.04, give figures showing how the results changed after running the same identical DeepXDE script file.  ENDNOTE]}
		
		%
		% define verbatim box in footnote
		\begin{verbbox}
			\footnotesize
			\verb+def set_random_seed(seed)+
		\end{verbbox}
		% define 2nd verbatim box in footnote
		\begin{verbboxT}
			% \verb|a__|
			\footnotesize
			\verb+dde.config.set_random_seed(42)+
		\end{verbboxT}
		%
		% DeepXDE 
		\DDET\ 
		allows for setting the random\_seed to a number, such as ``42'' (see \emph{fixed random grid} in Remark~\ref{rm:data-point-grids}), in its configuration to produce \emph{deterministically reproducible} results 
		% using the command \verb*|deepxde.config.set_random_seed(42)|, 
		at the cost of slowing down the computation, and therefore this setting should be used only for debugging.\footnote{ 
			For \DDET, see the function definition 
			% ``{def set\_random\_seed(seed)}''
			\verbBox 
			in \href{https://deepxde.readthedocs.io/en/latest/\_modules/deepxde/config.html?highlight=deepxde.config}{Source code for deepxde.config}, and the command 
			% ``dde.config.set\_random\_seed(42)'' 
			\verbBoxT
			is used for determinism.
		}
		As an example, using the same model with 4,482 parameters and 200,000 steps as in Figure~\ref{fig:23.8.16 R1a A-PP F(1)2a lr0.01 cy5-CA Nsteps200000 regular N51 W64 H2 He init} resulted in a 10\% increase in computational time.
		%
		% define verbatim box in footnote
		\begin{verbbox}
			% \verb|a__|
			\footnotesize
			\verb+os.environ["XLA_FLAGS"] = "--xla_gpu_deterministic_ops=true"+
		\end{verbbox}
		% define 2nd verbatim box in footnote
		\begin{verbboxT}
			% \verb|a__|
			\footnotesize
			\verb+key = jax.random.PRNGKey(42)+
		\end{verbboxT}
		See Remark~\ref{rm:JAX-deterministic-faster} regarding the opposite for our \JAX\ script, for which the deterministic mode is more efficient than the non-deterministic mode, an unexpected surprise.\footnote{
			For \JAX, the commands 
			% ``os.environ["XLA\_FLAGS"] = "--xla\_gpu\_deterministic\_ops=true"'' 
			% \verb+os.environ["XLA_FLAGS"] = "--xla_gpu_deterministic_ops=true"+ 
			\verbBox
			and 
			% ``key = jax.random.PRNGKey(42)'' 
			\verbBoxT
			are used for determinism.
		}
	}
	\hfill
	$\blacksquare$
\end{rem}

\begin{rem}
	% \label{rm:unstable-results}
	\label{rm:unstable-solution}
	% Unstable results.
	Unstable solution and method to fix.
	{\rm
		Sometimes, the solution could be unstable in the sense that it could not be recovered as if it occurred by chance such as that in Figure~\ref{fig:23.8.16 R1a A-PP F(1)2a lr0.01 cy5-CA Nsteps200000 regular N51 W64 H2 He init}, which could never again be reproduced. 
		Rerunning the exact same script again several times, from which Figure~\ref{fig:23.8.16 R1a A-PP F(1)2a lr0.01 cy5-CA Nsteps200000 regular N51 W64 H2 He init} was obtained, led to pre-static or static solutions.  
		That's \emph{stochastically irreproducible}.
		
		When such situation occurs, gradually reducing the initial learning rate with cyclic annealing, using Remark~\ref{rm:learning-rate-schedule-2} (VCA), and extending the learning-rate schedule, Remark~\ref{rm:extension-cycle-5} (ELRS), would usually help.  For the case of Figure~\ref{fig:23.8.16 R1a A-PP F(1)2a lr0.01 cy5-CA Nsteps200000 regular N51 W64 H2 He init} with lowest total loss of 2.12e-06 at Step 200,000, using VCA for Cycles 1-5 with init\_lr = [0.01, 0.009, 0.008, 0.007, 0.006] and NCA for Cycles 6-9, and N\_steps = 400,000 steps, yielded a \emph{quasi-perfect} solution with 
		% lowest total loss 1.92e-06 at Step 200,000 with total GPU time 401 sec, and
		lowest total loss 0.793e-06 at Step 400,000 with total GPU time 748 sec (RunID 23831R2c). This VCA-ELRS solution is \emph{stochastically reproducible}, i.e., \emph{stable} in the sense that rerunning the script several times led to stochastically equivalent solutions. 
		
		Hence with the above VCA-ELRS, a good trade-off would be at Step 200,000 with 4,482 parameters, lowest total loss 1.92e-06 at Step 200,000, total GPU time 401 sec (RunID 23831R2c), and 
		a \emph{quasi-perfect} solution not distinguishable from that in Figure~\ref{fig:random-grid-NO-static-solution}, which was obtained with 12,802 parameters, lowest total loss 0.537e-06 at Step 192,000, total GPU time 640 sec. 
	}
	\hfill$\blacksquare$
\end{rem}

%\subsection{Wave equation}
\subsection{Axial motion of elastic bar}
\label{sc:wave-equation}
\label{sc:axial-motion}

% 2023.06.28
% changed Form i to Form i+1 to be consistent with the code
% \subsubsection{PINN Form 1: Time shift and early stopping}
% \subsubsection{DDE-T, Form 1: Shift, amplification, early stopping, static solution}
\subsubsection{Form 1: Shift, amplification, early stopping, static solution}
\label{sc:axial-Form-1}
\label{sc:pinned-pinned-Form-1}
The numerical evidence presented below shows that our JAX script solves many problems encountered with DDE-T, specifically:
\begin{itemize}
	\item For the pinned-pinned bar, DDE-T produced shift and amplification in the solution, whereas our JAX script shows no such problems.  To solve these particular DDE-T problems (shift and amplification), we devised several methods (Form 2a, Form 3, Form 4) before writing our own JAX script.
	
	\item For the pinned-free bar, DDE-T produced a static-solution time history, i.e., not dynamic solution, regardless of the size of the learning rate or the network capacity, whereas our JAX script shows no such problem.  To solve this particular DDE-T problem (static solution), we devised several methods (barrier function, reduced network capacity, different learning-rate schedules) before writing our own JAX script.  
	% We devised several methods to solve this DDE-T problem before writing our own JAX script.
\end{itemize}

\noindent
{\bf Form 1.}
The Form 1 of the axial equation of motion of an elastic bar is a particular case of the Form 1 for the Kirchhoff-Love rod in Section~\ref{sc:Kirchhoff-rod-Form-1}, and 
was given above in Eq.~\eqref{eq:eom-euler-bernoulli-axial}, reproduced below for convenience:
\begin{align}
	\sldn \ubp{2} + \dfbs{\X} = \nddt \ub
	\ ,
	\tag{\ref{eq:eom-euler-bernoulli-axial}}
\end{align}
where 
%$\rtdn$ is the rotundness
$\sldn = 1$ is the slenderness 
of the bar, defined in Eq.~\eqref{eq:slenderness-rotundness},
$\dfbs{\X} = \frac12$, 
with the following two types of prescribed boundary conditions:

\noindent
\emph{Pinned-pinned bar,} using Eq.~\eqref{eq:pinned-left-end-bound-cond}$_1$ for $\Xb=0$ and a similar one for $\Xb=1$:
\begin{align}
	\ub (\Xb = 0, \tb) = \ub (\Xb = 1, \tb) = 0 
	\label{eq:axial-pinned-pinned-BCs}
\end{align}

\noindent
\emph{Pinned-free bar,} using Eq.~\eqref{eq:pinned-left-end-bound-cond}$_1$ and Eq.~\eqref{eq:roller-right-end-bound-cond}$_1$, and the extension $\extb$ from Eq.~\eqref{eq:finite-extension-2} and Eq.~\eqref{eq:extension-bar}:
\begin{align}
	\ub (\Xb = 0, \tb) = \widehat{\ub} = 0 
	\ , \quad
	\extb (\Xb = 1, \tb ) = \ubp{1} (1, \tb) = \frac{1}{\sldn} \widehat{\Nb} = 0
\end{align}

\noindent
The prescribed initial conditions, using Eq.~\eqref{eq:initial-conditions-1}, are:
\begin{align}
	\ub (\Xb , 0) = \widetilde{\ub}_0 (\Xb) = 0 , \quad
	\ubd (\Xb , 0) = \widetilde{\ubd}_0 (\Xb) = 0 \ .
	\label{eq:axial-ICs}
\end{align}
The reference solutions for the pinned-pinned bar and for the pinned-free bar are given in Figure~\ref{fig:axial-Mathematica-solutions}.

Since both the strain $\ubp{1}$ and the velocity $\ubd$ are not explicit outputs in Form 1, these boundary and initial conditions cannot be explicitly enforced using hard constraints, but have to be imposed after computing the derivatives of the outputs using backpropagation and then included as additional squared errors in the overall loss function (soft constraints).

\begin{rem}
	\label{rm:time-shift-Form-1}
	\label{rm:Form-1-shift-amplification}
	\label{rm:shift-amplification-Form-1}
	% Form-1 time shift.
	% Pinned-pinned bar, Form 1:
	% DDE-T vs JAX, 
	\DvJ,
	Form 1, pinned-pinned bar:
	Time shift, amplification, damping percentage.
	{\rm
		% DDE-T 
		\DDET\ with Form~1 produced a shift and amplification in the solution and larger errors in boundary conditions with derivatives (Neumann type).
		JAX does not exhibit these pathological problems. 
		
		For a pinned-pinned elastic bar subjected to a constant distributed axial load,  
		Figure~\ref{fig:23.7.22 R1 A-PP v1 midspan,shape100000} shows two vibration periods of 
		the deformed-shape time history and
		the midspan-displacement time history, shifted to the right, 
		at training Step 100,000 obtained from a \emph{regular} grid, whereas  
		Figure~\ref{fig:23.7.23 R1d A-PP v1 midspan,shape200000} shows the same output variables at Step 200,000 obtained from a \emph{fixed random} grid.
		% ALEX: removed paragraph
		Also for the same pinned-pinned bar, Figures~\ref{fig:DEBUG v1.9.3 seed42, 23.9.4 R1 A-P bcrP F1J(2a,3) lr0.001 Cy1-9-NCA Nsteps400000 regular N51 W64 H4 Glorot init-1}-\ref{fig:DEBUG 23.9.4 R1d A-P bcrP F1J(2a,3) lr0.001 Cy1-9-NCA Nsteps400000 regular N51 W64 H4 Glorot init-1} show the shift (time and vertical) and amplification parameters increase continuously with training step number, as the training progressed to the final Step 400,000.
		
		Using \JAX, Figure~\ref{fig:DEBUG 23.9.17 R1d JAX A-P bcrP F1H(2a,3) lr0.002 cy1-7-CA Nsteps300000 random N51 W64 H4 Uni init-1} shows no shift since the times of the local maxima were exactly where they should be, and so were the times for the local minima.  There was no amplification, 
		except for a small damping. 
		% ALEX: removed paragraph
		For results obtained with \JAX, the quality of a response curve is measured by the damping percentage (damping\% for short) defined by the ratio of the 1st local maximum over the 2nd local maximum minus one, with amplification\% = 1 - damping\% > 0.
		For example, the midspan displacement of a pinned-pinned bar in Figure~\ref{fig:DEBUG 23.9.17 R1d JAX A-P bcrP F1H(2a,3) lr0.002 cy1-7-CA Nsteps300000 random N51 W64 H4 Uni init-1} has a damping\% of 1.7\%, which is discernible.  The quality of the response is classified by the damping\% as follows:
		\begin{itemize}
			\item $\mid $ damping\% $\mid < 0.5\%$ : Quasi-perfect, damping/amplification is not discernible by naked eyes.
			
			\item $\mid $ damping\% $\mid \in \left[ 0.5\% , 1.5\% \right)$ : Very good, small damping, could be discernible by naked eyes.
			
			\item $\mid $ damping\% $\mid \in \left[ 1.5\% , 3.0\% \right)$ : Good, significant damping, discernible by naked eyes.
			
			\item $\mid $ damping\% $\mid \ge 3.0\%$ : High damping, easily visible.
		\end{itemize}
	}
	\hfill$\blacksquare$
\end{rem}

%		\red{[NOTE: 23.8.28, pinned-free bar in Form 1: Check whether there are time shift and amplification.  Static solution do exist for pinned-free bar, but NOT for pinned-pinned bar, both in Form 1.  ENDNOTE]}

%--------------------------------------------------------------Form 1
% 23.9.20, OLD Fig.25 to move to section Form 1
\begin{figure}[tph]
	\centering
	%_\includegraphics[width=0.479\textwidth]{Figures/23.7.21_R4_No_barrier_lr0.001_cycles=1_Nsteps=100000_regular_grid_Glorot.PNG}
	%
	% script file (new filename with cy3-NCA, OLD cycles=3 Nsteps=100000)
	% 23.7.22 R1: A-PP v1 NO barrier lr 0.001 cy3-NCA Nsteps=100000 - regular N51 W64 H4
	% https://colab.research.google.com/drive/1HNdRBm2GCIDnktwGJ1BIeV8eXxHXV_2q
	%_\includegraphics[width=0.49\textwidth]{Figures/23.7.21_R4_pinned-pinned_shape_2_bumps.PNG}
	\includegraphics[width=0.49\textwidth]{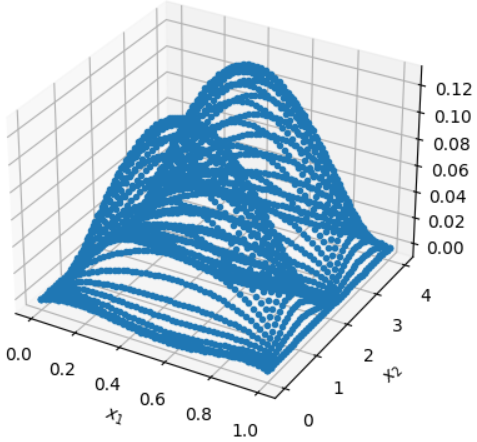}
	\includegraphics[width=0.479\textwidth]{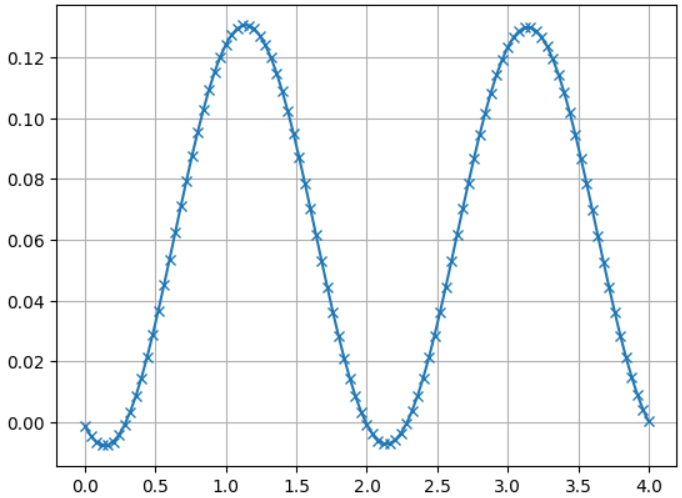}
	\caption{
		% \red{I AM HERE, 23.9.8. Rewrite below.}
		\DDET.
		\emph{Pinned-pinned bar}. 
		Shape (left),
		midspan displacement \emph{shifted to the right with small amplification} (right), Step 100,000.
		% Section~\ref{sc:wave-equation}.
		\red{$\star$}
		Section~\ref{sc:axial-Form-1}, \red{Form~1}. 
		% Two vibration periods. 
		\emph{Network:}
		Remarks~\ref{rm:parameter-names},~\ref{rm:data-point-grids},
		T=4,
		W=64,
		H=4, 
		n\_out=1, 
		Glorot-uniform initializer, 12,737 parameters,
		\emph{regular} grids.
		\emph{Training:}
		\red{$\star$}
		Remark~\ref{rm:learning-rate-schedule-3}, \red{LRS~3 (NCA)},
		init\_lr=0.001,
		n-cycles=3, N\_steps=100,000. 
		\red{$\star$}
		PDE loss 2.20e-06, Cycle 3 GPU time 154 sec. 
		% commented out figure number, 23.10.15
		% (cf. Figure~\ref{fig:23.7.23 R1d A-PP v1 midspan,shape100000} for PDE loss and GPU time). 
		{\footnotesize
			$\triangleright$
			Figure~\ref{fig:axial-Mathematica-solutions}, reference solution to compare.
		}
		{\scriptsize (23722R1-1)}
		% {\tiny (23722R1-1)} 
	}
	\label{fig:23.7.22 R1 A-PP v1 midspan,shape100000}
\end{figure}

% 23.9.20, OLD Fig.26 to move to section Form 1
\begin{figure}[tph]
	\includegraphics[width=0.49\textwidth]{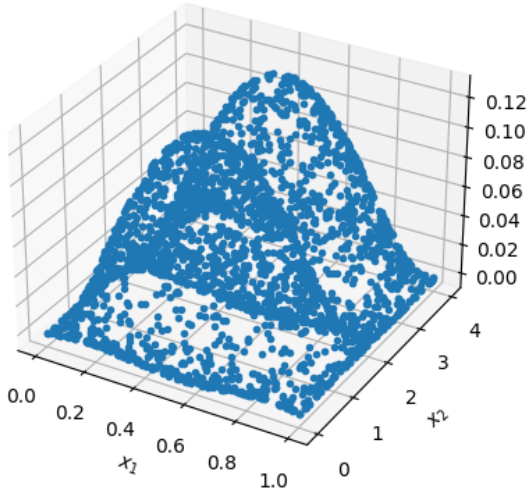}
	\includegraphics[width=0.49\textwidth]{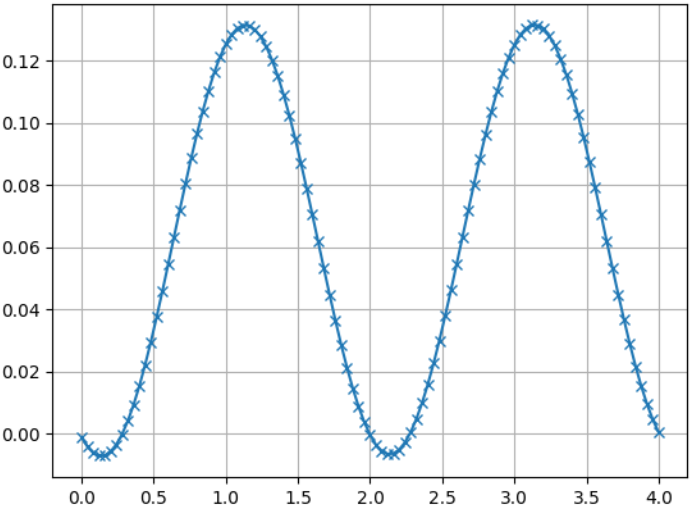}
	\caption{
		% \red{I AM HERE, 23.9.8. Rewrite below.}
		\DDET.
		\emph{Pinned-pinned bar}. 
		% Section~\ref{sc:wave-equation}.
		Shape (left),
		midspan displacement \emph{shifted to the right with small amplification} (right), Step 200,000.
		\red{$\star$}
		Section~\ref{sc:axial-Form-1}, \red{Form~1}.
		\emph{Network:}
		Remarks~\ref{rm:parameter-names},~\ref{rm:data-point-grids},
		T=4,
		H=4, W=64,
		n\_out=1, 
		Glorot-uniform initializer, 
		\red{12,737} parameters
		\emph{fixed random} grid. 
		\emph{Training:}
		\red{$\star$}
		Remark~\ref{rm:learning-rate-schedule-1}, LRS~1 (CA),
		init\_lr=0.001.
		\red{$\star$}
		Lowest total loss \red{4.31e-06}, Step 200,000 (sum of 5 losses).
		\red{$\star$}
		Total GPU time \red{621 sec}.
		$\bullet$
		{\footnotesize
			% See  
% commented out figure numbers, 23.10.15
%			Figures~\ref{fig:23.7.23 R1d A-PP v1 loss200000 shape25000} (Step 25,000), 
%			\ref{fig:23.7.23 R1d A-PP v1 midspan,shape50000} (Step 50,000),
%			\ref{fig:23.7.23 R1d A-PP v1 midspan,shape100000} (Step 100,000).
%			% \ref{fig:23.7.23 R1d A-PP v1 midspan,shape200000} (Step 200000).
%			% \ref{fig:23.07.04 R1 A-PP v1 midspan,shape200000 left unknown right} (Step 200000).
%			% Figure~\ref{fig:axial-Mathematica-solutions}, reference solution for axial motion and parameter values.
%			$\triangleright$
			Figures~\ref{fig:23.8.24 R1 A-PP F(1J,2a)3 lr0.04,0.03,0.02,0.01,0.01,0.005 cy1-6-VCA cy7-9-NCA Nsteps400000 random N51 W32 H2 He init-1}-\ref{fig:23.8.24 R1 A-PP F(1J,2a)3 lr0.04,0.03,0.02,0.01,0.01,0.005 cy1-6-VCA cy7-9-NCA Nsteps400000 random N51 W32 H2 He init-3}, \red{Form 3}, VCA-ELRS, 
			\red{1,251} parameters,
			\red{\emph{quasi-perfect}} solution (no shift),
			lowest total loss \red{1.25e-06},
			total GPU time \red{598 sec}.
			$\triangleright$
			Figure~\ref{fig:axial-Mathematica-solutions}, reference solution to compare.
		}
		{\scriptsize (23723R1d-4)} % \protect\footnotemark
	}
	\label{fig:23.7.23 R1d A-PP v1 midspan,shape200000}
\end{figure}
%\footnotetext{
	%	For Figure~\ref{fig:23.7.23 R1d A-PP v1 midspan,shape200000}, with  figure ID ``23723R1d-4,'' see Footnote~\ref{fn:23723R1d}.
	%}

% 23.9.20, OLD Fig.29 STAY HERE
\begin{figure}[tph]
	\includegraphics[width=0.484\textwidth]{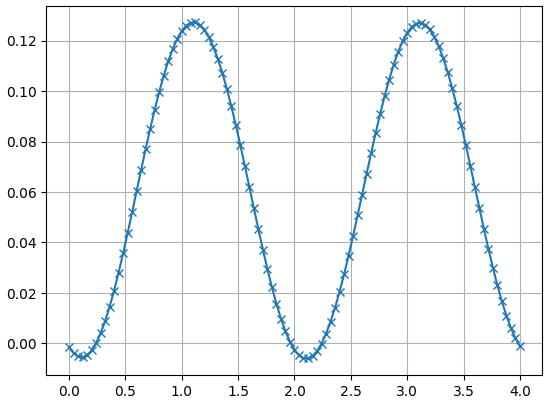}
	%_\includegraphics[width=0.496\textwidth]{Figures/DEBUG_v1.9.3_seed42,_23.9.4_R1_A-P_bcrP_F1J_2a,3__lr0.001_Cy1-9-NCA_Nsteps400000_regular_N51_W64_H4_Glorot_init_-_center_disp400000.png}
	%_\includegraphics[width=0.496\textwidth]{Figures/DEBUG_v1.9.3_seed42,_23.9.4_R1_A-P_bcrP_F1J_2a,3__lr0.001_Cy1-9-NCA_Nsteps400000_regular_N51_W64_H4_Glorot_init_-_shift_params400000.png}
	\includegraphics[width=0.496\textwidth]{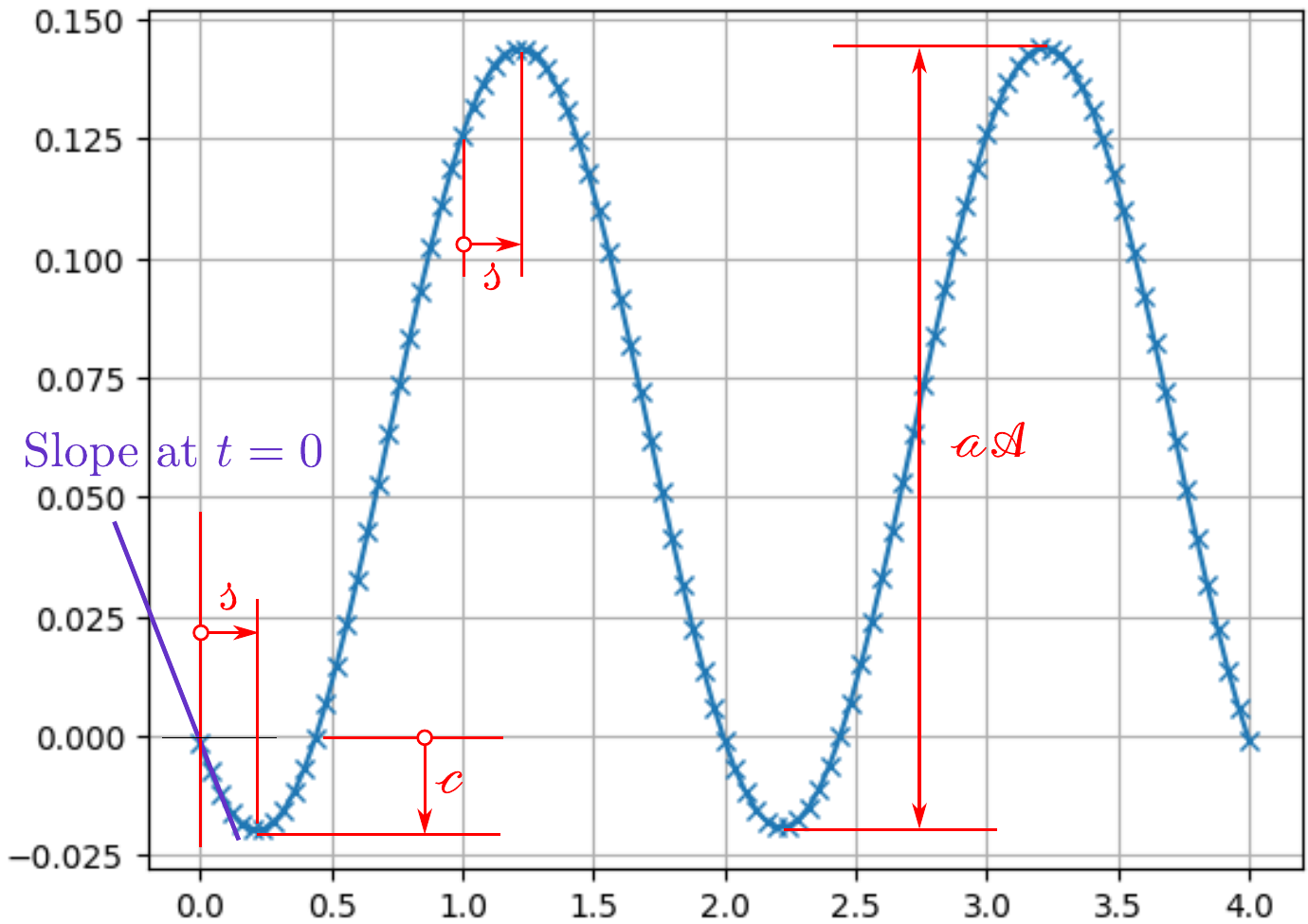}
	\caption{
		% \red{I AM HERE, 23.9.8. Rewrite below.}
		\DDET.
		\emph{Pinned-pinned bar, \red{continued shift \& amplification}}. 
		\red{$\star$}
		Midspan displacements: Step 50,000 (left), Step 400,000 (right). 
		% \emph{Left:} Midspan displacement, Step 50000.
		% \emph{Right:} Midspan displacement, Step 40000; hidden parameters for shift and amplification. 
		\red{$\star$}
		Section~\ref{sc:axial-Form-1}, \red{Form~1}:
		\emph{Hidden parameters} ($\shiftPs$, $\shiftPc$, $\shiftPa$), time shift $\shiftPs$, vertical shift $\shiftPc$, amplification $\shiftPa$, quasi-perfect peak-to-peak amplitude $\shiftPA$.
		\emph{Network:}
		Remarks~\ref{rm:parameter-names},~\ref{rm:data-point-grids},
		T=4
		W=64,
		H=4,
		n\_out=1, 
		He-uniform initializer, 
		12,737 parameters,
		regular grid.
		% {\color{red} REWRITE below.} 
		\emph{Training:}
		Remark~\ref{rm:learning-rate-schedule-1}, \red{LRS~3 (NCA)},
		init\_lr=0.005,
		N\_steps=400,000.
		$\bullet$
		{\footnotesize
			% Figure~\ref{fig:23.8.6 R2b A-PP F1(2a) lr0.03 cy3 Nsteps100000 regular N51 W64 H4 He init-2}, midspan displacement at Step 25,000 and zero midspan displacement at Step 100,000 (divergence).
			Figure~\ref{fig:DEBUG 23.9.4 R1d A-P bcrP F1J(2a,3) lr0.001 Cy1-9-NCA Nsteps400000 regular N51 W64 H4 Glorot init-1}, shift-amplification parameters increase with training step number.
		}
		{\scriptsize (2394R1a-1)} 
	}
	\label{fig:DEBUG v1.9.3 seed42, 23.9.4 R1 A-P bcrP F1J(2a,3) lr0.001 Cy1-9-NCA Nsteps400000 regular N51 W64 H4 Glorot init-1}
\end{figure}

\begin{figure}[tph]
	\centering
	% this image did not have clear lettering, because it was a screenshot taken
	% from Google Sheets, so I created the image below after downloading the svg
	% file from Google Sheets and use Inkscape
	%_\includegraphics[width=0.55\textwidth]{Figures/DEBUG_23.9.4_R1d_A-P_bcrP_F1J_2a,3__lr0.001_Cy1-9-NCA_Nsteps400000_regular_N51_W64_H4_Glorot_init_-_shift_params_vs_step_number.png}
	\includegraphics[width=0.55\textwidth]{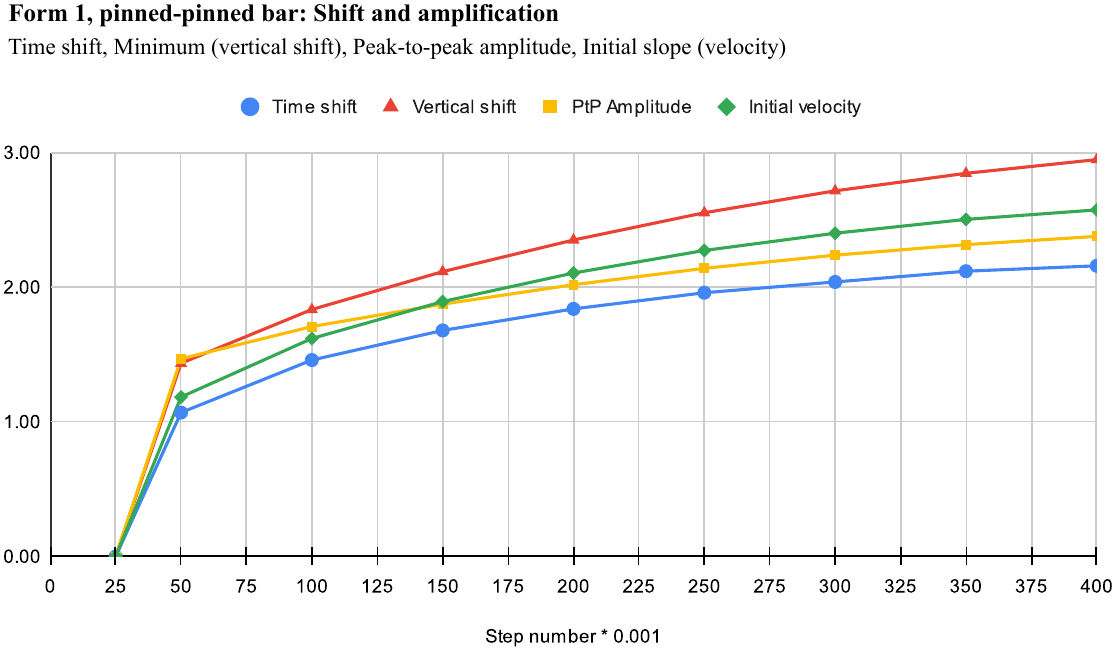}
	\includegraphics[width=0.40\textwidth]{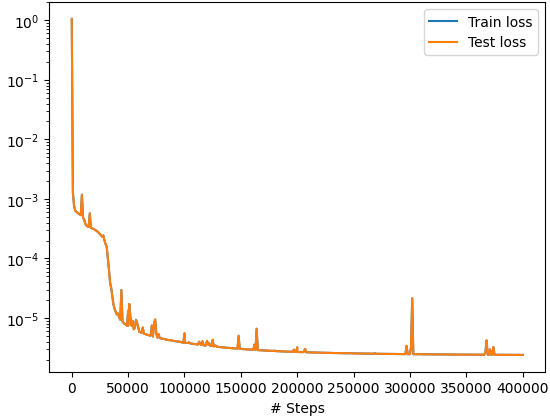}
	\caption{
		\DDET.
		\emph{Pinned-pinned bar.  Form 1}.
		Appendix~\ref{app:analysis-time-shift-amplification}, 
		% Section~\ref{sc:time-shift-amplification}, 
		\emph{Analysis}.
		Scaled
		shift-amplification parameters ($\shiftPs , \shiftPc , \shiftPa$) and initial velocity increase with training step number (left), while the
		Total loss continues to decrease (right).
		All parameters $p$, such as time shift $\shiftPs$, vertical shift $\shiftPc$, amplification factor $\shiftPa$ and thus peak-to-peak amplitude $\shiftPa \shiftPA$, and initial velocity, were scaled to fit in the chart.
		% using $p_{\text{scaled}} = \pm (p - p_{\text{25000}}) * k_p > 0$, where $p_{\text{scaled}}$ is the positive scaled parameter $p$, with $p_{\text{25000}}$ the value of $p$ at Step 25000, and $k_p$ a selected scale factor for $p$.
		$\bullet$
		{\footnotesize
			Table~\ref{tb:DEBUG 23.9.4 R1d A-P bcrP shift-amplification-values}, \emph{unscaled} values of $\shiftPs , \shiftPc , \shiftPa$ and computed initial velocity, scaling method.
			Figure~\ref{fig:DEBUG v1.9.3 seed42, 23.9.4 R1 A-P bcrP F1J(2a,3) lr0.001 Cy1-9-NCA Nsteps400000 regular N51 W64 H4 Glorot init-1}, midspan displacements, Steps 50000 \& 400000.
		}
		{\scriptsize (2394R1d-2)}
	}
	\label{fig:DEBUG 23.9.4 R1d A-P bcrP F1J(2a,3) lr0.001 Cy1-9-NCA Nsteps400000 regular N51 W64 H4 Glorot init-1}
\end{figure}

% 23.9.22
% image from 23.9.17 R1c not used, since reran in 23.9.17 R1d to get damping% and
% rating of quality of free-end disp 
% DEBUG 23.9.17 R1c JAX A-P bcrP F1H(2a,3) lr0.002 cy1-5-CA Nsteps200000 random N51 W64 H4 Uni init - loss200000.png
\begin{figure}[tph]
	\includegraphics[width=0.45\textwidth]{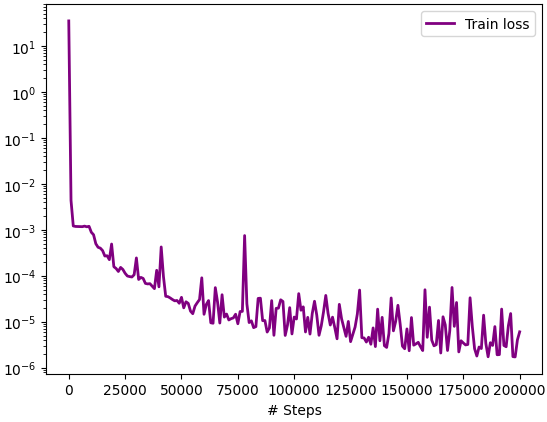}
	\includegraphics[width=0.53\textwidth]{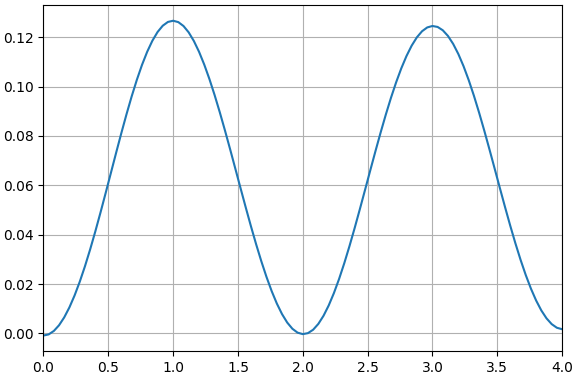}
	\caption{
		% \red{I AM HERE, 23.9.8. Rewrite below.}
		\JAX.
		\emph{Pinned-pinned bar, \red{No shift \& amplification}}, 
		Remark~\ref{rm:shift-amplification-Form-1}.
		\red{$\star$}
		\emph{Left:} Loss function.
		\emph{Right:}
		Step 198,000, lowest loss \red{1.724e-06}, midspan displacement,
		times at local maxima (1., 3.), times at local minima (0., 2., 4.),  damping\%=1.7\%, Quality: Good.  
		\red{$\star$}
		Section~\ref{sc:axial-Form-1}, \red{Form~1}.
		\emph{Network:}
		Remarks~\ref{rm:parameter-names},~\ref{rm:data-point-grids},
		T=4,
		W=64,
		H=4,
		% n\_out=1, 
		Uniform initializer, 
		12,737 parameters,
		random grid. 
		\emph{Training:}
		Remark~\ref{rm:learning-rate-schedule-4}, \red{LRS~4},
		init\_lr=0.002,
		factor\_lr=[0.9, 0.9, 0.9, 0.8, 0.8],
		n-cycles=5,
		N\_steps=200,000.
		Total GPU time \red{564 sec}.
		$\bullet$
		{\footnotesize
			Figure~\ref{fig:23.7.23 R1d A-PP v1 midspan,shape200000}, \DDET, shift and amplification.
			Figure~\ref{fig:DEBUG v1.9.3 seed42, 23.9.4 R1 A-P bcrP F1J(2a,3) lr0.001 Cy1-9-NCA Nsteps400000 regular N51 W64 H4 Glorot init-1}, \DDET, shift-amplification parameters.
			Figure~\ref{fig:DEBUG 23.9.4 R1d A-P bcrP F1J(2a,3) lr0.001 Cy1-9-NCA Nsteps400000 regular N51 W64 H4 Glorot init-1}, \DDET, shift-amplification parameters increase with training step number.
			$\triangleright$
			Figure~\ref{fig:random-grid-NO-static-solution}, \DDET, Form 2a, visually quasi-perfect midspan displacement.
		}
		{\scriptsize (23917R1d-1)} 
	}
	\label{fig:DEBUG 23.9.17 R1d JAX A-P bcrP F1H(2a,3) lr0.002 cy1-7-CA Nsteps300000 random N51 W64 H4 Uni init-1}
\end{figure}

\begin{figure}[tph]
	\includegraphics[width=0.49\textwidth]{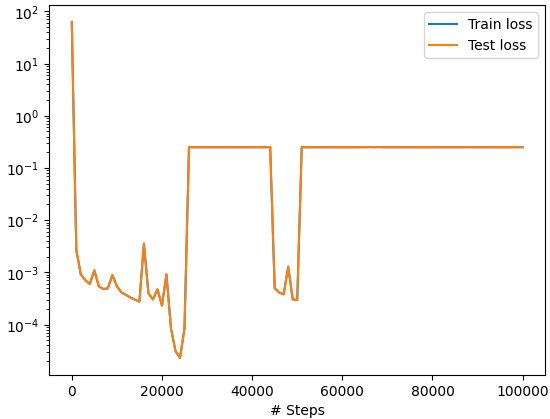}
	\includegraphics[width=0.49\textwidth]{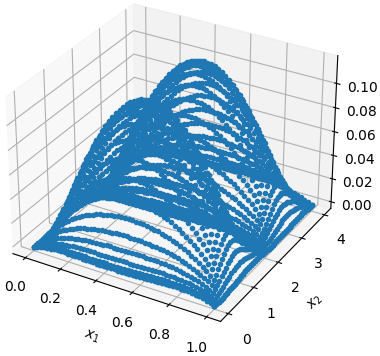}
	\caption{
		% \red{I AM HERE, 23.9.8. Rewrite below.}
		\DDET.
		\emph{Pinned-pinned bar, \red{early stopping}}. 
		% Section~\ref{sc:wave-equation}.
		\red{$\star$}
		% use \mbox{text} to prevent hyphenation of 2.30e-05 into "2.30e-" "05".
		\emph{Left:} Loss function, \red{\emph{lowest}} value \red{\mbox{2.30e-05}} % \newline 
		at Step \red{24,000},  value \red{0.25} at Step \red{100,000}, \red{\emph{divergence}}. 
		\emph{Right:} Shape, Step 24,000.
		\red{$\star$}
		Section~\ref{sc:axial-Form-1}, \red{Form 1}.
		\emph{Network:}
		Remarks~\ref{rm:parameter-names},~\ref{rm:data-point-grids}, 
		T=4,
		W=64,
		H=4,
		n\_out=1, 
		He-uniform initializer, 
		12,737 parameters,
		regular grid. 
		\emph{Training:}
		Remark~\ref{rm:learning-rate-schedule-1}, LRS~1 (CA),
		init\_lr=0.03.
		$\bullet$
		{\footnotesize
			Figure~\ref{fig:23.8.6 R2b A-PP F1(2a) lr0.03 cy3 Nsteps100000 regular N51 W64 H4 He init-2}, midspan displacement, Step 25,000, and zero midspan displacement, Step 100,000 (divergence).
			$\triangleright$
			Figure~\ref{fig:axial-Mathematica-solutions}, reference solution to compare.
		}
		{\scriptsize (2386R2b-1)} 
	}
	\label{fig:23.8.6 R2b A-PP F1(2a) lr0.03 cy3 Nsteps100000 regular N51 W64 H4 He init-1}
\end{figure}

\begin{figure}[tph]
	\includegraphics[width=0.49\textwidth]{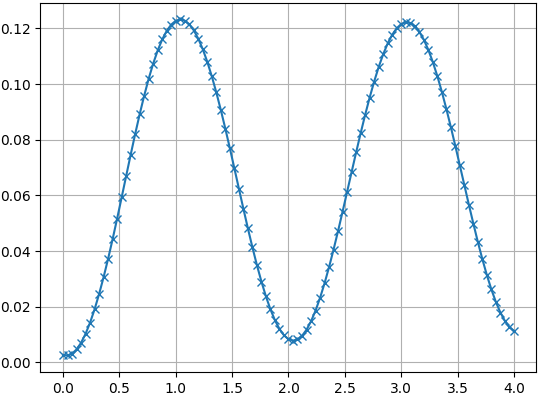}
	\includegraphics[width=0.49\textwidth]{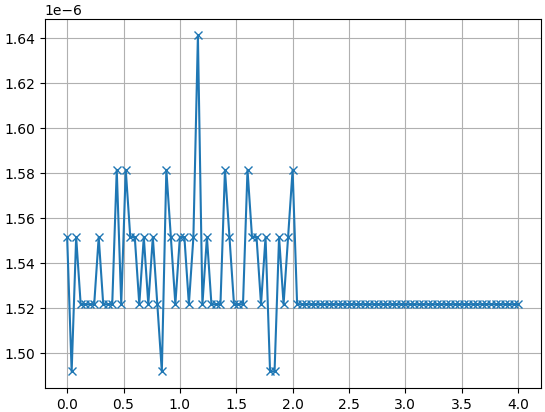}
	\caption{
		% \red{I AM HERE, 23.9.8. Rewrite below.}
		\DDET.
		\emph{Pinned-pinned bar}. 
		Midspan displacement, damping, Step 25,000 (left), 
		\emph{zero} midspan displacement, Step 100,000 (divergence, right).
		Section~\ref{sc:axial-Form-1}, \red{Form 1}.
		\emph{Network:}
		Remarks~\ref{rm:parameter-names},~\ref{rm:data-point-grids}, 
		T=4,
		W=64,
		H=4,
		n\_out=1, 
		He-uniform initializer, 
		12,737 parameters,
		regular grid.
		% {\color{red} REWRITE below.} 
		\emph{Training:}
		Remark~\ref{rm:learning-rate-schedule-1}, LRS~1 (CA),
		init\_lr=0.03.
		$\bullet$
		Figure~\ref{fig:23.8.6 R2b A-PP F1(2a) lr0.03 cy3 Nsteps100000 regular N51 W64 H4 He init-1}, loss function and shape, Step 24,000 with lowest loss.
		{\scriptsize (2386R2b-2)} 
	}
	\label{fig:23.8.6 R2b A-PP F1(2a) lr0.03 cy3 Nsteps100000 regular N51 W64 H4 He init-2}
\end{figure}

\begin{rem}
	\label{rm:early-stopping}
	% Pinned-pinned bar, Form 1:
	% DDE-T, 
	% \DDET,
	\DvJ.
	Form 1, pinned-pinned bar:
	Early stopping and divergence.
	{\rm
		After the lowest loss was reached earlier in the training process, and the loss function stopped to decrease after this lowest value, the training process can then be stopped earlier.  Depending on the behavior of the loss function after the lowest loss, the training could diverge.  See also ``early stopping'' in \cite{vuquoc2023deep}.
		
		Using \DDET, Figure~\ref{fig:23.8.6 R2b A-PP F1(2a) lr0.03 cy3 Nsteps100000 regular N51 W64 H4 He init-1} shows the loss function (left) with the lowest value of 2.30e-05 remaining at Step 24,000 as the training progressed to Step 100,000, with plateaus at loss value of 0.25, many orders of magnitude above the lowest loss, showing divergence.
		Once a plateau appears, there is no need to continue further, and stop the training process early, e.g., at Step 50,000 (or before).
		On the right subfigure is the shape time history at Step 24,000 (the lowest loss).
		
		Figure~\ref{fig:23.8.6 R2b A-PP F1(2a) lr0.03 cy3 Nsteps100000 regular N51 W64 H4 He init-2} (\DDET) shows a good midspan displacement with damping at Step 25,000, the last step at the end of Cycle 1 (left) and zero midspan displacement at Step 100,000 (divergence).
		
		Using \JAX, a similar behavior of ``early stopping'' was also observed.
	}
	\hfill$\blacksquare$
\end{rem}

% 23.8.6 R2b A-PP F1(2a) lr0.03 cy3 Nsteps100000 regular N51 W64 H4 He init - center disp50000.png

\begin{figure}[tph]
	\begin{subfigure}{0.24\textwidth}
		\includegraphics[width=\textwidth]{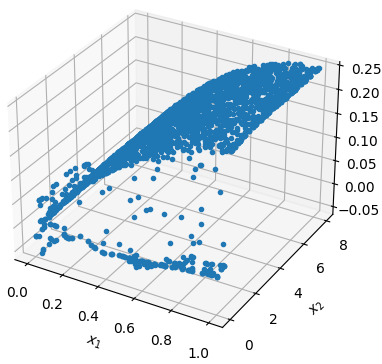}
		\caption{
			% Pinned-free bar, static solution.
			{\scriptsize 12,737 params, Shape 25000}
		}
		\label{fig:23.9.8 R2a shape25000}
	\end{subfigure}
	\begin{subfigure}{0.24\textwidth}
		\includegraphics[width=\textwidth]{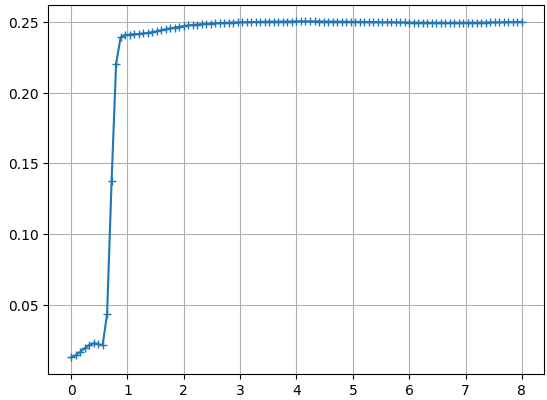}
		\caption{
			% Pinned-free bar, static solution.
			{\scriptsize 12,737 params, Free-end disp}
		}
		\label{fig:23.9.8 R2a free-end disp25000}
	\end{subfigure}
	\begin{subfigure}{0.24\textwidth}
		\includegraphics[width=\textwidth]{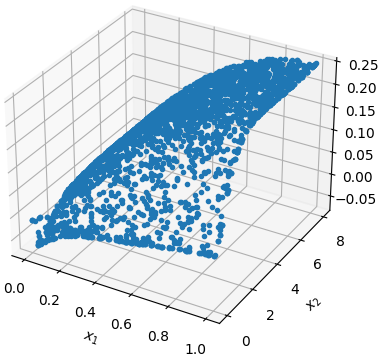}
		\caption{
			% Pinned-free bar, static solution.
			{\scriptsize 1,185 params, Shape 25000}
		}
		\label{fig:23.9.8 R2b shape25000}
	\end{subfigure}
	\begin{subfigure}{0.24\textwidth}
		\includegraphics[width=\textwidth]{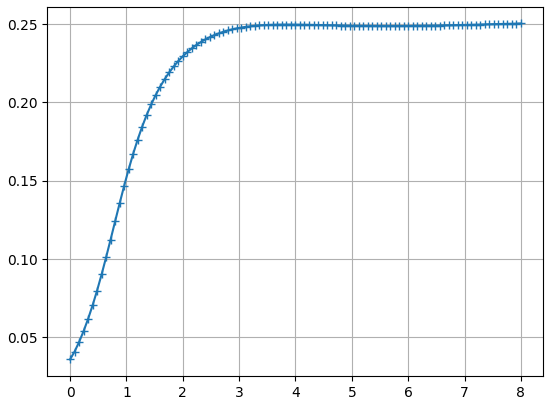}
		\caption{
			% Pinned-free bar, static solution.
			{\scriptsize 1,185 params, Free-end disp}
		}
		\label{fig:23.9.8 R2b free-end disp25000}
	\end{subfigure}
	\hfill
	\begin{subfigure}{0.24\textwidth}
		\includegraphics[width=\textwidth]{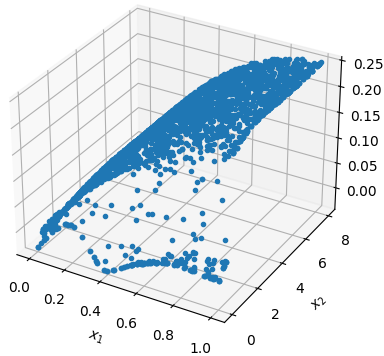}
		\caption{
			% Pinned-free bar, static solution.
			{\scriptsize 337 params, Shape 25000}
		}
		\label{fig:23.9.8 R2c shape25000}
	\end{subfigure}
	\begin{subfigure}{0.24\textwidth}
		\includegraphics[width=\textwidth]{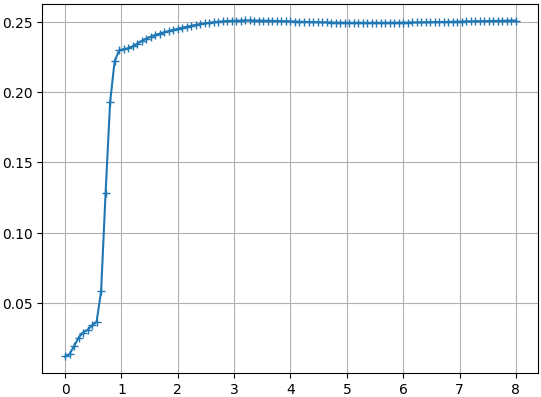}
		\caption{
			% Pinned-free bar, static solution.
			{\scriptsize 337 params, Free-end disp}
		}
		\label{fig:23.9.8 R2c free-end disp25000}
	\end{subfigure}
	\begin{subfigure}{0.24\textwidth}
		\includegraphics[width=\textwidth]{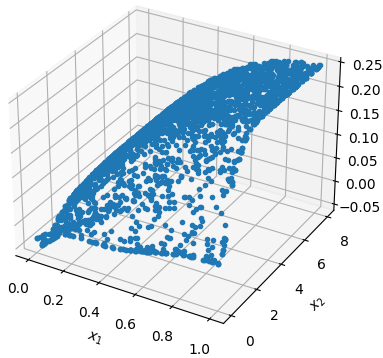}
		\caption{
			% Pinned-free bar, static solution.
			{\scriptsize 105 params, Shape 25000}
		}
		\label{fig:23.9.8 R2d shape25000}
	\end{subfigure}
	\begin{subfigure}{0.24\textwidth}
		\includegraphics[width=\textwidth]{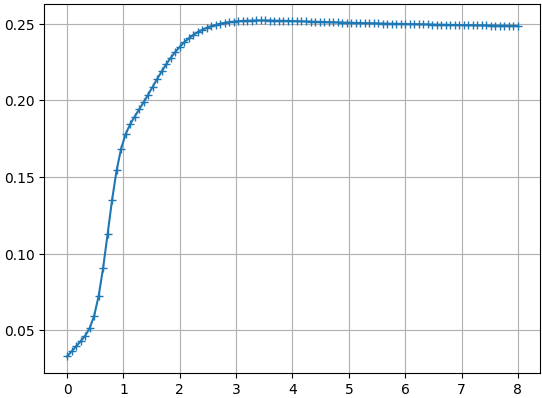}
		\caption{
			% Pinned-free bar, static solution.
			{\scriptsize 105 params, Free-end disp}
		}
		\label{fig:23.9.8 R2d free-end disp25000}
	\end{subfigure}
	\caption{
		% \DDET.
		\emph{\DDET, \red{Form 1}, pinned-free bar, \red{static} solutions}  
		(Remark~\ref{rm:static-solution-Form-1-pinned-free-2}, Section~\ref{sc:axial-Form-1}).
		$\bullet$
		SubFigs.~\ref{fig:23.9.8 R2a shape25000}-\ref{fig:23.9.8 R2a free-end disp25000}:
		% 12,737 parameters,
		Shape, free-end displacement, Step 25000; GPU time \red{655 sec} (Deterministic mode) for 200,000 steps.
		{\scriptsize (2398R2a)}.
		% Similarly for other image pairs:
		Shape, free-end disp, step number for other pairs: 
		(\ref{fig:23.9.8 R2b shape25000}-\ref{fig:23.9.8 R2b free-end disp25000}) 
		{\scriptsize (2398R2b)},
		(\ref{fig:23.9.8 R2c shape25000}-\ref{fig:23.9.8 R2c free-end disp25000})
		{\scriptsize (2398R2c)},
		(\ref{fig:23.9.8 R2d shape25000}-\ref{fig:23.9.8 R2d free-end disp25000})
		{\scriptsize (2398R2d)}.
		% SubFigs.~\ref{fig:23.9.8 R2b shape25000}-\ref{fig:23.9.8 R2b free-end disp25000}:
		% 1,185 params. 
		% Idem.
		% SubFigs.~\ref{fig:23.9.8 R2c shape25000}-\ref{fig:23.9.8 R2c free-end disp25000}:
		% 337 params.
		$\bullet$
		\emph{Network:} 
		Remarks~\ref{rm:parameter-names},~\ref{rm:data-point-grids}, 
		T=8,
		W=64,
		H=2,
		n\_out=1, 
		Glorot-uniform initializer, 
		% 12,737 parameters,
		random grid. 
		\emph{Training:}
		Remark~\ref{rm:learning-rate-schedule-1} LRS~1,
		init\_lr=0.001, N\_steps=25000.	
		$\bullet$
		{\footnotesize
			Figure~\ref{fig:JAX-Form-1-pinned-free-bar-NO-static-solutions}, using our \JAX\ script, no static solution.
			$\triangleright$
			Figure~\ref{fig:23.9.5 R3b A-P bcrF F(1J)2a(3) lr0.005 cy1-5-NCA Nsteps200000 random N51 W64 H4 He init-1}, \red{Form 2a}, GPU time \red{605 sec} for 200,000 steps.  
			% Rerun 23.9.5 R3b in [1] F2a Axial motion pinned-free.
			Figure~\ref{fig:23.9.5 R3d A-P bcrF F(1J,2a)3 lr0.005 cy1-5-NCA Nsteps200000 random N51 W64 H4 He init-1}, \red{Form~3}, GPU time \red{537 sec} for 200,000 steps.
			$\triangleright$
			Figure~\ref{fig:axial-Mathematica-solutions}, reference solution to compare.
			$\triangleright$
			Appendix~\ref{app:barrier-functions}, Figure~\ref{fig:barrier-function 23.10.5 R1c}, barrier-function effects.
		}	
	}
	\label{fig:DDE-T-Form-1-pinned-free-bar-static-solutions}
\end{figure}

\begin{figure}[tph]
	\begin{subfigure}{0.24\textwidth}
		%_\includegraphics[width=\textwidth]{Figures/23.9.15_R1_JAX_A-P_bcrF_F1H_2a,3__lr0.005_Nsteps100000_random_N51_W64_H4_uniform_init_-_shape100000.png}
		\includegraphics[width=\textwidth]{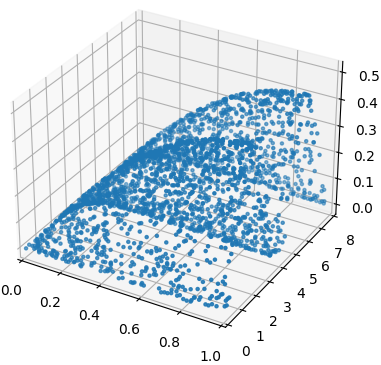}
		\caption{
			{\scriptsize 12,737 params, Shape 49000}
		}
		\label{fig:23.9.16 R1a.4 shape49000}
	\end{subfigure}
	\begin{subfigure}{0.24\textwidth}
		%_\includegraphics[width=\textwidth]{Figures/23.9.15_R1_JAX_A-P_bcrF_F1H_2a,3__lr0.005_Nsteps100000_random_N51_W64_H4_uniform_init_-_free-end_disp100000.png}
		\includegraphics[width=\textwidth]{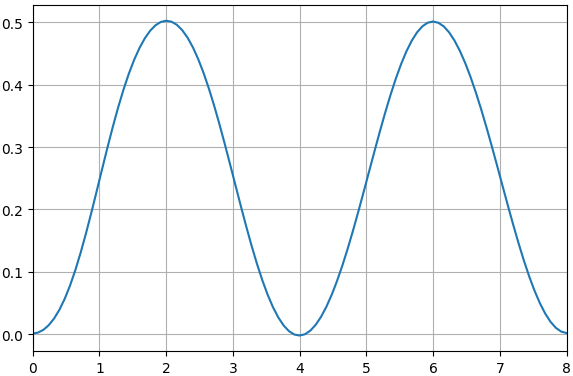}
		\caption{
			% {\scriptsize Free-end disp t.h.}
			% {\scriptsize 12,737 params, Free-end disp t.h.}
			{\scriptsize 12,737 params, Free-end disp}
		}
		\label{fig:23.9.16 R1a.4 free-end disp49000}
	\end{subfigure}
	\begin{subfigure}{0.24\textwidth}
		% replaced this result in 9.16 R1b by better result in 9.16 R1b.1.2
		%_\includegraphics[width=\textwidth]{Figures/DEBUG_23.9.16_R1b_JAX_A-P_bcrF_F1H_2a,3__lr0.005_cy1-3-CA_Nsteps100000_random_N51_W32_H2_Uni_init_-_shape100000.png}
		\includegraphics[width=\textwidth]{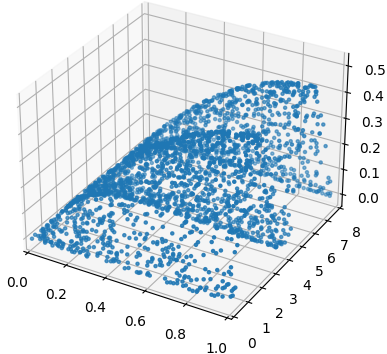}
		\caption{
			% {\scriptsize 1,185 params, Shape 100000}
			{\scriptsize 1,185 params, Shape 98000}
		}
		\label{fig:23.9.16 R1b.1 shape100000}
	\end{subfigure}
	\begin{subfigure}{0.24\textwidth}
		% replaced this result in 9.16 R1b by better result in 9.16 R1b.1.2
		%_\includegraphics[width=\textwidth]{Figures/DEBUG_23.9.16_R1b_JAX_A-P_bcrF_F1H_2a,3__lr0.005_cy1-3-CA_Nsteps100000_random_N51_W32_H2_Uni_init_-_free-end_disp100000.png}
		\includegraphics[width=\textwidth]{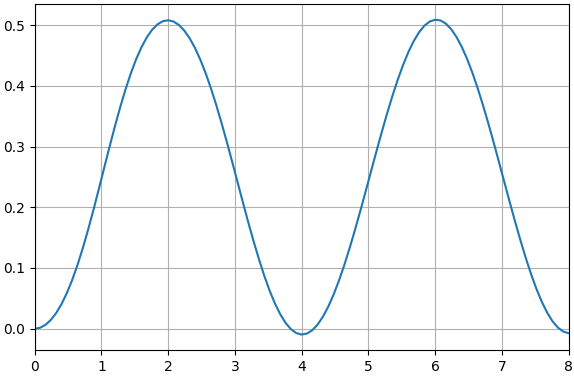}
		\caption{
			% {\scriptsize Free-end disp t.h.}
			{\scriptsize 1,185 params, Free-end disp}
		}
		\label{fig:23.9.16 R1b.1 free-end disp100000}
	\end{subfigure}
	\hfill
	\begin{subfigure}{0.24\textwidth}
		% replaced this result in 9.16 R1c.2 by better result in 9.16 Rc2.2
		%_\includegraphics[width=\textwidth]{Figures/DEBUG_23.9.16_R1c.2_JAX_A-P_bcrF_F1H_2a,3__lr0.002_cy1-9-CA_Nsteps400000_random_N51_W16_H2_Uni_init_-_shape400000.png}
		\includegraphics[width=\textwidth]{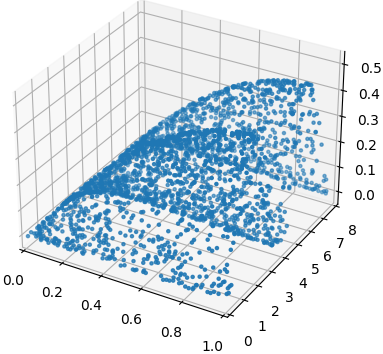}
		\caption{
			{\scriptsize 337 params, Shape 350000}
		}
		\label{fig:23.9.16 R1c.2 shape400000}
	\end{subfigure}
	\begin{subfigure}{0.24\textwidth}
		% replaced this result in 9.16 R1c.2 by better result in 9.16 Rc2.2
		%_\includegraphics[width=\textwidth]{Figures/DEBUG_23.9.16_R1c.2_JAX_A-P_bcrF_F1H_2a,3__lr0.002_cy1-9-CA_Nsteps400000_random_N51_W16_H2_Uni_init_-_free-end_disp400000.png}
		\includegraphics[width=\textwidth]{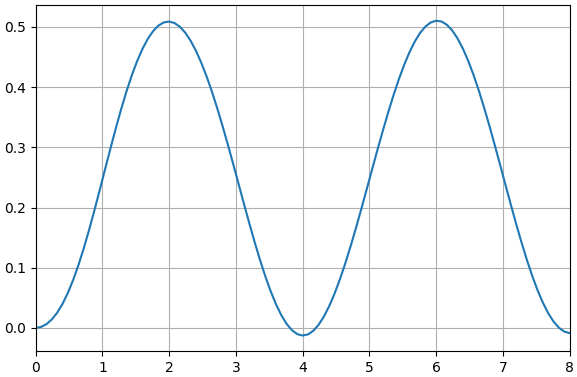}
		\caption{
			% {\scriptsize Free-end disp t.h.}
			{\scriptsize 337 params, Free-end disp}
		}
		\label{fig:23.9.16 R1c.2 free-end disp400000}
	\end{subfigure}
	\begin{subfigure}{0.24\textwidth}
		% replaced this result in 9.16 R1d.2 by better result in 9.16 R1d.2.3
		%_\includegraphics[width=\textwidth]{Figures/DEBUG_23.9.16_R1d.2_JAX_A-P_bcrF_F1H_2a,3__lr0.003_cy1-12-CA_Nsteps550000_random_N51_W8_H2_Uni_init_-_shape550000.png}
		\includegraphics[width=\textwidth]{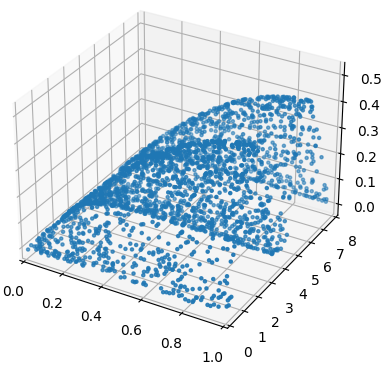}
		\caption{
			{\scriptsize 105 params, Shape 599000}
		}
		\label{fig:23.9.16 R1d.2.3 shape399000}
	\end{subfigure}
	\begin{subfigure}{0.24\textwidth}
		% replaced this result in 9.16 R1d.2 by better result in 9.16 R1d.2.3
		%_\includegraphics[width=\textwidth]{Figures/DEBUG_23.9.16_R1d.2_JAX_A-P_bcrF_F1H_2a,3__lr0.003_cy1-12-CA_Nsteps550000_random_N51_W8_H2_Uni_init_-_free-end_disp550000.png}
		\includegraphics[width=\textwidth]{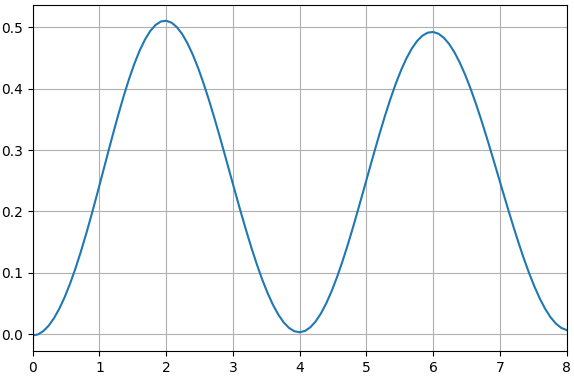}
		\caption{
			% {\scriptsize Free-end disp t.h.}
			{\scriptsize 105 params, Free-end disp}
		}
		\label{fig:23.9.16 R1d.2.3 free-end disp399000}
	\end{subfigure}
	\caption{
		\JAX.
		% \emph{\red{Form 1}, pinned-free bar, quasi-perfect to good solutions}.
		\emph{\red{Form 1}, pinned-free bar}.  
		% JAX.
		Remark~\ref{rm:learning-rate-schedule-4} (LRS~4).
		$\bullet$
		SubFigs~\ref{fig:23.9.16 R1a.4 shape49000}-\ref{fig:23.9.16 R1a.4 free-end disp49000}: 
		12,737 params, Step 49,000,
		Cycle 2 lowest loss 5.324e-06, 
		damping\%=0.25\%, \red{Quasi-perfect}, GPU time 161 sec.  
		\emph{Network:} 
		Remarks~\ref{rm:parameter-names},~\ref{rm:data-point-grids}, 
		T=8,
		W=64,
		H=2,
		n\_out=1, 
		Uniform initializer, 
		% 12,737 parameters,
		random grid. 
		\emph{Training:}
		Remark~\ref{rm:learning-rate-schedule-4} LRS~4,
		init\_lr=0.002,
		factor\_lr=[0.9, 0.9, 0.9, 0.8, 0.8, 1, 1, 1, 1].
		{\scriptsize (23916R1a.4)}. 
		% $\bullet$
		\red{$\star$}~{\color{OliveGreen} \mbox{\bf SubFigs}}~\ref{fig:23.9.16 R1b.1 shape100000}-\ref{fig:23.9.16 R1b.1 free-end disp100000}:
		% Cycle 3 lowest loss 3.864e-06, 
		% Damping\%=1.4\%, \red{very good}, GPU time 110 sec
		Cycle 3 lowest loss 4.381e-06, 
		damping\%=0.2\%, \red{Quasi-perfect}, GPU time 107 sec
		{\scriptsize (23916R1b.1)}.
		$\bullet$
		SubFigs~\ref{fig:23.9.16 R1c.2 shape400000}-\ref{fig:23.9.16 R1c.2 free-end disp400000}: 
		Cycle 8 lowest loss \mbox{5.445e-06}, 
		damping\%=\red{0.0014\%}, \red{Quasi-perfect}, GPU time 335 sec 
		{\scriptsize (23916R1c.2.2)}.
		$\bullet$
		SubFigs~\ref{fig:23.9.16 R1d.2.3 shape399000}-\ref{fig:23.9.16 R1d.2.3 free-end disp399000}:
		Cycle 13 lowest loss 2.749e-05, 
		damping\%=\red{3.7\%}, \red{High damping}, GPU time 505 sec
		{\scriptsize (23916R1d.2.3)}.
		$\bullet$
		{\footnotesize
			Figure~\ref{fig:DDE-T-Form-1-pinned-free-bar-static-solutions}, using DDE-T, static solutions.
			$\triangleright$
			Figure~\ref{fig:axial-Mathematica-solutions}, reference solution to compare.
		}
		% \red{23.9.24, TO WRITE.}
	}
	\label{fig:JAX-Form-1-pinned-free-bar-NO-static-solutions}
\end{figure}

\begin{rem}
	\label{rm:static-solution-Form-1-pinned-free-2}
	\label{rm:static-solution-Form-1-DDE-T-v-JAX}
	% {\color{Bittersweet} DDE-T}
	% {\color{RedOrange} DDE-T} vs {\color{purple} JAX}, 
	% \DDET, \JAX, 
	\DvJ, 
	Form 1, pinned-free bar:
	Static solution.
	{\rm
		Our DDE-T script for Form~1 produced a static-solution time history in every case.
		Our JAX script does not exhibit this pathological problem.
		% ALEX: removed paragraph
		For the pinned-free bar subjected to a constant distributed axial load, 
		DDE-T produced a static solution regardless of the size of the learning rate and the network capacity (Figure~\ref{fig:DDE-T-Form-1-pinned-free-bar-static-solutions}), whereas  
		JAX produced \red{no} static solution, but the expected waves with two vibration periods (Figure~\ref{fig:JAX-Form-1-pinned-free-bar-NO-static-solutions}), as obtained with Mathematica in Figure~\ref{fig:axial-Mathematica-solutions}.
		
		Using DDE-T, several different sets of parameters were considered with the number of network parameters varying from 12,737 down to 33, and with two initial learning rates, 0.001 and 0.0001.
		The results from four different network-parameter numbers (from 12,737 down to 105) and with the initial learning rate of 0.001 are given in
		Figure~\ref{fig:DDE-T-Form-1-pinned-free-bar-static-solutions}, in which all subfigures had these common parameters, which are not repeated in the figure caption: T=8, N\_out=1, Glorot-uniform initializer, random grid.
		Idential results (i.e., static solutions) were also obtained with the smaller initial learning rate 0.0001.
	}
	\hfill$\blacksquare$
\end{rem}

% pinned-free, DDE, static solution

% pinned-free bar, JAX, waves

\subsubsection{Form 2a: No time shift, static solution, efficiency}
\label{sc:axial-Form-2a}
\label{sc:pinned-pinned-Form-2a}
\noindent
The PINN Form 2a for the axial equation of motion of an elastic bar is a particular case of the PINN Form 2a for the Kirchhoff-Love rod in Section~\ref{sc:Kirchhoff-rod-Form-2a}.
\begin{align}
	\sldn 
	\ubp 2 
	+ \dfbs{\X} 
	= \pbsd{\X}
	\ , \quad
	\ndt \ub = \pbs{\X}
	\ ,
	\label{eq:wave-eq-form-2}
\end{align}
with the initial conditions:
\begin{align}
	\ub (\X , 0) = \ubs{0} 
	\ , \quad
	\pbs{\X} (\X , 0) = \ubds{0}
	\ .
\end{align}
In the numerical examples, we use homogeneous initial conditions:
\begin{align}
	\ub (\X , 0) = 0 
	\ , \quad
	\pbs{\X} (\X , 0) = 0
	\ .
	\label{eq:axial-ICs-Form-2a}
\end{align}

\begin{figure}[tph]
	\includegraphics[width=0.49\textwidth]{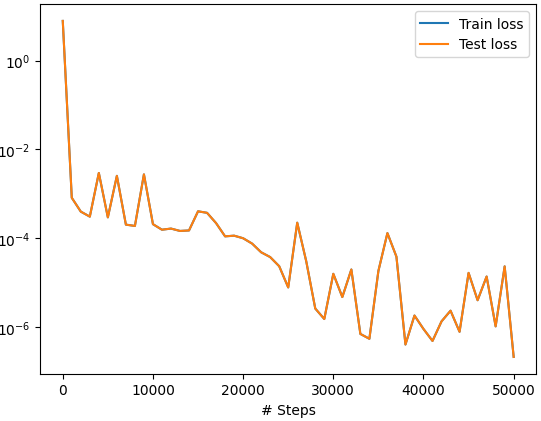}
	\includegraphics[width=0.49\textwidth]{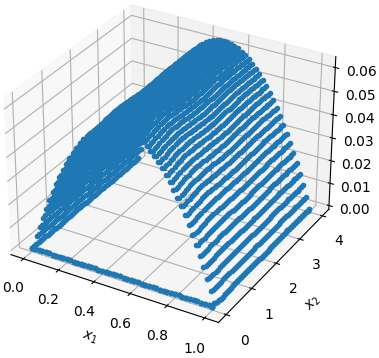}
	\caption{
		\DDET.
		\emph{Pinned-pinned bar, static solution}. 
		Loss function (left), static solution (right), Step 50,000, end of Cycle 2.
		\red{$\star$}
		Remark~\ref{rm:static-solution}, Static solution; Remark~\ref{rm:static-solution-avoid}, How to avoid.
		% Section~\ref{sc:networks-grids}.
		Section~\ref{sc:axial-Form-2a}, \red{Form 2a}.		 
		% Computational domain: $(x, t) \in [0, 1] \times [0, 4]$.
		\emph{Network:}
		Remarks~\ref{rm:parameter-names}, \ref{rm:data-point-grids},
		% n\_inp=2,
		T=4,
		W=64,
		H=4, 
		n\_out=2, 
		He-uniform initializer, 12,802 parameters,
		\red{$\star$}
		% Remark~\ref{rm:data-point-grids}, 
		\emph{regular} grid.
		\emph{Training:}
		Remark~\ref{rm:learning-rate-schedule-3}, LRS~3,
		init\_lr=0.01,
		n-cycles=2, N\_steps=50,000.
		$\bullet$
		{\footnotesize
			Figure~\ref{fig:pinned-pinned-static-solution-2}, velocity, midspan displacement.
			$\triangleright$
			Figure~\ref{fig:23.8.22 R3b A-PP F(1J)2a(3) lr0.001 cy5-CA cy6-NCA Nsteps250000 regular N51 W64 H4 He init}, regular grid, smaller init\_lr=0.001.
			Figure~\ref{fig:random-grid-NO-static-solution}, \emph{random} grid, no static solution (same LRS 3 and init\_lr).
		}
		{\scriptsize (23725R1-1)}
	}
	\label{fig:pinned-pinned-static-solution-1}
\end{figure}

\begin{figure}[tph]
	% NOTE: 2023.08.14, the OLD image filename said "cycles=5 NSteps=200000," but
	% actually it was "cy2-NCA NSteps50000" = 2 cycles with no cyclic annealing.
	% script file:
	% 23.7.25 R1 A-PP F2a lr0.01 cy2-NCA Nsteps=50000 regular N51 W64 H4 He init.ipynb
	% https://colab.research.google.com/drive/1aBEIjBdaki_oR3sWDAOBFPUuoMJI7kFX 
	%
	\includegraphics[width=0.49\textwidth]{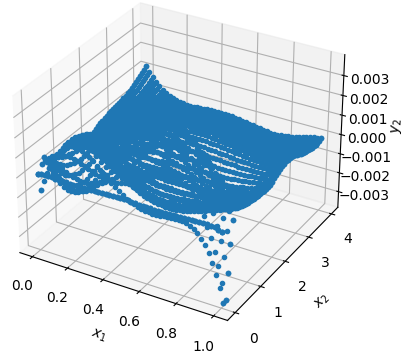}
	\includegraphics[width=0.49\textwidth]{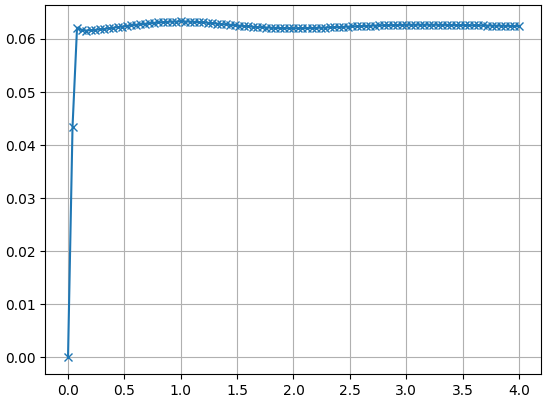}
	\caption{
		\DDET.
		% \red{I AM HERE, 23.9.8. Rewrite below.}
		\emph{Pinned-pinned bar, static solution}. 
		Essentially-zero velocity (left); \emph{static}  midspan displacement (right), Step 50,000, end of Cycle 2.
		\red{$\star$}
		Remark~\ref{rm:static-solution}, Static solution; Remark~\ref{rm:static-solution-avoid}, How to avoid.
		Section~\ref{sc:axial-Form-2a}, \red{Form~2a}.		 
		\emph{Network:}
		Remarks~\ref{rm:parameter-names}, \ref{rm:data-point-grids},
		% n\_inp=2,
		T=4,
		W=64,
		H=4, 
		n\_out=2, 
		He-uniform initializer, 12,802 parameters,
		\red{$\star$}
		% Remark~\ref{rm:data-point-grids}, 
		\emph{regular} grid.
		\emph{Training:}
		Remark~\ref{rm:learning-rate-schedule-3}, LRS~3,
		init\_lr=0.01,
		n-cycles=2, N\_steps=50,000.
		$\bullet$
		{\footnotesize
			Figure~\ref{fig:pinned-pinned-static-solution-1}, loss function, shape.
			$\triangleright$
			Figure~\ref{fig:23.8.22 R3b A-PP F(1J)2a(3) lr0.001 cy5-CA cy6-NCA Nsteps250000 regular N51 W64 H4 He init}, regular grid, smaller init\_lr.
			Figure~\ref{fig:random-grid-NO-static-solution}, \emph{random} grid, no static solution (same LRS 3 and init\_lr).
		}
		{\scriptsize (23725R1-2)}
	}
	\label{fig:pinned-pinned-static-solution-2}
\end{figure}

\begin{figure}[tph]
	\includegraphics[width=0.49\textwidth]{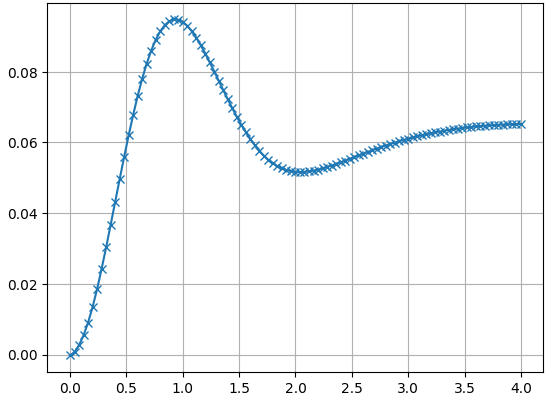}
	\includegraphics[width=0.49\textwidth]{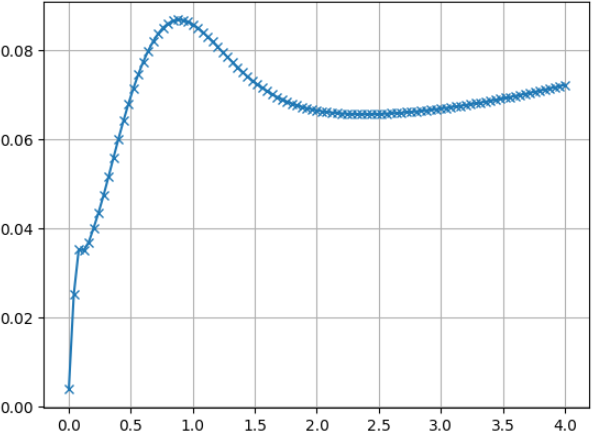}
	\caption{
		% \emph{Regular grid and pre-static time histories.} Remark~\ref{rm:data-point-grids}.  A regular grid would detect the static solution easier than a random grid.
		% \red{I AM HERE, 23.9.8. Rewrite below.}
		\DDET.
		\emph{Pinned-pinned bar}. 
		\emph{Pre-static} midspan-displacement time histories at Step 25,000 (end of Cycle 1).
		\red{$\star$}
		Remarks~\ref{rm:static-solution},~\ref{rm:static-solution-avoid}.
		% Section~\ref{sc:networks-grids}.
		Section~\ref{sc:axial-Form-2a}, \red{Form 2a}.
		\red{$\star$}
		\emph{Network:}
		Remark~\ref{rm:parameter-names},~\ref{rm:data-point-grids}, 
		T=4,
		W=64,
		H=2, 
		n\_out=2, 
		He-uniform initializer, 
		\red{4,482} parameters,
		\emph{regular} grid.
		\emph{Training:}
		Remark~\ref{rm:learning-rate-schedule-3}, LRS~3, init\_lr=0.03.
		\red{$\star$}
		\emph{Left:} n-cycles=5, 
		% no cyclic annealing, 
		N\_steps=200,000; Figures~\ref{fig:23.8.15 R3a A-PP F(1)2a lr0.03 cy5-NCA Nsteps200000 regular N51 W64 H2 He init-1}-\ref{fig:23.8.15 R3a A-PP F(1)2a lr0.03 cy5-NCA Nsteps200000 regular N51 W64 H2 He init-2}, Step 200,000.
		{\scriptsize (23815R3a-1)}. 
		\red{$\star$}
		\emph{Right:} 
		n-cycles=1,
		N\_steps=25,000, 
		jump at t=0.
		{\scriptsize (2389R2-1)}.
	}
	\label{fig:pre-static-solution-jump}
\end{figure}

\begin{figure}[tph]
	\includegraphics[width=0.49\textwidth]{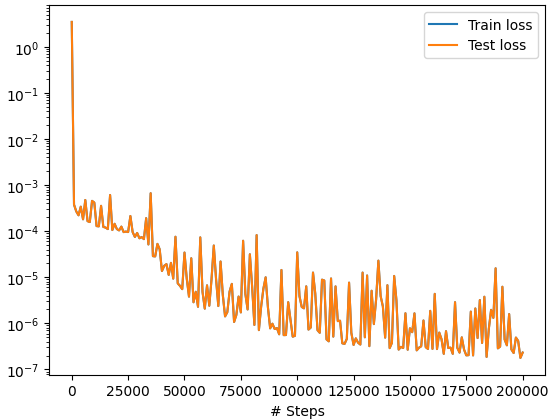}
	\includegraphics[width=0.49\textwidth]{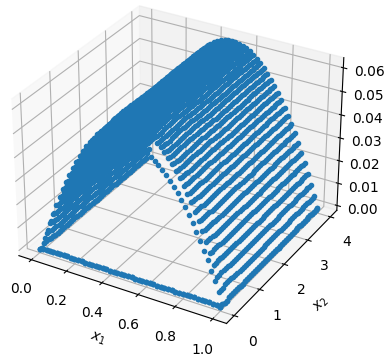}
	\caption{
		% NCA led to static solution.
		% Loss function and shape time history.
		% \red{I AM HERE, 23.9.8. Rewrite below.}
		\DDET.
		\emph{Pinned-pinned bar, no cyclic annealing} (NCA). 
		Loss function (left) and
		\emph{static-shape} (right), Step 200,000.
		\red{$\star$}
		Remarks~\ref{rm:static-solution},~\ref{rm:static-solution-avoid}.
		% Section~\ref{sc:networks-grids}.
		Section~\ref{sc:axial-Form-2a}, \red{Form 2a}.
		\emph{Network:}
		Remarks~\ref{rm:parameter-names},~\ref{rm:data-point-grids}, 
		T=4,
		W=64,
		H=2, 
		n\_out=2, 
		He-uniform initializer, 4,482 parameters,
		\emph{regular} grid.
		\red{$\star$}
		\emph{Training:}
		Remark~\ref{rm:learning-rate-schedule-3}, 
		LRS~3, 
		% no cyclic annealing (NCA),
		init\_lr=0.03.
		% n-cycles=5,  
		% N\_steps=200,000.
		% See
		$\bullet$
		{\footnotesize
			Figure~\ref{fig:pre-static-solution-jump}, left, midspan displacement, 
			Step 25,000. 
			Figure~\ref{fig:23.8.15 R3a A-PP F(1)2a lr0.03 cy5-NCA Nsteps200000 regular N51 W64 H2 He init-2}, velocity, midspan displacement, Step 200,000.
		}	
		{\scriptsize (23815R3a-2)}
	}
	\label{fig:23.8.15 R3a A-PP F(1)2a lr0.03 cy5-NCA Nsteps200000 regular N51 W64 H2 He init-1}
\end{figure}

\begin{figure}[tph]
	\includegraphics[width=0.49\textwidth]{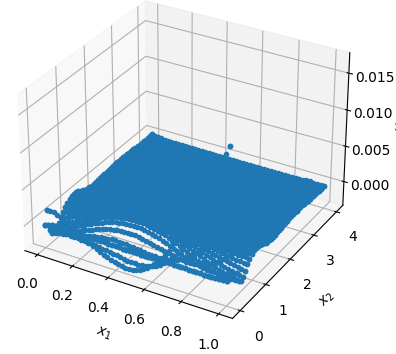}
	\includegraphics[width=0.49\textwidth]{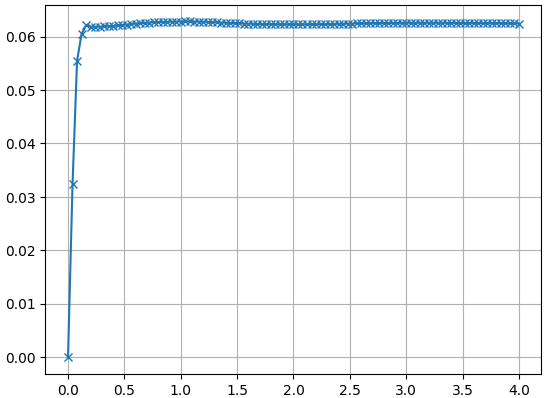}
	\caption{
		% NCA led to static solution.
		% Velocity and midspan displacement time histories.
		% \red{I AM HERE, 23.9.8. Rewrite below.}
		\DDET.
		\emph{Pinned-pinned bar, no cyclic annealing} (NCA). 
		\red{$\star$}
		Essentially-zero Velocity (left), 
		\red{static} midspan displacement (right), Step 200,000.
		\red{$\star$}
		Remarks~\ref{rm:static-solution},~\ref{rm:static-solution-avoid}.
		% Section~\ref{sc:networks-grids}.
		Section~\ref{sc:axial-Form-2a}, \red{Form~2a}.
		\emph{Network:}
		Remarks~\ref{rm:parameter-names},~\ref{rm:data-point-grids}, 
		T=4,
		W=64,
		H=2, 
		n\_out=2, 
		He-uniform initializer, 4,482 parameters,
		\emph{regular} grid.
		\emph{Training:}
		Remark~\ref{rm:learning-rate-schedule-3},
		LRS~3, 
		init\_lr=0.03.
		% n-cycles=5, no cyclic annealing, N\_steps=200,000.
		$\bullet$
		{\footnotesize
			Figure~\ref{fig:pre-static-solution-jump}, left, midspan displacement, Step 25,000. 
			Figure~\ref{fig:23.8.15 R3a A-PP F(1)2a lr0.03 cy5-NCA Nsteps200000 regular N51 W64 H2 He init-1}, loss function, static shape, Step 200,000.
		}
		{\scriptsize (23815R3a-3)}
	}
	\label{fig:23.8.15 R3a A-PP F(1)2a lr0.03 cy5-NCA Nsteps200000 regular N51 W64 H2 He init-2}
\end{figure}

\begin{figure}[tph]
	\includegraphics[width=0.49\textwidth]{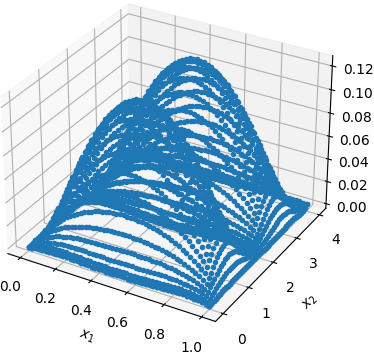}
	\includegraphics[width=0.49\textwidth]{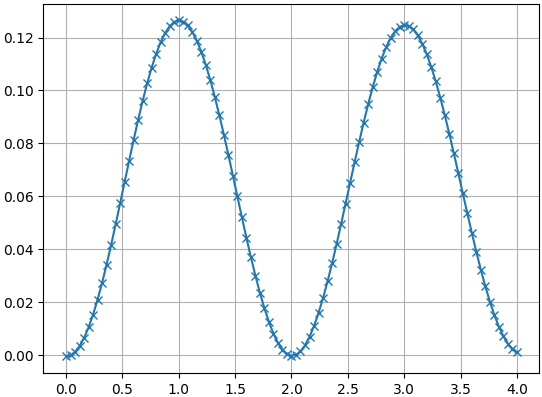}
	\caption{
		% 23.9.6, model figure caption
		% \red{quasi-perfect}, \red{good} solution
		\DDET.
		% v1.9.3.
		\emph{Pinned-pinned bar}. 
		Shape (left), midspan displacement (right), small damping, Step 250,000.
		\red{$\star$}
		% Remark~\ref{rm:static-solution-avoid}.
		Section~\ref{sc:axial-Form-2a}, \red{Form~2a}.
		\emph{Network:}
		Remarks~\ref{rm:parameter-names},~\ref{rm:data-point-grids},
		% n\_inp=2,
		T=4,
		W=64,
		H=4, 
		n\_out=2, 
		He-uniform initializer, \red{12,802} parameters,
		% Remark~\ref{rm:data-point-grids}, 
		\emph{regular} grid.		
		% \emph{Optimization:}
		\mbox{\emph{Training:}}
		\red{$\star$}
		Remarks~\ref{rm:learning-rate-schedule-1} (LRS~1),
		\ref{rm:extension-cycle-5} (ELRS),
		init\_lr=0.001,
		Cycles 1-5 (CA),
		Cycle 6 (NCA), 
		N\_steps=250,000.
		\red{$\star$}
		Lowest total loss \red{2.55e-06}, Step 250,000 (sum of 6 losses).
		$\bullet$
		{\footnotesize
			Figures~\ref{fig:pinned-pinned-static-solution-1}-\ref{fig:pinned-pinned-static-solution-2}, init\_lr=0.01, static solution.
			Figure~\ref{fig:23.8.16 R1a A-PP F(1)2a lr0.01 cy5-CA Nsteps200000 regular N51 W64 H2 He init}, 4,802 parameters, 
			good solution.
			$\triangleright$
			Figure~\ref{fig:axial-Mathematica-solutions}, reference solution to compare.
		}
		{\scriptsize (23822R3b-1)}
	}
	\label{fig:23.8.22 R3b A-PP F(1J)2a(3) lr0.001 cy5-CA cy6-NCA Nsteps250000 regular N51 W64 H4 He init}
\end{figure}

\begin{figure}[tph]
	% script file from which the figures below came:
	% 23.8.4 R2a A-PP F(1)2a lr0.01 cy1 Nsteps200000 fixed random N51 W64 H4 He init.ipynb
	% https://colab.research.google.com/drive/1PdZWHLNoJTo-jRQHUh6n6JOUgJbisp3V
	%
	% OLD DDE script, no damping%
%	\includegraphics[width=0.49\textwidth]{Figures/23.8.4_R2a_A-PP_F_1_2a_lr0.01_cy1_Nsteps200000_fixed_random_N51_W64_H4_He_init_-_shape200000.png}
%	\includegraphics[width=0.49\textwidth]{Figures/23.8.4_R2a_A-PP_F_1_2a_lr0.01_cy1_Nsteps200000_fixed_random_N51_W64_H4_He_init_-_center_disp200000.png}
	%
	\includegraphics[width=0.49\textwidth]{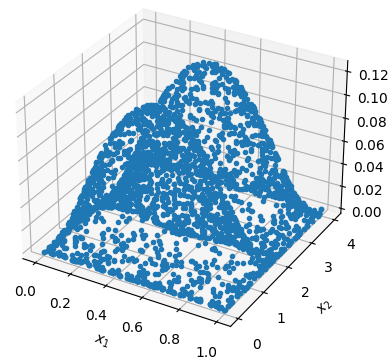}
	\includegraphics[width=0.49\textwidth]{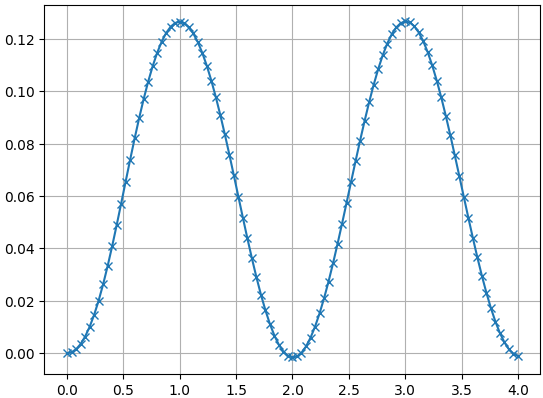}
	%
	% NEW DDE script, damping% and quality rating, 23.9.23
	%_\includegraphics[width=0.49\textwidth]{Figures/DEBUG_23.8.4_R2d_DDE_A-P_bcrP_F_1H_2a_3__lr0.01_cy1-5-CA_Nsteps200000_random_N51_W64_H4_He_init_-_shape146000.png}
	%_\includegraphics[width=0.49\textwidth]{Figures/DEBUG_23.8.4_R2d_DDE_A-P_bcrP_F_1H_2a_3__lr0.01_cy1-5-CA_Nsteps200000_random_N51_W64_H4_He_init_-_center_disp146000.png}
	\caption{
		% Random grid avoids static solution.
		\DDET.
		\emph{Pinned-pinned bar}. 
		Step 200,000:
		Shape (left), midspan displacement (right), 
		visually quasi-perfect.
		% NEW code
		% damping\%=0.3\%, 
		% \red{\emph{Quasi-perfect}}. 
		\red{$\star$}
		Remark~\ref{rm:optimal-capacity}, Optimal capacity.
		Remark~\ref{rm:unstable-solution}, Unstable solution.
		Remark~\ref{rm:static-solution-avoid}, Avoiding static solution.
		Section~\ref{sc:pinned-pinned-Form-2a}, {\color{red} Form~2a}.
		% Computational domain: $(x, t) \in [0, 1] \times [0, 4]$.
		\emph{Network:}
		\red{$\star$}
		Remarks~\ref{rm:parameter-names},~\ref{rm:data-point-grids}, 
		% n\_inp=2,
		T=4,
		W=64,
		H=4, 
		n\_out=2, 
		He-uniform initializer, 
		\red{12,802} parameters,
		\red{$\star$}
		\mbox{\red{\emph{random}}} grid.
		\emph{Training:}
		\red{$\star$}
		Remark~\ref{rm:learning-rate-schedule-3} 
		% (LRS~3),
		\red{LRS~3 (NCA)},
		init\_lr=0.01.
		n-cycles=1, N\_steps=200,000.
		\red{$\star$}
		% OLD DDE script, no damping%
		Lowest total loss {\color{red} \mbox{0.537e-06}}, Step 192,000 (sum of 6 losses).
		Total GPU time \red{640 sec}.
		%
		% NEW DDE script, damping% and quality rating, 23.9.23
		% Total loss {\color{red} \mbox{0.722e-06}}, Step 146,000 (sum of 6 losses).
		% Total GPU time \red{442 sec}. 
		% \protect\footnotemark\ 
		$\bullet$
		{\footnotesize
			Figure~\ref{fig:random-grid-NO-static-solution-2}, early stopping at Step 146,000, good trade-off.
			$\triangleright$
			Figures~\ref{fig:pinned-pinned-static-solution-1}-\ref{fig:pinned-pinned-static-solution-2}, \emph{regular} grid, static solution
			(same LRS 3 and init\_lr).
			Figure~\ref{fig:23.8.16 R1a A-PP F(1)2a lr0.01 cy5-CA Nsteps200000 regular N51 W64 H2 He init}, 4,802 parameters, good solution.
			Figure~\ref{fig:23.7.26 R5c A-PP F2a lr0.02x5 cycles=9 Nsteps=400000 regular N51 W32 H2 He init-3}, 1,218 parameters, \red{\emph{very-good}} solution. 
		}
		% OLD DDE script, no damping%
		{\scriptsize (2384R2a-1)}
		% NEW DDE script, damping% and quality rating, 23.9.23
		% {\scriptsize (2384R2d-1)}
	}
	\label{fig:random-grid-NO-static-solution}
\end{figure}

\begin{figure}[tph]
	% script file from which the figures below came:
	% 23.8.4 R2a A-PP F(1)2a lr0.01 cy1 Nsteps200000 fixed random N51 W64 H4 He init.ipynb
	% https://colab.research.google.com/drive/1PdZWHLNoJTo-jRQHUh6n6JOUgJbisp3V
	%
	% OLD DDE script, no damping%
	%_\includegraphics[width=0.49\textwidth]{Figures/23.8.4_R2a_A-PP_F_1_2a_lr0.01_cy1_Nsteps200000_fixed_random_N51_W64_H4_He_init_-_shape200000.png}
	%_\includegraphics[width=0.49\textwidth]{Figures/23.8.4_R2a_A-PP_F_1_2a_lr0.01_cy1_Nsteps200000_fixed_random_N51_W64_H4_He_init_-_center_disp200000.png}
	%
	% NEW DDE script, damping% and quality rating, 23.9.23
	%_\includegraphics[width=0.49\textwidth]{Figures/DEBUG_23.8.4_R2d_DDE_A-P_bcrP_F_1H_2a_3__lr0.01_cy1-5-CA_Nsteps200000_random_N51_W64_H4_He_init_-_shape146000.png}
	\includegraphics[width=0.47\textwidth]{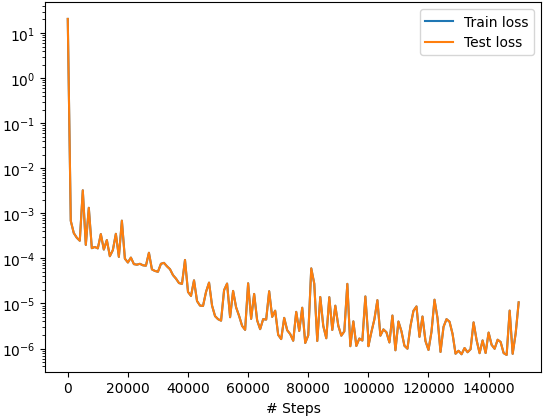}
	\includegraphics[width=0.51\textwidth]{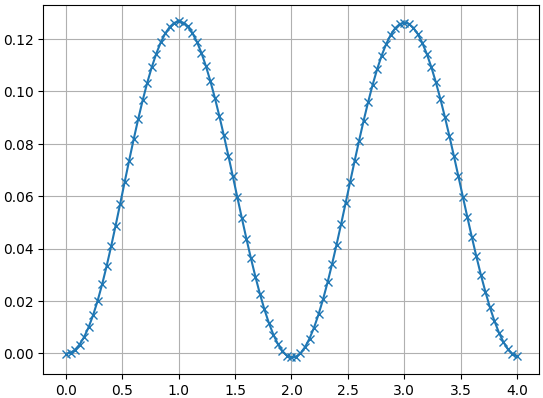}
	\caption{
		% Random grid avoids static solution.
		\DDET.
		\emph{Pinned-pinned bar}. 
		\emph{Left:} Loss function.
		\emph{Right:}
		Step \red{146,000},
		% Shape (left), 
		midspan displacement (right), damping\%=0.3\%, \red{\emph{Quasi-perfect}}. 
		\red{$\star$}
%		Remark~\ref{rm:optimal-capacity}, Optimal capacity.
%		Remark~\ref{rm:unstable-solution}, Unstable solution.
%		Remark~\ref{rm:static-solution-avoid}, Avoiding static solution.
		Section~\ref{sc:pinned-pinned-Form-2a}, {\color{red} Form~2a}.
		% Computational domain: $(x, t) \in [0, 1] \times [0, 4]$.
		\emph{Network:}
		\red{$\star$}
		Remarks~\ref{rm:parameter-names},~\ref{rm:data-point-grids}, 
		% n\_inp=2,
		T=4,
		W=64,
		H=4, 
		n\_out=2, 
		He-uniform initializer, 
		\red{12,802} parameters,
		\red{$\star$}
		\mbox{\red{\emph{random}}} grid.
		\emph{Training:}
		\red{$\star$}
		Remark~\ref{rm:learning-rate-schedule-3} 
		% (LRS~3),
		\red{LRS~3 (NCA)},
		init\_lr=0.01.
		% n-cycles=1, N\_steps=200,000.
		\red{$\star$}
		% OLD DDE script, no damping%
		% Lowest total loss {\color{red} \mbox{0.537e-06}}, Step 192,000 (sum of 6 losses).
		% Total GPU time \red{640 sec}.
		%
		% NEW DDE script, damping% and quality rating, 23.9.23
		Total loss {\color{red} \mbox{0.722e-06}}, Step 146,000 (sum of 6 losses).
		Total GPU time \red{442 sec}. 
		% \protect\footnotemark\ 
		$\bullet$
		{\footnotesize
			Figure~\ref{fig:random-grid-NO-static-solution}, running through to Step 200,000.
			$\triangleright$
			Figures~\ref{fig:pinned-pinned-static-solution-1}-\ref{fig:pinned-pinned-static-solution-2}, \emph{regular} grid, static solution
			(same LRS 3 and init\_lr).
			Figure~\ref{fig:23.8.16 R1a A-PP F(1)2a lr0.01 cy5-CA Nsteps200000 regular N51 W64 H2 He init}, 4,802 parameters, good solution.
			Figure~\ref{fig:23.7.26 R5c A-PP F2a lr0.02x5 cycles=9 Nsteps=400000 regular N51 W32 H2 He init-3}, 1,218 parameters, \red{\emph{very-good}} solution.
			$\triangleright$
			Figure~\ref{fig:axial-Mathematica-solutions}, reference solution to compare.
		}
		% OLD DDE script, no damping%
		% {\scriptsize (2384R2a-1)}
		% NEW DDE script, damping% and quality rating, 23.9.23
		{\scriptsize (2384R2d-1)}
	}
	\label{fig:random-grid-NO-static-solution-2}
\end{figure}

\begin{figure}[tph]
	% 23.9.29: replaced images from 9.5 R3b by images from 9.5 R3c 
	%_\includegraphics[width=0.49\textwidth]{Figures/23.9.5_R3b_A-P_bcrF_F_1J_2a_3__lr0.005_cy1-5-NCA_Nsteps200000_random_N51_W64_H4_He_init_-_shape200000.png}
	\includegraphics[width=0.49\textwidth]{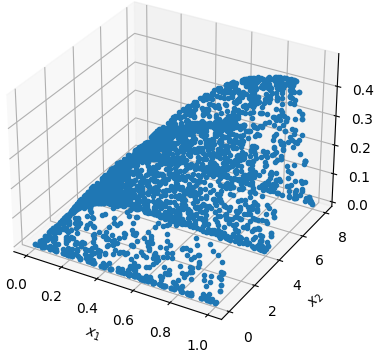}
	%
	%_\includegraphics[width=0.49\textwidth]{Figures/23.9.5_R3b_A-P_bcrF_F_1J_2a_3__lr0.005_cy1-5-NCA_Nsteps200000_random_N51_W64_H4_He_init_-_free-end_disp200000.png}
	\includegraphics[width=0.49\textwidth]{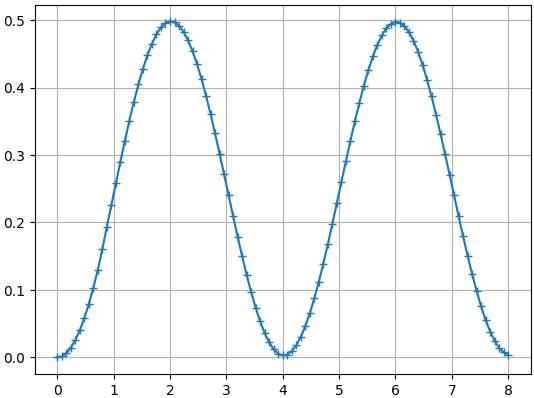}
	\caption{
		\DDET.
		% \red{I AM HERE, 23.9.8. Rewrite below.}
		\emph{Pinned-free bar. {He-uniform} initializer}. 
		\red{$\star$}
		Step 100,000:
		Shape (left), free-end displacement (right), damping\%=0.14\%, \red{Quasi-perfect}.
		Section~\ref{sc:axial-Form-2a}, \red{Form~2a}.
		\emph{Network:}
		Remarks~\ref{rm:parameter-names},~\ref{rm:data-point-grids}, 
		T=8,
		W=64,
		H=4, 
		n\_out=2, 
		He-uniform initializer, 12,802 parameters,
		\red{\emph{random}} grid.
		\emph{Training:}
		Remark~\ref{rm:learning-rate-schedule-3},
		LRS~3, 
		init\_lr=0.005,
		Cycle 3 lowest loss 2.472e-06, 
		GPU time 318 sec.
		{\scriptsize (2395R3b-1)}
		% (Cycle 5 lowest loss 1.268e-06, GPU time 605 sec.)
		$\bullet$
		{\footnotesize
			% \red{Change figures here.}
			(Cycle 5 lowest loss 1.268e-06, GPU time \red{605 sec} (Deterministic mode). {\scriptsize (2395R3b.3)})
			Figure~\ref{fig:DEBUG v1.9.3 seed42, 23.9.5 R3a A-P bcrF F(1J)2a(3) lr0.005 cy1-3-NCA Nsteps100000 regular N51 W64 H4 He init-1}, same model, LRS~3, init\_lr=0.005, \red{regular} grid, \emph{static} solution.
			$\triangleright$
			Figure~\ref{fig:DDE-T-Form-1-pinned-free-bar-static-solutions} \red{Form 1}, \red{static} solutions.
			Figure~\ref{fig:23.9.5 R3d A-P bcrF F(1J,2a)3 lr0.005 cy1-5-NCA Nsteps200000 random N51 W64 H4 He init-1}, \red{Form 3}, static and \red{quasi-perfect} solutions.
			$\triangleright$
			Figure~\ref{fig:axial-Mathematica-solutions}, reference solution to compare.
		}
	}
	\label{fig:23.9.5 R3b A-P bcrF F(1J)2a(3) lr0.005 cy1-5-NCA Nsteps200000 random N51 W64 H4 He init-1}
\end{figure}

\begin{figure}[tph]
	\begin{subfigure}{0.24\textwidth}
		\includegraphics[width=\textwidth]{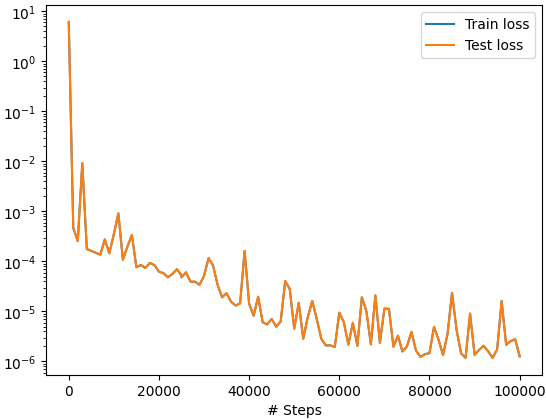}
		\caption{
			% Pinned-pinned bar, trade-off
			{\scriptsize Loss function}
		}
		\label{fig:23.8.4 R2e.1 loss88000.png}
	\end{subfigure}
	\begin{subfigure}{0.24\textwidth}
		\includegraphics[width=\textwidth]{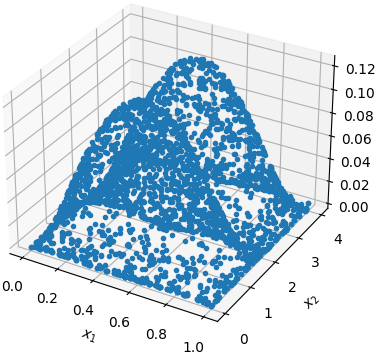}
		\caption{
			% Pinned-pinned bar, trade-off
			{\scriptsize Shape t.h.}
		}
		\label{fig:23.8.4 R2e.1 shape88000.png}
	\end{subfigure}
	\begin{subfigure}{0.24\textwidth}
		\includegraphics[width=\textwidth]{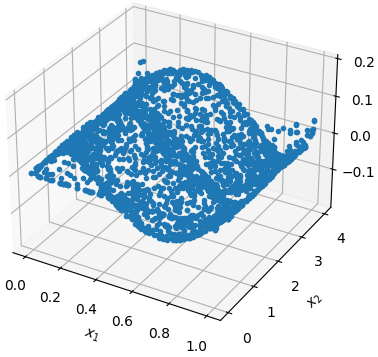}
		\caption{
			% Pinned-pinned bar, trade-off
			{\scriptsize Velocity t.h.}
		}
		\label{fig:23.8.4 R2e.1 velocity88000.png}
	\end{subfigure}
	\begin{subfigure}{0.24\textwidth}
		\includegraphics[width=\textwidth]{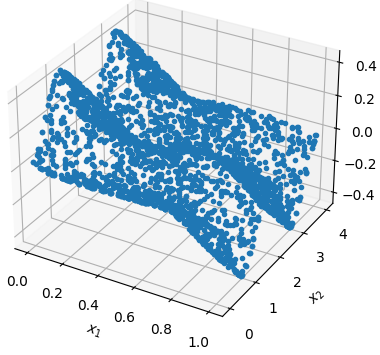}
		\caption{
			% Pinned-pinned bar, trade-off
			{\scriptsize Slope t.h.}
		}
		\label{fig:23.8.4 R2e.1 slope88000.png}
	\end{subfigure}
	\\
	\centering
	\begin{subfigure}{0.24\textwidth}
		\includegraphics[width=\textwidth]{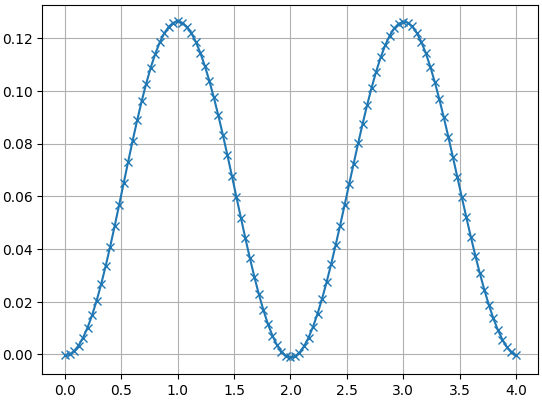}
		\caption{
			% Pinned-pinned bar, trade-off
			{\scriptsize Midspan disp t.h.}
		}
		\label{fig:23.8.4 R2e.1 center disp88000.png}
	\end{subfigure}
	\begin{subfigure}{0.24\textwidth}
		\includegraphics[width=\textwidth]{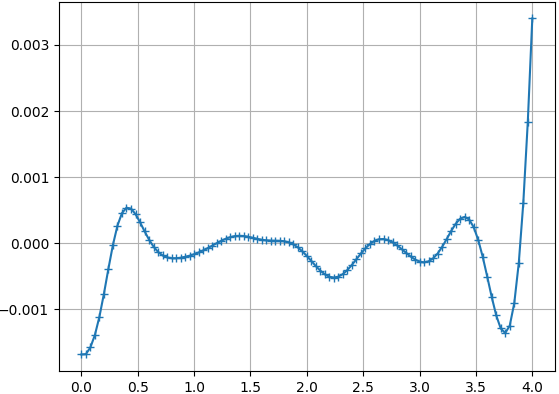}
		\caption{
			% Pinned-pinned bar, trade-off
			{\scriptsize Right-end disp t.h.}
		}
		\label{fig:23.8.4 R2e.1 right-end disp88000.png}
	\end{subfigure}
	\begin{subfigure}{0.24\textwidth}
		\includegraphics[width=\textwidth]{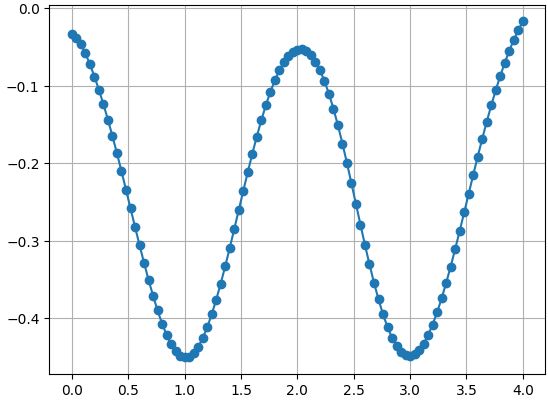}
		\caption{
			% Pinned-pinned bar, trade-off
			{\scriptsize Right-end slope t.h.}
		}
		\label{fig:23.8.4 R2e.1 right-end slope88000.png}
	\end{subfigure}
	\caption{
		\DDET.
		\emph{Pinned-pinned bar}.
		Step \red{88,000}: 
		``t.h.'' = time history.
		SubFig.~\ref{fig:23.8.4 R2e.1 center disp88000.png},
		midspan displacement t.h., damping\%=0.1\%, \red{\emph{Quasi-perfect}}. 
		\red{$\star$}
		Section~\ref{sc:pinned-pinned-Form-2a}, {\color{red} Form~3}.
		\emph{Network:}
		\red{$\star$}
		Remarks~\ref{rm:parameter-names},~\ref{rm:data-point-grids}, 
		T=4,
		W=64,
		H=4, 
		n\_out=3, 
		He-uniform initializer, 
		\red{12,867} parameters,
		\red{$\star$}
		\mbox{\red{\emph{random}}} grid.
		\emph{Training:}
		\red{$\star$}
		Remark~\ref{rm:learning-rate-schedule-3} 
		\red{LRS~3 (NCA)},
		init\_lr=0.01.
		%n-cycles=1, N\_steps=200,000.
		\red{$\star$}
		Total loss {\color{red} \mbox{1.179e-06}}, 
		% Step 88,000 
		(sum of 7 losses).
		Total GPU time \red{278 sec}. 
		$\bullet$
		{\footnotesize
			The above is a good trade-off with the results in
			Figures~\ref{fig:random-grid-NO-static-solution}-\ref{fig:random-grid-NO-static-solution-2}.
			$\triangleright$
			Figure~\ref{fig:axial-Mathematica-solutions}, reference solution to compare.
		}
		{\scriptsize (2384R2e-1)}
	}
	\label{fig:23.8.4 R2e.1 good trade-off}
\end{figure}

\begin{rem}
	\label{rm:static-solution-avoid}
	Static solution, methods to avoid.
	{\rm
		\emph{Pinned-pinned bar.}
		Using \DDET\ Form 2a (time-derivative splitting) and a \emph{regular} grid,
		Figure~\ref{fig:pinned-pinned-static-solution-1} show the loss function and a static-shape time history (or ``static solution'' for short, Remark~\ref{rm:static-solution}), with the velocity and midspan-displa\-cement time histories shown in
		Figure~\ref{fig:pinned-pinned-static-solution-2}, all figures were from Step 50,000 of the same run.
		A sharp jump with a flat plateau along the time axis in the time history, such as in the right subfigure of Figure~\ref{fig:pinned-pinned-static-solution-2}, is a telltale sign of a static solution.
		
		Because of the perfect alignment of the data points along the time axis, \emph{regular} grids would lead to \emph{static} solutions easier than random grids. 
		
		A \emph{pre-static} time history is a solution at an early optimization step (e.g., 25,000) having the characteristic (a jump with a wavy plateau) that portend a static time history at subsequent optimization steps, and thus is a sign to stop the optimization process early (no need to continue further). 
		Two pre-static midspan-displacement time histories at Step 25,000 are shown in
		Figure~\ref{fig:pre-static-solution-jump}, both with init\_lr=0.03. The optimization for the left subfigure continues to Step 200,000 to exhibit a clear static solution shown in Figures~\ref{fig:23.8.15 R3a A-PP F(1)2a lr0.03 cy5-NCA Nsteps200000 regular N51 W64 H2 He init-1}-\ref{fig:23.8.15 R3a A-PP F(1)2a lr0.03 cy5-NCA Nsteps200000 regular N51 W64 H2 He init-2}, whereas the optimization for the right subfigure, where the jump at time $t=0$ is already sure sign of a static solution, stopped at Step 25,000. 
		
		Several methods below (or a combination of these) could be used to avoid the static solution: 
		(1) Use \emph{random} grid (Figure~\ref{fig:random-grid-NO-static-solution}) instead of regular grid (Figures~\ref{fig:pinned-pinned-static-solution-1}-\ref{fig:pinned-pinned-static-solution-2}), 
		(2) Reduce the initial learning rate init\_lr (see below),
		% {\color{red} Figures as examples.}
		(3) Reduce network capacity (number of network parameters) (see below),
		% {\color{red} Figures as examples.}
		(4) Use barrier function (see Appendix~\ref{app:barrier-functions}).
		% {\color{red} Figures as examples.}
		
		Using the same Form 2a and network with 12,802 parameters as in Figures~\ref{fig:pinned-pinned-static-solution-1}, and reduce init\_lr from 0.01 to various values below.  
		Reducing the initial learning rate by 100 times to init\_lr=0.0001 with NCA produced un-accentuated waves with high damping at Step 200,000.  On the other hand, reducing by 10 times to init\_lr=0.001 with CA, and extending Step 5 to 250000 with NCA (Remark~\ref{rm:extension-cycle-5}), resulted in a % \emph{quasi-perfect} 
		good
		solution with small damping, Figure~\ref{fig:23.8.22 R3b A-PP F(1J)2a(3) lr0.001 cy5-CA cy6-NCA Nsteps250000 regular N51 W64 H4 He init}. 
		See also Remark~\ref{rm:unstable-solution}, \emph{Unstable solution}.
		
		Reducing the network capacity would allow using a much larger learning rate for similar results.
		Figure~\ref{fig:23.8.17 R10a A-PP F(1)2a lr0.07 cy2-CA Nsteps50000 regular N51 W32 H2 He init} shows a static solution using a network with only 1,218 parameters (ten times smaller) with init\_lr=0.07 (seven times larger) compared to  Figures~\ref{fig:pinned-pinned-static-solution-1}-\ref{fig:pinned-pinned-static-solution-2} using a network with 12,802 parameters and init\_lr=0.01.
		% ALEX: removed paragraph
		Figure~\ref{fig:23.8.4 R2e.1 good trade-off} shows a solution obtained with \DDET\ Form 3 using a network with 12,867 parameters, init\_lr=0.01, NCA, resulting in the Cycle 3 lowest loss of 1.179e-06 at Step 88,000, and a \red{quasi-perfect} midspan displacement of a pinned-pinned bar.  
		
		\emph{Pinned-free bar.}  The static solutions in \DDET\ Form 1 in Figure~\ref{fig:DDE-T-Form-1-pinned-free-bar-static-solutions} are totally avoided by using \DDET\ Form 2a (Figure~\ref{fig:23.9.5 R3b A-P bcrF F(1J)2a(3) lr0.005 cy1-5-NCA Nsteps200000 random N51 W64 H4 He init-1}), and of course \DDET\ Form 3 ({\color{OliveGreen} \mbox{\bf SubFigs}}~\ref{fig:23.9.5 R3d.4 A-PF loss200000}-\ref{fig:23.9.5 R3d.4 A-PF slope200000}).
		% (\red{NEED FIGURE of DDE-T FORM 3}).  
		% \red{TO COMPLETE 23.10.6}
	}
	\phantom{blank}\hfill$\blacksquare$
\end{rem}

%--------------------------------------------------------------Form 2a

% 23.9.20, OLD Figure 8, moved to section on Form 2a
% not used to provide shape with regular grid
% 23.8.22 R5a A-PP F(1J)2a(3) lr0.0001 cy5-NCA Nsteps200000 regular N51 W64 H4 He init - center disp200000.png

% 23.9.20, OLD Figure 9, moved to section on Form 2a
%\footnotetext{
	%	Two reruns of the same exact script that produced Figure~\ref{fig:random-grid-NO-static-solution} yielded the total loss of 
	%	0.626e-06 (RunID 2384R2b) and 0.518e-06 (RunID 2384R2c), respectively.
	%}

%\noindent
%\red{[NOTE: Compare Figure~\ref{fig:random-grid-NO-static-solution} to the results in 23.9.16 R1c.3 lr=0.002 400000 steps quasi-perfect.  Need to quantify what it means by quasi-perfect.  Use damping / amplification percentage.  Less than 1\% will not be visible. ENDNOTE]}

\begin{figure}[tph]
	\includegraphics[width=0.49\textwidth]{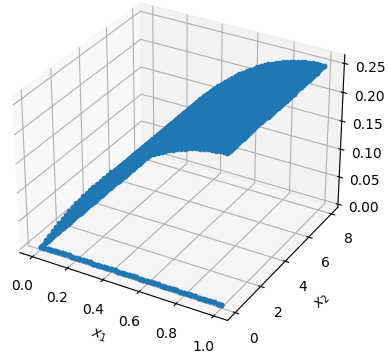}
	\includegraphics[width=0.49\textwidth]{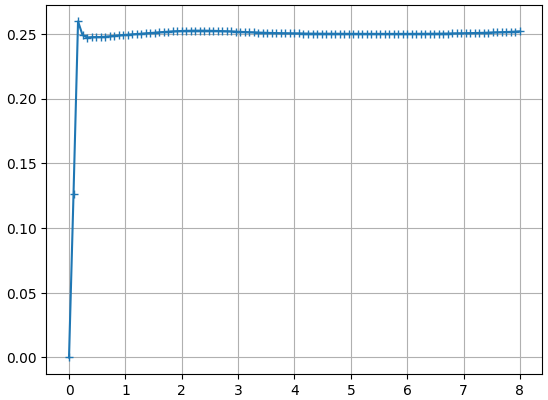}
	\caption{
		% \red{I AM HERE, 23.9.8. Rewrite below.}
		\DDET.
		\emph{Pinned-free bar. {He-uniform} initializer}. 
		% \emph{Pinned-free bar}.
		\red{$\star$}
		\emph{Static} shape (left), free-end displacement (right), Step 25,000, 
		end of Cycle~1.
		Section~\ref{sc:axial-Form-2a}, \red{Form 2a}.
		\emph{Network:}
		Remarks~\ref{rm:parameter-names},~\ref{rm:data-point-grids}, 
		T=8,
		W=64,
		H=4, 
		n\_out=2, 
		\red{\emph{He-uniform} initializer}, 
		12,802 parameters,
		\red{\emph{regular}} grid.
		\emph{Training:}
		Remark~\ref{rm:learning-rate-schedule-3},
		LRS~3, 
		init\_lr=0.005.
		$\bullet$
		{\footnotesize
			% \red{Change figures here.}
			% Figure~\ref{fig:DEBUG v1.9.3 seed42, 23.9.5 R3a A-P bcrF F(1J)2a(3) lr0.005 cy1-3-NCA Nsteps100000 regular N51 W64 H4 He init-1}, same model, LRS~3, init\_lr=0.005, \red{regular} grid, \emph{static} solution.
			Figure~\ref{fig:23.9.5 R3b A-P bcrF F(1J)2a(3) lr0.005 cy1-5-NCA Nsteps200000 random N51 W64 H4 He init-1}, same model, LRS~3, init\_lr=0.005, \red{random} grid, \red{\emph{quasi-perfect}} solution. 
		}
		{\scriptsize (2395R3a-1)}
	}
	\label{fig:DEBUG v1.9.3 seed42, 23.9.5 R3a A-P bcrF F(1J)2a(3) lr0.005 cy1-3-NCA Nsteps100000 regular N51 W64 H4 He init-1}
\end{figure}

\subsubsection{Form 2b: Same issues as Form 1}
\label{sc:axial-Form-2b}
Form 2b for the axial equation of motion of an elastic bar is a particular case of Form 2b for the Kirchhoff-Love rod in Section~\ref{sc:Kirchhoff-rod-Form-2b}:
\begin{align}
	%	\sldn 
	%	\left[ \slpbsp{\du}{1} + \eslpbp{1} \sin \eslpb \right]
	%	\bullet
	\sldn 
	%	\ubp 2
	\slpbsp{\du}{1} 
	+ \dfbs{\X} 
	= \nddt \ub
	\ , \quad
	\ubp{1} = \slpbs{\du} 
	% \ , \quad
	% \ndt \ub = \pbs{\X}
	\ ,
	\label{eq:wave-eq-form-2b}
\end{align}
with the initial conditions
% \red{I AM HERE 23.10.6.}
\begin{align}
	\ub (\X , 0) = \ubs{0}  
	\ , \quad
	% \pbs{\X} (\X , 0) = \ubds{0}
	\ndt{\ub} (\X , 0) = \ubds{0}
	\ , \quad
	\slpbs{\ub} (\X , 0) = \slpbs{\ub}_{0}
	\ .
	\label{eq:wave-eq-form-2b-initial-conditions-1}
\end{align}
In the numerical examples, we use homogeneous initial conditions:
\begin{align}
	\ub (\X , 0) = 0 
	\ , \quad
	% \pbs{\X} (\X , 0) = 0
	\ndt{\ub} (\X , 0) = 0
	\ , \quad
	\slpbs{\ub} (\X , 0) = 0 
	\ .
	\label{eq:wave-eq-form-2b-initial-conditions-2}
\end{align}

Form 2b (splitting the space derivative operator) does not resolve the shift and amplification like Form 2a (splitting the time derivative operator).  For example, using the same parameters as in Figure~\ref{fig:23.7.22 R1 A-PP v1 midspan,shape100000} for the pinned-pinned bar, shift and amplification, albeit with a different magnitude, still persisted.

Form 2b also yielded the same static solutions as those in Figure~\ref{fig:DDE-T-Form-1-pinned-free-bar-static-solutions}.

% 2023.06.28
% changed Form i to Form i+1 to be consistent with the code
\subsubsection{Form 3 (or 4): No space/time shift, static solution, efficiency}
\label{sc:axial-Form-3}
\label{sc:pinned-pinned-Form-3}
\noindent
Form 3 (or 4) for the axial equation of motion of an elastic bar is a particular case of Form 3 (or 4) for the Kirchhoff-Love rod in Section~\ref{sc:Kirchhoff-rod-Form-3} (or \ref{sc:Kirchhoff-rod-Form-4}):
\begin{align}
%	\sldn 
%	\left[ \slpbsp{\du}{1} + \eslpbp{1} \sin \eslpb \right]
%	\bullet
	\sldn 
%	\ubp 2
	\slpbsp{\du}{1} 
	+ \dfbs{\X} 
	= \pbsd{\X}
	\ , \quad
	\ubp{1} = \slpbs{\du} 
	\ , \quad
	\ndt \ub = \pbs{\X}
	\ ,
	\label{eq:wave-eq-form-3}
\end{align}
with the initial conditions
\begin{align}
	\ub (\X , 0) = \ubs{0} 
	\ , \quad
	\slpbs{\ub} (\X , 0) = \slpbs{\ub}_{0} 
	\ , \quad
	\pbs{\X} (\X , 0) = \ubds{0}
	\ .
\end{align}
In the numerical examples, we use homogeneous initial conditions:
\begin{align}
	\ub (\X , 0) = 0 
	\ , \quad
	\slpbs{\ub} (\X , 0) = 0 
	\ , \quad
	\pbs{\X} (\X , 0) = 0
	\ .
\end{align}

%--------------------------------------------------------------Form 3

% 23.9.20, OLD Fig.14 to move to section Form 3
\begin{figure}[tph]
	%_\includegraphics[width=0.49\textwidth]{Figures/23.8.21_R1c_A-PP_F_1J,2a_3_lr0.04,0.03,0.02,0.01,0.005_cy5-VCA_Nsteps200000_regular_N51_W32_H2_He_init_-_loss200000.png}
%	\includegraphics[width=0.49\textwidth]{Figures/23.8.21_R1c_A-PP_F_1J,2a_3_lr0.04,0.03,0.02,0.01,0.005_cy5-VCA_Nsteps200000_regular_N51_W32_H2_He_init_-_shape200000.png}
%	\includegraphics[width=0.49\textwidth]{Figures/23.8.21_R1c_A-PP_F_1J,2a_3_lr0.04,0.03,0.02,0.01,0.005_cy5-VCA_Nsteps200000_regular_N51_W32_H2_He_init_-_center_disp200000.png}
	%
	\includegraphics[width=0.49\textwidth]{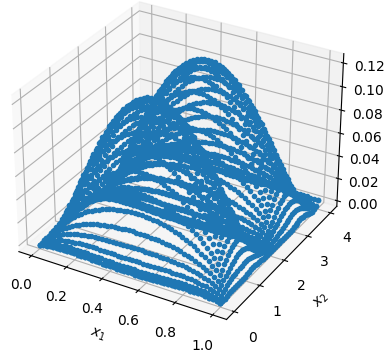}
	\includegraphics[width=0.49\textwidth]{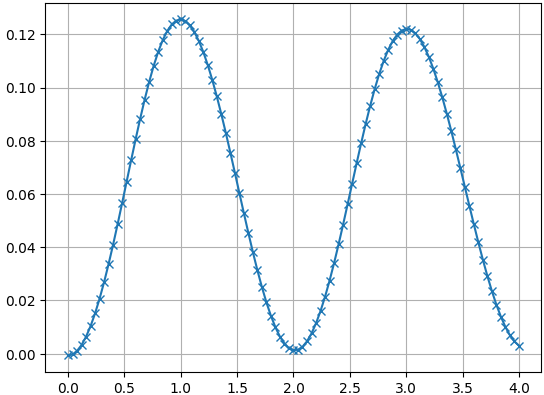}
	\caption{
		% CA led to waves with damping.
		% In 23.8.13 R4ab, 1218 params, NCA led to waves with damping, whereas CA led to the static solution.  In 23.8.15 R3ab, 4482 params, it's the opposite: NCA led to the static solution, whereas CA led to waves with damping ! 
		% 23.8.21 R1c: Form F3
		%
		%		\emph{Pinned-pinned bar, Varying init\_lr cyclic annealing} (VCA). 
		% \red{I AM HERE, 23.9.8. Rewrite below.}
		\DDET.
		\emph{Pinned-pinned bar, VCA}.
		% Section~\ref{sc:optimization-learning-rate-scheduling}.
		% Loss function (left), 
		Shape (left), 
		midspan displacement, Step 200,000 (right), \emph{waves, small damping}.
		\red{$\star$}
		Section~\ref{sc:pinned-pinned-Form-3}, {\color{red} Form~3}.
		% \red{$\star$}
		Remark~\ref{rm:damping-low-capacity-models}, \emph{Damping}.
		% {\color{red} 2023.08.21, Write the section for Form 3.}
		\emph{Network:}
		Remarks~\ref{rm:parameter-names},~\ref{rm:data-point-grids}, 
		T=4,
		W=32,
		H=2, 
		n\_out=2, 
		He-uniform initializer, \red{1,251} parameters,
		\red{$\star$}
		\emph{regular} grid.
		\emph{Training:} 
		\red{$\star$}
		Remark~\ref{rm:learning-rate-schedule-2},
		\red{LRS~2 (VCA)},
		init\_lr=[0.04, 0.03, 0.02, 0.01, 0.005].
		\red{$\star$}
		Lowest total loss \red{4.31e-06}, Step 200,000 (sum of 7 losses).
		\red{$\star$}
		Total GPU time \red{309 sec}.
		$\bullet$
		{\footnotesize
			Figure~\ref{fig:23.8.21 R1c A-PP F(1J,2a)3 lr0.04,0.03,0.02,0.01,0.005 cy5-VCA Nsteps200000 regular N51 W32 H2 He init-2}, velocity, slope.
			$\triangleright$
			Figure~\ref{fig:23.8.13 R4b A-PP F(1)2a lr0.04 cy5-CA Nsteps200000 regular N51 W32 H2 He init}, CA, Form 2a, same model and init\_lr=0.04, static solution.
			Figure~\ref{fig:23.8.21 R1 A-PP F(1)2a lr0.04,0.03,0.02,0.01,0.005 cy5-VCA Nsteps200000 regular N51 W32 H2 He init}, VCA, Form 2a, same model and init\_lr sequence, waves, small damping.
		}
		{\scriptsize (23821R1c-1)}
	}
	\label{fig:23.8.21 R1c A-PP F(1J,2a)3 lr0.04,0.03,0.02,0.01,0.005 cy5-VCA Nsteps200000 regular N51 W32 H2 He init-1}
\end{figure}

% 23.9.20, OLD Fig.15 to move to section Form 3
\begin{figure}[tph]
	\includegraphics[width=0.49\textwidth]{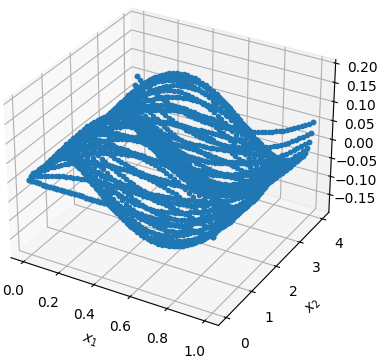}
	\includegraphics[width=0.49\textwidth]{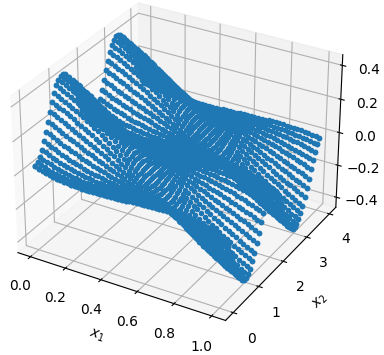}
	\caption{
		% CA led to waves with damping.
		% In 23.8.13 R4ab, 1218 params, NCA led to waves with damping, whereas CA led to the static solution.  In 23.8.15 R3ab, 4482 params, it's the opposite: NCA led to the static solution, whereas CA led to waves with damping ! 
		% 23.8.21 R1c: Form F3
		%
		% \red{I AM HERE, 23.9.8. Rewrite below.}
		\DDET.
		\emph{Pinned-pinned bar, VCA}. 
		%		Remark~\ref{rm:learning-rate-schedule-2},
		%		Section~\ref{sc:optimization-learning-rate-scheduling}.
		% Loss function (left), 
		Velocity (left), slope (right), Step 200,000, \emph{waves, small damping}.
		\red{$\star$}
		Section~\ref{sc:pinned-pinned-Form-3}, \red{Form~3}.
		% {\color{red} 2023.08.21, Write the section for Form 3.}
		\emph{Network:}
		\red{$\star$}
		Remarks~\ref{rm:parameter-names},~\ref{rm:data-point-grids}, 
		T=4,
		W=32,
		H=2, 
		n\_out=2, 
		He-uniform initializer, 
		\red{1,251} parameters,
		\red{$\star$}
		\emph{regular} grid.
		\emph{Training:}
		\red{$\star$}
		Remark~\ref{rm:learning-rate-schedule-2}, 
		\red{LRS~2 (VCA)},
		init\_lr=[0.04, 0.03, 0.02, 0.01, 0.005].
		$\bullet$
		{\footnotesize
			Figure~\ref{fig:23.8.21 R1c A-PP F(1J,2a)3 lr0.04,0.03,0.02,0.01,0.005 cy5-VCA Nsteps200000 regular N51 W32 H2 He init-1}, VCA, Form 3, same model and init\_lr sequence, waves with small damping.
			$\triangleright$
			Figure~\ref{fig:23.8.13 R4b A-PP F(1)2a lr0.04 cy5-CA Nsteps200000 regular N51 W32 H2 He init}, CA, Form 2a, same model, init\_lr=0.04, static solution.
			Figure~\ref{fig:23.8.21 R1 A-PP F(1)2a lr0.04,0.03,0.02,0.01,0.005 cy5-VCA Nsteps200000 regular N51 W32 H2 He init}, VCA, Form 2a, waves with small damping. 
		} 
		{\scriptsize (23821R1c-2)}
	}
	\label{fig:23.8.21 R1c A-PP F(1J,2a)3 lr0.04,0.03,0.02,0.01,0.005 cy5-VCA Nsteps200000 regular N51 W32 H2 He init-2}
\end{figure}

Figures~\ref{fig:23.8.21 R1c A-PP F(1J,2a)3 lr0.04,0.03,0.02,0.01,0.005 cy5-VCA Nsteps200000 regular N51 W32 H2 He init-1}-\ref{fig:23.8.21 R1c A-PP F(1J,2a)3 lr0.04,0.03,0.02,0.01,0.005 cy5-VCA Nsteps200000 regular N51 W32 H2 He init-2} for the axial motion of a \emph{pinned-pinned bar} were obtained using Form 3, VCA (Remark~\ref{rm:learning-rate-schedule-2}), network with 1,251 parameters, regular grid, for 200,000 steps, and resulted in the lowest total loss of 4.31e-06 and GPU time of 309 sec.
%{\color{red} NOTE: 2023.08.27, Figures~\ref{fig:23.8.21 R1c A-PP F(1J,2a)3 lr0.04,0.03,0.02,0.01,0.005 cy5-VCA Nsteps200000 regular N51 W32 H2 He init-1}-\ref{fig:23.8.21 R1c A-PP F(1J,2a)3 lr0.04,0.03,0.02,0.01,0.005 cy5-VCA Nsteps200000 regular N51 W32 H2 He init-2} may be removed since VCA and ELRS are used in   Figures~\ref{fig:23.8.24 R1 A-PP F(1J,2a)3 lr0.04,0.03,0.02,0.01,0.01,0.005 cy1-6-VCA cy7-9-NCA Nsteps400000 random N51 W32 H2 He init-1},
%\ref{fig:23.8.24 R1 A-PP F(1J,2a)3 lr0.04,0.03,0.02,0.01,0.01,0.005 cy1-6-VCA cy7-9-NCA Nsteps400000 random N51 W32 H2 He init-2},
%\ref{fig:23.8.24 R1 A-PP F(1J,2a)3 lr0.04,0.03,0.02,0.01,0.01,0.005 cy1-6-VCA cy7-9-NCA Nsteps400000 random N51 W32 H2 He init-3}. ENDNOTE}

For the same pinned-pinned bar, Form 3, network with 1,251 parameters, but VCA with ELRS (Remark~\ref{rm:extension-cycle-5}) and random grid for 400,000 steps,
Figures~\ref{fig:23.8.24 R1 A-PP F(1J,2a)3 lr0.04,0.03,0.02,0.01,0.01,0.005 cy1-6-VCA cy7-9-NCA Nsteps400000 random N51 W32 H2 He init-1},
\ref{fig:23.8.24 R1 A-PP F(1J,2a)3 lr0.04,0.03,0.02,0.01,0.01,0.005 cy1-6-VCA cy7-9-NCA Nsteps400000 random N51 W32 H2 He init-2},
\ref{fig:23.8.24 R1 A-PP F(1J,2a)3 lr0.04,0.03,0.02,0.01,0.01,0.005 cy1-6-VCA cy7-9-NCA Nsteps400000 random N51 W32 H2 He init-3},
yielded a solution with a total loss of 1.25e-06 (3 times smaller compared to 4.31e-06 in Figures~\ref{fig:23.8.21 R1c A-PP F(1J,2a)3 lr0.04,0.03,0.02,0.01,0.005 cy5-VCA Nsteps200000 regular N51 W32 H2 He init-1}-\ref{fig:23.8.21 R1c A-PP F(1J,2a)3 lr0.04,0.03,0.02,0.01,0.005 cy5-VCA Nsteps200000 regular N51 W32 H2 He init-2}) and total GPU time of 598 sec (less than 2 times longer than 309 sec in Figures~\ref{fig:23.8.21 R1c A-PP F(1J,2a)3 lr0.04,0.03,0.02,0.01,0.005 cy5-VCA Nsteps200000 regular N51 W32 H2 He init-1}-\ref{fig:23.8.21 R1c A-PP F(1J,2a)3 lr0.04,0.03,0.02,0.01,0.005 cy5-VCA Nsteps200000 regular N51 W32 H2 He init-2}).
See Remark~\ref{rm:data-point-grids} regarding the lowest total loss achieved, \red{0.518e-06}, for axial motion of a pinned-pinned bar.

While \emph{symmetry} can be observed in the \emph{quasi-perfect} midspan-displacement time history in Figure~\ref{fig:23.8.24 R1 A-PP F(1J,2a)3 lr0.04,0.03,0.02,0.01,0.01,0.005 cy1-6-VCA cy7-9-NCA Nsteps400000 random N51 W32 H2 He init-1}, a lack of symmetry was revealed in the computed slope (space derivative) time history of the right pinned-end in Figure~\ref{fig:23.8.24 R1 A-PP F(1J,2a)3 lr0.04,0.03,0.02,0.01,0.01,0.005 cy1-6-VCA cy7-9-NCA Nsteps400000 random N51 W32 H2 He init-3}, which showed a shift downward by 0.03 (since the slope must be zero at time $t = 0$) and a clear departure from symmetry near time $t = 4$.

Compared to Form 1 used in Figure~\ref{fig:23.7.23 R1d A-PP v1 midspan,shape200000}, with a network of 12,737 parameters, CA (Remark~\ref{rm:learning-rate-schedule-1}) over 200,000 steps,
with a total GPU time of 621 sec, Form 3 used in Figure~\ref{fig:23.8.24 R1 A-PP F(1J,2a)3 lr0.04,0.03,0.02,0.01,0.01,0.005 cy1-6-VCA cy7-9-NCA Nsteps400000 random N51 W32 H2 He init-1} with 1,251 parameters and a total GPU time of 598 sec is still 4\% more efficient, despite running over twice the number of steps to 400,000.

\begin{figure}[tph]
	\includegraphics[width=0.49\textwidth]{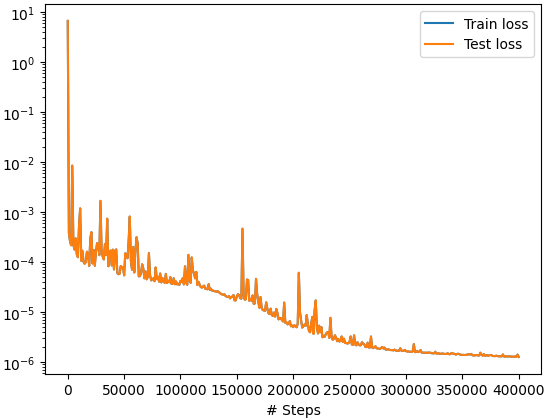}
	\includegraphics[width=0.49\textwidth]{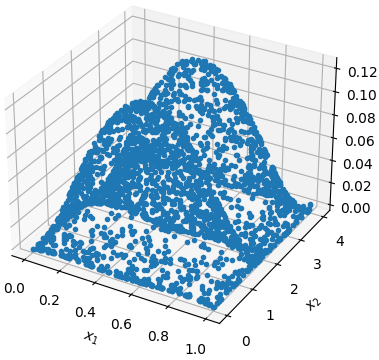}
	\caption{
		% \red{I AM HERE, 23.9.9. Rewrite below.}
		\DDET.
		\emph{Pinned-pinned bar}. 
		% Remark~\ref{rm:learning-rate-schedule-2}, \emph{Varying init\_lr cyclic annealing} (VCA). 
		% Remark~\ref{rm:extension-cycle-5}, \emph{Extended learning-rate schedule} (ELRS).
		Loss function (left), shape, Step 400,000 (right).
		\red{$\star$}
		Remark~\ref{rm:damping-low-capacity-models}, Damping.
		Section~\ref{sc:axial-Form-3}, {\color{red} Form~3}.
		\emph{Network:}
		Remarks~\ref{rm:parameter-names},~\ref{rm:data-point-grids}, 
		T=4,
		W=32,
		H=2, 
		n\_out=3, 
		He-uniform initializer, 
		\red{1,251} parameters,
		\emph{random} grid.
		\emph{Training:}
		\red{$\star$}
		Remarks~\ref{rm:learning-rate-schedule-2} (LRS~2), 
		\ref{rm:extension-cycle-5} (ELRS),
		init\_lr=[0.04, 0.03, 0.02, 0.01, 0.01, 0.005] Cycles 1-6 (VCA);
		Cycles 7-9 (NCA).
		\red{$\star$}
		Lowest total loss {\color{red} 1.25e-06}, Step 400,000 (sum  of 7 losses).
		\red{$\star$}
		Total GPU time {\color{red} 598 sec}.
		% {\color{red} REWRITE below.}
		$\bullet$
		{\footnotesize 
			Figure~\ref{fig:23.8.24 R1 A-PP F(1J,2a)3 lr0.04,0.03,0.02,0.01,0.01,0.005 cy1-6-VCA cy7-9-NCA Nsteps400000 random N51 W32 H2 He init-2}, velocity, midspan displacement, 
			% time histories at 
			Step 400,000.
			Figure~\ref{fig:23.8.24 R1 A-PP F(1J,2a)3 lr0.04,0.03,0.02,0.01,0.01,0.005 cy1-6-VCA cy7-9-NCA Nsteps400000 random N51 W32 H2 He init-3}, slope, 
			% time histories at 
			Step 400,000.
			$\triangleright$
			% Figure~\ref{fig:23.7.23 R1d A-PP v1 midspan,shape200000}, 
			Figure~\ref{fig:random-grid-NO-static-solution},
			total loss \red{0.537e-06}.
%			Figure~\ref{fig:23.7.23 R1d A-PP v1 midspan,shape200000}, Form 1, 12,737 parameters, CA over 200,000 steps, total GPU time of 621 sec.
			$\triangleright$
			Figure~\ref{fig:axial-Mathematica-solutions}, reference solution to compare.
		}
		{\scriptsize (23824R1-1)}
	}
	\label{fig:23.8.24 R1 A-PP F(1J,2a)3 lr0.04,0.03,0.02,0.01,0.01,0.005 cy1-6-VCA cy7-9-NCA Nsteps400000 random N51 W32 H2 He init-1}
\end{figure}

\begin{figure}[tph]
	\includegraphics[width=0.49\textwidth]{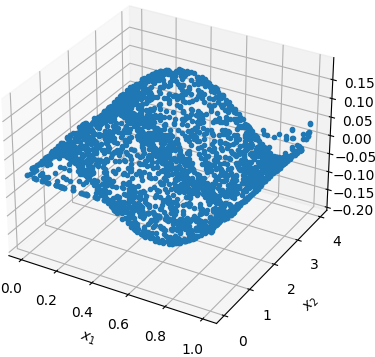}
	\includegraphics[width=0.49\textwidth]{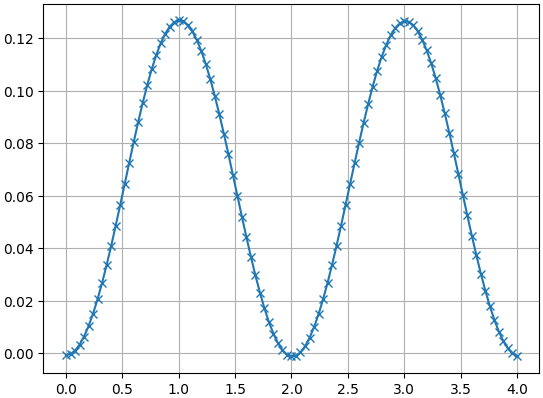}
	\caption{
		% \red{I AM HERE, 23.9.9. Rewrite below.}
		\DDET.
		\emph{Pinned-pinned bar}. 
%		Remark~\ref{rm:learning-rate-schedule-2}, \emph{Varying init\_lr cyclic annealing} (VCA). 
%		Remark~\ref{rm:extension-cycle-5}, \emph{Extended learning-rate schedule} (ELRS).
		\red{$\star$}
		Velocity (left), visually \emph{quasi-perfect} midspan displacement, Step 400,000 (right).
		\red{$\star$}
		Section~\ref{sc:axial-Form-3}, \red{Form~3}.
		\emph{Network:}
		Remarks~\ref{rm:parameter-names},~\ref{rm:data-point-grids}, 
		T=4,
		W=32,
		H=2, 
		n\_out=3, 
		He-uniform initializer, 
		1,251 parameters,
		\emph{random} grid.
		\emph{Training:}
		Remark~\ref{rm:learning-rate-schedule-2} (LRS~2), 
		Remark~\ref{rm:extension-cycle-5} (ELRS),
		init\_lr=[0.04, 0.03, 0.02, 0.01, 0.01, 0.005] Cycles 1-6 (VCA);
		Cycles 7-9 (NCA).
		% Total GPU time 598 sec.
		% {\color{red} REWRITE below.}
		$\bullet$
		{\footnotesize 
			Figure~\ref{fig:23.8.24 R1 A-PP F(1J,2a)3 lr0.04,0.03,0.02,0.01,0.01,0.005 cy1-6-VCA cy7-9-NCA Nsteps400000 random N51 W32 H2 He init-1}, loss function, shape, Step 400,000, lowest total loss, GPU time.
			Figure~\ref{fig:23.8.24 R1 A-PP F(1J,2a)3 lr0.04,0.03,0.02,0.01,0.01,0.005 cy1-6-VCA cy7-9-NCA Nsteps400000 random N51 W32 H2 He init-3}, slope, Step 400,000.
		}
		{\scriptsize (23824R1-2)}
	}
	\label{fig:23.8.24 R1 A-PP F(1J,2a)3 lr0.04,0.03,0.02,0.01,0.01,0.005 cy1-6-VCA cy7-9-NCA Nsteps400000 random N51 W32 H2 He init-2}
\end{figure}

\begin{figure}[tph]
	\includegraphics[width=0.49\textwidth]{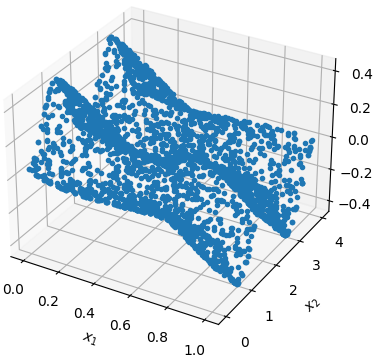}
	\includegraphics[width=0.49\textwidth]{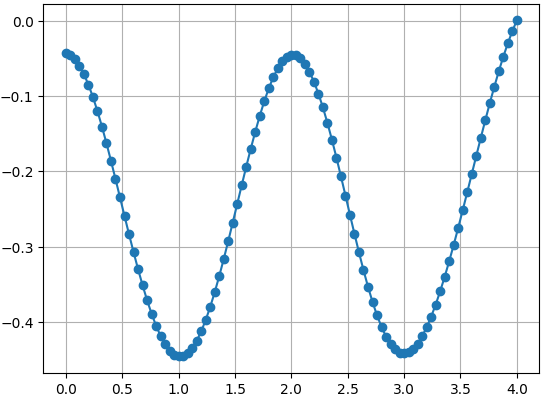}
	\caption{
		% \red{I AM HERE, 23.9.9. Rewrite below.}
		\DDET.
		\emph{Pinned-pinned bar}. 
%		Remark~\ref{rm:learning-rate-schedule-2}, \emph{Varying init\_lr cyclic annealing} (VCA). 
%		Remark~\ref{rm:extension-cycle-5}, \emph{Extended learning-rate schedule} (ELRS).
		\red{$\star$}
		Slope (left), right-pinned-end slope, Step 400,000 (right).
		\red{$\star$}
		Section~\ref{sc:axial-Form-3}, \red{Form~3}.
		\emph{Network:}
		Remark~\ref{rm:parameter-names},~\ref{rm:data-point-grids}, 
		T=4,
		W=32,
		H=2, 
		n\_out=3, 
		He-uniform initializer, 
		1,251 parameters,
		\emph{random} grid.
		\emph{Training:}
		\red{$\star$}
		Remark~\ref{rm:learning-rate-schedule-2} (LRS~2), 
		\ref{rm:extension-cycle-5} (ELRS),
		init\_lr=[0.04, 0.03, 0.02, 0.01, 0.01, 0.005] Cycles 1-6 (VCA);
		Cycles 7-9 (NCA).
		% Total GPU time 598 sec.
		% {\color{red} REWRITE below.}
		$\bullet$
		{\footnotesize
			Figure~\ref{fig:23.8.24 R1 A-PP F(1J,2a)3 lr0.04,0.03,0.02,0.01,0.01,0.005 cy1-6-VCA cy7-9-NCA Nsteps400000 random N51 W32 H2 He init-1}, loss function, shape, Step 400,000, lowest total loss, GPU time.
			Figure~\ref{fig:23.8.24 R1 A-PP F(1J,2a)3 lr0.04,0.03,0.02,0.01,0.01,0.005 cy1-6-VCA cy7-9-NCA Nsteps400000 random N51 W32 H2 He init-2}, velocity, midspan displacement, Step 400,000.
		}
		{\scriptsize (23824R1-3)}
	}
	\label{fig:23.8.24 R1 A-PP F(1J,2a)3 lr0.04,0.03,0.02,0.01,0.01,0.005 cy1-6-VCA cy7-9-NCA Nsteps400000 random N51 W32 H2 He init-3}
\end{figure}

\begin{figure}[tph]
	\centering
	\begin{subfigure}{0.24\textwidth}
		\includegraphics[width=1.\textwidth]{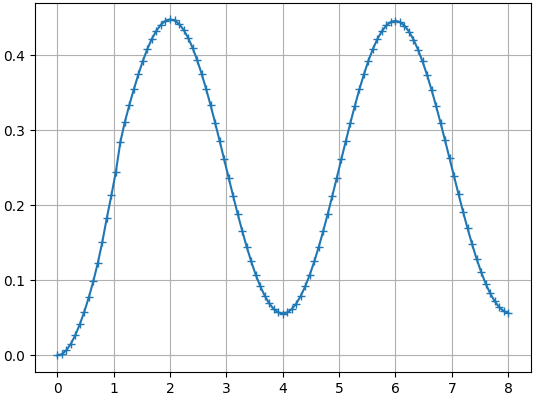}
		\caption{
			\scriptsize
			% Free-end disp, jump 197,000.
			Init\_lr=0.005, He, pre-static.
		}
		\label{fig:23.9.5 R3d A-PF free-end disp197000}
	\end{subfigure}
	\begin{subfigure}{0.24\textwidth}
		\includegraphics[width=1.\textwidth]{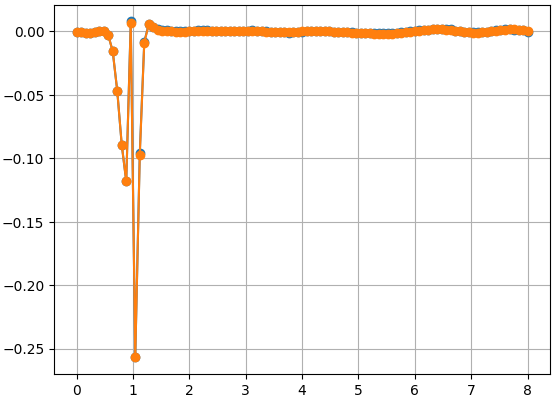}
		\caption{
			\scriptsize
			% Free-end slope, jump.
			Free-end slope, 197000.
		}
		\label{fig:23.9.5 R3d A-PF free-end slope197000}
	\end{subfigure}
	\begin{subfigure}{0.24\textwidth}
		\includegraphics[width=1.\textwidth]{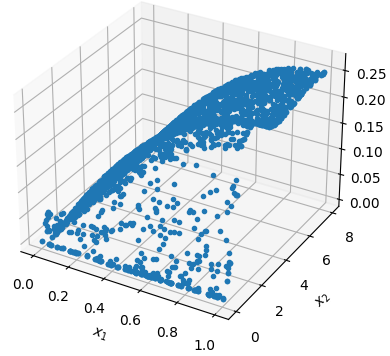}
		\caption{
			\scriptsize
			% Shape, static 50,000.
			Init\_lr=0.01, He, static.
		}
		\label{fig:23.9.5 R3d.2 A-PF shape50000}
	\end{subfigure}
	\begin{subfigure}{0.24\textwidth}
		\includegraphics[width=1.\textwidth]{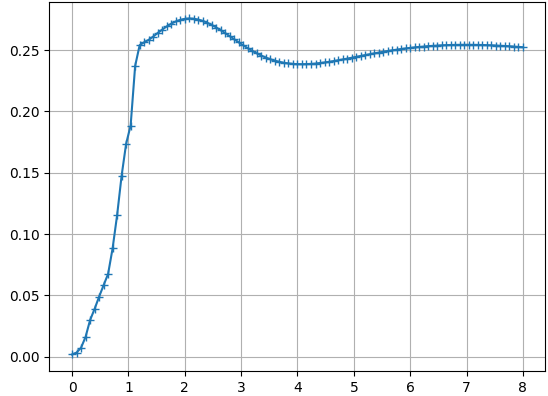}
		\caption{
			\scriptsize
			% Free-end disp, static.
			Free-end disp, 50,000.
		}
		\label{fig:23.9.5 R3d.2 A-PF free-end disp50000}
	\end{subfigure}
	\\
	\begin{subfigure}{0.24\textwidth}
		% changed shape to free-end disp
		%_\includegraphics[width=1.\textwidth]{Figures/DEBUG_23.9.5_R3d.5_A-P_bcrF_F_1J,2a_3_lr0.001_cy1-5-NCA_Nsteps200000_random_N51_W64_H4_He_init_-_shape200000.png}
		\includegraphics[width=1.\textwidth]{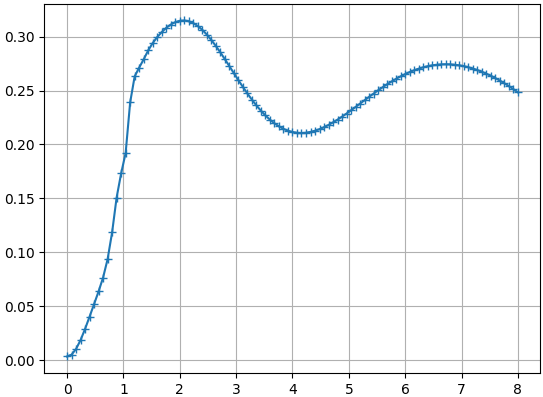}
		\caption{
			\scriptsize
			% Shape, pre-static 200,000.
			% Free-end disp, pre-static 200,000.
			Init\_lr=0.001, He, pre-static.
		}
		% \label{fig:23.9.5 R3d.5 A-PF shape200000}
		\label{fig:23.9.5 R3d.5 A-PF free-end disp200000}
	\end{subfigure}
	\begin{subfigure}{0.24\textwidth}
		% changed free-end disp to free-end slope
		%_\includegraphics[width=1.\textwidth]{Figures/DEBUG_23.9.5_R3d.5_A-P_bcrF_F_1J,2a_3_lr0.001_cy1-5-NCA_Nsteps200000_random_N51_W64_H4_He_init_-_free-end_disp200000.png}
		\includegraphics[width=1.\textwidth]{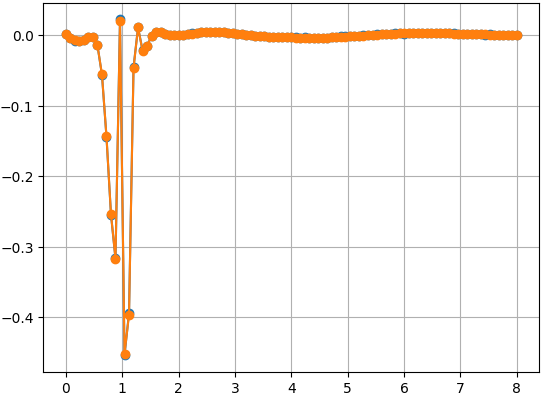}
		\caption{
			\scriptsize
			% Free-end disp, static.
			% Free-end slope, jump.
			Free-end slope, 200,000.
		}
		% \label{fig:23.9.5 R3d.5 A-PF free-end disp200000}
		\label{fig:23.9.5 R3d.5 A-PF free-end slope200000}
	\end{subfigure}
	\begin{subfigure}{0.24\textwidth}
		\includegraphics[width=1.\textwidth]{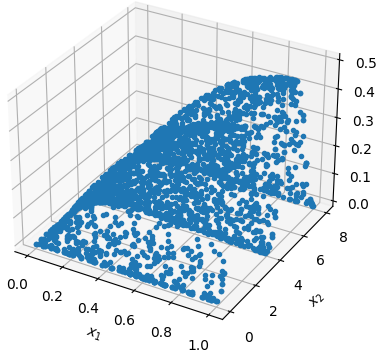}
		\caption{
			\scriptsize
			% Shape 184,000.
			Init\_lr=0.005, \red{Glorot}.
		}
		\label{fig:23.9.5 R3d.3 A-PF shape184000}
	\end{subfigure}
	\begin{subfigure}{0.24\textwidth}
		\includegraphics[width=1.\textwidth]{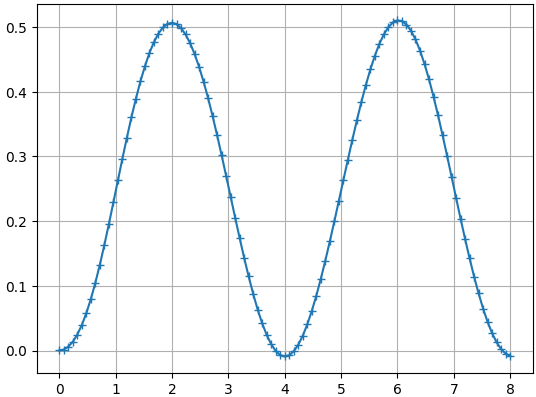}
		\caption{
			\scriptsize
			Free-end disp, \red{Very good}.
		}
		\label{fig:23.9.5 R3d.3 A-PF free-end disp184000}
	\end{subfigure}
	\\
	\centering
	\begin{subfigure}{0.24\textwidth}
		\includegraphics[width=1.\textwidth]{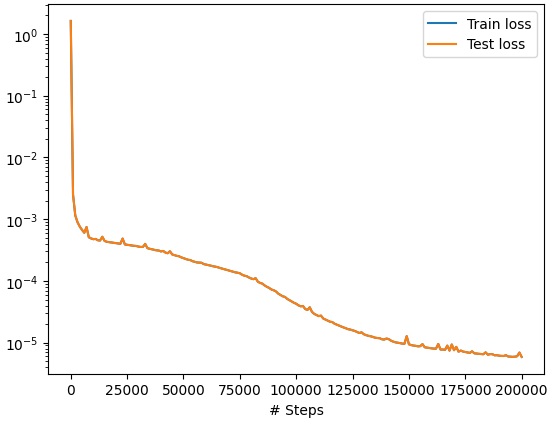}
		\caption{
			\scriptsize
			% Lowest loss 200,000.
			Init\_lr=0.001, loss, \red{Glorot}.
		}
		\label{fig:23.9.5 R3d.4 A-PF loss200000}
	\end{subfigure}
	\begin{subfigure}{0.24\textwidth}
		\includegraphics[width=1.\textwidth]{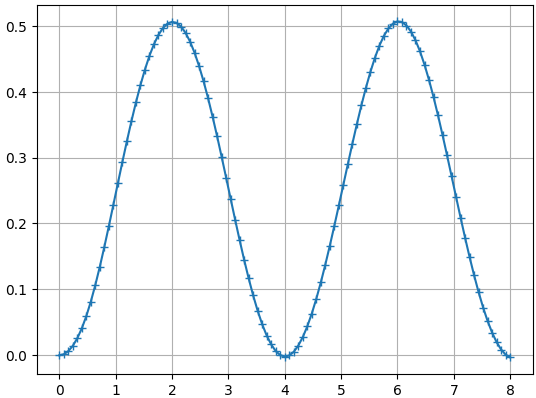}
		\caption{
			\scriptsize
			Free-end disp, \red{Quasi-perfect}.
		}
		\label{fig:23.9.5 R3d.4 A-PF free-end disp200000}
	\end{subfigure}
	\begin{subfigure}{0.24\textwidth}
		\includegraphics[width=1.\textwidth]{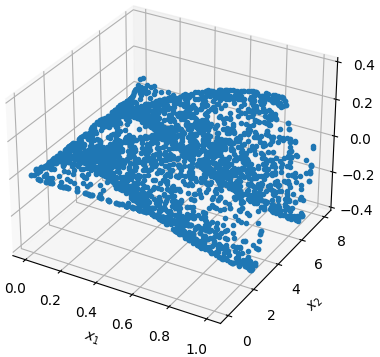}
		\caption{
			\scriptsize
			Velocity time history.
		}
		\label{fig:23.9.5 R3d.4 A-PF velocity200000}
	\end{subfigure}
	\begin{subfigure}{0.24\textwidth}
		\includegraphics[width=1.\textwidth]{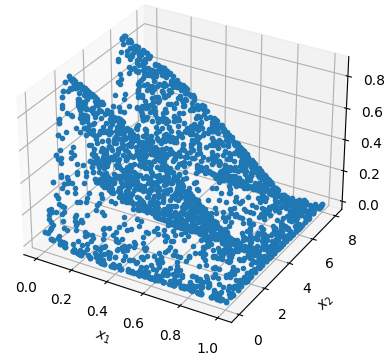}
		\caption{
			\scriptsize
			Slope t.h. 200,000.
		}
		\label{fig:23.9.5 R3d.4 A-PF slope200000}
	\end{subfigure}
	\caption{
		\DDET.
		% \emph{Pinned-free bar, pre-static solution}. 
		% \emph{Pinned-free bar}.
		\emph{Pinned-free bar, \red{Form 3}}
		(Section~\ref{sc:axial-Form-3}).
		Remark~\ref{rm:switch-Form-2a-to-Form-3}, Switch Form 2a to Form 3.
		$\bullet$
		SubFigs~\ref{fig:23.9.5 R3d A-PF free-end disp197000}-\ref{fig:23.9.5 R3d A-PF free-end slope197000},
		\emph{Pre-static solution.}
		Step 197,000:
		(\ref{fig:23.9.5 R3d A-PF free-end disp197000})
		Free-end displacement with velocity discontinuity at t=1 and upward shift.
		(\ref{fig:23.9.5 R3d A-PF free-end slope197000}) 
		Free-end slope jump with large amplitude 
		% two-period, 
		oscillations (negative space derivative) starting before t=1.
		% \red{Rewrite below 23.9.30.}
		\red{$\star$}
		Section~\ref{sc:axial-Form-3}, \red{Form~3}.
		\emph{Network:}
		Remark~\ref{rm:parameter-names},~\ref{rm:data-point-grids}, 
		T=8,
		W=64,
		H=4, 
		n\_out=3, 
		\red{He-uniform} initializer, 
		12,867 parameters,
		\emph{random} grid.
		\emph{Training:}
		\red{$\star$}
		Remark~\ref{rm:learning-rate-schedule-3} (LRS~3, NCA),
		init\_lr=0.005,
		% N\_steps=200,000, 
		% Cycle 5 lowest loss xxx
		GPU time 537 sec.
		{\scriptsize (2395R3d-1)}
		% {\color{red} REWRITE below.}
		$\bullet$
		SubFigs~\ref{fig:23.9.5 R3d.2 A-PF shape50000}-\ref{fig:23.9.5 R3d.2 A-PF free-end disp50000}, \emph{\red{Static}} shape, free-end disp, Step 50,000, init\_lr=0.01, He-uniform. 
		{\scriptsize (2395R3d.2)}
		$\bullet$
		SubFigs~\ref{fig:23.9.5 R3d.5 A-PF free-end disp200000}-\ref{fig:23.9.5 R3d.5 A-PF free-end slope200000}, \red{pre-static} free-end disp, slope (space derivative), Step 200,000, init\_lr=0.001, He-uniform.
		{\scriptsize (2395R3d.5)}
		$\bullet$
		SubFigs~\ref{fig:23.9.5 R3d.3 A-PF shape184000}-\ref{fig:23.9.5 R3d.3 A-PF free-end disp184000}, lowest loss \mbox{1.851e-06}, Step 187,000, shape, \red{Very good} free-end disp, \red{Glorot-uniform} initializer (Remark~\ref{rm:He-vs-Glorot}), init\_lr=0.005, damping\%=-0.7\%.
		{\scriptsize (2395R3d.3)}
		\red{$\star$}
		{\color{OliveGreen} \mbox{\bf SubFigs}}~\ref{fig:23.9.5 R3d.4 A-PF loss200000}-\ref{fig:23.9.5 R3d.4 A-PF slope200000}, Step 200,000: (\ref{fig:23.9.5 R3d.4 A-PF loss200000}) Lowest loss 5.877e-06, 
		init\_lr=0.001, 
		\red{Glorot-uniform},
		% Cycle 5 lowest loss 5.877e-06, 
		damping\%=\mbox{-0.13\%}, 
		(\ref{fig:23.9.5 R3d.4 A-PF free-end disp200000})
		\emph{\red{Quasi-perfect}} free-end disp, 
		(\ref{fig:23.9.5 R3d.4 A-PF velocity200000})
		Velocity time history (t.h.),
		(\ref{fig:23.9.5 R3d.4 A-PF velocity200000})
		Slope (space derivative) t.h.,
		GPU \red{\mbox{533 sec}} (Deterministic mode). 
		{\scriptsize (2395R3d.4)}
		$\bullet$
		{\footnotesize
			Figure~\ref{fig:23.9.5 R3b A-P bcrF F(1J)2a(3) lr0.005 cy1-5-NCA Nsteps200000 random N51 W64 H4 He init-1}, \red{Form 2a}, same model parameters as SubFigs~\ref{fig:23.9.5 R3d A-PF free-end disp197000}-\ref{fig:23.9.5 R3d A-PF free-end slope197000}, \red{Quasi-perfect} free-end disp.
			$\triangleright$
			Figure~\ref{fig:axial-Mathematica-solutions}, reference solution to compare.
		}
	}
	\label{fig:23.9.5 R3d A-P bcrF F(1J,2a)3 lr0.005 cy1-5-NCA Nsteps200000 random N51 W64 H4 He init-1}
\end{figure}

\begin{figure}[tph]
	\begin{subfigure}{0.32\textwidth}
		\includegraphics[width=\textwidth]{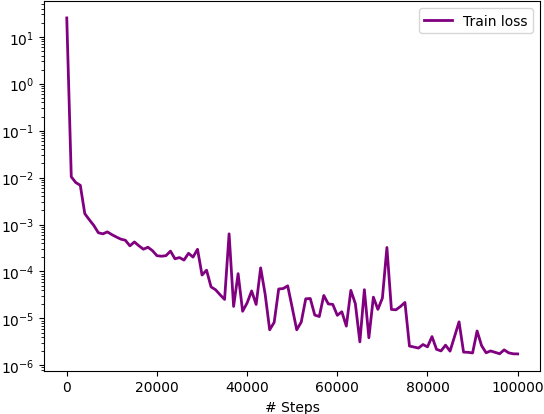}
		\caption{
			\scriptsize
			\JAX\ Form 1, loss 100,000.
		}
		\label{fig:23.9.16 R1a.8 loss100000}
	\end{subfigure}
	\begin{subfigure}{0.32\textwidth}
		\includegraphics[width=\textwidth]{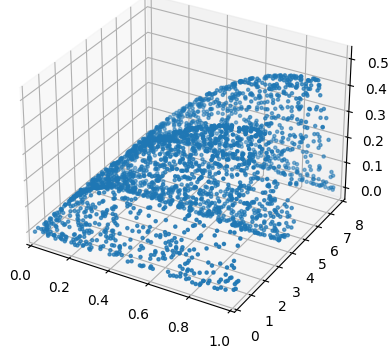}
		\caption{
			\scriptsize
			Shape time history.
		}
		\label{fig:23.9.16 R1a.8 shape100000}
	\end{subfigure}
	\begin{subfigure}{0.32\textwidth}
		\includegraphics[width=\textwidth]{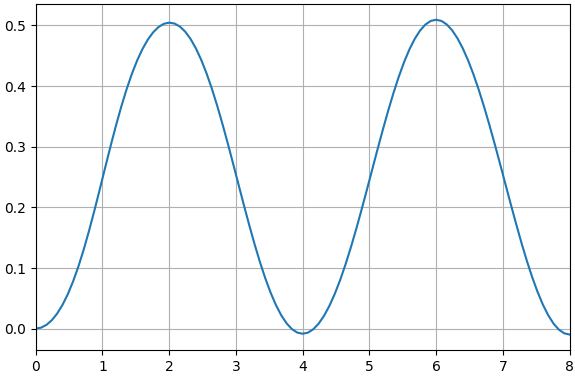}
		\caption{
			\scriptsize
			Free-end disp, \red{Very good}.
		}
		\label{fig:23.9.16 R1a.8 free-end disp100000}
	\end{subfigure}
	\\
	\begin{subfigure}{0.24\textwidth}
		\includegraphics[width=\textwidth]{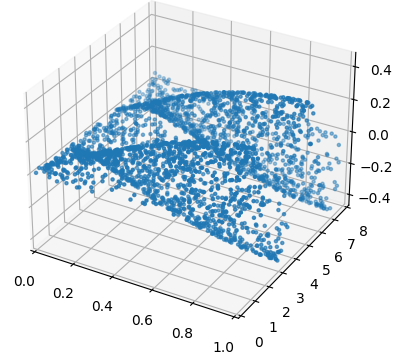}
		\caption{
			\scriptsize
			\JAX\ Form 3, velocity 145,000.
		}
		\label{fig:23.9.16 R1a.7.2 velocity145000}
	\end{subfigure}
	\begin{subfigure}{0.24\textwidth}
		\includegraphics[width=\textwidth]{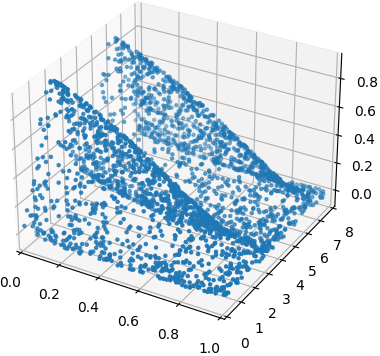}
		\caption{
			\scriptsize
			Slope time history.
		}
		\label{fig:23.9.16 R1a.7.2 slope145000}
	\end{subfigure}
	\begin{subfigure}{0.24\textwidth}
		\includegraphics[width=\textwidth]{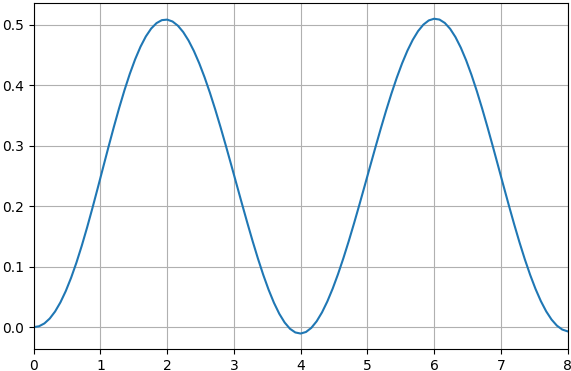}
		\caption{
			\scriptsize
			Free-end disp, \red{Quasi-perfect}.
		}
		\label{fig:23.9.16 R1a.7.2 free-end disp145000}
	\end{subfigure}
	\begin{subfigure}{0.24\textwidth}
		\includegraphics[width=\textwidth]{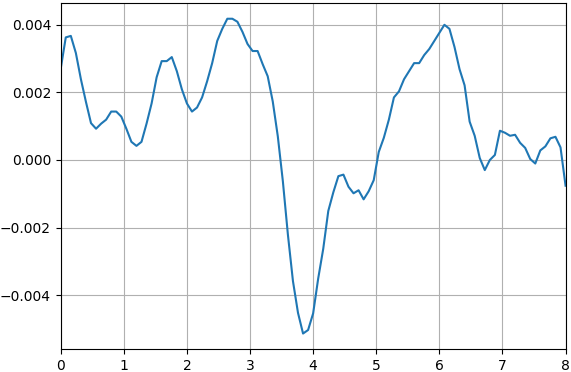}
		\caption{
			\scriptsize
			Free-end slope. 
			% (space derivative).
		}
		\label{fig:23.9.16 R1a.7.2 free-end slope145000}
	\end{subfigure}
	\caption{
		\JAX.
		\emph{Pinned-free bar.}
		SubFigs~\ref{fig:23.9.16 R1a.8 loss100000}-\ref{fig:23.9.16 R1a.8 free-end disp100000}, 
		\red{Form 1:} 
		(\ref{fig:23.9.16 R1a.8 loss100000}) 
		Step 100,000, Total loss \mbox{1.749e-06} (sum of 5 losses), Average loss 0.350e-06. 
		\emph{Network:}
		Remark~\ref{rm:parameter-names},~\ref{rm:data-point-grids}, 
		T=8,
		W=64,
		H=4, 
		n\_out=1, 
		\red{Uniform} initializer, 
		12,737 parameters,
		\emph{random} grid.
		\emph{Training:}
		\red{$\star$}
		Remark~\ref{rm:learning-rate-schedule-3} (LRS~4),
		init\_lr=0.002,
		factor\_lr = [0.2, 0.2, 0.2, 0.2, 0.2, 1, 1, 1, 1].
		(\ref{fig:23.9.16 R1a.8 shape100000}) Shape time history (t.h.). Damping\%=-0.94\%, 
		(\ref{fig:23.9.16 R1a.8 free-end disp100000}) \red{Very good} free-end displacement.
		% PLOTS: Step no. 100000, Train loss TOTAL = 1.7487657260062406e-06, Form F1 Train loss AVERAGE = 3.4975314520124814e-07
		{\scriptsize (23916R1a.8)}
		$\bullet$ 
		\red{Form~3:}
		(\ref{fig:23.9.16 R1a.7.2 velocity145000})
		Step \red{145,000}, velocity t.h., Total loss \red{4.063e-06}, Average loss 0.580e-06.
		\emph{Network:}
		Same as the above, except
		% for SubFigs~\ref{fig:23.9.16 R1a.8 loss100000}-\ref{fig:23.9.16 R1a.8 free-end disp100000}, except 
		n\_out=3, 
		% \red{Uniform} initializer, 
		12,867 parameters.
		% \emph{random} grid.
		\emph{Training:}
		Same as the above. 
		%for SubFigs~\ref{fig:23.9.16 R1a.8 loss100000}-\ref{fig:23.9.16 R1a.8 free-end disp100000}.
		(\ref{fig:23.9.16 R1a.7.2 slope145000})
		Slope (space derivative) t.h.
		Damping\%=-0.3\%,
		(\ref{fig:23.9.16 R1a.7.2 free-end disp145000})
		\red{Quasi-perfect} free-end disp t.h.. 
		(\ref{fig:23.9.16 R1a.7.2 free-end slope145000})
		Free-end slope (space derivative) t.h.,
		peak-to-peak amplitude 8.658e-03.
		{\scriptsize (23916R1a.7.2)}
	}
	\label{fig:23.9.16 R1a.7.2}
\end{figure}

\begin{rem}
	\label{rm:switch-Form-2a-to-Form-3}
	Switching from Form 2a to Form 3,  
	% \DvJ.
	\DDET.
	{\rm
		To obtain results using Form 3 in \DDET, the starting point was to use the same parameters (with init\_lr=0.005, He-uniform initializer) used for Form 2a in \DDET\ that resulted in Figure~\ref{fig:23.9.5 R3b A-P bcrF F(1J)2a(3) lr0.005 cy1-5-NCA Nsteps200000 random N51 W64 H4 He init-1}, with only Form 2a changed to Form 3; the result was a pre-static solution shown in SubFigs~\ref{fig:23.9.5 R3d A-PF free-end disp197000}-\ref{fig:23.9.5 R3d A-PF free-end slope197000}.
		% ALEX: removed paragraph
		Indeed,
		increasing init\_lr to 0.01 resulted in a static solution shown in SubFigs~\ref{fig:23.9.5 R3d.2 A-PF shape50000}-\ref{fig:23.9.5 R3d.2 A-PF free-end disp50000}, and reducing init\_lr to 0.001 was not enough to avoid a pre-static solution, SubFigs~\ref{fig:23.9.5 R3d.5 A-PF free-end disp200000}-\ref{fig:23.9.5 R3d.5 A-PF free-end slope200000}.
		% ALEX: removed paragraph
		Returning init\_lr to 0.005 and 
		switching from He-uniform to Glorot-uniform initializer resulted in a very-good free-end displacement with a Cycle 5 lowest loss of 1.851e-06 at Step 187,000, as shown in  SubFigs~\ref{fig:23.9.5 R3d.3 A-PF shape184000}-\ref{fig:23.9.5 R3d.3 A-PF free-end disp184000}.
		% ALEX: removed paragraph
		Decreasing init\_lr to 0.001, maintaining the Glorot-uniform initialization, then resulted in a quasi-perfect free-end displacement, albeit with a higher Cycle 5 lowest loss of 5.877e-06 at Step 200,000, as shown in {\color{OliveGreen} \mbox{\bf SubFigs}}~\ref{fig:23.9.5 R3d.4 A-PF loss200000}-\ref{fig:23.9.5 R3d.4 A-PF free-end disp200000}.
		
		It is important to note that the above results do not absolutely imply the superiority of the Glorot initializer over the He initializer, but as noted before in Remark~\ref{rm:He-vs-Glorot}, the Glorot initializer is equivalent to a smaller learning rate compared to the He initializer.
	}
	\hfill$\blacksquare$
\end{rem}

\begin{rem}
	\label{rm:efficiency-DDE-T-vs-JAX}
	Efficiency: Form 1 vs Form 2a vs Form 3.  \DvJ.
	{\rm
		First we used our \DDET\ script.
		Based on our numerical experiments on the axial motion of a bar, Form 3 is computationally more efficient than Form 2a, which is itself computationally more efficient than Form 1.  
		
		% 23.10.1
		% Googld doc: "[1] F3 Axial motion pinned-pinned"
		% Section "GPU time average F1H 624 sec, F2aH 561 sec, F3 531 sec" 
		For the \emph{pinned-pinned bar},
		going from Form 3 
		(531 sec) 
		to Form 2a 
		% (638 sec) 
		(561 sec)
		requires 
		% 20\% 
		6\%
		more computational time.  Going from Form 2a to Form 1 
		% (692 sec)
		(624 sec) 
		requires 
		% 8\% 
		11\%
		more computational time.  Going from Form 3 to Form 1 requires 
		% 30\% 
		18\%
		more computational time.\footnote{
			For each form, three runs were carried out to obtain the average GPU time.  The additional GPU time required for moving from one form to another was based on these average GPU times. 
			For the pinned-pinned bar, see RunIDs 2388R2a,b,c for Form 3 (531 sec), RunIDs 2388R6a,b,c for Form 2a (561 sec), and RunIDs 23101R3a,b,c for Form 1 with Hessian (624 sec). 
		}
		For the \emph{pinned-pinned} bar, Form 1 using the Hessian to compute the 2nd derivative has an average GPU time of 624 sec (RunIDs 2388R5a,b,c).
		Switching from a Hessian 
		%
		% 24.7.16, code snippets
%		(Listing~\ref{lst:DDE-T time derivative by hessian}) 
		to two Jacobians 
		%
		% 24.7.16, code snippets
%		(Listing~\ref{lst:DDE-T time derivative by two jacobians}) 
		increased the average GPU time to 693 sec, or 11\%.  With the Form 3 GPU time remaining the same at 531 sec, going from Form 3 to Form 1 using two Jacobians requires 31\% more computational time.
		
		For the \emph{pinned-free bar}, the results (obtained in Deterministic mode) in Figure~\ref{fig:DDE-T-Form-1-pinned-free-bar-static-solutions} (Form 1, 655 sec), Figure~\ref{fig:23.9.5 R3b A-P bcrF F(1J)2a(3) lr0.005 cy1-5-NCA Nsteps200000 random N51 W64 H4 He init-1} (Form 2a, 605 sec), and Figure~\ref{fig:23.9.5 R3d A-P bcrF F(1J,2a)3 lr0.005 cy1-5-NCA Nsteps200000 random N51 W64 H4 He init-1} (Form 3, 533 sec), indicate that going from From 3 to Form 2a requires 14\% more computational time,
		going from Form 2a to Form 1 requires 8\% more computational time,
		going from Form 3 to Form 1 requires 23\% more computational time.
		The \DDET\ results in Non-deterministic mode are tabulated in Table~\ref{tb:pinned-free-bar-DDE-T-vs-JAX}, which is explained below.
		% ALEX: removed paragraph
		These percentages would change with the number of derivatives on a PDE system, with the axial motion having the lowest number of derivatives considered here, with two space derivatives and two time derivatives.
		
		After obtaining the above results for the pinned-free bar with our \DDET\ script, we ran the same cases with our \JAX\ script to compare Form 1 to Form 3.  
		Using JAX Form 3, the total loss of, 5.135e-06 (average loss 0.734e-06) at Step 97000 after 357 sec (at Step 100,000, which is before Step 145,000 the results at which are shown in SubFigs~\ref{fig:23.9.16 R1a.7.2 velocity145000}-\ref{fig:23.9.16 R1a.7.2 free-end slope145000}) is a good trade-off for the loss 5.877e-06 at Step 200,000 (end of Cycle 5) after 533 sec using \DDET\ Form 3 in \mbox{SubFigs}~\ref{fig:23.9.5 R3d.4 A-PF loss200000}-\ref{fig:23.9.5 R3d.4 A-PF slope200000}.
		In SubFigs~\ref{fig:23.9.16 R1a.7.2 velocity145000}-\ref{fig:23.9.16 R1a.7.2 free-end slope145000}, it took 529 sec to train 150,000 steps of \JAX\ Form~3 to get a loss of 4.104e-06 and Cycle 4 lowest loss of 4.063e-06 at Step 145000. 
		
		All results depend crucially on the learning-rate schedule and the initializer used, meaning in principle one could conceivably design a learning-rate schedule with an appropriate initializer for DDE-T Form 3 to reach a total loss below 4.104e-06 in 150000 steps or less. Hence in the end, the GPU time for a fixed number of steps is important and more meaningful for comparing between two different implementations.
		
		The GPU times for \DDET\ and for \JAX\ for 200,000 training steps in each form in Non-deterministic mode are tabulated in Table~\ref{tb:pinned-free-bar-DDE-T-vs-JAX}, which shows that
		JAX Form 1 is faster than JAX Form 2a, which is faster than JAX Form 3, the \emph{opposite} of DDE-T.\footnote{
			Table~\ref{tb:pinned-free-bar-DDE-T-vs-JAX}.
			\label{fn:DDE-T vs JAX 3 runs}
			For \DDET\ \red{Form 1}, the average GPU time of three runs was 529 sec: RunIDs
			2398R2a.3 (538 s), 2398R2a.4 (523 s), 2398R2a.5 (525 s).
			DDE-T \red{Form 2a}, average GPU time 498 sec:
			2398R2a.6 (503 s), 2398R2a.6 (504 s), 2398R2a.6 (487 s).
			DDE-T \red{Form 3}, average GPU time 488 sec:
			2398R2a.9 (499 s), 2398R2a.10 (477 s), 2398R2a.11 (488 s). 
			$\bullet$
			For \JAX\ \red{Form 1}, average GPU time 575 sec: 
			23916R1a8 (588 s), 23.9.16 R1a.8.3 (554 s), 23.9.16 R1a.8.4 (583 s).  
			JAX \red{Form 2a}, average GPU time 
			% four 
			% three runs 
			644 sec: 
			% RunIDs 
			% 231016R1 (642 sec), 
			231024R2c (618 s), 231025R2b (640 s), 231025R2c (675 s). 
			JAX \red{Form 3}, average GPU time 733 sec:  
			23916R1a7 (701 s), 23916R1a7.2 (751 s), 23916R1a7.3 (747 s).
		}    
		The results in Table~\ref{tb:pinned-free-bar-DDE-T-vs-JAX} show that our \DDET\ script is faster than our \JAX\ script.
		Using DDE-T GPU time as reference, 
		DDE-T Form 1 (529 s) is more efficient than JAX Form 1 (575 s) by 
		% (575 - 529) / 529 = 9\%,
		9\%, 
		whereas 
		DDE-T Form 2a (498 s) is more efficient than JAX Form 2a  
		(644 s)
		by
		% (644 - 498) / 498 = 29\%
		29\%, and
		DDE-T Form 3 (488 sec) is more efficient than JAX Form 3 (733 sec) by 
		% (733 - 488) / 488 =  50\%.
		50\%.
	}
	% \phantom{Blank}
	\hfill$\blacksquare$
\end{rem}

\begin{table}[tph]
	\centering
	\caption{
		\emph{Pinned-free bar.}  
		Remark~\ref{rm:efficiency-DDE-T-vs-JAX}, Efficiency, \DvJ.
		GPU time for 200,000 training steps.  All ran in \emph{Non-deterministic} mode.
		Each number represents the average GPU time of three runs (see Footnote~\ref{fn:DDE-T vs JAX 3 runs}).
		\emph{Network:}
		Remark~\ref{rm:parameter-names},~\ref{rm:data-point-grids}, 
		T=8,
		W=64,
		H=4.
%		$\bullet$
%		{\footnotesize
%			Figure~\ref{fig:DEBUG 23.9.4 R1d A-P bcrP F1J(2a,3) lr0.001 Cy1-9-NCA Nsteps400000 regular N51 W64 H4 Glorot init-1}, HELLO
%		}
	}
	\label{tb:pinned-free-bar-DDE-T-vs-JAX}
	\vspace{3mm}
	\begin{tabular}{r| 
			>{\columncolor[HTML]{FFFFFF}}r| r|r|r|r}
		\multicolumn{1}{r|}{\phantom{\DDET}} 	& 
		\multicolumn{1}{r|}{Form 1} 			& 
		\multicolumn{1}{r|}{Form 2a} 			& 
		\multicolumn{1}{r|}{Form 3} 
		% \vspace{-0.5mm}                      
		\\	\hline\hline
		% \\[-4mm]                        
		\DDET & 529 sec & 498 sec  	& 488 sec  
		\\ \hline
		\JAX  & 575 sec & 644 sec  	& 733 sec  
	\end{tabular}
\end{table}

\begin{rem}
	\label{rm:JAX-deterministic-faster}
	\JAX\ Deterministic slightly faster than Non-deterministic.
	{\rm
		Unlike \DDET\ for which running in deterministic mode carries an average penalty of about 4\% in GPU time compared to non-deterministic mode (or 10\% for a single case as noted in Remark~\ref{rm:stochastic-reproducibility}), it is surprising to note that \JAX\ in deterministic mode is more efficient than in non-deterministic mode.
		Re-running each \JAX\ non-deterministic case in Table~\ref{tb:pinned-free-bar-DDE-T-vs-JAX} three times in deterministic mode, the corresponding average GPU times are: 
		550 sec (JAX Form 1, or 4\% less than non-deterministic), 
		% (575 - 550) / 575
		628 sec (JAX Form 2a, or 2.5\% less), 
		% (644 - 628) / 644
		689 sec (JAX Form 3, or 6\% less).
		% (733 - 689) / 733   
	}
	\phantom{blank}\hfill$\blacksquare$
\end{rem}

% {\color{red} NOTE: 2023.08.21, Write Form 3 here. ENDNOTE}

%\subsubsection{Space derivative constraint}
%\subsubsection{Time derivative constraint}
%\subsubsection{Error at end of time interval}

%\red{I AM HERE 2023.09.15.}
% \red{I AM HERE 2023.10.10.}

% \subsection{Euler-Bernoulli beam}

% \subsection{Kirchhoff-Love rod}

% stochastic gradient descent
\section{Closure}
\noindent
Following up our review of deep learning applied to computational mechanics \cite{vuquoc2023deep}, we developed 
%\red{novel} 
novel
PINN formulations for nonlinear dynamics governed by partial-differential-algebraic equations (PDAEs), exemplified by the equations of motion of a geometrically-exact Kirchhoff rod in higher-level forms, initially developed to address the pathological problems of shift and amplification encountered with \DDET, which is likely the most well-documented among PINN frameworks.

For the static-solution problem, we consider the use of barrier functions, which themselves represent a 
%\red{novel} 
novel 
application in PINN, to prevent the training iterate to get close to a static solution.  Even though the barrier function is an elegant method, it still has room for improvement, while simpler methods such as reducing the model capacity (using lower number of parameters), use a different initializer, and smaller learning rate led to satisfactory results.

Another 
% \red{novelty}
novelty 
is to apply PINN directly to the highest-level momentum form (balance of momenta) as it would save much time and effort to hand derive from the highest-level from to the lowest-level form through to computational formulation and implementation as done in {traditional} approach.

In parallel with using \DDET, we developed a script based on \JAX, which does not exhibit the pathological problems of \DDET, 
% and is faster than \DDET\ for the lower-level Form 1 of axial motion of an elastic bar, but much slower for the highest-level Form 3 (or 4).
it is slower than \DDET\ in all forms tested. 

To help readers reproduce our results, we 
%\red{normalize/standardize the training process}, 
normalize/standardize the training process,
which is important to compare different training strategies due to the large number of parameters involved.  Each figure is accompanied with a caption that details all parameters used to obtain the results.
%
% 24.7.16, code snippets  
%In addition, we provide extensively annotated code snippets of how we implemented the gradient computation, in both \DDET\ and \JAX, as the computed gradients likely hold the key for the behavior (pathological problems and computational speed) observed.
The computed gradients in both \DDET\ and \JAX likely hold the key for the behavior (pathological problems and computational speed) observed.

%Forward and inverse problems and results involving the Kirchhoff rod will be presented in a follow-up publication. 
%\red{Large computational efficiency with the proposed PINN formulations recently obtained for the transverse motion of the Euler-Bernoulli beam (more than 400\%) and for the geometrically-exact Kirchhoff rod (more than 1000\%) will be reported in a follow-up paper.}
Large computational efficiency with the proposed PINN formulations recently obtained for the transverse motion of the Euler-Bernoulli beam and for the geometrically-exact Kirchhoff rod will be reported in a follow-up paper.

% 2020.02.17
% add References to table of contents
%https://tex.stackexchange.com/questions/8458/making-the-bibliography-appear-in-the-table-of-contents
\addcontentsline{toc}{section}{References}
% 2019.02.10, need to put all bib files on one line for TexStudio to automatically 
% list the bib entry labels to select the desired label, without having to remember
% the desired label or having to look up the desired label in the bib file.  
% splitting the bib files into several lines took away this useful feature.
\bibliography{data-driven-computation-mechanic,data-driven-computing-or-computation,deep-learning_alex,deep-learning-comput-mechanic,deep-learning-misc,ghaboussi,mechanics,miscellaneous,neural-network-mechanics,vuquoc-DL}
\label{sc:references}

\newpage
% 2022.06.08 - USE the appendices environment for CMES style
\noindent
{\Large \bf Appendices}
\addcontentsline{toc}{section}{Appendices}
\begin{appendices}
	
	% Form 1, time shift, error in slope of free end
	% \input{Time shift}
	
	% analysis of time shift and amplification
	% \subsubsection{Analysis of time shift and amplification}
\section{Analysis of time shift and amplification}
\label{sc:time-shift-amplification}
\label{app:analysis-time-shift-amplification}
First, we analyze the computed solution of Form 1 of the axial motion of a pinned-pinned elastic bar with distributed axial force (Section~\ref{sc:axial-Form-1}), with $\xb \approx \Xb$:
\begin{align}
	\text{PDE:} \quad
	&
	% \sldn 
	\ubp{2} + \dfbs{\x} = \nddt \ub
	\ ,
	\tag{\ref{eq:eom-euler-bernoulli-axial}}
	\\
	\text{BCs:} \quad
	&
	\ub (0, \tb) = \ub (1, \tb) = 0
	\ ,
	\tag{\ref{eq:axial-pinned-pinned-BCs}}
	\\
	\text{ICs:} \quad
	&
	\ub (\xb , 0) = 0 , \quad
	\ubd (\xb , 0) = 0 
	\ .
	\tag{\ref{eq:axial-ICs}}
\end{align}
To alleviate the notation, we omit the overbar designating non-dimensionalization on $(x,t)$, while keeping the overbar in $\ub$ and $\dfbs{\x}$.
Let the computed solution $\uc (x,t)$ be related to the exact solution $\ub (x,t)$ by (see Figures~\ref{fig:DEBUG v1.9.3 seed42, 23.9.4 R1 A-P bcrP F1J(2a,3) lr0.001 Cy1-9-NCA Nsteps400000 regular N51 W64 H4 Glorot init-1}-\ref{fig:DEBUG 23.9.4 R1d A-P bcrP F1J(2a,3) lr0.001 Cy1-9-NCA Nsteps400000 regular N51 W64 H4 Glorot init-1}):
\begin{align}
	\uc (x,t) = \shiftPa \ub (x, t - \shiftPs) + c
	\Rightarrow
	\ub (x,t) = \frac{1}{\shiftPa} \left[\uc (x, t + \shiftPs) - \shiftPc \right]
	\ ,
	\label{eq:computed-solution-exact-solution}
\end{align}
where $\shiftPa$ is the amplification parameter, $\shiftPs$ the time shift parameter, and $\shiftPc$ the vertical shift parameter.  
These are the three \emph{hidden} parameters, in addition to the network parameters, for the training to adjust to reduce the loss.  Figure~\ref{fig:DEBUG 23.9.4 R1d A-P bcrP F1J(2a,3) lr0.001 Cy1-9-NCA Nsteps400000 regular N51 W64 H4 Glorot init-1} shows that as the shift and amplification (scaled from their true values in Table~\ref{tb:DEBUG 23.9.4 R1d A-P bcrP shift-amplification-values}) increase, the total loss continues to decrease during training.

% table obtained using https://www.tablesgenerator.com/#
% just copy and paste from Google Doc
% 
% Please add the following required packages to your document preamble:
% \usepackage[table,xcdraw]{xcolor}
% If you use beamer only pass "xcolor=table" option, i.e. \documentclass[xcolor=table]{beamer}
\begin{table}[tph]
	\centering
	\caption{
		\emph{Pinned-pinned bar.  Form 1.}  
		Section~\ref{sc:time-shift-amplification}, \emph{Analysis}.
		\emph{Unscaled} shift-amplification parameters ($\shiftPs , \shiftPc$, $\shiftPa$), with $\shiftPA$ being constant, initial velocity, left-end displacement versus checkpoint training step number.
		$\bullet$
		{\footnotesize
			Figure~\ref{fig:DEBUG 23.9.4 R1d A-P bcrP F1J(2a,3) lr0.001 Cy1-9-NCA Nsteps400000 regular N51 W64 H4 Glorot init-1}, \emph{scaled} shift-amplification parameters and initial velocity increase with training, using 
			$p_{\text{scaled}} = \pm (p - p_{\text{25000}}) * k_p > 0$, where $p_{\text{scaled}}$ is the positive scaled parameter $p$, with $p_{\text{25000}}$ the value of $p$ at Step 25000, and $k_p$ a selected scale factor for $p$:
			% the scaled factors: 
			$k_\shiftPs = 10$,
			$k_\shiftPc = 100$,
			$k_{(\shiftPa \shiftPA)} = 10$,
			$k_{\text{Init vel}} = 10$.
		}
	}
	\label{tb:DEBUG 23.9.4 R1d A-P bcrP shift-amplification-values}
	\vspace{3mm}
	\begin{tabular}{r| 
			>{\columncolor[HTML]{FFFFFF}}r| r|r|r|r}
		\multicolumn{1}{r|}{Step} & \multicolumn{1}{r|}{Time shift $\shiftPs$} & \multicolumn{1}{r|}{Vertical shift $\shiftPc$} & \multicolumn{1}{r|}{Amplitude $\shiftPa \shiftPA$} & \multicolumn{1}{r|}{Initial velocity} &
		\multicolumn{1}{r}{Left-end disp}
		% \vspace{-0.5mm}                      
		\\	\hline\hline
		% \\[-4mm]                        
		25000 & 0 & {0.00973727}  & {0.08445985} & \cellcolor[HTML]{FFFFFF}{\color[HTML]{212121} 0.01150370} &
		0.00973727 \\
		50000                    & 0.107                      & \cellcolor[HTML]{FFFFFF}{\color[HTML]{212121} -0.00462352} & \cellcolor[HTML]{FFFFFF}{\color[HTML]{212121} 0.13343337} & \cellcolor[HTML]{FFFFFF}{\color[HTML]{212121} -0.06753206} & 
		-0.00094142 \\
		100000                   & {0.146}                      & \cellcolor[HTML]{FFFFFF}{\color[HTML]{212121} -0.00861757} & \cellcolor[HTML]{FFFFFF}{\color[HTML]{212121} 0.14141834} & \cellcolor[HTML]{FFFFFF}{\color[HTML]{212121} -0.09649991} &
		-0.00139512 \\
		150000                   & {0.168}                      & \cellcolor[HTML]{FFFFFF}{\color[HTML]{212121} -0.01143161} & \cellcolor[HTML]{FFFFFF}{\color[HTML]{212121} 0.14696970}  & \cellcolor[HTML]{FFFFFF}{\color[HTML]{212121} -0.11479854} &
		-0.00152922 \\
		200000                   & {0.184}                      & \cellcolor[HTML]{FFFFFF}{\color[HTML]{212121} -0.01377921} & \cellcolor[HTML]{FFFFFF}{\color[HTML]{212121} 0.15177655} & \cellcolor[HTML]{FFFFFF}{\color[HTML]{212121} -0.12892485} &
		-0.00155308 \\
		250000                   & {0.196}                      & \cellcolor[HTML]{FFFFFF}{\color[HTML]{212121} -0.01580278} & \cellcolor[HTML]{FFFFFF}{\color[HTML]{212121} 0.15583168} & \cellcolor[HTML]{FFFFFF}{\color[HTML]{212121} -0.14010072} &
		-0.00157554 \\
		300000                   & {0.204}                      & \cellcolor[HTML]{FFFFFF}{\color[HTML]{212121} -0.01743147} & \cellcolor[HTML]{FFFFFF}{\color[HTML]{212121} 0.15909590}  & \cellcolor[HTML]{FFFFFF}{\color[HTML]{212121} -0.14862418} &
		-0.00158554 \\
		350000                   & {0.212}                      & \cellcolor[HTML]{FFFFFF}{\color[HTML]{212121} -0.01873142} & \cellcolor[HTML]{FFFFFF}{\color[HTML]{212121} 0.16170362} & \cellcolor[HTML]{FFFFFF}{\color[HTML]{212121} -0.15547872} &
		-0.00158999 \\
		400000                   & {0.216}                      & \cellcolor[HTML]{FFFFFF}{\color[HTML]{212121} -0.01975145} & \cellcolor[HTML]{FFFFFF}{\color[HTML]{212121} 0.16374922} & \cellcolor[HTML]{FFFFFF}{\color[HTML]{212121} -0.16018747} &
		-0.00159215
	\end{tabular}
\end{table}

The peak-to-peak amplitude of $\ub$ (if known), or of a computed quasi-perfect solution, used as a reference is denoted by $\shiftPA$, so that $\shiftPa \shiftPA$ is the peak-to-peak amplitude of the computed solution $\uc$ (Figure~\ref{fig:DEBUG v1.9.3 seed42, 23.9.4 R1 A-P bcrP F1J(2a,3) lr0.001 Cy1-9-NCA Nsteps400000 regular N51 W64 H4 Glorot init-1}).

The space and time derivatives of $\uc (x, t)$ are
\begin{align}
	&
	\uc_{xx} (x, t) = \shiftPa \uc_{xx} (x, t - \shiftPs)
	\Rightarrow
	\uc_{xx} (x, t) = \frac{1}{\shiftPa} \uc_{xx} (x, t + \shiftPs)
	\ , 
	\\
	&
	\uc_{t} (x, t) = \frac{1}{\shiftPa} \uc_{t} (x, t + \shiftPs)
	\ , \quad
	\uc_{tt} (x, t) = \frac{1}{\shiftPa} \uc_{tt} (x, t + \shiftPs)
	\ ,
	\label{eq:computed-solution-derivatives}
\end{align}
which when substituted in Eq.~\eqref{eq:eom-euler-bernoulli-axial} yield
\begin{align}
	\uc_{xx} (x, t + \shiftPs) + \shiftPa \dfbs{\x} (x,t) = \uc_{tt} (x , t + \shiftPs)
	\Rightarrow
	\uc_{xx} (x, t) + \shiftPa \dfbs{\x} (x,t - \shiftPs) = \uc_{tt} (x , t)
	\ ,
	\label{eq:eom-u-check-shifted}
\end{align} 
with boundary conditions at $x_b = 0$ and at $x_b = 1$:
\begin{align}
	% &
	\ub (x_b ,t) = 0 = \frac{1}{\shiftPa} \left[ \uc (x_b , t + \shiftPs) - c \right]
	\Rightarrow
	\uc (x_b , t + \shiftPs) = c
	\ ,
	% \\
	% &
\end{align}
and initial conditions:
\begin{align}
	&
	\ub (x ,0) = 0 = \frac{1}{\shiftPa} \left[ \uc (x , \shiftPs) - c \right]
	\Rightarrow
	\uc (x , \shiftPs) = c
	\ ,
	\\
	&
	\ub_t (x ,0) = 0 = \frac{1}{\shiftPa} \left[ \uc_t (x , \shiftPs) \right]
	\Rightarrow
	\uc_t (x , \shiftPs) = 0
	\ .
\end{align}

Due to the computed displacement being enforced by the Dirichlet boundary condition, i.e.,
\begin{align}
	\uc (x,0) = \epsilon \approx 0 = \shiftPa \ub (x, - \shiftPs) + c
	\Rightarrow
	c = - \shiftPa \ub (x, - \shiftPs) + \epsilon
	\ ,
\end{align}
where $\epsilon$ is a small number or the order of $10^{-3}$, the ``essentially zero'' for this type of network training, as shown in Table~\ref{tb:DEBUG 23.9.4 R1d A-P bcrP shift-amplification-values}.  
So if $\shiftPs > 0$ (right time shift), since $\ub (x, -\shiftPs) > 0$ and $\shiftPa > 0$, it follows that $c < 0$ if $\shiftPa \ub (x, - \shiftPs) > \epsilon$, which is a condition met in the example in Table~\ref{tb:DEBUG 23.9.4 R1d A-P bcrP shift-amplification-values}.

The Neumann boundary condition was used in DeepXDE for the initial velocity, but did not enforce the time derivative of the computed solution to the ``essentially zero'' $\epsilon$.  Instead, the initial velocity was allowed to be non-zero:
\begin{align}
	\uc_t (x,0) = \shiftPa \ub_t (x , -\shiftPs) \ne 0
	\ , 
\end{align}  
and thus if $\shiftPs > 0$, since $\ub_t (x , -\shiftPs) < 0$, it follows that $\uc_t (x, 0) < 0$, which is observed in Figure~\ref{fig:DEBUG v1.9.3 seed42, 23.9.4 R1 A-P bcrP F1J(2a,3) lr0.001 Cy1-9-NCA Nsteps400000 regular N51 W64 H4 Glorot init-1}. 

Second, we show that Form 2a of the pinned-pinned elastic bar (Section~\ref{sc:axial-Form-2a})
\begin{align}
	\text{PDE-1:} \quad
	&
	% \sldn 
	\ubp{2} + \dfbs{\x} = \pbsd{\x}
	\ , 
	% \quad
	% \ndt \ub = \pbs{\x}
	% \ ,
	\tag{\ref{eq:wave-eq-form-2}}
	\\
	\text{PDE-2:} \quad
	&
	\ndt \ub = \pbs{\x}
	\ ,
	\tag{\ref{eq:wave-eq-form-2}}
	\\
	\text{BCs:} \quad
	&
	\ub (0, \tb) = \ub (1, \tb) = 0
	\ ,
	\tag{\ref{eq:axial-pinned-pinned-BCs}}
	\\
	\text{ICs:} \quad
	&
	\ub (\xb , 0) = 0 , \quad
	% \ubd (\xb , 0) = 0 
	\pbs{\x} (\x , 0) = 0
	\ .
	\tag{\ref{eq:axial-ICs-Form-2a}}
\end{align}
have neither shift nor amplification.  The analysis for Form 3 is the same, and is not repeated.  The explicit Dirichlet boundary condition imposed an ``essentially zero'' $\epsilon$ on the initial velocity of the computed solution:
\begin{align}
	\uc_t (x,0) = \epsilon = \shiftPa \ub_t (x, 0 - \shiftPs)
	\Rightarrow
	\ub_t (x, -\shiftPs) = \epsilon / \shiftPa \approx 0
	\ ,
\end{align}
with $\shiftPa \ne 0$ and $\shiftPs$ in the vicinity of zero.  But the initial velocity of the exact solution is $\ub_t (x,0) = 0$.  Hence, $\uc_t (x, -\shiftPs) - \epsilon / a = \ub_t (x,0) = 0$, and thus $s \approx 0$ (``essentially zero'').
At this point $\shiftPa > 0$ can be anything.
The explicit Dirichlet boundary condition also imposes an ``essentially zero'' $\epsilon$ on the initial displacement of the computed solution:
\begin{align}
	\uc (x,0) = \epsilon = \shiftPa \ub (x,0 - \shiftPs) + \shiftPc =
	\shiftPa \ub (x,0) + \shiftPc
	\Rightarrow
	\shiftPc = \epsilon \approx 0
	\ ,
\end{align}
which is also an ``essentially zero'' (see the left-end displacement in Table~\ref{tb:DEBUG 23.9.4 R1d A-P bcrP shift-amplification-values}).  Now substituting $\ub = \uc / \shiftPa$ into the above Form 2a of the PDE system to obtain:
\begin{align}
	\text{PDE-1:} \quad
	&
	% \sldn 
	\uc_{xx} / \shiftPa + \dfbs{\x} = \pbsd{\x}
	\Rightarrow
	\uc_{xx} + \shiftPa \dfbs{\x} = \shiftPa  \pbsd{\x}
	\ , 
	\\
	\text{PDE-2:} \quad
	&
	\uc_t / \shiftPa = \pbs{\x}
	\Rightarrow
	\uc_t = \shiftPa \pbs{\x}
	\ ,
	\\
	\text{BCs:} \quad
	&
	\uc (0, t) = \uc (1, t) = 0
	\ ,
	\tag{\ref{eq:axial-pinned-pinned-BCs}}
	\\
	\text{ICs:} \quad
	&
	\uc (\x , 0) = 0 , \quad
	\pbs{\x} (\x , 0) = 0
	\ .
	\tag{\ref{eq:axial-ICs-Form-2a}}
\end{align}
Thus $\uc$ would be the solution of the above PDE system only if we had amplified the applied force $\dfbs{\x}$ by the amplification factor $\shiftPa$.  In other words, if we did not amplify the applied force, i.e., $\shiftPa = 1$, then $\uc = \ub$, or in practice, $\uc \approx \ub$.

	% barrier functions to filter static solutions
	\begin{figure}[tph]
	\centering
	\includegraphics[width=0.40\textwidth]{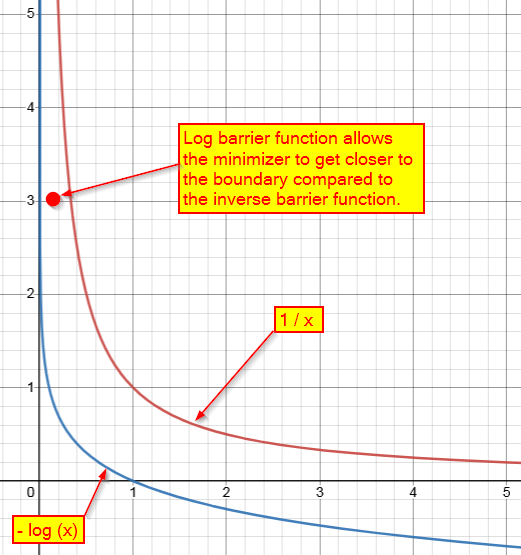}
	\includegraphics[width=0.40\textwidth]{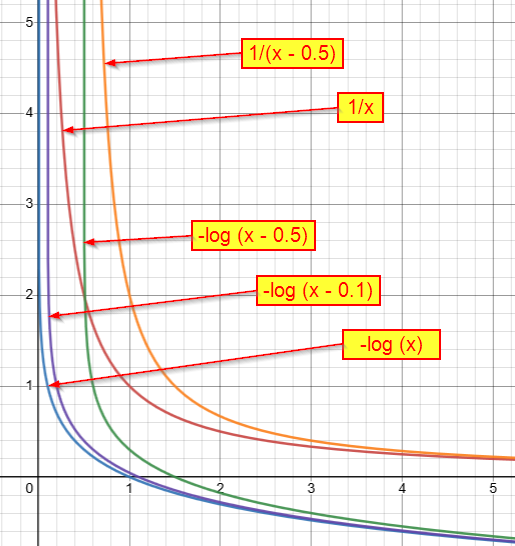}
	\caption{
		\emph{Inverse barrier function vs log barrier function.} Appendix~\ref{app:barrier-functions}. 
		% (Section~\ref{sc:filter-static-solution})
		Shifted barrier,  
		buffer-zone depth $\bde$.
	}
	\label{fig:inverse-log-barrier functions}
\end{figure}

\section{Filtering static solutions, barrier functions}
\label{sc:filter-static-solution}
\label{app:barrier-functions}
\noindent
A static solution is a solution of the governing equations of motion, such as Eqs.~\eqref{eq:eom-u-bar}-\eqref{eq:eom-v-bar}, with zero inertia force on the right hand side, abstractly written based on the dynamic nonlinear differential operator $\dpdesp{i}{(k)}$ defined in Eq.~\eqref{eq:generic-PDE-aux-conditions} as:
\begin{align}
	%&
	\left. \dpdesp{i}{(k)} \right|_{n_q^{(k)} = 0} 
	=
	&
	\spdesp{i}{(k)} (	\{\partial_\X^p u_j , \partial_\X^{p-1} u_j , \ldots, \partial_\X^0 u_j\} ) 
	= \spdesp{i}{(k)} (\busta (\X))
	% = \spdes{i} (\X)
	= 0
	\ , %\quad
	\nonumber
	\\
	&
	k= 0, \ldots , 3
	\ , \quad i = 1, \ldots, n_i^{(k)}
	\ ,
	\label{eq:generic-static-pde}
\end{align} 
where $\spdesp{i}{(k)}$ is a \emph{static} nonlinear differential operator, which is a static part of the balance of linear momenta in PINN Form $k$, with $n_q^{(k)} = 0$ in Eq.~\eqref{eq:generic-PDE-aux-conditions}, and 
$\busta = (\usta, \vsta)$ the exact \emph{static} solutions for $\bub = (\ub, \vb)$.  

\noindent
For convenience, the following static-operator array is defined based on the dynamic-operator array defined in Eq.~\eqref{eq:array-operators-dynamics}:
\begin{align}
	\bspdep{(k)}
	= \left\{ \spdesp{i}{(k)} , \ i=1, \ldots , n_i^{(k)} \right\}
	\ , \quad
	k = 0, \ldots, 3
	\ .
\end{align}

\noindent
To prevent the minimization iterate during the minimization process from reaching the static solution $\busta$, the loss function $\loss$ can be augmented by a barrier (or penalty) function $\barf{\bub}$ that would drastically increase the loss, thus creating a barrier (penalty), when the minimization iterate comes close to these solutions.

\noindent
For a generic constrained nonlinear optimization problem of the form:
\begin{align}
	\min_x f(x) \text{ such that } 0 < x
	\ ,
\end{align}
the nonlinear objective function $f(x)$ can be augmented by a barrier (penalty) function $\pf(x)$:
\begin{align}
	&
	J^{(\pf)} = f(x) + \bwe \, \pf(x)
	\ , \text{ with }
	\\
	&
	\pf (x)  
	= \pfp{\text{(inv)}} (x) 
	= \pfp{\text{(i)}} (x)
	= \frac{1}{x - \bde} , \text{ or }
	\label{eq:inverse-barrier-function-1}
	\\
	&
	\pf (x)  
	= \pfp{\text{(log)}} (x)
	= \pfp{\text{(l)}} (x) 
	= \log\left(\frac{1}{x - \bde}\right) = - \log (x - \bde)
	\label{eq:log-barrier-function-1}
	\ .
\end{align}
where 
$\pfp{\text{(inv)}} = \pfp{\text{(i)}}$ is the \emph{inverse} barrier function,
$\pfp{\text{(log)}} = \pfp{\text{(l)}}$ the \emph{logarithmic} (log) barrier function, and $\bde$ the shift of the origin to create a buffer zone that prevents the minimization iterate from getting too close to the feasible-domain boundary, as shown in Figure~\ref{fig:inverse-log-barrier functions}.  
% As the minimization iterate gets closer to the feasible-domain boundary, the value of the inverse barrier function $1/x$ is much higher than that of the log barrier function $-\log (x)$, thus preventing the iterate to reach a minimizer that is very close to this boundary.  The shifted barrier functions $1 / (x - \bde)$ and $-\log (x - \bde)$, with application in nonlinear structural dynamics, provide a buffer zone that prevents the iterate from getting too close to the boundary.   
Shown in the right subfigure are the cases for $\bde = 0.1$ (for log barrier only) and $\bde = 0.5$.  While the buffer-zone depth of the inverse function gradually decreases toward the boundary, the buffer-zone depth for the shifted log function remains seemingly constant, with much slower decrease rate.

The shorter abbreviations defined as ``i = inv = inverse'' and ``l = log = logarithm'' are to be used in an expanded supercripts $(mn)$ that determine the character of the barrier function  $\pfp{\text{(mn)}}$, with $m \in \{\text{i, l}\}$ (already defined above) and $n \in \{\text{s, p, a}\}$ (to be defined further below in the context of structural dynamics).

%\noindent
%References to cite: \cite{Nocedal.2006} p.~417, \cite{Fiacco.1990} p.~7, \cite{frisch1955logarithmic}.

\noindent
Let $\busta$ be a static solution satisfying the static PDE, i.e., the left-hand side of Eq.~\eqref{eq:generic-PDE-aux-conditions}:
\begin{align}
	\spdes{i} (\busta (\X , \ti)) = 0 \ , \text{ for } i = 1,\ldots, n_i
	\ .
	\tag{\ref{eq:generic-static-pde}}
\end{align}
\begin{align}
	% J^{\mathcal{b}}
	J^{(\pf)} (\bub) = J(\bub) + \bwe \, \psf{\bub}
	\ .
\end{align} 

\noindent
Barrier functions based on a \emph{static solution}, with superscripts ``i = inverse,'' ``l = log,'' and ``s = statics'':
\begin{align}
	\pfp{\text{(is)}} (\bub) := \frac{1}{|\bub - \busta| - \bde } \ , \quad
	\pfp{\text{(ls)}} (\bub) := - \log \left(| \bub - \busta | - \bde\right)
	\ .
	\label{eq:barrier-static-solution}
\end{align}

\noindent
Barrier functions based on satisfying the \emph{static nonlinear operator}, with additional superscripts ``p = pde,'' and ``$k$ = Form $k$'' (other than the above):
\begin{align}
	\pfp{\text{(ip,$k$)}} (\bub) := \frac{1}{\parallel \bspdep{(k)} \parallel - \bde } \ , \quad
	\pfp{\text{(lp,$k$)}} (\bub) := - \log \left(\parallel \bspdep{(k)} \parallel - \bde\right)
	\ , \quad
	\parallel \bspdep{(k)} \parallel = 
	\left[ \sum_{i=1}^{i=n_i^{(k)}} \left( \spdesp{i}{(k)} \right)^2  \right]^{1/2}
	\ .
	\label{eq:barrier-static-operator}
\end{align}

\noindent
Barrier functions based on the system \emph{acceleration} (or time derivative of linear momenta), with an additional superscript ``a = acceleration'' (other than the above):
\begin{align}
	\pfp{\text{(ia)}} (\bub) := \frac{1}{\parallel \ndt{\bpb} \parallel - \bde } \ , \quad
	\pfp{\text{(la)}} (\bub) := - \log \left(\parallel \ndt{\bpb} \parallel - \bde\right) 
	\ , \quad
	\parallel \bpb \parallel = 
	\left[ (\pbsd{1})^2 + (\pbsd{2})^2 \right]^{1/2}
	\ .
	\label{eq:barrier-acceleration}
\end{align}

\begin{rem}
	\label{rm:barrier-acceleration}
	Acceleration barrier function.
	{\rm
		For PDEs without an analytical static solution, 
		the static-solution barrier function in Eq.~\eqref{eq:barrier-static-solution} is clearly not convenient to use. 
		For PDEs with complex, nonlinear static operators, the barrier function in Eq.~\eqref{eq:barrier-static-operator} aiming at satisfying such a static operator is computationally inefficient. 
		On the contrary, since the system acceleration is readily available, the acceleration barrier function in Eq.~\eqref{eq:barrier-acceleration} is easily implemented, and computationally efficient.
%		{\color{red} [NOTE: Refer to figures using the acceleration barrier function.  ENDNOTE]}
	}
	\hfill$\blacksquare$
\end{rem}

\begin{figure}[tph]
	\centering
	\includegraphics[width=0.340\textwidth]{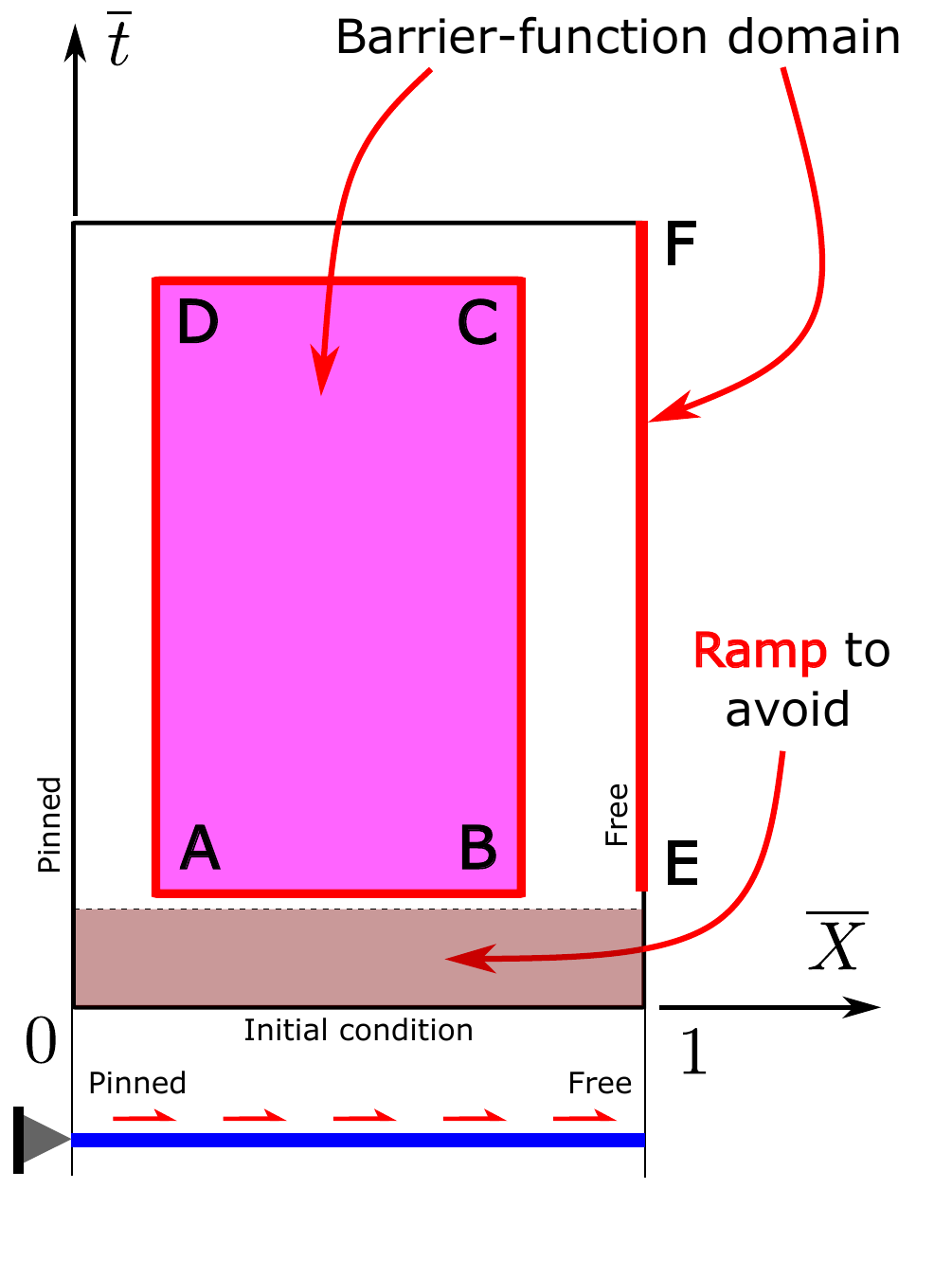}
	\includegraphics[width=0.410\textwidth]{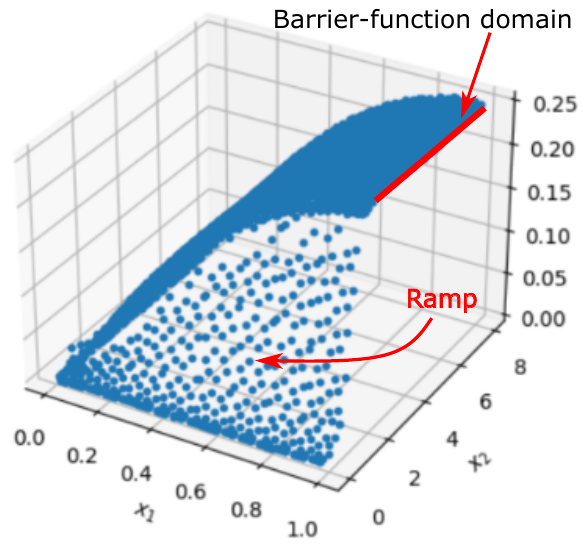}
	\caption{
		\emph{Barrier-function domain}.  Section~\ref{sc:filter-static-solution}.
		The barrier function domain (left), a subset of the computational domain bounded by $\Xb \in \left[ 0 , 1 \right]$ and $\tb \in \left[0, \tb_{\text{max}} \right]$ could be the rectangle ABCD, or just the segment EF, with both avoiding the ramp zone that occurs for certain distribution of nodes in the computational domain and the network architecture, such as a random distribution of nodes (in a static solution) with a neural network having 4 hidden layers, and 32 neurons per layer for a pinned-free bar (right).
%		{\color{red} [NOTE: 2023.06.28, Need to report the number of nodes and the number of parameters. Find input file. ENDNOTE]}
	}
	\label{fig:barrier-1}
\end{figure}

\begin{figure}[tph]
	\centering
	\includegraphics[width=0.45\textwidth]{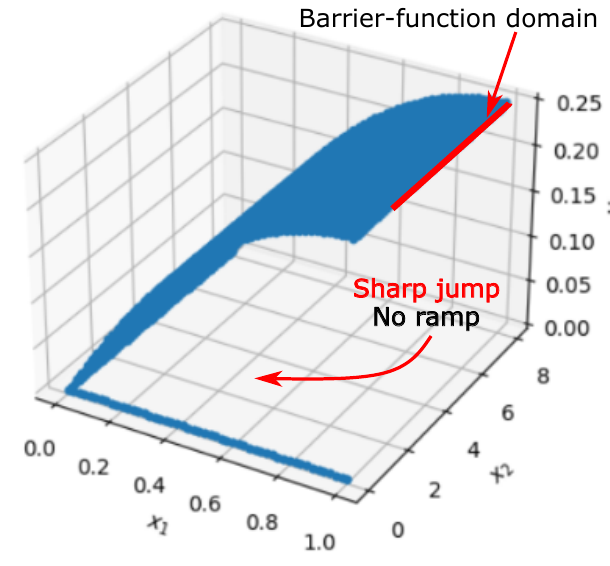}
	\includegraphics[width=0.45\textwidth]{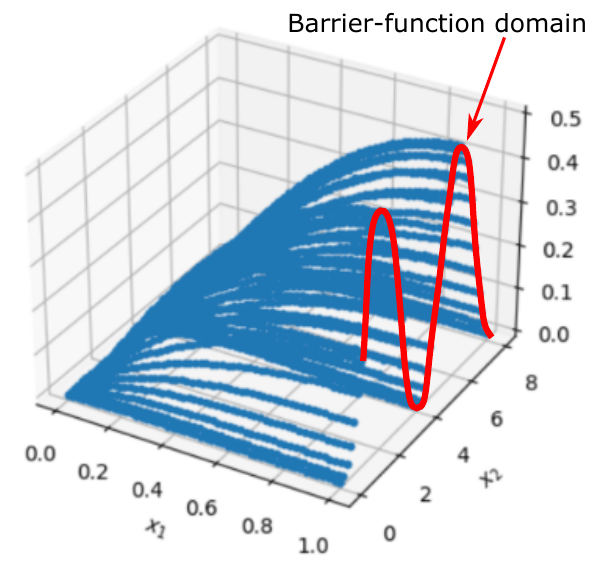}
	\caption{
		\emph{Expected effects of barrier function}.  Section~\ref{sc:filter-static-solution}.
		Pinned-free bar, axial motion.
		$\bullet$
		\emph{Left:} NO barrier,
		\emph{static} solution with 4 hidden layers, 128 neurons per payer.  
		$\bullet$
		\emph{Right:} With barrier,
		\emph{dynamic} solution with two periods, 2 hidden layers, 32 neurons per layer.
%		$\bullet$
%		{\color{red} [NOTE: 2023.08.14, Need to add and to refer to figures that actually use a barrier function. ENDNOTE]}
%		$\bullet$
%		{\color{red} [NOTE: 2023.06.28, Need to report the number of nodes and the number of parameters. Find script files, which were not systematically perserved for these early results.  RECREATE these script files.   ENDNOTE]}
	}
	\label{fig:barrier-2}
\end{figure}

\begin{rem}
	History of logarithmic and inverse barrier functions.
	{\rm
		In structural dynamics, as soon as we observed through our numerical experiments that the PINN minimization process could lead to a static solution, we immediately thought about additing an inverse penalty function to the PINN loss function to prevent the  minimization iterate from reaching the static solution.
		By taking the logarigthm of the inverse barrier function, we obtain the log barrier function.  This was done before we checked the optimization literature, such as \cite{Nocedal.2006} p.~417, from where we learned of \cite{Fiacco.1990} and \cite{frisch1955logarithmic}. 
		It was surprising that, 
		according to the classic book by Fiacco and McCormick \cite{Fiacco.1990} p.~7, the years of appearance of these barrier functions were reversed, with the logarithmic barrier function  introduced by Frisch in 1954 \cite{frisch1954principles} and 1955 \cite{frisch1955logarithmic}, and the inverse barrier function  introduced by Carroll a few years later in 1959 \cite{Carroll.1959.dissertation} and 1961 \cite{Carroll.1961}. 
		%		References to cite:
		%		\cite{Nocedal.2006} p.~417, \cite{Fiacco.1990} p.~7, \cite{frisch1955logarithmic}.
	}
	\phantom{blank}\hfill$\blacksquare$
\end{rem}

\begin{figure}[tph]
	\begin{subfigure}{0.24\textwidth}
		\includegraphics[width=\textwidth]{Figures/23.9.8_R2b_A-P_bcrF_F1_J,2a,3__lr0.001_Nsteps25000_random_N51_W32_H2_Glorot_init_-_shape25000.png}
		\caption{
			% Pinned-free bar, static solution.
			% {\scriptsize 1,185 params, Shape 25000}
			{\tiny (\ref{fig:23.9.8 R2b shape25000}): 1,185 params, Shape 25000}
		}
		% \tag{\ref{fig:23.9.8 R2b shape25000}}
		\label{fig:23.9.8 R2b shape25000-2}
	\end{subfigure}
	\begin{subfigure}{0.24\textwidth}
		\includegraphics[width=\textwidth]{Figures/23.9.8_R2b_A-P_bcrF_F1_J,2a,3__lr0.001_Nsteps25000_random_N51_W32_H2_Glorot_init_-_free-end_disp25000.png}
		\caption{
			% Pinned-free bar, static solution.
			% {\scriptsize 1,185 params, Free-end disp}
			{\tiny (\ref{fig:23.9.8 R2b free-end disp25000}): 1,185 params, Free-end disp}
		}
		% \tag{\ref{fig:23.9.8 R2b free-end disp25000}}
		\label{fig:23.9.8 R2b free-end disp25000-2}
	\end{subfigure}
	\begin{subfigure}{0.24\textwidth}
		\includegraphics[width=1.\textwidth]{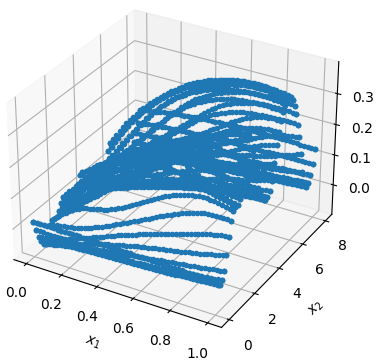}
		\caption{
			{\scriptsize \red{Barrier effects}, Shape}
		}
		\label{fig:23.10.5 R1c shape25000}
	\end{subfigure}
		\begin{subfigure}{0.24\textwidth}
		\includegraphics[width=1.\textwidth]{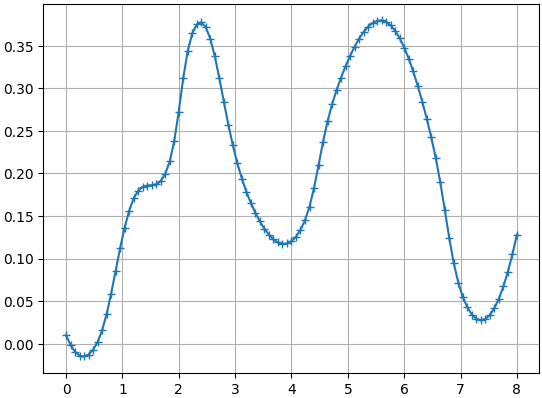}
		\caption{
			{\scriptsize Free-end disp, \red{Very good}}
		}
		\label{fig:23.10.5 R1c free-end disp25000}
	\end{subfigure}
	\\
	\begin{subfigure}{0.32\textwidth}
		\includegraphics[width=\textwidth]{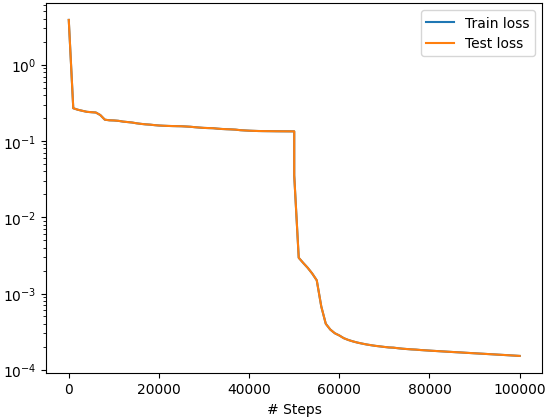}
		\caption{
			{\scriptsize \red{Barrier off} at 50000}
		}
		\label{fig:23.10.5 R1c.2 loss100000}
	\end{subfigure}
	\begin{subfigure}{0.32\textwidth}
		\includegraphics[width=1.\textwidth]{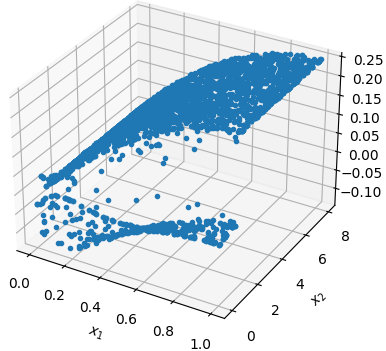}
		\caption{
			% {\scriptsize No barrier, Shape 100000}
			{\scriptsize No barrier, \red{Static} 100000}
		}
		\label{fig:23.10.5 R1c.2 shape100000}
	\end{subfigure}
	\begin{subfigure}{0.32\textwidth}
		\includegraphics[width=1.\textwidth]{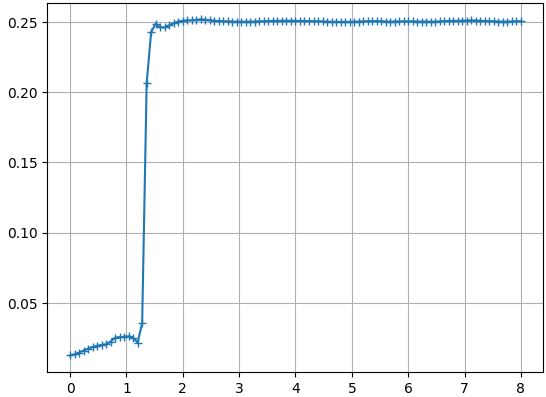}
		\caption{
			{\scriptsize No barrier, Free-end disp}
		}
		\label{fig:23.10.5 R1c.2 free-end disp100000}
	\end{subfigure}
	\caption{
		\DDET.
		\emph{Pinned-free bar, \red{static} solutions, \red{barrier} functions.}
		Appendix~\ref{sc:filter-static-solution}.
		SubFigs~\ref{fig:23.9.8 R2b shape25000-2}-\ref{fig:23.9.8 R2b free-end disp25000-2} =
		SubFigs~\ref{fig:23.9.8 R2b shape25000}-\ref{fig:23.9.8 R2b free-end disp25000}, \red{No barrier}, \red{static} solution, recalled for convenience.
		$\bullet$
		SubFigs~\ref{fig:23.10.5 R1c shape25000}-\ref{fig:23.10.5 R1c free-end disp25000}: Step 25,000, Shape, Free-end disp. 
		\red{$\star$}
		\red{Barrier} depth $\bde = 1$, weight $\bwe = 1$. 
		\emph{Network:}
		Remarks~\ref{rm:parameter-names}, \ref{rm:data-point-grids},
		% n\_inp=2,
		T=8,
		W=32,
		H=2, 
		n\_out=1, 
		Glorot-uniform initializer, 
		1,185 parameters,
		\red{$\star$} 
		\emph{regular} grid.
		\emph{Training:}
		Remark~\ref{rm:learning-rate-schedule-1} (LRS~1),
		init\_lr=0.001.
		Damping\%=-0.6\%, \red{Very good}.
		Local maxima = \{0.378, 0.379\}.
		Time at local max = \{2.4, 5.6\}.
		$\triangleright$
		Step 50,000,
		damping\%=\mbox{-0.3}\%, \red{Quasi-perfect} free-end disp.
		{\scriptsize (23105R1c)}
		$\bullet$
		SubFigs~\ref{fig:23.10.5 R1c.2 loss100000}-\ref{fig:23.10.5 R1c.2 free-end disp100000}:
		\red{Barrier turned off} at Step 50,000.
		% Step 100,000,
		\red{Static solution came back}, even with random grid.
		{\scriptsize (23105R1c.2)}
		$\bullet$
		{\footnotesize
			Figure~\ref{fig:DDE-T-Form-1-pinned-free-bar-static-solutions}, \DDET\ static solutions.
			Figure~\ref{fig:JAX-Form-1-pinned-free-bar-NO-static-solutions}, using our \JAX\ script, no static solution.
			$\triangleright$
			Figure~\ref{fig:23.9.5 R3b A-P bcrF F(1J)2a(3) lr0.005 cy1-5-NCA Nsteps200000 random N51 W64 H4 He init-1}, \red{Form 2a}, GPU time \red{605 sec} for 200,000 steps.  
			% Rerun 23.9.5 R3b in [1] F2a Axial motion pinned-free.
			Figure~\ref{fig:23.9.5 R3d A-P bcrF F(1J,2a)3 lr0.005 cy1-5-NCA Nsteps200000 random N51 W64 H4 He init-1}, \red{Form~3}, GPU time \red{537 sec} for 200,000 steps.
			$\triangleright$
			Figure~\ref{fig:axial-Mathematica-solutions}, reference solution to compare.
			$\triangleright$
			Figure~\ref{fig:barrier-function 23.10.6 R1a.2}, barrier-weight effects.
		}
	}
	\label{fig:barrier-function 23.10.5 R1c}
\end{figure}

As an example of how the barrier function works, consider the static solution in SubFigs~\ref{fig:23.9.8 R2b shape25000}-\ref{fig:23.9.8 R2b free-end disp25000}, obtained with \DDET\ Form 1, using a network with N51 W32 H2 and 1,185 parameters.
In Figure~\ref{fig:barrier-function 23.10.5 R1c}, SubFigs~\ref{fig:23.10.5 R1c shape25000}-\ref{fig:23.10.5 R1c free-end disp25000}, an acceleration barrier (Appendix~\ref{app:barrier-functions}, Remark~\ref{rm:barrier-acceleration}) was used, and effectively helped prevent the static solution, while achieving a ``Very good'' free-end displacement at Step 25,000, with a damping percentage of -0.6\%, which measures the ratio between two successvive local maxima (peaks) of free-end displacement at \{0.378, 0.379\}, respectively, below the more accurate peak amplitude of about 0.5 shown in Figure~\ref{fig:JAX-Form-1-pinned-free-bar-NO-static-solutions}. 
At Step 50,000, this ratio yielded a damping percentage of -0.3\%, and was thus even classified as quasi-perfect.
Even though the shape of the response was far from being recognized as good, the vibration period ressembled the accurate vibration period in Figure~\ref{fig:JAX-Form-1-pinned-free-bar-NO-static-solutions}  (or in Figure~\ref{fig:axial-Mathematica-solutions}), with the times for the local maxima being \{2.4, 5.6\}, instead of \{2., 6.\}. 

\begin{figure}[tph]
	\begin{subfigure}{0.32\textwidth}
		%_\includegraphics[width=1.\textwidth]{Figures/DEBUG_23.10.5_R1b.2_A-PF_F1H_2ab,3__barrier_d1_w0.1_lr0.001_Nsteps50000_random_N51_W32_H2_Glorot_init_-_free-end_disp50000.png}
		\includegraphics[width=1.\textwidth]{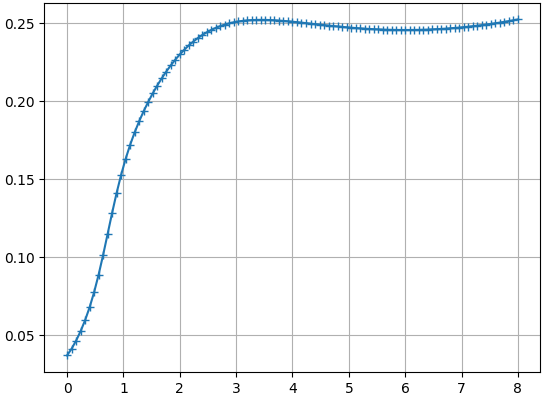}
		\caption{
			{\scriptsize $\bwe=0.1$, too small.}
		}
		\label{fig:23.10.5 R1b.2 free-end disp25000}
	\end{subfigure}
	\begin{subfigure}{0.32\textwidth}
		%_\includegraphics[width=1.\textwidth]{Figures/DEBUG_23.10.6_R1a.2_A-PF_F1H_2ab,3__barrier_d1_w0.5_lr0.001_Nsteps50000_random_N51_W32_H2_Glorot_init_-_free-end_disp50000.png}
		\includegraphics[width=1.\textwidth]{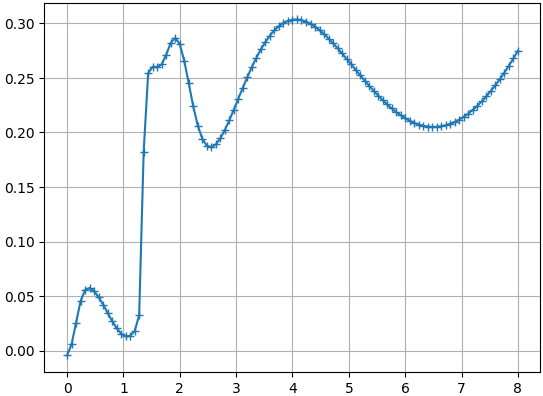}
		\caption{
			{\scriptsize $\bwe=0.5$, better.}
		}
		\label{fig:23.10.6 R1a.2 free-end disp25000}
	\end{subfigure}
	\begin{subfigure}{0.32\textwidth}
		%_\includegraphics[width=1.\textwidth]{Figures/DEBUG_23.10.5_R1d.2_A-PF_F1H_2ab,3__barrier_d1_w0.7_lr0.001_Nsteps50000_random_N51_W32_H2_Glorot_init_-_free-end_disp50000.png}
		\includegraphics[width=1.\textwidth]{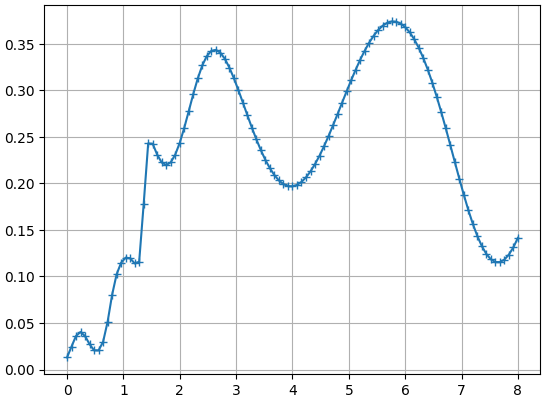}
		\caption{
			{\scriptsize $\bwe=0.7$, better.}
		}
		\label{fig:23.10.5 R1d.2 free-end disp25000}
	\end{subfigure}
	\caption{
		\DDET.
		\emph{Pinned-free bar, \red{barrier weight} $\bwe$ effects.}
		Appendix~\ref{sc:filter-static-solution}.
		All parameters were the same as in Figure~\ref{fig:barrier-function 23.10.5 R1c}, except for the barrier weight $\bwe$. 
		$\bullet$
		SubFig~\ref{fig:23.10.5 R1b.2 free-end disp25000}:
		% (\ref{fig:23.10.5 R1b.2 free-end disp25000}):
		$\bwe = 0.1$, smooth, flat pateau began to curve.
		$\bullet$
		SubFig~\ref{fig:23.10.6 R1a.2 free-end disp25000}:
		% (\ref{fig:23.10.6 R1a.2 free-end disp25000}):
		$\bwe = 0.5$, waves, bigger jump after $t=1$.
		$\bullet$
		SubFig~\ref{fig:23.10.5 R1d.2 free-end disp25000}:
		% (\ref{fig:23.10.5 R1d.2 free-end disp25000}): 
		$\bwe = 0.7$, waves, smaller jump after $t=1$, with two peaks clearly formed.
		$\bullet$
		SubFig~\ref{fig:23.10.5 R1c free-end disp25000}, $\bwe = 1.0$, smooth, no jump, with two distinct peaks.
	}
	\label{fig:barrier-function 23.10.6 R1a.2}
\end{figure}

As for the selection of the barrier depth and weight, see Figure~\ref{fig:barrier-function 23.10.6 R1a.2} for the influence of the barrier weight $\bwe$ on the solution, keeping the barrier depth $\bde=1$.  It is also possible to vary the barrier depth to illustrate its influence on the solution.

Upon removal of the barrier after Step 50,000, the static solution sprung back quickly, as shown by a precipitous drop in the loss function, the shape, and the characteristic free-end displacement with a flat plateau in SubFigs~\ref{fig:23.10.5 R1c.2 loss100000}-\ref{fig:23.10.5 R1c.2 free-end disp100000}, demonstrating that the barrier function was the key player that effectively held back the static solution.  

	% training error, generalization
	\section{Error relative to exact solution, generalization}
\label{app:exact-solution-generalization}
\noindent
We provide here (1) a comparison of PINN results to the exact solution for the \emph{axial} motion $u(x,t)$ of an elastic bar Eq.~\eqref{eq:eom-euler-bernoulli-axial}, under a constant, uniform distributed load $\dfbs{\X}$, with \emph{pinned-pinned} boundary conditions Eq~\eqref{eq:axial-pinned-pinned-BCs} and zero initial conditions Eq.~\eqref{eq:axial-ICs}, of the form
\begin{align}
	&
	u(x,t)  
	=
	u_{stat} (x)
	\left[1 - \cos(\pi t) \right]
	\ ,
	\\
	&
	u_{stat} (x) = a x^2 + b x + c
	\ , \quad
	a = \frac{\dfbs{\X}}{2 \sldn}
	\ , \quad
	b = - a
	\ , \quad
	c = 0
	\ ,
\end{align}
with $u_{stat} (x)$ being the static solution,
(2) revealing the essentially-zero value in the transverse displacement $v(x,t)$ and the pinned boundary conditions, and (3) a generalization of the motion in time beyond the space-time domain used for training.

Compared to the exact solution, the error of the displacement time history is measured by the square-root of the mean squared error (MSE), abbreviated by SMSE for short, during the time interval of interest.  The relative error (RE) is obtained by dividing the SMSE by the peak-to-peak amplitude (0.25) of the exact solution (0.125) at the center of the pinned-pinned elastic bar. 

\begin{figure}[tph]
	\begin{subfigure}{0.32\textwidth}
		\includegraphics[width=\textwidth]{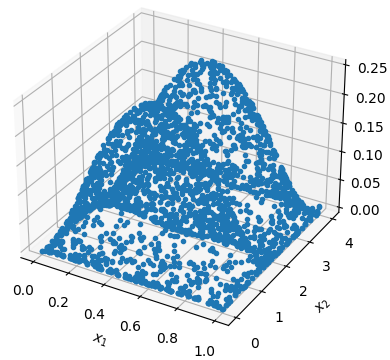}
		\caption{
			\scriptsize
			Shape $u$ time history.
		}
		\label{fig:24.1.19 R1a.1b.3 shape u 200000}
	\end{subfigure}
	\begin{subfigure}{0.32\textwidth}
		\includegraphics[width=\textwidth]{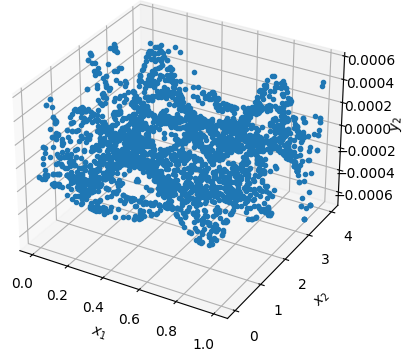}
		\caption{
			\scriptsize
			Shape $v$ time history.
		}
		\label{fig:24.1.19 R1a.1b.3 shape v 200000}
	\end{subfigure}
	\begin{subfigure}{0.32\textwidth}
		\includegraphics[width=\textwidth]{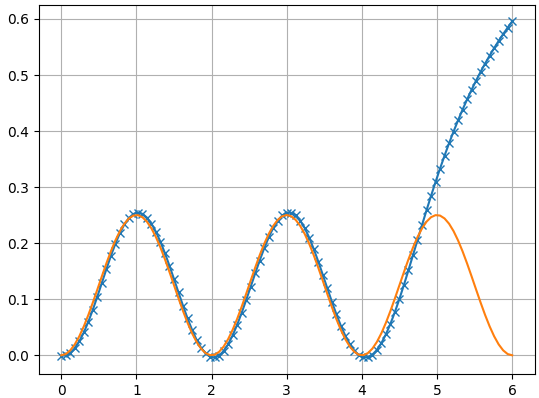}
		\caption{
			\scriptsize
			Center disp $u$, generalized.
		}
		\label{fig:24.1.19 R1a.1b.3 center u 200000}
	\end{subfigure}
	\\
	\begin{subfigure}{0.32\textwidth}
		\includegraphics[width=\textwidth]{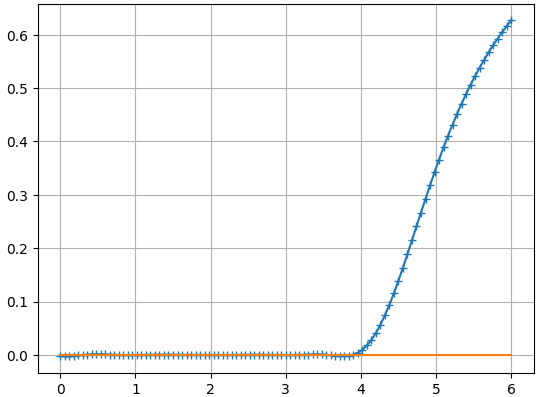}
		\caption{
			\scriptsize
			Right-pinned-end $u$.
		}
		\label{fig:24.1.19 R1a.1b.3 right-pinned-end u 200000}
	\end{subfigure}
	\begin{subfigure}{0.32\textwidth}
		\includegraphics[width=\textwidth]{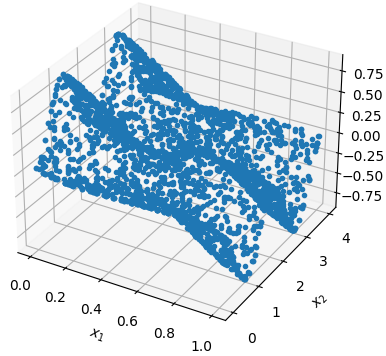}
		\caption{
			\scriptsize
			Space-slope time history.
		}
		\label{fig:24.1.19 R1a.1b.3 slope u_x shape 200000}
	\end{subfigure}
	\begin{subfigure}{0.32\textwidth}
		\includegraphics[width=\textwidth]{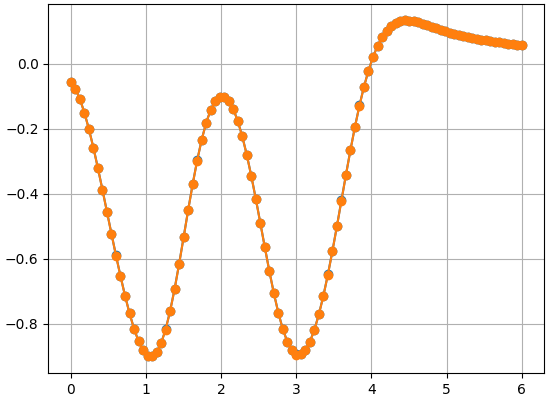}
		\caption{
			\scriptsize
			Right-pinned-end slope.
		}
		\label{fig:24.1.19 R1a.1b.3 right-pinned-end u_x 200000}
	\end{subfigure}
	\caption{
		\DDET.
		\emph{Pinned-pinned bar.}
		Appendix~\ref{app:exact-solution-generalization}.
		\red{Form 3.2:}
		Step 200,000.  
		\emph{Network:}
		Remark~\ref{rm:parameter-names},~\ref{rm:data-point-grids}, 
		T=4,
		W=64,
		H=4, 
		n\_out=10, 
		{Glorot} initializer, 
		18,554 parameters,
		\emph{random} grid.
		\emph{Training:}
		\red{$\star$}
		Remark~\ref{rm:learning-rate-schedule-3} (LRS~3, NCA),
		init\_lr=0.01.
		SubFig.~\ref{fig:24.1.19 R1a.1b.3 shape u 200000}, shape time history (t.h.).
		(\ref{fig:24.1.19 R1a.1b.3 shape v 200000}) Transverse displacement t.h., essentially-zero value of order 1e-3.
		(\ref{fig:24.1.19 R1a.1b.3 center u 200000}) Axial displacement at bar center (blue), exact solution (orange), training SMSE = 8e-3, RE=3\%.
		(\ref{fig:24.1.19 R1a.1b.3 right-pinned-end u 200000}) Right-pinned-end disp $u$ (blue) cannot generalize (should be zero, orange).
		(\ref{fig:24.1.19 R1a.1b.3 slope u_x shape 200000}) Bar space slope t.h.
		(\ref{fig:24.1.19 R1a.1b.3 right-pinned-end u_x 200000}) Right-pinned end space slope t.h.
		{\scriptsize (24.1.19 R1a.1b.3)}
	}
	\label{fig:24.1.19 R1a.1b.3}
\end{figure}

Figure~\ref{fig:24.1.19 R1a.1b.3} demonstrates the limitation of the generalization---i.e., predicting the solution beyond the space-time domain used to train the network---that PINN could provide.  It could however predict the displacement curving up to follow closely the upswing from the end of the training time domain ($\tb = T = 4$) at which the displacement and the velocity were zero as shown in SubFigure~\ref{fig:24.1.19 R1a.1b.3 center u 200000} until close to half a period beyond $\tb = T = 4$, i.e., the interval $[ 4, 5 ]$, with the MSE = 5.47e-4 compared to the training MSE = 0.66e-4.

\begin{rem}
	PINN vs traditional FEM for generalization.
	{\rm
		Extending the time domain to obtain the numerical solution in traditional FEM using the implicit Newmark time-integration method requires additional discrete-system solution steps, which could significantly increase the computational cost.  With PINN, there is no need to do optimization again, as the generalization to extend the solution beyond the training time domain $[0, T]$ is by performing the low-cost prediction (function evaluation) using exactly the same neural-network minimizers (weights and biases) obtained during the training optimization process on the previous (unextended) space-time domain.  No further computational-intensive optimization is required for PINN.
	}
	\hfill
	$\blacksquare$
\end{rem}

The theoretically-zero \emph{transverse} displacement $v(x,t)$ in SubFigure~\ref{fig:24.1.19 R1a.1b.3 shape v 200000} shows the \emph{essentially-zero} value of the order of 1e-3, typical of PINN results.  The zero displacement in the pinned boundary condition is not exactly satisfied, but \emph{essentially zero} (i.e., of the order 1e-3).  From this standpoint, the PINN displacement time history in the training interval $[0, 4]$, with its SMSE = 8e-3 and RE = 3\%, agreed well the exact solution $u(x,t)$, with \emph{essentially-zero} difference.  

\begin{rem}
	Uncertainty in real boundary and initial condition.
	{\rm
		Traditional FEM can satisfy zero boundary and initial (BC-IC) conditions.  In 
		practice, exactly-zero BC-IC conditions are, however, an idealization, and thus \emph{essentially-zero} BC-IC conditions could represent uncertainty in the real BC-IC conditions.
	}
	\hfill $\blacksquare$
\end{rem}

The zero right-end-pinned boundary condition did not generalize, i.e., remains at zero in the time interval $[4, 6]$, as shown in SubFigure~\ref{fig:24.1.19 R1a.1b.3 right-pinned-end u 200000}.  
SubFigure~\ref{fig:24.1.19 R1a.1b.3 slope u_x shape 200000} shows the time history ($\tb \in [0, 4]$) of the space slope along the bar length $\xb \in [0, 1]$, whereas SubFigure~\ref{fig:24.1.19 R1a.1b.3 right-pinned-end u_x 200000} shows the right-pinned-end space-slope time history (at $\xb = 1$), which is negative in the training time interval $[0, 4]$, indicating compressive state (as it should be), and is positive during generalization, corresponding to the rising curve in SubFigure~\ref{fig:24.1.19 R1a.1b.3 right-pinned-end u 200000} in the generalization time interval $[4, 6]$. 
In general, PINN cannot generalize to have periodicity, which is longer than the half period $[4, 5]$ beyond $T = 4$.

	% code snippets
	\section{Python-script snippets}
\label{app:code-snippets}
\label{app:python-script-snippets}
\noindent
Since the speed of execution depends crucially on how the derivatives (gradients) were computed, and since there are several ways to implement to computation of the derivatives, we provide below snippets of our Python scripts (for \DDET\ and for \JAX) that were executed on Colab so readers (and \DDET\ developers) are better informed of how we obtained our results for potential future improvements, since the computation of these derivatives is likely the key for the pathological problems encountered in both the \DDET\ script (shift, amplification, static solution) and the \JAX\ script (slowness in Form 3) described in previous sections.

\begin{lstlisting}[
	language=Python, 
	caption={
		\DDET\ and \JAX. Import packages and modules common to both \DDET\ and \JAX.
	}, 
	label={lst:DDE-T import module}
]	
	import numpy as np
	import matplotlib.pyplot as plt
	from scipy.signal import argrelextrema
	import pickle
	import timeit
	from timeit import default_timer
\end{lstlisting}

% See sample code in Listing~\ref{lst:listing example}.
% 
% example from Overleaf
% https://www.overleaf.com/learn/latex/Code_listing#Code_styles_and_colours
% for listing label, see the documentation of The Listings Package
% https://texdoc.org/serve/listings.pdf/0
%\begin{lstlisting}[language=Python, caption={Python example from Overleaf}, label={lst:listing example}]
%	import numpy as np
%	
%	def incmatrix(genl1,genl2):
%	m = len(genl1)
%	n = len(genl2)
%	M = None #to become the incidence matrix
%	VT = np.zeros((n*m,1), int)  #dummy variable
%	
%	#compute the bitwise xor matrix
%	M1 = bitxormatrix(genl1)
%	M2 = np.triu(bitxormatrix(genl2),1) 
%	
%	for i in range(m-1):
%	for j in range(i+1, m):
%	[r,c] = np.where(M2 == M1[i,j])
%	for k in range(len(r)):
%	VT[(i)*n + r[k]] = 1;
%	VT[(i)*n + c[k]] = 1;
%	VT[(j)*n + r[k]] = 1;
%	VT[(j)*n + c[k]] = 1;
%	
%	if M is None:
%	M = np.copy(VT)
%	else:
%	M = np.concatenate((M, VT), 1)
%	
%	VT = np.zeros((n*m,1), int)
%	
%	return M
%\end{lstlisting}

% \subsection{\DDET\ code snippets}
\subsection{DDE-T script snippets}
\label{sc:DDE-T script snippets}
\noindent
For the axial motion of an elastic bar, the acceleration can be computed in two ways: (a) Using the Hessian operator (Listing~\ref{lst:DDE-T time derivative by hessian}), or (b) Using two successive Jacobian operators (Listing~\ref{lst:DDE-T time derivative by two jacobians}). 

\begin{lstlisting}[
	language=Python, 
	caption={\DDET. Imported packages and modules.}, 
	label={lst:DDE-T import module}
]
	!pip install deepxde
	import os
	os.environ["DDE_BACKEND"] = "tensorflow"
	import deepxde as dde
	from deepxde.geometry import *
\end{lstlisting}

\begin{lstlisting}[
	language=Python, 
	escapechar=|,
	caption={
		\DDET. \emph{Axial motion of elastic bar. \red{Forms 1, 2b.}  Second time derivative} by Hessian.  
		$\bullet$ 
		See Listing~\ref{lst:DDE-T time derivative by two jacobians} for the use of two Jacobians.
		$\blacktriangleright$
		See \JAX\ \emph{\red{Form 1}} in Listing~\ref{lst:JAX Form 1}.
	}, 
	label={lst:DDE-T time derivative by hessian}
]
	# time derivative
	if (Form == "F1"):
		du_tt = dde.grad.hessian(y, x, i=1, j=1) |\label{ln:Form 1 acceleration du-tt}|
	elif (Form == "F2b"):
		du_tt = dde.grad.hessian(y, x, i=1, j=1, component=0)
\end{lstlisting}

\begin{lstlisting}[
	language=Python, 
	escapechar=|,
	caption={
		\DDET. \emph{Axial motion of elastic bar. \red{Forms 1, 2b.} First and second time derivatives} (velocity and acceleration) by a Jacobian.  Remark~\ref{rm:efficiency-DDE-T-vs-JAX}, increased GPU time when using two Jacobians compared to a Hessian (Listing~\ref{lst:DDE-T time derivative by hessian}).
		$\blacktriangleright$
		See \JAX\ \emph{\red{Form 1}} in Listing~\ref{lst:JAX Form 1}.
	}, 
	label={lst:DDE-T time derivative by two jacobians}
]	
	if (Form == "F1"):
		du_t  = dde.grad.jacobian(y, x, i=0, j=1)
		du_tt = dde.grad.jacobian(du_t, x, i=0, j=1) |\label{ln:Form 1 acceleration du-tt}|
	elif (Form == "F2b"):
		du_t  = dde.grad.jacobian(y, x, i=0, j=1, component=0)
		du_tt = dde.grad.jacobian(du_t, x, i=0, j=1, component=0) |\label{ln:Form 2b acceleration du-tt}|
\end{lstlisting}

\begin{lstlisting}[
	language=Python, 
	escapechar=|,
	caption={
		\DDET. \emph{Axial motion of elastic bar. \red{Forms 2a, 3.}}  Time derivatives.
	}, 
	label={lst:DDE-T F2a F3 time derivatives}
]
	elif ((Form == "F2a") or (Form == "F3")):
		# split 2nd time derivative
		# momentum
		py    = y[:, 1:2] |\label{ln:Form 2a or Form 3 momentum}|
		# velocity as time derivative of displacement
		dy_t  = dde.grad.jacobian(y, x, i=0, j=1) |\label{ln:Form 2a or Form 3 velocity}|
		# acceleration as time derivative of momentum
		py_t  = dde.grad.jacobian(py, x, i=0, j=1) |\label{ln:Form 2a acceleration py-t}|
\end{lstlisting}

\begin{lstlisting}[
	language=Python, 
	escapechar=|,
	caption={
		\DDET. \emph{Axial motion of elastic bar. \red{Forms 1, 2a.}  Second space-derivative of displacement} by a Hessian or by two Jacobians.  For 2nd time derivative of displacement (acceleration), see Listing~\ref{lst:DDE-T time derivative by hessian} for the use of a Hessian, and Listing~\ref{lst:DDE-T time derivative by two jacobians} for the use of two Jacobians.  For Form 2a, Line~\ref{ln:2nd space derivative du-xx component=0},\protect\footnotemark 
		component=0 (col index 0 or 1st col) must be specified since the y array has two columns.
	}, 
	label={lst:DDE-T 1st and 2nd space derivatives}
]
	# space derivative
	if ((Form == "F1") or (Form == "F2a")):
		if (hessian_jacobian == "hessian"):
			if (Form == "F1"):
				# method 1, use hessian function
				du_xx = dde.grad.hessian(y, x, i=0, j=0, component=None) |\label{ln:Form 1 2nd deriv du-xx}|
			if (Form == "F2a"):
				# method 1, use hessian function
				du_xx = dde.grad.hessian(y, x, i=0, j=0, component=0) |\label{ln:2nd space derivative du-xx component=0}|
		elif (hessian_jacobian == "jacobian"):
			# method 2, use 2 jacobian functions
			y0 = y[:, 0:1]
			# 1st space derivative of disp y0
			# method 1
			# du_x    = dde.grad.jacobian(y  , x, i=0, j=0)
			# method 2
			du_x      = dde.grad.jacobian(y0, x, i=0, j=0) |\label{ln:compute slope du-x}|
			du_xx     = dde.grad.jacobian(du_x, x, i=0, j=0) |\label{ln:Form 2a 2nd deriv by Jacobian du-xx}|
\end{lstlisting}
\footnotetext{A line number not preceded by a listing number indicates that the said line is in the current listing.  When referring to a code line in a different listing, see, e.g., the caption of Listing~\ref{lst:DDE-T, Form 2b, space derivative splitting} and the abbreviation ``L-Line~\ref{lst:DDE-T, all forms, losses}.\ref{ln:Form 2b losses, with barrier}'' for ``Listing~\ref{lst:DDE-T, all forms, losses} Line~\ref{ln:Form 2b losses, with barrier}.''}

\begin{lstlisting}[
	language=Python, 
	escapechar=|,
	caption={
		\DDET. \emph{Axial motion of elastic bar.  \red{Form 2b}.  Second space-derivative splitting.} 
		Variable (Var) y0 = \unexpanded{y[:, 0:1]} is the disp stored in the 1st col of the y array (col index 0).
		Var du\_x = space derivative of disp y0 (slope) computed using the Jacobian of y0 (col index i=0), wrt space coordinate $x$ stored in 1st col of x array (index j=0).
		Var y1 = \unexpanded{y[:, 1:2]} is the slope (space derivative of the displacement) stored in the 2nd col of the y array (col index 1).  
		Var y1\_x is the 1st space derivative of y1, and thus 
		the 2nd space derivative of the displacement, computed using the Jacobian of y1 = \unexpanded{y[:, 1:2]}, i.e., the 2nd col of y array with col index i=1, with respect to the space coordinate $\x$ stored in \unexpanded{x[:, 0]}, i.e., the 1st col of the space-time x array (col index j=0).
		The constraint (or loss) on the slope is (y1 - du\_x), which is the negative of Eq.~\eqref{eq:wave-eq-form-2b-initial-conditions-2}$_2$ and coded as shown in Listing~\ref{lst:DDE-T, all forms, losses} Line~\ref{ln:Form 2b losses, with barrier} (abbreviated as L-Line~\ref{lst:DDE-T, all forms, losses}.\ref{ln:Form 2b losses, with barrier}) and L-Line~\ref{lst:DDE-T, all forms, losses}.\ref{ln:Form 2b losses, with NO barrier}, with du\_x computed as coded in 
		% Listing~\ref{lst:DDE-T 1st and 2nd space derivatives} Line~\ref{ln:compute slope du-x}
		L-Line~\ref{lst:DDE-T 1st and 2nd space derivatives}.\ref{ln:compute slope du-x}.
	}, 
	label={lst:DDE-T, Form 2b, space derivative splitting}
]
	elif (Form == "F2b"):
		# split space derivative
		y0 = y[:, 0:1]
		# 1st space derivative of disp y0
		# method 1
		# du_x    = dde.grad.jacobian(y  , x, i=0, j=0)
		# method 2
		du_x      = dde.grad.jacobian(y0, x, i=0, j=0) |\label{ln:Form 2b slope du-x}|
		#
		# slope in 2nd col, index 1
		y1 = y[:, 1:2] |\label{ln:Form 2b y1 slope in 2nd col of y}|
		# constraint: slope y1 = du_x, 1st space derivative of disp y0
		#
		# space derivative of slope y1 = 2nd space derivative of disp y0
		y1_x = dde.grad.jacobian(y, x, i=1, j=0) |\label{ln:Form 2b 2nd space deriv y1_x = du_xx}|
\end{lstlisting}

\begin{lstlisting}[
	language=Python, 
	escapechar=|,
	caption={
		\DDET. \emph{Axial motion of elastic bar.  \red{Form 3}.  Split 2nd derivatives wrt both space and time.} 
		Var y0 = \unexpanded{y[:, 0:1]} (1st col of y array, col index 0) = displacement.
		Var du\_x = derivative of disp y0 (index i=0) wrt x (index j=0).
		Var y1 = \unexpanded{y[:, 2:3]} (3rd col of y array, col index 2) = slope (space derivative of disp).
		Var y1\_x = space derivative of the slope stored in 3rd col of y array (index i=2, i.e., 2nd derivative of disp), using the Jacobian wrt to $x$ coordinate stored in 1st col of x array (index j=0).
	}, 
	label={lst:DDE-T, Form 3, space derivative splitting}
]
	elif (Form == "F3"):
		# split space derivative
		y0 = y[:, 0:1]
		# 1st space derivative of disp y0
		# method 1
		# du_x    = dde.grad.jacobian(y  , x, i=0, j=0)
		# method 2
		du_x      = dde.grad.jacobian(y0, x, i=0, j=0)
		#
		# slope in 3rd col, index 2
		y1 = y[:, 2:3]
		# constraint: slope y1 = du_x, 1st space derivative of disp y0
		#
		# space derivative of slope y1 = 2nd space derivative of disp y0
		y1_x = dde.grad.jacobian(y, x, i=2, j=0)
\end{lstlisting}

\begin{lstlisting}[
	language=Python, 
	escapechar=|,
	caption={
		\DDET. \emph{Axial motion of elastic bar.  \red{Forms 1, 2a, 2b, 3}.}  
		Var f\_ext = distributed external axial force.
		In case of zero barrier weight $\bwe$ (barrier\_w), Form 1 has one PDE loss, Form 2 has two, Form 3 has three.
		In case of non-zero barrier weight $\bwe$, add one barrier loss (residue, barrier\_r) in each form.
		$\bullet$
		\emph{\red{Form 1:}} 
		Var du\_xx = 2nd deriv of disp computed in 
		% Listing~\ref{lst:DDE-T 1st and 2nd space derivatives} Line~\ref{ln:Form 1 2nd deriv du-xx}
		\mbox{L-Line}~\ref{lst:DDE-T 1st and 2nd space derivatives}.\ref{ln:Form 1 2nd deriv du-xx}.
		Var du\_tt = acceleration computed in 
		% Listing~\ref{lst:DDE-T time derivative by hessian} Line~\ref{ln:Form 1 acceleration du-tt} 
		L-Line~\ref{lst:DDE-T time derivative by hessian}.\ref{ln:Form 1 acceleration du-tt}
		or 
		% Listing~\ref{lst:DDE-T time derivative by two jacobians} Line~\ref{ln:Form 1 acceleration du-tt}
		L-Line~\ref{lst:DDE-T time derivative by two jacobians}.\ref{ln:Form 1 acceleration du-tt}.
		PDE loss = (du\_xx + f\_ext - du\_tt) $\Leftarrow$  
		% Eq.~\eqref{eq:eom-euler-bernoulli-axial} 
		Eq.~\eqref{eq:wave-eq-form-2}$_1$.
		$\bullet$
		\emph{\red{Form 2a:} Split 2nd time-derivative.} 
		Var py = momentum in 2nd col of y array (col index 1), 
		% Listing~\ref{lst:DDE-T F2a F3 time derivatives} Line~\ref{ln:Form 2a or Form 3 momentum}
		L-Line~\ref{lst:DDE-T F2a F3 time derivatives}.\ref{ln:Form 2a or Form 3 momentum}.  
		Var dy\_t = velocity computed in 
		% Listing~\ref{lst:DDE-T F2a F3 time derivatives} Line~\ref{ln:Form 2a or Form 3 velocity}
		L-Line~\ref{lst:DDE-T F2a F3 time derivatives}.\ref{ln:Form 2a or Form 3 velocity}.
		Momentum loss = (dy\_t - py) $\Leftarrow$ Eq.~\eqref{eq:wave-eq-form-2}$_2$.  
		Var du\_xx = 2nd deriv of disp computed in 
		% Listing~\ref{lst:DDE-T 1st and 2nd space derivatives} Line~\ref{ln:2nd space derivative du-xx component=0} 
		\mbox{L-Line}~\ref{lst:DDE-T 1st and 2nd space derivatives}.\ref{ln:2nd space derivative du-xx component=0}
		or 
		% Line~\ref{ln:Form 2a 2nd deriv by Jacobian du-xx}
		\mbox{L-Line}~\ref{lst:DDE-T 1st and 2nd space derivatives}.\ref{ln:Form 2a 2nd deriv by Jacobian du-xx}.
		Var py\_t = time derivative of momentum (acceleration) computed in 
		% Listing~\ref{lst:DDE-T F2a F3 time derivatives} Line~\ref{ln:Form 2a acceleration py-t}.
		L-Line~\ref{lst:DDE-T F2a F3 time derivatives}.\ref{ln:Form 2a acceleration py-t}.
		PDE loss = (du\_xx + f\_ext - py\_t) $\Leftarrow$ Eq.~\eqref{eq:wave-eq-form-2}$_1$.
		$\bullet$
		\emph{\red{Form 2b:} Split 2nd space-derivative.}
		Var y1 = slope in 2nd col of y array (col index 1) coded in 
		% Listing~\ref{lst:DDE-T, Form 2b, space derivative splitting} Line~\ref{ln:Form 2b y1 slope in 2nd col of y}.
		L-Line~\ref{lst:DDE-T, Form 2b, space derivative splitting}.\ref{ln:Form 2b y1 slope in 2nd col of y}. 
		Var du\_x = space derivative of displacement (slope) computed in 
		% Listing~\ref{lst:DDE-T, Form 2b, space derivative splitting} Line~\ref{ln:Form 2b slope du-x}.
		L-Line~\ref{lst:DDE-T, Form 2b, space derivative splitting}.\ref{ln:Form 2b slope du-x}.
		Slope loss = (y1 - du\_x) $\Leftarrow$ Eq.~\eqref{eq:wave-eq-form-2b}$_2$.
		Var y1\_x = 2nd deriv of disp computed in 
		% Listing~\ref{lst:DDE-T, Form 2b, space derivative splitting} Line~\ref{ln:Form 2b 2nd space deriv y1_x = du_xx}.
		L-Line~\ref{lst:DDE-T, Form 2b, space derivative splitting}.\ref{ln:Form 2b 2nd space deriv y1_x = du_xx}.
		Var du\_tt = acceleration computed in 
		% Listing~\ref{lst:DDE-T time derivative by two jacobians} Line~\ref{ln:Form 2b acceleration du-tt}.
		L-Line~\ref{lst:DDE-T time derivative by two jacobians}.\ref{ln:Form 2b acceleration du-tt}. 
		PDE loss = (y1\_x + f\_ext - du\_tt) $\Leftarrow$ Eq.~\eqref{eq:wave-eq-form-2b}$_1$.
		$\bullet$
		\emph{\red{Form 3:} Split 2nd derivatives wrt both space and time.}
		A combination of Form 2a and Form 2b.
		Listing~\ref{lst:DDE-T F2a F3 time derivatives} for py, dy\_t, and py\_t.
		Listing~\ref{lst:DDE-T, Form 3, space derivative splitting} for du\_x and y1\_x.
		PDE loss = (y1\_x + f\_ext - py\_t) $\Leftarrow$ Eq.~\eqref{eq:wave-eq-form-3}$_1$.
		Momentum loss = (dy\_t - py) $\Leftarrow$ Eq.~\eqref{eq:wave-eq-form-3}$_3$.
		Slope loss = (y1 - du\_x) $\Leftarrow$ Eq.~\eqref{eq:wave-eq-form-3}$_2$.
		% $\blacksquare$
		% \red{I AM HERE 23.10.9.}
	}, 
	label={lst:DDE-T, all forms, losses}
]
	if (Form == "F1"):
		if (barrier_w != 0):
			return [du_xx + f_ext - du_tt, barrier_r]
		else:
			return [du_xx + f_ext - du_tt]

	elif (Form == "F2a"):
		if (barrier_w != 0):
			return [du_xx + f_ext - py_t, dy_t - py, barrier_r]
		else:
			return [du_xx + f_ext - py_t, dy_t - py]

	elif (Form == "F2b"):
		if (barrier_w != 0):
			return [y1_x + f_ext - du_tt, y1 - du_x, barrier_r] |\label{ln:Form 2b losses, with barrier}|
		else:
			return [y1_x + f_ext - du_tt, y1 - du_x] |\label{ln:Form 2b losses, with NO barrier}|

	elif (Form == "F3"):
		if (barrier_w != 0):
			return [y1_x + f_ext - py_t, dy_t - py, y1 - du_x, barrier_r]
		else:
			return [y1_x + f_ext - py_t, dy_t - py, y1 - du_x]
\end{lstlisting}

\begin{lstlisting}[
	language=Python, 
	caption={
		\DDET. \emph{Axial motion of elastic bar.  \red{Forms 1, 2a, 2b, 3}.}  
		\emph{Left}-end boundary condition, \emph{\red{pinned}}, using DirichletBC operator.
	}, 
	label={lst:DDE-T, left-end boundary condition}
]
	# BOUNDARY CONDITIONS
	# left end, x = 0
	if (Form == "F1"):
		bc1 = dde.icbc.DirichletBC(geomtime, lambda x: 0, boundary_l)
	elif ((Form == "F2a") or (Form == "F2b") or (Form == "F3")):
		bc1 = dde.icbc.DirichletBC(geomtime, lambda x: 0, boundary_l, component=0)
\end{lstlisting}

\begin{lstlisting}[
	language=Python, 
	escapechar=|,
	caption={
		\DDET. \emph{Axial motion of elastic bar.  \red{Forms 1, 2a, 2b, 3}.}  
		\emph{Right}-end boundary condition, pinned or free, using DirichletBC or NeumannBC operator.
		\emph{\red{Form 1:}} 
		\emph{Pinned}, use DirichletBC, Line~\ref{ln:Form 1 pinned right end DirichletBC};
		\emph{Free}, use NeumannBC, Line~\ref{ln:Form 1 free right end NeumannBC}.
		\emph{\red{Form 2a:}}  
		\emph{Pinned}, use DirichletBC with component=0 (1st col of y array), Line~\ref{ln:Form 2a pinned right end DirichletBC};
		\emph{Free}, use NeumannBC with component=0 (1st col of y array), Line~\ref{ln:Form 2a free right end NeumannBC}.
		\emph{\red{Form 2b:}} 
		\emph{Pinned}, use DirichletBC with component=0 (1st col of y array, storing disp), Line~\ref{ln:Form 2b PINNED boundary condition};
		\emph{Free}, use DirichletBC with component=1 (2nd col of y array, storing slope), Line~\ref{ln:Form 2b FREE boundary condition}. 
		\emph{\red{Form 3:}}
		\emph{Pinned}, use DirichletBC with component=0 (1st col of y array, storing disp), Line~\ref{ln:Form 3 PINNED boundary condition};
		\emph{Free}, use DirichletBC with component=2 (3rd col of y array, storing slope), Line~\ref{ln:Form 3 FREE boundary condition}.
	}, 
	label={lst:DDE-T, right-end boundary condition}
]
	# right end, x = 1
	# Free right end, Newmann BC
	# bc2 = dde.icbc.NeumannBC(geomtime, lambda x: 0, boundary_r)
	# Pinned right end, Dirichlet BC
	if (Form == "F1"):
		# right-end boundary condition
		if (right_end_bc == "pinned"):
			# pinned-pinned bar, right-end pinned
			bc2 = dde.icbc.DirichletBC(geomtime, lambda x: 0, boundary_r) |\label{ln:Form 1 pinned right end DirichletBC}|
		elif (right_end_bc == "free"):
			# pinned-free bar, right-end free
			bc2 = dde.icbc.NeumannBC(geomtime, lambda x: 0, boundary_r) |\label{ln:Form 1 free right end NeumannBC}|
	
	elif (Form == "F2a"):
		# right-end boundary condition
		if (right_end_bc == "pinned"):
			# pinned-pinned bar, right-end pinned
			bc2 = dde.icbc.DirichletBC(geomtime, lambda x: 0, boundary_r, component=0) |\label{ln:Form 2a pinned right end DirichletBC}|
		elif (right_end_bc == "free"):
			# pinned-free bar, right-end free
			bc2 = dde.icbc.NeumannBC(geomtime, lambda x: 0, boundary_r, component=0) |\label{ln:Form 2a free right end NeumannBC}|
	
	elif (Form == "F2b"):
		# right-end boundary condition
		if (right_end_bc == "pinned"):
			# pinned-pinned bar, right-end pinned
			bc2 = dde.icbc.DirichletBC(geomtime, lambda x: 0, boundary_r, component=0) |\label{ln:Form 2b PINNED boundary condition}|
		elif (right_end_bc == "free"):
			# pinned-free bar, right-end free
			bc2 = dde.icbc.DirichletBC(geomtime, lambda x: 0, boundary_r, component=1) |\label{ln:Form 2b FREE boundary condition}|
	
	elif (Form == "F3"):
		# right-end boundary condition
		if (right_end_bc == "pinned"):
			# pinned-pinned bar, right-end pinned
			bc2 = dde.icbc.DirichletBC(geomtime, lambda x: 0, boundary_r, component=0) |\label{ln:Form 3 PINNED boundary condition}|
		elif (right_end_bc == "free"):
			# pinned-free bar, right-end free
			bc2 = dde.icbc.DirichletBC(geomtime, lambda x: 0, boundary_r, component=2) |\label{ln:Form 3 FREE boundary condition}|
\end{lstlisting}

\begin{lstlisting}[
	language=Python, 
	escapechar=|,
	caption={
		\DDET. \emph{Axial motion of elastic bar.  \red{Forms 1, 2a, 2b, 3}.}  
		Initial conditions. 
		$\bullet$
		\emph{Zero initial displacement.} \red{\emph{Form 1}}, Line~\ref{ln:Form 1 initial displacement}.  \emph{\red{Forms 2a, 2b, 3}}, with component=0, Line~\ref{ln:Forms 2a, 2b, 3 initial displacement}, since disp is stored in the 1st col in the output y array in all these forms.  
		For example, for \emph{\red{Form 2a}}, see L-Line~\ref{lst:DDE-T F2a F3 time derivatives}.\ref{ln:Form 2a or Form 3 momentum}, where the momentum py was placed in the 2nd col, i.e., \unexpanded{py = y[:, 1:2]} after the disp in the 1st col \unexpanded{y[:, 0]}, and the computation of the Hessian (or 2nd space derivative) of the disp in L-Line~\ref{lst:DDE-T F2a F3 time derivatives}.\ref{ln:2nd space derivative du-xx component=0} must indicate that the disp is component=0 of the y array. 
		$\bullet$
		\emph{Zero initial velocity.}
		\emph{\red{Forms 1, 2b:}} 
		Since the 2nd time derivative is not split in these forms, to impose the initial velocity, use the function dde.icbc.OperatorBC to define the boundary condition (BC) operator with
		the Jacobian to compute the derivative of disp (col index i=0 in y array) wrt time (col index 1 in x array), i.e., velocity, which is implicitly equated to zero, Line~\ref{ln:Form 1 or Form 2b initial condition jacobian}. 
		\emph{\red{\mbox{Forms 2a, 3:}}} 
		Since the 2nd time derivative (acceleration) is split in these forms, the velocity is readily available in the 2nd col of the y array (col index 1, thus component=1, L-Line~\ref{lst:DDE-T F2a F3 time derivatives}.\ref{ln:Form 2a or Form 3 momentum}), and can be directly used to impose the initial velocity without having to take the derivative, Line~\ref{ln:Form 2a or Form 3 initial velocity component=1}.
		% \red{I AM HERE 23.10.10.}
	}, 
	label={lst:DDE-T, initial conditions}
]
	# INITIAL CONDITIONS
	if (Form == "F1"):
		ic1 = dde.icbc.IC(geomtime, lambda x: 0, lambda _, on_initial: on_initial) |\label{ln:Form 1 initial displacement}|
	elif ((Form == "F2a") or (Form == "F2b") or (Form == "F3")):
		ic1 = dde.icbc.IC(geomtime, lambda x: 0, lambda _, on_initial: on_initial, component=0) |\label{ln:Forms 2a, 2b, 3 initial displacement}|

	if ((Form == "F1") or (Form == "F2b")):
		# the code below is for prescribed velocity
		ic2 = dde.icbc.OperatorBC(
			geomtime,
			# Zero initial velocity
			lambda x, y, _: dde.grad.jacobian(y, x, i=0, j=1), |\label{ln:Form 1 or Form 2b initial condition jacobian}|
			#
			lambda _, on_initial: on_initial,
		)
	
	elif ((Form == "F2a") or (Form == "F3")):
		# Dirichlet BC for velocity
		ic2 = dde.icbc.IC(geomtime, lambda x: 0, lambda _, on_initial: on_initial, component=1) |\label{ln:Form 2a or Form 3 initial velocity component=1}|
\end{lstlisting}

% \subsection{\JAX\ code snippets}
\subsection{JAX script snippets}
\noindent
In what follows, the quintessential parts of our \JAX\ script, which does not rely on any \DDET\ features, are to be outlined. 
Its documentation describes \JAX\ as \textit{``\href{https://github.com/hips/autograd}{Autograd} and \href{https://www.tensorflow.org/xla}{XLA} [Accelerated Linear Algebra], brought together for high-performance numerical computing''}, i.e., \JAX' strengths lie on automatic differentation and parallel computation. 
As opposed to well-known machine-learning frameworks as \href{https://www.tensorflow.org}{TensorFlow} and \href{https://pytorch.org}{PyTorch}, \JAX\ adopts a functional programming paradigm, which is already visible in the very first code snippet 
%(Listing \ref{lst:JAX code 1}) 
defining the PDE loss of the axial bar.

\begin{lstlisting}[
	language=Python,
	escapechar=|, 
	caption={
		\JAX. Import
		$\bullet$
		jax (\href{https://jax.readthedocs.io/en/latest/notebooks/quickstart.html}{JAX Quickstart})
		$\bullet$
		\href{https://jax.readthedocs.io/en/latest/jax.numpy.html}{jax.numpy} (JAX Accelerated NumPy) as `jnp'
		$\bullet$
		\href{https://jax.readthedocs.io/en/latest/jax.nn.html}{jax.nn} (JAX neural network module) as `nn'
		$\bullet$
		\href{https://jax.readthedocs.io/en/latest/jax.example_libraries.optimizers.html}{jax.example\_libraries.optimizers} module  (Examples of how to write optimizers with JAX) as `optimizer'
		$\bullet$
		\href{https://jax.readthedocs.io/en/latest/_autosummary/jax.vmap.html}{vmap} (\href{https://jax.readthedocs.io/en/latest/jax.html\#vectorization-vmap}{Vectorizing map}) submodule of jax module
		$\bullet$
		\href{https://jax.readthedocs.io/en/latest/jax-101/02-jitting.html}{jit} (Just In Time Compilation) submodule
		$\bullet$
		\href{https://jax.readthedocs.io/en/latest/jax-101/01-jax-basics.html\#jax-first-transformation-grad}{grad} submodule (takes a numerical function written in Python and returns a new Python function that computes the gradient of the original function)
		$\bullet$
		\href{https://jax.readthedocs.io/en/latest/\_autosummary/jax.vjp.html\#jax.vjp}{vjp} (\href{https://jax.readthedocs.io/en/latest/notebooks/autodiff\_cookbook.html\#vector-jacobian-product}{vector-Jacobian product}) submodule \unexpanded{[Compute a (reverse-mode) vjp of a function]}
		$\bullet$
		\href{https://optax.readthedocs.io/en/latest/}{optax} module (gradient processing and optimization library for JAX).	
	}, 
	label={lst:JAX code 1}
]
import jax
import jax.numpy as jnp
import jax.nn as nn
import jax.example_libraries.optimizers as optimizers
from jax import vmap, jit, grad, vjp
import optax
\end{lstlisting}

%\noindent
%\red{[NOTE: 23.10.15 Check ``??'' to comment out due to figures in commented-out appendix. Add tag \DDET\ and \JAX\ to all code snippets. ENDNOTE]}

\begin{lstlisting}[
	language=Python,
	escapechar=|,
	caption={
		\JAX.
		$\bullet$
		Line~\ref{ln:JAX init-network-params}: Initialize network parameters using function init\_network\_params (see
		\href{https://jax.readthedocs.io/en/latest/notebooks/neural_network_with_tfds_data.html\#hyperparameters}{Hyperparameters} in
		\href{https://jax.readthedocs.io/en/latest/notebooks/neural_network_with_tfds_data.html}{Training a Simple Neural Network} in
		\href{https://jax.readthedocs.io/en/latest/advanced_guide.html}{JAX Advanced Tutorials}).
		$\bullet$
		Lines~\ref{ln:JAX func pinn BEGIN}-\ref{ln:JAX func pinn END}: Define function `pinn' to predict the output given the input as a single pair of space-time coordinates $(x, t)$, a collocation point.  
		% The input is not the whole x array that contains all pairs of $(x, t)$ points (in a vertical stack) so that vectorization using `vmap' can be used to process all points in parallel.  
		The \JAX\ vectorization function 
		% vmap 
		\href{https://jax.readthedocs.io/en/latest/_autosummary/jax.vmap.html}{vmap}
		processes simultaneously all collocation points $(x, t)$, arranged in a vertical stack in the space-time x array. 
		See function `predict' in  \href{URLhttps://jax.readthedocs.io/en/latest/notebooks/neural_network_with_tfds_data.html\#auto-batching-predictions}{Auto-batching predictions}.
		$\triangleright$
		Line~\ref{ln:JAX for w, b loop}: Loop over the layer pairs, each represented by (w, b), except the last layer pair.
		$\circ$
		Line~\ref{ln:JAX outputs}: For each layer pair, compute the outputs from inputs.
		$\circ$
		Line~\ref{ln:JAX activations}: Apply tanh activation function to outputs.
		$\triangleright$
		Lines~\ref{ln:JAX final_w}-\ref{ln:JAX output}: Compute network output from last layer pair without activation function.
		$\bullet$
		Line~\ref{ln:JAX func pinn END}: Return network output.
		$\bullet$
		Line~\ref{ln:JAX y_fun_param}: define function y\_fun\_param with two arguments `params' and `x' and `jit'-compile this function, i.e., subject the function to \href{https://jax.readthedocs.io/en/latest/_autosummary/jax.jit.html}{Just-In-Time} compilation.
		$\bullet$
		Line~\ref{ln:JAX y_fun}: Define function y\_fun to have only the space-time `x' array as argument to simplify taking the  partial differentiation wrt `x'.
		% \red{[NOTE: 23.10.17, Tue, added after our discussion.  To document. ENDNOTE]}
	},
	label={lst:JAX define y-fun}
]
	layer_size = [n_inp] + [W] * H + [n_out]
	params = init_network_params(layer_size, key) |\label{ln:JAX init-network-params}|
	
	def pinn(params, input):	|\label{ln:JAX func pinn BEGIN}|
		activations = input.reshape(-1)
		for w, b in params[:-1]:					|\label{ln:JAX for w, b loop}|
			outputs = jnp.dot(w, activations) + b	|\label{ln:JAX outputs}|
			activations = nn.tanh(outputs)			|\label{ln:JAX activations}|
		
		final_w, final_b = params[-1]				|\label{ln:JAX final_w}|
		output = jnp.dot(final_w, activations) + final_b	|\label{ln:JAX output}|
	return output  				|\label{ln:JAX func pinn END}|
	
	y_fun_param = jit(lambda params, x: pinn(params, x))  	|\label{ln:JAX y_fun_param}|
	y_fun = lambda x: y_fun_param(params, x)				|\label{ln:JAX y_fun}|
\end{lstlisting}
% grad_y_fun = grad(lambda params, x: y_fun_param(params, x)[0])

\begin{lstlisting}[
	language=Python,
	escapechar=|, 
	caption={
		\JAX.  \emph{Axial motion of elastic bar}. \emph{\red{Form 1}}.
		$\bullet$
		Line~\ref{ln:JAX F1 u}: Define function u as the 1st element of the output (col index 0) of y\_fun applied on the 
		% coordinate pair 
		collocation point
		$(x, t)$ in x array. 
		$\bullet$
		Line~\ref{ln:JAX F1 grad(u)}: Set du as gradient of u.
		$\bullet$
		Line~\ref{ln:JAX F1 u_x}: Define u\_x as 1st elmt (col index 0) of the output of du evaluated at $(x, t)$ in x array.
		$\bullet$
		Line~\ref{ln:JAX F1 u_t}: Define u\_t as 2nd elmt (col index 1) of the output of du evaluated at $(x, t)$ in x array.
		$\bullet$
		Line~\ref{ln:JAX F1 du_x}: Set du\_x as grad of u\_x.
		$\bullet$
		Line~\ref{ln:JAX F1 du_t}: Set du\_t as grad of u\_t.
		$\bullet$
		Line~\ref{ln:JAX F1 du_xx}: Define du\_xx, 2nd space deriv of u, as 1st elmt (col index 0) of the output of du\_x, evaluated at $(x, t)$ in x array.
		$\bullet$
		Line~\ref{ln:JAX F1 du_tt}: Define du\_tt, 2nd time deriv of u, as 2nd elmt (col index 1) of the output of du\_t, evaluated at $(x, t)$ in x array.
		$\bullet$
		Line~\ref{ln:JAX F1 pde}: Define pde-loss func.
		$\bullet$
		Line~\ref{ln:JAX F1 return jnp array}: Return as jnp array.
		$\bullet$
		Line~\ref{ln:JAX F1 pde at x}: The return object is pde-loss func evaluated at x array. The vectorized computation of the pde loss of \emph{\red{Form 1}} is in L-Line~\ref{lst:JAX Form 3}.\ref{ln:JAX F3 pde_vmap}.
		$\blacktriangleright$
		See~\DDET\ \mbox{\emph{\red{Form 1}}} Listings~\ref{lst:DDE-T time derivative by hessian}-\ref{lst:DDE-T time derivative by two jacobians} for time derivatives, and Listing~\ref{lst:DDE-T 1st and 2nd space derivatives} for space derivatives.  
	}, 
	label={lst:JAX Form 1}
]		
if (Form == "F1"):
	def pde(x, y_fun):				|\label{ln:JAX F1 def pde}|
		u = lambda x: y_fun(x)[0] 	|\label{ln:JAX F1 u}|
		du = grad(u)				|\label{ln:JAX F1 grad(u)}|
		u_x = lambda x: du(x)[0]	|\label{ln:JAX F1 u_x}|
		u_t = lambda x: du(x)[1]	|\label{ln:JAX F1 u_t}|
		
		du_x = grad(u_x)			|\label{ln:JAX F1 du_x}|
		du_t = grad(u_t)			|\label{ln:JAX F1 du_t}|
		u_xx = lambda x: du_x(x)[0] |\label{ln:JAX F1 du_xx}|
		u_tt = lambda x: du_t(x)[1] |\label{ln:JAX F1 du_tt}|
		
		pde = lambda x: u_xx(x) + f_ext - u_tt(x) |\label{ln:JAX F1 pde}|
		
		return jnp.asarray((		|\label{ln:JAX F1 return jnp array}|
			pde(x)					|\label{ln:JAX F1 pde at x}|
		))
\end{lstlisting}

\begin{rem}
	\label{rm:JAX Form 1 snippet}
	Listing~\ref{lst:JAX Form 1}: JAX \red{Form 1}.
	{\rm
		The argument 
		% \texttt{y\_fun} 
		y\_fun of function pde in L-Line~\ref{lst:JAX Form 1}.\ref{ln:JAX F1 def pde} is a function object that represents the neural network, which can be evaluated at the collocation points 
		% \texttt{x} 
		$(x, t)$ in the x array.
		% in the computational domain.
		L-Line~\ref{lst:JAX Form 1}.\ref{ln:JAX F1 u} wraps the first (and, for \emph{\red{Form 1}}, sole) output (index 0) of the neural network into a function 
		% \texttt{u},
		u,  
		the (reverse-mode) gradient of which (Var 
		% \texttt{du})
		du) 
		is computed using \JAX's 
		% \texttt{grad} 
		\href{https://jax.readthedocs.io/en/latest/jax-101/01-jax-basics.html#jax-first-transformation-grad}{grad}
		function in \mbox{L-Line}~\ref{lst:JAX Form 1}.\ref{ln:JAX F1 grad(u)}.
		Var 
		% \texttt{du} 
		du
		is a function object comprising both spatial and temporal derivatives, which are isolated as functions % \texttt{u\_x} and \texttt{u\_t}, 
		u\_x and u\_t, respectively (L-Line~\ref{lst:JAX Form 1}.\ref{ln:JAX F1 u_x} and 
		L-Line~\ref{lst:JAX Form 1}.\ref{ln:JAX F1 u_t}).
		Functions 
		% \texttt{u\_xx} and \texttt{u\_tt} 
		u\_xx and u\_tt
		representing the second derivatives are set up analogously. 
		L-Line~\ref{lst:JAX Form 1}.\ref{ln:JAX F1 pde} sets up a function 
		% \texttt{pde\_fun} 
		pde
		% that evaluates the PDE loss for a given set of space-time points 
		that evaluates the PDE loss for a given collocation point 
		% \texttt{x}.
		x.
		Var 
		% \texttt{f\_ext}
		f\_ext 
		is the (constant) distributed load.
		The function 
		% \texttt{pde\_fun} 
		pde
		is eventually evaluated in the 
		% \texttt{return} 
		return
		statement 
		(L-Line~\ref{lst:JAX Form 1}.\ref{ln:JAX F1 return jnp array}), 
		which returns a \JAX\ array containing the PDE loss at 
		the collocation space-time point $(x, t)$ in 
		% which returns a \JAX\ array containing the PDE losses at 
		% all collocation space-time points $(x, t)$ in 
		% \texttt{x}.
		x array.
	}
	\hfill$\blacksquare$
\end{rem}

In our implementation of PINNs, we used \JAX's 
% \texttt{grad} 
\href{https://jax.readthedocs.io/en/latest/jax-101/01-jax-basics.html#jax-first-transformation-grad}{grad}
function to compute derivatives, relying on the reverse-mode automatic differentiation, commonly used when training neural networks by means of backpropagation. 
Note that \JAX\ does have more possiblities to take derivatives, including forward-mode differentiation (see, e.g., the function 
% \texttt{jacfwd}). 
\href{https://jax.readthedocs.io/en/latest/_autosummary/jax.jacfwd.html}{jacfwd}),
with the \JAX\ function 
% \texttt{grad} 
\href{https://jax.readthedocs.io/en/latest/jax-101/01-jax-basics.html#jax-first-transformation-grad}{grad}
offering the simplest interface, even though it is not necessarily the most efficient way to compute and evaluate derivatives. 

\begin{lstlisting}[
	language=Python,
	escapechar=|, 
	caption={
		\JAX. \emph{Axial motion of elastic bar.  \red{Form 3}. Vectorized computation of PDE loss.}
		$\bullet$
		Line~\ref{ln:JAX F3 u}: Define func u as the 1st elmt (col index 0) of the output of y\_fun applied on $(x, t)$ in x array.
		% $\bullet$
		% Line~\ref{ln:JAX F1 grad(u)}: Define func du as gradient of u.
		$\bullet$
		Line~\ref{ln:JAX F3 u_x}: Define u\_x as the 2nd elmt (col index 1) of the output of y\_fun applied on $(x, t)$.
		$\bullet$
		Line~\ref{ln:JAX F3 u_t}: Define u\_t as the 3rd elmt (col index 2) of the output of y\_fun applied on $(x, t)$.
		$\bullet$
		Line~\ref{ln:JAX F3 du_x}: Set du\_x as grad of u\_x.
		$\bullet$
		Line~\ref{ln:JAX F3 du_t}: Set du\_t as grad of u\_t.
		$\bullet$
		Line~\ref{ln:JAX F3 u_xx}: Define u\_xx as 1st elmt (col index 0) of the output of du\_x applied on $(x, t)$.
		$\bullet$
		Line~\ref{ln:JAX F3 u_tt}: Define u\_tt as 2nd elmt (col index 1) of the output of du\_t at $(x, t)$.
		$\bullet$
		Line~\ref{ln:JAX F3 grad(u)}: Set du as grad of u.
		$\bullet$
		Line~\ref{ln:JAX F3 u_x_}: Define u\_x\_ as 1st elmt (col index 0) of the output of du at $(x, t)$.
		$\bullet$
		Line~\ref{ln:JAX F3 u_t_}: Define u\_t\_ as 2nd elmt (col index 1) of du at $(x, t)$.
		$\bullet$
		Line~\ref{ln:JAX F3 pde1}: Define pde1 as pde loss.
		$\bullet$
		Line~\ref{ln:JAX F3 pde2}: Define pde2 as space-slope loss.
		$\bullet$
		Line~\ref{ln:JAX F3 pde3}: Define pde3 as time-slope loss.
		$\bullet$
		Line~\ref{ln:JAX F3 return jnp array}: Return output as jnp array.
		$\bullet$
		Line~\ref{ln:JAX F3 pde1 at x}: Return pde loss at $(x, t)$.
		$\bullet$
		Line~\ref{ln:JAX F3 pde2 at x}: Return space-slope loss at $(x, t)$.
		$\bullet$
		Line~\ref{ln:JAX F3 pde3 at x}: Return time-slope loss at $(x, t)$.
		$\bullet$
		Line~\ref{ln:JAX F3 pde_vmap}: Vectorized computation of pde loss (for either \emph{\red{Form 1}} or \emph{\red{Form 3}}) with vmap and its argument (0, None), then jit compile the vectorized pde loss with static\_argnums=1 (see explanation below).
		% \red{NOTE: 23.11.01, need to explain "(0, None)" and "static\_argnums=1".}
		$\blacktriangleright$
		See \DDET\ \mbox{\emph{\red{Form 3}}} in Listing~\ref{lst:DDE-T time derivative by hessian}.
	}, 
	label={lst:JAX Form 3}
	]	
	if (Form == "F3"):
		def pde(x, y_fun):				|\label{ln:JAX F3 def pde BEGIN}|
			u = lambda x: y_fun(x)[0]	|\label{ln:JAX F3 u}|
			u_x = lambda x: y_fun(x)[1] |\label{ln:JAX F3 u_x}|
			u_t = lambda x: y_fun(x)[2] |\label{ln:JAX F3 u_t}|
			
			du_x = grad(u_x)			|\label{ln:JAX F3 du_x}|
			du_t = grad(u_t)			|\label{ln:JAX F3 du_t}|
			u_xx = lambda x: du_x(x)[0]	|\label{ln:JAX F3 u_xx}|
			u_tt = lambda x: du_t(x)[1]	|\label{ln:JAX F3 u_tt}|
			
			du = grad(u)				|\label{ln:JAX F3 grad(u)}|
			u_x_ = lambda x: du(x)[0]	|\label{ln:JAX F3 u_x_}|
			u_t_ = lambda x: du(x)[1]	|\label{ln:JAX F3 u_t_}|
			
			pde1 = lambda x: u_xx(x) + f_ext - u_tt(x)	|\label{ln:JAX F3 pde1}|
			pde2 = lambda x: u_x(x) - u_x_(x)			|\label{ln:JAX F3 pde2}|
			pde3 = lambda x: u_t(x) - u_t_(x)			|\label{ln:JAX F3 pde3}|
			
			return jnp.asarray((	|\label{ln:JAX F3 return jnp array}|
				pde1(x),			|\label{ln:JAX F3 pde1 at x}|
				pde2(x),			|\label{ln:JAX F3 pde2 at x}|
				pde3(x)				|\label{ln:JAX F3 pde3 at x}|
			))							|\label{ln:JAX F3 def pde END}|
	
	pde_vmap_jit = jit(vmap(pde, (0, None)), static_argnums=1) |\label{ln:JAX F3 pde_vmap}|
\end{lstlisting}

For \JAX\ \emph{\red{Form 3}} in Listing~\ref{lst:JAX Form 3}, the neural network has three outputs, i.e., the displacement 
% (\texttt{u}) 
u
and its space derivative 
% (\texttt{u\_x}) 
u\_x
and its time derivative 
% (\texttt{u\_t}), 
u\_t,
as indicated in L-Lines~\ref{lst:JAX Form 3}.\ref{ln:JAX F3 u}-\ref{ln:JAX F3 u_t}.
To set up functions for evaluating the second derivatives, L-Lines~\ref{lst:JAX Form 3}.\ref{ln:JAX F3 du_x}-\ref{ln:JAX F3 du_t} employ
automatic differentiation of the first derivatives u\_x and u\_t via the
%\texttt{grad} 
\href{https://jax.readthedocs.io/en/latest/jax-101/01-jax-basics.html#jax-first-transformation-grad}{grad}
function, which 
%when applied to the computed displacement u provides 
is also used to provide
the  
first derivatives of the displacement u wrt 
space (u\_x\_ in L-Line~\ref{lst:JAX Form 3}.\ref{ln:JAX F3 u_x_}, 1st element of du(x), index 0) and wrt
time  (u\_t\_ in L-Line~\ref{lst:JAX Form 3}.\ref{ln:JAX F3 u_t_}, 2nd element of du(x), index 1).
The loss of the 1st space derivative (i.e., space-slope loss in L-Line~\ref{lst:JAX Form 3}.\ref{ln:JAX F3 pde2}) is the difference between the independent variable u\_x defined in L-Line~\ref{lst:JAX Form 3}.\ref{ln:JAX F3 u_x}) and the 1st space derivative u\_x\_ in L-Line~\ref{lst:JAX Form 3}.\ref{ln:JAX F3 u_x_}, both evaluated at the space-time point $(x, t)$ of the x array.
Similarly for the loss of the 1st time derivative in L-Line~\ref{lst:JAX Form 3}.\ref{ln:JAX F3 pde3}.
The total loss is the sum of the pde loss, pde1, the space-slope loss, pde2, and the time-slope loss, pde3; see Eq.~\eqref{eq:wave-eq-form-3}.

In L-Line~\ref{lst:JAX Form 3}.\ref{ln:JAX F3 pde_vmap}, the
% vmap 
\href{https://jax.readthedocs.io/en/latest/_autosummary/jax.vmap.html}{vmap} operation (vectorization)
is applied first on the pde loss function defined in \mbox{L-Lines}~\ref{lst:JAX Form 3}.\ref{ln:JAX F3 def pde BEGIN}-\ref{ln:JAX F3 def pde END}, followed by a
% jit.  
\href{https://jax.readthedocs.io/en/latest/_autosummary/jax.jit.html}{jit} compilation (Just-In-Time).
The arguments ``pde, (0, None)'' of 
% vmap 
\href{https://jax.readthedocs.io/en/latest/_autosummary/jax.vmap.html}{vmap}
has the pde loss function and the tuple ``(0, None)'' (which is the value of the in\_axes field of 
% vmap) 
\href{https://jax.readthedocs.io/en/latest/_autosummary/jax.vmap.html}{vmap})
corresponding to the arguments ``(x, y\_fun),'' respectively, of the pde loss function defined in L-Line~\ref{lst:JAX Form 3}.\ref{ln:JAX F3 def pde BEGIN}, with ``0'' meaning vectorizing row-wise the `x' array, i.e., process the rows (collocation points) in the x array simultaneously, and with ``None'' meaning no vectorization for the function y\_fun.
The second argument ``static\_argnums=1'' of \href{https://jax.readthedocs.io/en/latest/_autosummary/jax.jit.html}{jit} tells the jit compilation to treat the second argument y\_fun of the `pde' loss function as static. 
Since the composition of the
\href{https://jax.readthedocs.io/en/latest/_autosummary/jax.vmap.html}{vmap} operation and the
\href{https://jax.readthedocs.io/en/latest/_autosummary/jax.jit.html}{jit} compilation can be 
written in two lines, and can also be
composed in reverse order, 
Listing~\ref{lst:alternatives for pde_vmap_jit} provides three alternative coding methods for L-Line~\ref{lst:JAX Form 3}.\ref{ln:JAX F3 pde_vmap}.

\begin{lstlisting}[
	language=Python,
	escapechar=|, 
	caption={
		\JAX. 
		Alternatives to pde\_vmap\_jit in L-Line~\ref{lst:JAX Form 3}.\ref{ln:JAX F3 pde_vmap}.
		Either Method 2, or \mbox{Method 3}, or Method 4 can replace Method 1 with no GPU-time penalty.
		$\bullet$
		Line~\ref{ln:JAX F3 pde_vmap_jit} $\equiv$ L-Line~\ref{lst:JAX Form 3}.\ref{ln:JAX F3 pde_vmap}: Method 1, vmap first, jit after, combined in 1 line,
		GPU time 547 s {\scriptsize (23916R1a91c)}.
		$\bullet$
		Lines~\ref{ln:JAX F3 pde_vmap_jit BEGIN}-\ref{ln:JAX F3 pde_vmap_jit END}: Method 2, vmap first, jit after in 2 lines, 
		% with no GPU-time penalty, i.e., approximately the same GPU time as for Method 1.
		GPU time 547 s {\scriptsize (23119R1b)}.
		$\bullet$
		Line~\ref{ln:JAX F3 pde_jit_vmap}: Method 3, reverse composition, jit first, vmap after, in 1 line, 
		% with no GPU-time penalty.
		GPU time 538 s {\scriptsize (23119R1c)}.
		$\bullet$
		Lines~\ref{ln:JAX F3 pde_jit_vmap BEGIN}-\ref{ln:JAX F3 pde_jit_vmap END}: Method 4, jit first, vmap after in 2 lines, GPU time 539 s {\scriptsize (23119R1d)}. 
		% no GPU-time penalty.  
		All 4 methods required very close GPU times to complete 200,000 steps, with the same network and training parameters as in SubFigs~\ref{fig:23.9.16 R1a.8 loss100000}-\ref{fig:23.9.16 R1a.8 free-end disp100000} for \emph{\red{Form 1}}.
	}, 
	label={lst:alternatives for pde_vmap_jit}
]
	#-------------------------------------------------------------------------------
	# Method 1: vmap first, jit after combined in 1 line
	pde_vmap_jit = jit(vmap(pde, (0, None)), static_argnums=1)		|\label{ln:JAX F3 pde_vmap_jit}|
	#-------------------------------------------------------------------------------
	# Method 2: vmap first, jit after in 2 lines, no GPU-time penalty
	pde_vmap = vmap(pde, (0, None))									|\label{ln:JAX F3 pde_vmap_jit BEGIN}|
	pde_vmap_jit = jit(pde_vmap, static_argnums=1)					|\label{ln:JAX F3 pde_vmap_jit END}|
	#-------------------------------------------------------------------------------
	# Method 3: jit first, vmap after combined in 1 line, with no GPU-time penalty
	pde_jit_vmap = vmap(jit(pde, static_argnums=1), (0, None))		|\label{ln:JAX F3 pde_jit_vmap}|
	#-------------------------------------------------------------------------------
	# Method 4: jit first, vmap after in 2 lines, no GPU-time penalty
	pde_jit = jit(pde, static_argnums=1)							|\label{ln:JAX F3 pde_jit_vmap BEGIN}|
	pde_jit_vmap = vmap(pde_jit, (0, None))							|\label{ln:JAX F3 pde_jit_vmap END}|
\end{lstlisting}

\begin{lstlisting}[
	language=Python,
	escapechar=|, 
	caption={
		\JAX. 
		\emph{Axial motion of elastic bar.  Boundary conditions}  (BCs).  \emph{\red{\mbox{Form 1, 2a, or 3}}}.
		$\bullet$
		Lines~\ref{ln:JAX bc_left pinned BEGIN}-\ref{ln:JAX bc_left pinned END}: Define left-end \emph{pinned} BC.
		$\bullet$
		Lines~\ref{ln:JAX bc_right pinned BEGIN}-\ref{ln:JAX bc_right pinned END}: Define right-end \emph{pinned} BC.
		$\bullet$
		Lines~\ref{ln:JAX bc_right free BEGIN}-\ref{ln:JAX bc_right free END}: Define right-end \emph{free} BC for either \emph{\red{\mbox{Form 1, or 2a, or 3}}}.
	}, 
	label={lst:JAX BCs}
]	
	def bc_left(xt_left, y_fun):		|\label{ln:JAX bc_left pinned BEGIN}|
		u = lambda x: y_fun(x)[0]
		return u(xt_left)
	bc_left_vmap = jit(vmap(bc_left, (0, None)), static_argnums=1)	|\label{ln:JAX bc_left pinned END}|
	
	if (right_end_bc == "pinned"):		|\label{ln:JAX bc_right pinned BEGIN}|
		def bc_right(xt_right, y_fun):
			u = lambda x: y_fun(x)[0]
			return u(xt_right)
		bc_right_vmap = jit(vmap(bc_right, (0, None)), static_argnums=1)	|\label{ln:JAX bc_right pinned END}|
			
	if (right_end_bc == "free"):										|\label{ln:JAX bc_right free BEGIN}|
		def bc_right(xt_right, y_fun):
			if (Form == "F1" or Form == "F2a"):
				u = lambda x: y_fun(x)[0]
				du = grad(u)
				u_x = lambda x: du(x)[0]
			
			if (Form == "F3"):
				u_x = lambda x: y_fun(x)[1]
			
			return jnp.asarray((
				u_x(xt_right)
			))
		
		bc_right_vmap = jit(vmap(bc_right, (0, None)), static_argnums=1) |\label{ln:JAX bc_right free END}|
\end{lstlisting}

\begin{lstlisting}[
language=Python,
escapechar=|, 
caption={
	\JAX. 
	\emph{Axial motion of elastic bar.  Initial conditions}  (ICs).  \emph{\red{Forms 1, 2a, or 3}}.
	$\bullet$
	Lines~\ref{ln:JAX F1 IC BEGIN}-\ref{ln:JAX F1 IC END}: Initial conditions for \emph{\red{Form 1}}.
	$\bullet$
	Lines~\ref{ln:JAX F2a IC BEGIN}-\ref{ln:JAX F2a IC END}: Initial conditions for \emph{\red{Form 2a}}.
	$\bullet$
	Lines~\ref{ln:JAX F3 IC BEGIN}-\ref{ln:JAX F3 IC END}: Initial conditions for \emph{\red{Form 3}}.
}, 
label={lst:JAX ICs}
]
def ic(xt_ic, y_fun):
	u = lambda x: y_fun(x)[0]
	
	# Form 1											|\label{ln:JAX F1 IC BEGIN}|
	if (Form == "F1"):
		# find time derivative u_t using grad function
		du = grad(u)
		u_t = lambda x: du(x)[1]						|\label{ln:JAX F1 IC END}|
	
	# Form 2a											|\label{ln:JAX F2a IC BEGIN}|
	if (Form == "F2a"):
		# set u_t as 2nd col (index 1) of y_fun
		u_t = lambda x: y_fun(x)[1]						|\label{ln:JAX F2a IC END}|
	
	# Form 3											|\label{ln:JAX F3 IC BEGIN}|	
	if (Form == "F3"):
		# set u_t as 3rd col (index 2) of y_fun
		u_t = lambda x: y_fun(x)[2]						|\label{ln:JAX F3 IC END}|
	
	return jnp.asarray((
		u(xt_ic),
		u_t(xt_ic)
	)).T

ic_vmap = jit(vmap(ic, (0, None)), static_argnums=1)
\end{lstlisting}

	% optimization and iterate evolution - NOT written yet
	% \input{Optimization and iteration evolution}
	
	% trial and error for non-static solution
	% capacity and static solution - NOT written yet
	% \input{Capacity and static solution}
	
\end{appendices}

\end{document}